%% file: Diverses.tex
\documentclass[11pt,a4paper,german]{book}

\input{MathPreamble2.tex}

\usepackage{ifthen}
\usepackage{calc}
\usepackage{bibgerm}   
\usepackage{fancyhdr}
\usepackage{verbatim}  
\usepackage{mathrsfs}  

\usepackage{makeidx}
\usepackage{index}
\usepackage{lmodern}

\usepackage{hyperref}

\usepackage[all]{xy}

\newboolean{finalVersion}
\setboolean{finalVersion}{true}

\newcommand{\textVersion}{2.5.2009}

\newboolean{insertdComments}
\newboolean{insertdcite}
\newboolean{showVersion}  
\newboolean{schreibidx}   
\setboolean{schreibidx}{true}

\ifthenelse{\boolean{finalVersion}}
{\setboolean{insertdComments}{false}
 \setboolean{insertdcite}{false}
 \setboolean{showVersion}{true}
 \setboolean{showVersion}{false}  
}
{\setboolean{insertdComments}{true}
 \setboolean{insertdcite}{true}
  \setboolean{showVersion}{true}
}

\newlength{\lengthtempa}  

\newcommand{\defemph}[1]{\textbf{#1\/}}   
\newcommand{\defemphi}[1]{\textbf{#1\/}\index[B]{#1}}  
\newcommand{\bewitemph}[1]{\emph{#1\/}}   
\newcommand{\bewitemphq}[1]{\bewitemph{\glqq#1\grqq:}}
\newcommand{\bewitemphqq}[1]{\bewitemph{\glqq#1\grqq}}

\newcommand{\dheisst}{d.\,h.\xspace}  
\newcommand{\iallg}{i.\,a.\xspace}   

\newcommand{\abbsa}[1]{{#1}^\star}%
\newcommand{\absBig}[1]{\Bigl\lvert#1\Bigr\rvert}
\newcommand{\abslr}[1]{\left\lvert#1\right\rvert}
\newcommand{\Adjlco}[1]{\mathscr{A}_{l,s}^{C_0}(#1)}
\newcommand{\Adjhm}[1]{\mathbb{B}(#1)}
\newcommand{\Adjhmi}[2]{\mathbb{B}_{#1}(#2)}
\newcommand{\Adjloru}[1]{\mathscr{A}_l^{C_0}(#1)}
\newcommand{\Adjlor}[1]{\mathscr{A}_l(#1)}
\newcommand{\adjor}[1]{{#1}^\star}
\newcommand{\alg}{\operatorname{alg}}
\newcommand{\arrowbij}{\rightarrowtail\!\!\!\!\!\rightarrow}
\newcommand{\cohg}{$C_0$"=Halbgruppe\xspace}%
\newcommand{\cohgn}{$C_0$"=Halbgruppen\xspace}%
\newcommand{\cketbra}[2]{|#1\rangle\langle#2|}
\newcommand{\cmGamma}{\mathcal{M}}%
\newcommand{\cmiZ}{i.\,Z.\xspace}
\newcommand{\cmkug}{B}
\newcommand{\cmlin}{\operatorname{lin}}%
\newcommand{\cmlmult}{\tilde{L}}
\newcommand{\cmpmatrix}[1]{\begin{pmatrix}#1\end{pmatrix}}%
\newcommand{\cms}{ss\xspace}
\newcommand{\cmsa}{s.\,a.\xspace}%
\newcommand{\cmsmallpmatrix}[1]{\begin{smallpmatrix}#1\end{smallpmatrix}}%
\newcommand{\cmtimes}{\times}%
\newcommand{\cmvb}[1]{vollst"andig beschr"ankt#1\xspace}%
\newcommand{\cmvk}[1]{vollst"andig kontraktiv#1\xspace}%
\newcommand{\cmvp}[1]{vollst"andig positiv#1\xspace}%
\newcommand{\cmzB}{z.\,B.\xspace}%
\newcommand{\csalgebra}{$C^*$"=Algebra\xspace}%
\newcommand{\csalgebren}{$C^*$"=Algebren\xspace}%
\newcommand{\cstern}{$C^*$-}%
\newcommand{\dass}{da\cms}%
\newcommand{\diag}{\operatorname{diag}}
\newcommand{\dualr}[1]{{#1}^*}
\newcommand{\einlemph}[1]{\textbf{#1}}
\newcommand{\erzl}[1]{[#1]}
\newcommand{\erzur}[2]{{#1}_{\restring{#2}}}
\newcommand{\ftskalpr}[2]{\langle #1,#2\rangle}
\newcommand{\geltransN}{\mathscr{G}}
\newcommand{\Graph}{\mathcal{G}}
\newcommand{\hcsmodul}{Hilbert"=$C^*$"=Modul\xspace}%
\newcommand{\hcsmoduln}{Hilbert"=$C^*$"=Moduln\xspace}%
\renewcommand{\Id}{\operatorname{id}}
\renewcommand{\Im}{\operatorname{Im}}
\newcommand{\indHG}{{\mbbR_{\geq 0}}}
\newcommand{\indtHG}{{t \in \mbbR_{\geq 0}}}
\newcommand{\indtGr}{{t \in \mbbR}}
\newcommand{\injh}[1]{I(#1)}
\newcommand{\isoIMl}{\Theta}
\newcommand{\jbstripel}{$JB^*$-Tripel}%
\newcommand{\kptOp}{\mcalK}
\newcommand{\kptOphm}{\mathbb{K}}
\newcommand{\LinkAlg}[1]{\mcalL(#1)}
\newcommand{\matnull}[1]{\eta_{#1}}%
\newcommand{\matoplus}{\oplus}%
\newcommand{\matoplusc}[2]{#1 \oplus #2}%
\newcommand{\mCB}{\operatorname{CB}}%
\newcommand{\Mensch}[1]{\textsc{#1}}%
\newcommand{\mfrakA}{\mathfrak{A}}%
\newcommand{\mfrakB}{\mathfrak{B}}%
\newcommand{\mfrakC}{\mathfrak{C}}%
\newcommand{\mfrakD}{\mathfrak{D}}%
\newcommand{\mfrakT}{\mathfrak{T}}%
\newcommand{\mHerm}[1]{{H(#1)}}
\newcommand{\mIMl}{\mcalI{}M_\ell}
\newcommand{\mIMls}{\mcalI{}M_\ell^*}
\newcommand{\mLinStet}{L}
\newcommand{\mPosEl}[1]{{#1}_{\geq 0}}
\newcommand{\mre}{\mathrm{e}}
\newcommand{\mri}{\mathrm{i}}
\newcommand{\multis}{\cdot}
\newcommand{\Multcs}[1]{\mathcal{M}(#1)}
\newcommand{\Multlco}[1]{\mathscr{M}_l^{C_0}(#1)}%
\newcommand{\Multlcs}[1]{\mathcal{M}_l(#1)}%
\newcommand{\Multlor}[1]{{\mathscr{M}_l(#1)}}
\newcommand{\Multrcs}[1]{\mathcal{M}_r(#1)}%
\newcommand{\Multror}[1]{\mathscr{M}_r(#1)}
\newcommand{\Multwor}[1]{\mcalR(#1)}
\newcommand{\normbig}[1]{\bigl\lVert#1\bigr\rVert}
\newcommand{\normBig}[1]{\Bigl\lVert#1\Bigr\rVert}
\newcommand{\normBigg}[1]{\Biggl\lVert#1\Biggr\rVert}
\newcommand{\normcb}[1]{\norm{#1}_\text{cb}}
\newcommand{\normcblr}[1]{\normlr{#1}_\text{cb}}
\newcommand{\normi}[2]{{\lVert#1\rVert}_{#2}}
\newcommand{\normlr}[1]{\left\lVert#1\right\rVert}
\newcommand{\normMl}[2]{\norm{#1}_{\Multlor{#2}}}
\newcommand{\oplushm}{\oplus_2}
\newcommand{\pr}{\operatorname{pr}}
\newcommand{\refb}[2]{#1 oder #2}%
\renewcommand{\Re}{\operatorname{Re}}
\newcommand{\skalpr}[2]{\langle #1,#2\rangle}
\newcommand{\skalprb}[1]{\skalpr{#1}{#1}}
\newcommand{\skalpri}[3]{{\langle #1,#2\rangle}_{#3}}
\newcommand{\skalpril}[3]{\rule{0pt}{1pt}_{#3}\langle #1,#2\rangle}
\newcommand{\skalprd}[2]{\langle\!\langle #1,#2\rangle\!\rangle}
\newcommand{\skiptext}{\bigskip}
\newcommand{\spa}{\,\,}
\newcommand{\spradN}{r}%
\newcommand{\sprad}[1]{\spradN(#1)}
\newcommand{\sterns}{$\,\!^*$"~}
\newcommand{\sterndarst}{\sterns{}Dar\-stel\-lung\xspace}%
\newcommand{\sectiona}[1]{\section{\,\,\,#1}}%
\newcommand{\ternh}[1]{\mathcal{T}(#1)}
\newcommand{\ternsh}[1]{\mathcal{T}^*(#1)}
\newcommand{\troerzi}[2]{\mathop{TRO}_{#2}(#1)}
\newcommand{\UMLan}{\operatorname{UM}}
\newcommand{\vrerz}[1]{\operatorname{lin}\{#1\}}  

\newcommand{\dtext}[1]{{\footnotesize[#1]}}%
\newcommand{\dvorauss}[1]{}%

\ifthenelse{\boolean{insertdComments}}
{\newcommand{\dremark}[1]{\dtext{#1}}
 \newcommand{\dremarkm}[1]{[#1]}
 \newcommand{\dmarginpar}[1]{\marginpar{\footnotesize #1}}%
 \newcommand{\dmarginfz}{\marginpar{\textbf{?}}}
 \newcommand{\dused}[1]{\dtext{$\mbbD\,$#1}}
 \newcommand{\dremww}[1]{\dtext{#1}}
}
{\newcommand{\dremark}[1]{}%
 \newcommand{\dremarkm}[1]{}%
 \newcommand{\dmarginpar}[1]{}%
 \newcommand{\dmarginfz}{}%
 \newcommand{\dused}[1]{}%
 \newcommand{\dremww}[1]{}%
}
\newcommand{\dremc}[1]{}

\newcommand{\dliter}[1]{\dremark{#1}}
\newcommand{\dfrage}[1]{\dmarginpar{Frage}\dremark{\textbf{Frage:} #1}}%

\ifthenelse{\boolean{insertdcite}}
{\newcommand{\dcitem}[1]{\dtext{#1}} }   
{\newcommand{\dcitem}[1]{} }

\newcommand{\newlinef}{\newline}
\renewcommand{\qedsymbol}{$\qedsymbolF$}  

\makeatletter \@addtoreset {equation}{section}
\makeatother
\renewcommand {\theequation}{\arabic{chapter}.\arabic{section}.\arabic{equation}}

\ifthenelse{\boolean{showVersion}}
{\rfoot{{\scriptsize\textVersion}}}
{}

\newenvironment{enumaequiv}
{\begin{enumerate}[(a)]}%
{\end{enumerate}}
\newenvironment{enumaufz}
{\begin{enumerate}[(i)]}%
{\end{enumerate}}
\newenvironment{enumaufzB}
{\begin{enumerate}[(I)]}%
{\end{enumerate}}
\newenvironment{enumbsp}
{\begin{enumerate}[(i)]}%
{\end{enumerate}}
\newenvironment{smallpmatrix}%
{\left( \begin{smallmatrix}}%
{\end{smallmatrix} \right)}

\ifthenelse{\boolean{schreibidx}}%
{\newindex{B}{bdx}{bnd}{Stichwortverzeichnis}%
 \newindex{S}{sdx}{snd}{Symbolverzeichnis}%
}
{}

\newtheoremstyle{plainM}{3ex plus 3pt minus 5pt}{3ex plus 3pt minus 5pt}{\itshape}{}{\bf}{.}{.5em}{}
\newtheoremstyle{definitionM}{3ex plus 3pt minus 5pt}{3ex plus 3pt minus 5pt}{}{}{\bf}{.}{.5em}{}

\theoremstyle{plainM}
\newtheorem{satz}{Satz}[chapter]

\newtheorem{anmerkung}[satz]{Anmerkung}
\newtheorem{beispiel}[satz]{Beispiel}

\newtheorem{bemerkung}[satz]{Proposition}

\newtheorem{erinnerung}[satz]{Erinnerung}
\newtheorem{folgerung}[satz]{Folgerung}
\newtheorem{lemma}[satz]{Lemma}
\newtheorem{korollar}[satz]{Korollar}
\newtheorem{satzDefinition}[satz]{Satz/Definition}

\theoremstyle{definitionM}

  \newtheorem{defBemerkung}[satz]{Definitions-Proposition}
  \newtheorem{defSatz}[satz]{Definitions-Satz}
  \newtheorem{definitn}[satz]{Definition}
  \newtheorem{generalvor}[satz]{Generalvoraussetzung}

  \newtheorem{notation}[satz]{Notation}
  \newtheorem{warnung}[satz]{Warnung}

\begin{document}
\input{MathCommands.tex}
\renewcommand{\schlietmceMathbb}[1]{\ensuremath{\mathbb{#1}}}
\renewcommand{\setfdg}{\,;\;}

\swapnumbers

\makeatletter
\renewcommand\cleardoublepage{
\clearpage \if@twoside \ifodd \c@page \else \hbox {}\thispagestyle{empty}\newpage
\if@twocolumn \hbox {}\thispagestyle{empty}\newpage  \fi \fi \fi
}
\makeatother

\frontmatter

\ifthenelse{\boolean{finalVersion}}{

\thispagestyle{empty}
\begin{center}

  \textbf{Mathematik}\\
  \vspace*{\stretch{3}}
  \textbf{\huge
     Unbeschr"ankte Multiplikatoren auf Operatorr"aumen\\}

  \vspace*{\stretch{6}}

    Inaugural-Dissertation\\
    zur Erlangung des Doktorgrades\\
    der Naturwissenschaften im Fachbereich\\
    Mathematik und Informatik\\
    der Mathematisch-Naturwissenschaftlichen Fakultät\\
    der Westfälischen Wilhelms-Universität Münster

  \vspace*{\stretch{4}}%
    vorgelegt von\\
    \Mensch{Hendrik Schlieter}\\
    aus Kiel\\
    -- 2009 --
\end{center}

\clearpage
\thispagestyle{empty}

 \vspace*{16cm}
 \begin{tabbing}
\normalsize\textrm{{Dekan:}}\qquad\qquad\qquad\qquad\qquad\qquad \=
 \normalsize\textrm{{\Mensch{Prof.\,Dr.\,Dr.\,h.c.\,Joachim\,Cuntz}}}\\
\normalsize\textrm{{Erster Gutachter:}}\>
\normalsize\textrm{{\Mensch{apl.\,Prof.\,Dr.\,Wend\,Werner}}}\\
\normalsize\textrm{{Zweiter Gutachter:}}
\>\normalsize\textrm{{\Mensch{Prof.\,Dr.\,Siegfried\,Echterhoff}}}\\
\\
\normalsize\textrm{{Tag der mündlichen Prüfung:}}
\>\normalsize\textrm{{8. Juli 2009}}\\

\\\normalsize\textrm{{Tag der Promotion:}}
\>\normalsize\textrm{{8. Juli 2009}}
\end{tabbing}

\cleardoublepage
\thispagestyle{empty}

\pagestyle {fancyplain}
\pagenumbering{roman}
\setcounter{page}{3}  

}{}

\tableofcontents


\chapter{Einleitung}

\pagenumbering{arabic}
\setcounter{page}{1}

\markboth{EINLEITUNG}{}
\dremark{Kernpunkte herausstellen}%
\dremark{Klaus: Wozu wird dies gemacht? Ideen schildern.}%

Unbeschr"ankte Multiplikatoren spielen eine wichtige Rolle in der Mathematik,
beispielsweise im Spektralsatz f"ur unbeschr"ankte, selbstadjungierte Operatoren
oder in Verbindung mit Differentialoperatoren.\dremark{letzteres funktioniert nur gut auf $L^2$}
Auch bei Berechnungen in der Quantenphysik treten oft unbeschränkte
Größen auf, die durch unbeschränkte Operatoren dargestellt werden.

Unter unbeschr"ankten Multiplikatoren auf \hcsmoduln versteht man
"ublicherweise die sogenannten regul"aren Operatoren.
Die Definition eines solchen Operators geht auf A. Connes zur"uck,
der in \cite{ConnesThomIsom}\dremark{Definition 7} selbstadjungierte regul"are Operatoren
auf \csalgebren einf"uhrt und
diese Operatoren selbstadjungierte unbeschr"ankte Multiplikatoren nennt.

S. L. Woronowicz definiert in \cite{Woronowicz91} allgemeine regul"are Operatoren
auf \csalgebren.
Diese Operatoren kann man auch auf \hcsmoduln betrachten,
wie beispielsweise von E. C. Lance (\cite{LanceHmod}) oder
J. Kustermans (\cite{Kustermans97FunctionalCalcRegOp}) durchgef"uhrt.\\

Die regul"aren Operatoren sind wichtig in der Theorie der nicht"=kompakten
Quantengruppen (siehe zum Beispiel \cite{Woronowicz91}, \cite{WoronowiczNapi92OpThInCsAlg}
und \cite{Woronowicz95CsalgGenByUnbEl}).
So werden verschiedene nicht-kompakte Quantengruppen von
regul"aren Operatoren der zugrundeliegenden \csalgebra erzeugt,
und auch die Komultiplikation l"a\cms{}t sich mit Hilfe dieser Erzeuger beschreiben.
Au"serdem kann man in Kasparovs bivarianter $KK$-Theorie $KK(A,B)$
mit Hilfe unbeschr"ankter Kasparov-Moduln definieren und benutzt daf"ur
regul"are Operatoren (siehe \cite{BaajJulg83TheorieBivariante}).\\

Einen weiteren Zugang zu unbeschr"ankten Multiplikatoren auf \csalgebren
studieren A. J. Lazar und D. C. Taylor in \cite{LazarTaylorMultPedId}.
Sie betrachten  die Multiplikatoralgebra (Algebra der Bizentralisatoren)
des Pedersen-Ideals.
Wie C.~Webster in \cite{Webster04UnboundedOp}\dremark{Theorem 3.1} feststellt,
umfassen die regul"aren Operatoren die Multiplikatoren des Pedersen-Ideals.\\

Das Ziel dieser Arbeit ist es, die Theorie der unbeschr"ankten Multiplikatoren
von \hcsmoduln auf die wesentlich gr"o"sere Klasse der Operatorr"aume
fortzusetzen.
Bekanntlich tr"agt jeder \hcsmodul eine kanonische Operatorraumstruktur.
Da ein Operatorraum weniger Struktur als ein \hcsmodul besitzt,
ist nicht ohne weiteres ersichtlich, wie man die Definition eines
regul"aren Operators auf Operatorr"aume "ubertragen kann.

Es gelingt uns, unbeschr"ankte Multiplikatoren
auf Operatorr"aumen so zu definieren,
\dass diese Multiplikatoren auf \hcsmoduln, versehen mit der kanonischen Operatorraumstruktur,
mit den regul"aren Operatoren "ubereinstimmen,
also eine Verallgemeinerung der regul"aren Operatoren auf Operatorr"aume sind.
Au"serdem enthalten die unbeschr"ankten Multiplikatoren auf einem Operatorraum
die beschr"ankten, von links adjungierbaren Multiplikatoren.
Des weiteren werden verschiedene Eigenschaften unbeschr"ankter Multiplikatoren
untersucht und in\-trin\-sische Charakterisierungen sowie Anwendungen angegeben.\\

Im folgenden schildern wir die Struktur der Arbeit und
den Inhalt der einzelnen Kapitel.\\

Im ersten Abschnitt von \einlemph{Kapitel~1} werden
Grundlagen der Theorie der regul"aren Operatoren auf \hcsmoduln wiederholt.
Neu ist hier eine Variante des Satzes von Stone f"ur \hcsmoduln,
die im zweiten Abschnitt bewiesen wird.
Es stellt sich heraus, \dass man die schiefadjungierten regul"aren Operatoren
als Erzeuger bestimmter $C_0$-Gruppen charakterisieren kann.
Zusammenh"ange zwischen regul"aren Operatoren auf $E$ und auf $E \oplus E$
werden im darauf\/folgenden Abschnitt untersucht.
Zusammen mit dem hier bewiesenen Satz von Stone wird dies im dritten Kapitel verwendet,
um zu zeigen, \dass die unbeschr"ankten Multiplikatoren auf einem
\hcsmodul mit den regul"aren Operatoren "ubereinstimmen.
Im letzten Abschnitt werden die Multiplikatoren des Pedersen-Ideals behandelt.\\

In \einlemph{Kapitel~2} wird zun"achst an Grundlagen der Theorie der Operatorr"aume erinnert.
So werden mehrere Beispiele f"ur Operatorr"aume vorgestellt und
unitale Operatorsysteme behandelt.
Ferner werden wichtige Eigenschaften zu injektiven Operatorr"aumen und tern"aren
Ringen von Operatoren zusammengetragen.
Im vorletzten Abschnitt wird auf selbstadjungierte Operatorr"aume eingegangen.
Neu ist hier das Studium der injektiven H"ulle solcher Operatorr"aume im letzten Abschnitt.\\

\pagebreak
Das \einlemph{dritte Kapitel} ist der zentrale Teil dieser Arbeit.
Hier werden verschiedene Begriffe von unbeschr"ankten Multiplikatoren
auf einem Operatorraum jeweils mit Hilfe von $C_0$"=Halbgruppen
(bzw. $C_0$-Gruppen) definiert.
Es werden verschiedene Eigenschaften dieser Multiplikatoren untersucht
und in\-trin\-sische Charakterisierungen bewiesen.
Beispielsweise ist $A$ genau dann ein unbeschr"ankter Multiplikator,
wenn eine unbeschr"ankter Operator $B$ in $X$ derart existiert, \dass
$\cmsmallpmatrix{ 0 & \mri A & 0 & 0 \\ \mri B & 0 & 0 & 0 \\
                  0 & 0 & 0 & 0 \\ 0 & 0 & 0 & 0 }$
eine vollst"andig kontraktive $C_0$-Gruppe auf dem Spaltenoperatorraum $C_4(X)$ erzeugt.

Somit ist von Interesse, wann ein Operator Erzeuger einer vollst"andig kontraktiven
$C_0$-Gruppe ist.
Eine Charakterisierung hierf"ur wird im vierten Abschnitt gegeben,
in dem die S"atze von Hille-Yosida und von Lumer-Phillips
auf Operatorr"aume "ubertragen werden.
Ferner werden unbeschr"ankte Multiplikatoren von einem Operatorraum $X$
auf die tern"are H"ulle von $X$ und auf einen Hilbertraum "ubertragen.
Des weiteren wird eine Reihe von Anwendungen angegeben.
Zum Beispiel wird ein Resultat aus der St"orungstheorie regul"arer Operatoren
von S. Damaville (\cite{Damaville04RegulariteDesOp})
auf Operatorr"aume verallgemeinert.\\

In \einlemph{Anhang A} wird an beschr"ankte Multiplikatoren auf Operatorr"aumen
erinnert, genauer an Linksmultiplikatoren und von links adjungierbare Multiplikatoren.\\

Grundlagen der Theorie der $C_0$-Halbgruppen werden in \einlemph{Anhang B} behandelt.
\section*{Danksagung}

Ich danke allen, die zum Gelingen dieser Arbeit beigetragen haben.
Mein ganz besonderer Dank gilt apl. Prof. Dr. Wend Werner f"ur die Betreuung meiner Promotion.
Sein stetes Interesse und viele anregende Diskussionen haben mich angespornt.

Au"serdem m"ochte ich mich bei den Mitgliedern der Arbeitsgruppe
\emph{Funktionalanalysis,  Operatoralgebren und Nichtkommutative Geometrie} bedanken,
insbesondere bei Prof. Dr. Dr. h.\,c. Joachim Cuntz und Prof. Dr. Siegfried Echterhoff.
Des weiteren danke ich Prof. Dr. Reimar Wulkenhaar.
Ferner gilt mein Dank allen Kollegen aus dem SFB
\emph{Geometrische Strukturen in der Mathematik}.
Insbesondere m"ochte ich Ralf Kasprowitz, Erik M"uller und Dr. Walther Paravicini hervorheben,
des weiteren Klaus Loerke, Katharina Neumann, Dr. Thomas Timmermann und Moritz Weber.
Weiter danke ich Kyung-huy Moon, Alexander Ullmann und den Kollegen aus dem Blauen Pavillon,
insbesondere Dr. Steve Br"uske, Dr. Daniel Epping, Christian Kappen und Dr. Christian Wahle.
\dremark{Anja Wenning, Bj"orn Hille}%
Bei Dennis Bohle m"ochte ich mich f"ur viele anregende Diskussionen bedanken.
F"ur fachliche und finanzielle Unterst"utzung danke ich der Deutschen Forschungsgemeinschaft.
Ferner gilt mein Dank meinen Eltern für ihre vielfältige Unterstützung und allen
anderen, die mir im Verlauf der Promotion besonderen Rückhalt gegeben haben.

\chapter{Notation}

\markboth{NOTATION}{}

Es sei $\mbbN$ die Menge der nat"urlichen Zahlen ohne $0$ und $\mbbN_0 := \mbbN \cup \{0\}$.
\index[S]{N@$\mbbN$ ($:= \{1,2,3,\dots\}$)}%
\index[S]{N0@$\mbbN_0$ ($:= \{0,1,2,\dots\})$}%
F"ur alle $n \in \mbbN$ sei $\haken{n} := \{1, \dots, n\}$.
\index[S]{nh@$\haken{n}$ ($:= \{1, \dots, n\}$)}%
\dremark{F"ur eine beliebige komplexe Zahl $\lambda$ wird mit $\Re(\lambda)$
der Realteil und mit $\Im(\lambda)$ der Imagin"arteil von $\lambda$ bezeichnet.}%
Mit $\mathbbm{1}$ wird die Funktion auf $\mbbR$, die konstant $1$ ist, bezeichnet.
\index[S]{1@$\mathbbm{1}$}%

Sei $V$ ein $\mbbC$-Vektorraum und $W \subseteq V$.
Der komplexe Vektorraumaufspann von $W$ in $V$ wird mit $\cmlin W$ notiert.
\index[S]{lin@$\cmlin M$}%
Ist $V$ ein normierter Raum, so setze $\overline{\cmlin}\, W := \overline{\cmlin W}$.
\index[S]{linc@$\operatorname{\overline{\cmlin}} M$}%

Sei $\mfrakA$ eine Algebra und $G \subseteq \mfrakA$ eine Teilmenge.
Als $\alg(G)$ schreiben wir die von $G$ in $\mfrakA$ erzeugte Algebra.
\index[S]{alg@$\alg M$}
Das Einselement einer unitalen Algebra $\mfrakA$ wird mit $e_\mfrakA$ bezeichnet.
\index[S]{ekleinA@$e_\mfrakA$ (Einselement)}%
\skiptext

Seien $X$, $Y$, $Z$ normierte R"aume, und es sei eine Paarung $X \times Y \to Z$
gegeben, die wir als Abbildung $(x,y) \mapsto xy$ notieren.
Seien $M \subseteq X$ und $N \subseteq Y$.
Es sei definiert
\[ [MN] := \cmlin \{ xy \setfdg x \in M, y \in N \}  \quad\text{und}\quad
MN := \overline{[MN]}. \]
\index[S]{MN@$MN$}%
\index[S]{[MN]@$[MN]$}%
F"ur alle $x \in X$ und $r \in \mbbR_{\geq 0}$ sei mit
$\cmkug_X(x,r)$ (bzw. $\overline{\cmkug}_X(x,r)$) die offene (bzw. abgeschlossene) Kugel
in $X$ um $x$ mit Radius $r$ bezeichnet.
\index[S]{KXxr@$\cmkug_X(x,r)$ (offene Kugel)}
\index[S]{KXxra@$\overline{\cmkug}_X(x,r)$ (abgeschlossene Kugel)}
Die Menge der stetigen, linearen Abbildungen von $X$ nach $Y$ wird mit $L(X,Y)$ abgek"urzt,
die Menge der kompakten Operatoren von $X$ nach $Y$ mit $\kptOp(X,Y)$.
Setze $L(X) := L(X,X)$ und $\kptOp(X) := \kptOp(X,X)$.
\index[S]{LXY@$L(X,Y)$}%
\index[S]{KcalXY@$\kptOp(X,Y)$}%
\index[S]{LX@$L(X)$}%
\index[S]{KcalX@$\kptOp(X)$ (kompakte Operatoren)}%

F"ur eine Abbildung $A$ wird mit $D(A)$ der Definitionsbereich von $A$ notiert.
\index[S]{DA@$D(A)$ (Definitionsbereich)}%
\skiptext

Sei $\mfrakA$ eine \csalgebra.
Mit $\Multlcs{\mfrakA}$ wird die Menge der
Linksmultiplikatoren auf $\mfrakA$ bezeichnet und
mit $\mPosEl{\mfrakA}$ die Menge der positiven Elemente von $\mfrakA$.
\index[B]{Linksmultiplikator}%
\index[S]{Mcall@$\Multlcs{\mfrakA}$ (Linksmultiplikatoren)}%
\index[S]{A0@$\mPosEl{\mfrakA}$ (positive Elemente)}%

Unter einem \hcsmodul verstehen wir ein Rechts-\hcsmodul,
dessen Skalarprodukt linear in der zweiten Variablen ist.
Ebenso sei das Skalarprodukt eines Hilbertraumes
linear in der zweiten Variablen.

Seien $E$, $F$ \hcsmodul{}n "uber $\mfrakA$.
Die Menge der adjungierbaren Abbildungen von $E$ nach $F$ wird als $\Adjhm{E,F}$ notiert.
Setze $\Adjhm{E} := \Adjhm{E,E}$.
\index[S]{BEF@$\Adjhm{E,F}$}%
\index[S]{BE@$\Adjhm{E}$}%
F"ur alle $x \in F$ und $y \in E$ wird durch
\[ \theta_{x,y} : E \to F, z \mapsto x \skalpri{y}{z}{E}, \]
eine Abbildung in $\Adjhm{E,F}$ definiert.
Den abgeschlossenen, linearen Aufspann von
$\{ \theta_{x,y} \setfdg x \in F, y \in E \}$ in $\Adjhm{E,F}$ schreiben wir als $\kptOphm(E,F)$.
\index[S]{KEF@$\kptOphm(E,F)$}%
Setze $\kptOphm(E) := \kptOphm(E,E)$.
\index[S]{KE@$\kptOphm(E)$}%

Sei $\Omega$ ein topologischer Raum.
Mit $C(\Omega)$ wird die Menge der stetigen Funktionen von $\Omega$ nach $\mbbC$ bezeichnet,
mit $C_b(\Omega)$ die Teilmenge der beschr"ankten, stetigen Funktionen von $C(\Omega)$
und mit $C_c(\Omega)$ die Teilmenge der stetigen Funktionen
mit kompaktem Tr"ager von $C(\Omega)$.
Ist $\Omega$ ein lokalkompakter Hausdorffraum, so wird die Menge der
im Unendlichen verschwindenden, stetigen Funktionen von $\Omega$ nach $\mbbC$
als $C_0(\Omega)$ notiert.
\dremark{Vorschlag von Dennis: $L_g$ (Linksmultiplikation) auf einer Gruppen einf"uhren,
  \cmzB $L^G_g$.}%
\dremark{$\dualr{X}$ ist noch nicht definiert.}%

\mainmatter
\setcounter{page}{7}  

\chapter{Unbeschr"ankte Multiplikatoren auf \hcsmoduln}

\dremark{Insgesamt: Kernpunkte herausstellen.}%
In diesem Kapitel werden zun"achst Grundlagen der Theorie der regul"aren Operatoren auf
Hilbert-\cstern{}Moduln behandelt.
Anschlie"send wird eine Variante des Satzes von Stone f"ur \hcsmodul{}n bewiesen.
Au"serdem wird untersucht, wie regul"are Operatoren auf einem beliebigen \hcsmodul $E$
und auf $E \oplus E$ zusammenh"angen.
Zusammen mit der hier gezeigten Variante des Satzes von Stone wird in Kapitel 3 bewiesen,
\dass die unbeschr"ankten schiefadjungierten Multiplikatoren auf $E$
gleich den schiefadjungierten regul"aren Operatoren auf $E$ sind.
Des weiteren wird an unbeschr"ankte Multiplikatoren des Pedersen-Ideals erinnert.

\begin{generalvor}
Sei $\mfrakA$ eine \csalgebra.
\end{generalvor}

\section{Regul"are Operatoren}

Wir f"uhren in diesem Abschnitt den Begriff des regul"aren Operators auf
Hilbert-\cstern{}Moduln ein.
Regul"are Operatoren sind eine Verallgemeinerung dicht
definierter, abgeschlossener Operatoren auf Hilbertr"aumen.
Es ist eine wesentliche Eigenschaft eines auf einem Hilbertraum $H$
dicht definierten, abgeschlossenen Operators $T$,
\dass die Adjungierte $T^*$ dicht definiert ist und \dass $1+T^*T$ das Inverse
eines beschr"ankten Operators auf $H$ ist,
wobei mit~$1$ die Identit"atsabbildung auf $H$ bezeichnet wird.
Diese Eigenschaft besitzen dicht definierte, abgeschlossene
Operatoren auf \hcsmoduln \iallg{} nicht,
da nicht einmal jeder beschr"ankte Operator adjungierbar ist
(\cite{LanceHmod}, S. 8).\dmarginpar{Bsp. eintippen}
Daher verlangt man hier zus"atzlich, \dass $T^*$ dicht definiert ist und
\dass $1+T^*T$ dichtes Bild hat,
um eine reichhaltige\dremark{evtl. besser: brachbare, engl. viable} Theorie zu erhalten.

\begin{defBemerkung}[\cite{LanceHmod}, S. 95]
Seien $E,F$ Hilbert-\cstern{}Mo\-duln "uber $\mfrakA$,
sei $T : D(T) \subseteq E \to F$ dicht definiert und $\mfrakA$-linear.
Durch
\[ D(T^*) := \bigl\{ y \in F \setfdg
   \exists z_y \in E\,\, \forall x \in D(T): \skalpr{Tx}{y} = \skalpr{x}{z_y} \bigr\} \]
wird ein Untermodul von $F$ definiert.
F"ur alle $y \in D(T^*)$ ist das Element $z_y$ in der obigen Formel eindeutig bestimmt und
wird mit $T^*y$ bezeichnet.
Somit wird eine $\mfrakA$-lineare, abgeschlossene Abbildung $T^* : D(T^*) \subseteq F \to E$,
\index[S]{Ts@$T^*$}%
genannt die \defemphi{Adjungierte} von $T$, definiert, die folgendes erf"ullt:
\dremark{Idee: Definiere $T^*$ wie in der Hilbertraumsituation.}%
\[ \skalpr{x}{T^*y} = \skalpr{Tx}{y} \]
f"ur alle $x \in D(T)$ und $y \in D(T^*)$.
\dliter{$T^*$ abgeschlossen: \cite{LanceHmod}, S. 95 unten}%
\end{defBemerkung}

\begin{definitn}\label{defRegAbb}
Seien $E,F$ Hilbert-\cstern{}Moduln "uber $\mfrakA$.
\begin{enumaufz}
\item Ein \defemphi{regulärer Operator} von $E$ nach $F$ ist eine dicht definierte, abgeschlossene,
$\mfrakA$-lineare Abbildung $T : D(T) \subseteq E \rightarrow F$ so,
\dass $T^*$ dicht definiert ist und $1+T^*T$ dichtes Bild hat.
\dremark{Im Gegensatz zum Hilbertraumfall mu\cms{} man fordern, \dass $T^*$ dicht definiert ist,
  denn sogar beschr"ankte Operatoren m"ussen nicht adjungierbar sein.
  Beispiel: \cite{Wegge-OlsenKTheory}, Example 15.E}%
\item Die Menge der regul"aren Operatoren von $E$ nach $F$ wird mit $\Multwor{E,F}$ bezeichnet.
  Setze $\Multwor{E} := \Multwor{E,E}$.
\index[S]{REF@$\Multwor{E,F}$ (reguläre Operatoren)}%
\index[S]{RE@$\Multwor{E}$ (reguläre Operatoren)}%
\dremark{Weitere Literatur: \cite{Kustermans97FunctionalCalcRegOp}
  (hier werden viele Resultate genannt), \cite{Hilsum89FonctorialiteEnKth}
  (einige interessante Aussagen), \cite{Damaville04RegulariteDesOp}}%
\dremark{Motivation f"ur das Betrachten von regul"aren Operatoren:
  \cite{Kustermans97FunctionalCalcRegOp}, S. 1/2}%
\dliter{\cite{LanceHmod}, S. 96}%
\end{enumaufz}
\end{definitn}

F"ur jeden regul"aren Operator $T : D(T) \subseteq E \to F$ gilt:
$\Bild(1+T^*T) = E$ (\cite{LanceHmod}, siehe den Beweis von Theorem 9.3).
\skiptext

Im Hilbertraum-Fall hat jeder dicht definierte, abgeschlossene Operator $T$
eine dicht definierte Adjungierte $T^*$.
In einem Hilbert-\cstern{}Modul $E$ mu\cms{} man dies hingegen
als Voraussetzung an den Operator stellen.
\dremark{da nicht einmal jeder beschr"ankte Operator adjungierbar ist
(\cite{LanceHmod}, S. 8).}%
Au"serdem gilt f"ur einen dicht definierten, abgeschlossenen Operator $T$ auf $E$, dessen Adjungierte
dicht definiert ist, \iallg nicht, \dass $1+T^*T$ dichtes Bild hat
(siehe beispielsweise \cite{LanceHmod}, S. 102--104).
Daher wird auch dies in die Definition eines regul"aren Operators aufgenommen.

Da im Hilbertraum f"ur alle dicht definierten, abgeschlossenen Operatoren $T$
stets $1+T^*T$ dichtes Bild hat, sind regul"are Operatoren eine Verallgemeinerung
solcher Operatoren auf Hilbert-\cstern{}Moduln.\dremark{\ref{imHRAequivRegulaer}}
Beispielsweise sind gewisse Differentialoperatoren
dicht definierte, abgeschlossene Operatoren auf einem Hilbertraum
(siehe zum Beispiel \cite{DunfordSchwartzLinOp2}\dremark{S. 1278/1279 (Anfang Ch. XIII),
  Def. XIII.2.7.8 (S. 1291), Remark S. 1294 ($T_1(\tau)$ ist dicht def. und abgeschlossen)}).
\dremark{Sei $T : D(T) \subseteq E \to F$ dicht definiert, abgeschlossen,
  $\mfrakA$-linear mit $T^*$ dicht definiert.
  Dann hat \iallg{} $1+T^*T$ kein dichtes Bild, siehe \cite{LanceHmod}, S. 102--104.}%

Die regul"aren Operatoren treten in der Theorie der nicht"=kompakten
Quantengruppen auf (siehe zum Beispiel \cite{Woronowicz91}, \cite{WoronowiczNapi92OpThInCsAlg}
und \cite{Woronowicz95CsalgGenByUnbEl}).
Eine Verallgemeinerung von regul"aren Operatoren, sogenannte semiregul"are Operatoren,
werden in \cite{Pal99RegOpOnHilbertCsMod} untersucht.
\dmarginpar{pr}\dremark{Sind die semiregul"aren Operatoren Multiplikatoren?}%

\dremark{In Einleitung:
Die regul"aren Operatoren spielen eine wichtige Rolle in der Theorie der nicht-kompakten
Quantengruppen (siehe zum Beispiel \cite{Woronowicz91}, \cite{WoronowiczNapi92OpThInCsAlg}
und \cite{Woronowicz95CsalgGenByUnbEl}).
So werden verschiedene nicht-kompakte Quantengruppen von
regul"aren Operatoren erzeugt, und die Komultiplikation l"a\cms{}t sich
mit Hilfe dieser Erzeuger beschreiben.
Au"serdem kann man in Kasparovs bivarianter $KK$-Theorie $KK(A,B)$
mit Hilfe unbeschr"ankter Kasparov-Moduln definieren und benutzt daf"ur
regul"are Operatoren (siehe \cite{BaajJulg83TheorieBivariante}).
\dremark{Siehe auch \cite{BlackadarKTheory2ed}, §17.11, S. 163ff.}%
\dremark{Weitere Anwendung: Die unbeschr"ankten Multiplikatoren des Pedersen-Ideals
  lassen fallen unter die unbeschr"ankten Multiplikatoren nach Woronowicz,
  siehe \cite{Webster04UnboundedOp}.}
}%

Die Definition eines regul"aren Operators geht auf A. Connes zur"uck,
der diese Operatoren selbstadjungierte unbeschr"ankte Multiplikatoren
nennt und wie folgt auf \csalgebren{} definiert:

\begin{definitn}[Vgl. \cite{ConnesThomIsom}, Definition 7]\label{defMultConnes}
Ein \defemph{selbstadjungierter unbeschr"ankter Multiplikator} auf $\mfrakA$
ist ein dicht definierter, abgeschlossener Operator $T : D(T) \subseteq \mfrakA \to \mfrakA$
mit den Eigenschaften:
\begin{enumaufz}
\item $\exists \lambda \in \mbbR_{>0}: -\lambda \mri, \lambda \mri \in \rho(T)$,
\item $y^* T(x) = T(y)^* x$ f"ur alle $x,y \in D(T)$.
\end{enumaufz}
\end{definitn}

Wie im Hilbertraum definiert man:

\begin{definitn}
Sei $E$ ein Hilbert-\cstern{}Modul "uber $\mfrakA$.
Eine dicht definierte, $\mfrakA$-lineare Abbildung $T : D(T) \subseteq E \to E$ hei"st
\begin{enumaufz}
\item \defemphi{symmetrisch}, falls $T^* \supseteq T$ gilt,
\item \defemphi{selbstadjungiert}, falls $T^*=T$ ist, und
\item \defemphi{schiefadjungiert}, falls $T^*=-T$ gilt.
\dliter{\cite{LanceHmod}, S. 101}%
\end{enumaufz}
\end{definitn}

Auf jeder \csalgebra $\mfrakB$ wird durch
\[ \skalpri{x}{y}{\mfrakB} := x^* y \]
f"ur alle $x,y \in \mfrakB$ ein inneres Produkt derart definiert,
\dass $\mfrakB$ zugleich die Struktur eines \hcsmodul{}s tr"agt.
\skiptext

Die selbstadjungierten unbeschr"ankten Multiplikatoren von Connes
fallen f"ur $\lambda = 1$ in Definition~\ref{defMultConnes}.(i)
unter den Begriff des regul"aren Operators:

\begin{bemerkung}[\cite{Baaj80Multiplicateurs}, S. 1]\label{charMultConnes}
Sei $T : D(T) \subseteq \mfrakA \to \mfrakA$ linear, dicht definiert und abgeschlossen.
Dann sind die folgenden Aussagen "aquivalent:
\dremark{M"u\cms{}te auch allgemeiner f"ur Hilbert-\cstern{}Modul gelten.}%
\begin{enumaequiv}
\item $T$ erf"ullt (i) und (ii) aus Definition~\ref{defMultConnes} mit $\lambda = 1$.
\item
\begin{enumaufz}
\item $\Bild(1+T^2)$ ist dicht in $\mfrakA$,
\item $T=T^*$.
\end{enumaufz}
\end{enumaequiv}
\end{bemerkung}

\dremark{Multiplikatoren nach Connes: von $\lambda > 0$ zu $\lambda=1$:

Sei $\mfrakA$ eine \csalgebra, $T$ ein selbstadjungierter unbeschr"ankter
Multiplikator von $\mfrakA$ nach Connes bzgl. $\lambda$.
Dann ist $\hat{T} := \frac{1}{\lambda}T$ ein selbstadjungierter unbeschr"ankter
Multiplikator nach Connes mit $D(T\hat{T}) = D(T)$ und
$\rho(\hat{T}) \overset{\dremarkm{\ref{rhoUnterSkalar}}}{=} \frac{1}{\lambda}\rho(T)$.
\dremark{Z. 13.3.'08/2}}%

Da bei S. Baaj kein Beweis f"ur die obige Proposition aufgef"uhrt ist,
geben wir einen Beweis an.
Dieser verwendet die folgende Proposition, die man analog zum Hilbertraum-Fall beweist:

\dremww{
\begin{bemerkung}\label{KernTPlusIEqBild}
Sei $E$ ein $\mfrakA$-Hilbert-\cstern{}Modul,
$T : D(T) \subseteq E \to E$ linear und dicht definiert.
Dann gilt: $\Kern(T^*\pm \mri) = \Bild(T\mp \mri)^\bot$.
\end{bemerkung}

\begin{proof}
\dremark{Analog zum Beweis von \cite{WernerFunkana3}, Lemma VII.2.7}%
\dremark{Beachte: $(T + \mri)^* = T^* - \mri$.}%
\bewitemph{\glqq$\subseteq$\grqq:}
Sei $y \in \Kern(T^* \pm \mri)$.
Es folgt:
\[ 0
=  \skalpr{z}{(T^* \pm \mri)y}
=  \skalpr{z}{(T \mp \mri)^* y}
=  \skalpr{(T \mp \mri)z}{y} \]
f"ur alle $z \in D(T)$,
also $y \in \Bild(T \mp \mri)^\perp$.

\bewitemph{\glqq$\supseteq$\grqq:}
Sei $y \in \Bild(T \mp \mri)^\perp$.
Dann gilt
$0
=  \skalpr{(T \mp \mri)z}{y}
=  \skalpr{Tz}{y} \pm \skalpr{z}{\mri y}$
f"ur alle $z \in D(T)$,
also $y \in D(T^*)$ und
$0 = \skalpr{z}{(T \mp \mri)^*y} = \skalpr{z}{(T^* \pm \mri)y}$
f"ur alle $z \in D(T)$.
Da $D(T)$ dicht in $E$ liegt, erh"alt man: $y \in \Kern(T^* \pm \mri)$.
\dremark{Z. 13.3.'08/1}%
\end{proof}
}%

\begin{bemerkung}\label{BildTPlusIEqEFolgtTsa}
Sei $E$ ein Hilbert-\cstern{}Modul und
$T : D(T) \subseteq E \to E$ linear, dicht definiert und symmetrisch.
Ist $\Bild(T-\mri) = E = \Bild(T+\mri)$, so ist $T$ selbstadjungiert.
\dremark{Im Hilbertraum gilt auch \glqq$\Rightarrow$\grqq{} sowie eine weitere "Aquivalenz.}
\end{bemerkung}

\dremww{
\begin{proof}
\dremark{Analog zum Beweis von \cite{WernerFunkana3}, Satz VII.2.8}%
Wegen $T \subseteq T^*$ ist nur $D(T^*) \subseteq D(T)$ zu zeigen.
Sei $y \in D(T^*)$.
Dann findet man ein $x \in D(T)$ mit $(T^* - \mri)y = (T - \mri)x$.
Wegen $T \subseteq T^*$ ist $(T^* - \mri)y = (T^* - \mri)x$.
Nach \ref{KernTPlusIEqBild} ist\dremark{$(*)$: \cite{Wegge-OlsenKTheory}, Lemma 15.3.4}
$\Kern(T^* - \mri) = \Bild(T + \mri)^\perp = E^\perp \overset{\dremarkm{(*)}}{=} \{0\}$.
Somit folgt $y = x \in D(T)$.
\dremark{Z. 13.3.'08/1}%
\end{proof}
}%

\begin{proof}[Beweis von Proposition~\ref{charMultConnes}]
\bewitemphq{(a)$\Rightarrow$(b)}
Es gelte \bewitemph{(a)}.
Dann gilt: $\skalpri{y}{Tx}{\mfrakA} = y^* T(x) = T(y)^* x = \skalpri{Ty}{x}{\mfrakA}$
f"ur alle $x,y \in D(T)$.
Somit folgt $T \subseteq T^*$, und nach Proposition~\ref{BildTPlusIEqEFolgtTsa}
ist $T$ selbstadjungiert.

Da $-\mri + T$ und $\mri + T$ surjektiv auf $\mfrakA$ sind,
ist auch $1 + T^2 = (-\mri + T)(\mri + T)$ surjektiv und
hat insbesondere dichtes Bild.
\dremark{Sei $x \in \mfrakA$.
  Man findet $y \in D(T)$ mit $(-\mri + T)y = x$ und
  $z \in D(T)$ mit $(\mri + T)z = y$.
  Also gilt: $x = (\mri + T)y = (-\mri + T)(\mri + T)z = (1+T^2)z$.}%

\bewitemphq{(b)$\Rightarrow$(a)}
Es gelte \bewitemph{(b)}.
Da $T$ wegen $T = T^*$ regul"ar ist,
sind $T+\mri$ und $T-\mri$ nach \cite{LanceHmod}, Lemma 9.7,
injektiv und nach \cite{LanceHmod}, Lemma 9.8, surjektiv auf $\mfrakA$.
Da $T+\mri$ und $T-\mri$\dremark{nach \ref{AabgAMinusIdIstAbg}} abgeschlossen sind,
folgt mit \cite{WernerFunkana6}, Satz IV.4.4:
$-\mri,\mri \in \rho(T)$.
Weiter gilt:
$  y^*T(x)
=  \skalpri{y}{Tx}{\mfrakA}
\dremarkm{=  \skalpri{y}{T^*x}{\mfrakA}}
=  \skalpri{Ty}{x}{\mfrakA}
=  T(y)^*x$
f"ur alle $x,y \in D(T)$.
\dremark{Z. 10.3.'08/2}%
\end{proof}

\dremww{
\begin{bemerkung}\label{imHRAequivRegulaer}
Sei $H$ ein Hilbertraum, $S : D(S) \subseteq H \to H$ ein Operator.
Dann sind die folgenden Aussagen "aquivalent:
\begin{enumaequiv}
\item $S$ ist regul"ar.
\item $S$ ist linear, dicht definiert und abgeschlossen.
\end{enumaequiv}
\end{bemerkung}

\begin{proof}
\bewitemph{\glqq$\Rightarrow$\grqq} ist klar.

\bewitemph{\glqq$\Leftarrow$\grqq:}
Nach \cite{RudinFuncAna:73}, 13.12 und 13.13.(a), ist $S^*$ dicht definiert
und $\Id_H + S^*S$ surjektiv auf $H$.
Also ist $S$ regul"ar.
\end{proof}

\begin{definitn}
Ein regul"arer Operator $m : D(m) \subseteq A \to A$ auf dem $A$-Hilbert-\cstern{}Modul $A$
hei"st \defemph{unbeschränkter Multiplikator} von $A$.
Die Menge der unbeschr"ankten Multiplikatoren von $A$ wird mit $\UMLan(A)$ bezeichnet.
\dremark{Nach \cite{LanceHmod}, S. 120, tauchte die Definition f"ur einen unbeschr"ankten
  Multiplikator zuerst bei \cite{ConnesThomIsom} auf.}%
\dliter{\cite{LanceHmod}, S. 117}%
\end{definitn}
}%

Man kann regul"are Operatoren wie folgt charakterisieren:

\begin{satz}[\cite{Kustermans97FunctionalCalcRegOp}, Theorem 1.5]
Seien $E$, $F$ Hilbert-\cstern{}Moduln "uber $\mfrakA$.
Sei $T : D(T) \subseteq E \to F$ dicht definiert.
Dann ist $T$ genau dann regul"ar,
wenn ein $z \in \Adjhm{E, F}$ mit $\norm{z} \leq 1$ so existiert,
\dass $D(T) = (1-z^*z)^{1/2}E$ und $T((1-z^*z)^{1/2}x) = zx$ f"ur alle $x \in E$ gilt.
\dremark{Idee: $T$ wirkt als unbeschr"ankter Multiplikator.}%
\end{satz}

\dremark{Somit ist im Falle $E = F = \mfrakA$ der Operator $T$ genau dann regul"ar,
wenn $T$ zu $\mfrakA$ affiliiert ist.}%

\begin{definitn}\label{defzTransf}
Die Abbildung $z$ im obigen Satz ist eindeutig definiert, wird mit $z_T$ bezeichnet und
\defemph{$z$-Transformierte} von $T$ genannt.
\index[B]{zTransformierte@$z$-Transformierte}%
\index[S]{zT@$z_T$ ($z$-Transformierte)}%
\end{definitn}

\dremark{Sei $T \in \Multwor{\mfrakA}$.
Dann ist $z_T \in \Adjhm{\mfrakA} \cong \Multcs{\mfrakA}$,
also ein Multiplikator, und man kann $T$ als eine Art unbeschr"ankten Multiplikator auf\/fassen.}%

Im folgenden notieren wir einige Beispiele regul"arer Operatoren:

\begin{beispiel}\label{bspRegOp}\label{bspRegOpC0Omega}
\begin{enumaufz}
\item Seien $E$, $F$ \hcsmoduln "uber $\mfrakA$.
  Dann gilt:
\[ \Adjhm{E,F} = \{ T \in \Multwor{E,F} \setfdg D(T) = E \}. \]

\item Ist $\mfrakA$ unital, so gilt: $\Multwor{\mfrakA} \cong \mfrakA$.

\item Sei $\Omega$ ein lokalkompakter Hausdorffraum.
  Es gilt: $\Multwor{C_0(\Omega)} \cong C(\Omega)$.
  Genauer hat man:
  Sei
\[ M_f : D(M_f) \subseteq C_0(\Omega) \to C_0(\Omega), g \mapsto fg, \]
  f"ur alle $f \in C(\Omega)$.
  Dann ist die Abbildung
\[ \Phi : C(\Omega) \arrowbij \Multwor{C_0(\Omega)}, f \mapsto M_f, \]
  bijektiv auf $\Multwor{C_0(\Omega)}$
  mit $\Phi(f)^* = \Phi\left( \overline{f} \right)$ f"ur alle $f \in C(\Omega)$.
\dremark{I.\,a. gilt nicht: $\Phi(f+g) = \Phi(f) + \Phi(g)$.
  Gegenbeispiel:
  $\Omega := \mbbR$, $f := \Id_\mbbR$, $g := -f$.
  Dann $\Phi(f+g) = \Phi(0) = M_0 = 0$, also $D(M_0) = C_0(\mbbR)$.
  Weiter $\Phi(f) + \Phi(g) = M_f + M_{-f}$,
  $D(M_f) \neq C_0(\mbbR)$ \iallg.}%
\dliter{\cite{Woronowicz91}, S. 404; knapper: \cite{LanceHmod}, S. 117;
  \cite{VaesEtAl01LocCptQuantumGroups} Examples 2.1.20}%

\item Sei $H$ ein Hilbertraum.
  Dann gibt es eine bijektive Abbildung $\Psi$ von $\Multwor{\kptOp(H)}$
  auf die Menge der dicht definierten, abgeschlossenen Operatoren auf $H$
  mit der Eigenschaft: $\Psi(T^*) = \Psi(T)^*$ f"ur alle $T \in \Multwor{\kptOp(H)}$.
\dremark{I.\,a. gilt nicht: $\Psi(S+T) = \Psi(S) + \Psi(T)$.
  Gegenbeispiel wie in (iii):
  $D(\Psi(T + (-T))) = \kptOp(H)$, $D(\Psi(T)) \neq \kptOp(H)$ \iallg.}%
\dremark{Was kann man f"ur ein \hcsmodul $E$ und $\Multwor{\kptOphm(E)}$ sagen?}%
\dliter{\cite{LanceHmod}, S. 117}%
\dremark{\item Es gilt: $\Multcs{\mfrakA} \subseteq \Multwor{\mfrakA}$.
Genauer hat man: F"ur alle $a \in \mfrakA$ ist
die Abbildung $L_a : \mfrakA \to \mfrakA, x \mapsto ax$, aus $\Multwor{\mfrakA}$.
\dliter{Vgl. \cite{LanceHmod}, S. 117}%
\dremark{Richtig: $a \in \mfrakA$?}%
}%
\end{enumaufz}
\end{beispiel}

\begin{proof}
\bewitemph{(i):} \cite{Kustermans97FunctionalCalcRegOp}, Result 1.7.

\bewitemph{(ii):} \cite{LanceHmod}, S. 117.\dremark{$\Multcs{\mfrakA} \cong \mfrakA$}

\bewitemph{(iii):} \cite{Woronowicz91}, Example~2.

\bewitemph{(iv):} \cite{Woronowicz91}, Example 3.
\end{proof}

Wir erinnern an die folgende

\begin{definitn}\label{defBizentralisator}
Sei $\mfrakB$ eine Algebra.
\begin{enumaufz}
\item Unter einem \defemphi{Multiplikator} (bzw.\ \defemphi{Bizentralisator}) auf $\mfrakB$ versteht man
ein Paar $(S,T)$, wobei $S$ und $T$ Funktionen von $\mfrakB$ nach $\mfrakB$ sind, derart,
\dass $aS(b) = T(a)b$ f"ur alle $a,b \in \mfrakB$ gilt.
\dremark{$S$ und $T$ sind linear.}%
\dremark{$\mfrakB$ als Algebra: $\mfrakB$ kann das Pedersen-Ideal sein.}%
\dvorauss{$\mfrakB$ $K$-Algebra}%
\item Die Menge aller Multiplikatoren auf $\mfrakB$ wird mit $\cmGamma(\mfrakB)$ bezeichnet.
\index[S]{Mcal@$\cmGamma(\mfrakA)$ (Multiplikatoren)}%
\end{enumaufz}
\end{definitn}

\begin{bemerkung}[\cite{LazarTaylorMultPedId}, S. 5]
Sei $\mfrakB$ eine Algebra.
Dann ist $\cmGamma(\mfrakB)$, versehen mit der Multiplikation
$(S,T)\cdot(U,V) = (S \circ U, V \circ T)$, eine Algebra.
Falls $\mfrakB$ eine \sterns{}Algebra ist, so wird durch
$(S,T) \mapsto (T^*,S^*)$, wobei $S^*(a) := S(a^*)^*$ f"ur alle $a \in \mfrakB$ ist,
eine Involution auf $\cmGamma(\mfrakB)$ definiert.
\dvorauss{$\mfrakB$ $K$-Algebra}%
\end{bemerkung}

Nach Beispiel~\ref{bspRegOp}.(i) gilt
$\Multcs{\mfrakA} \cong \Adjhm{\mfrakA} \subseteq \Multwor{\mfrakA}$,
also kann man einen regul"aren Operator als einen unbeschr"ankten Multiplikator auf\/fassen.
\skiptext

\dremark{
Beweis.
\bewitemph{(i):}
Sei $m \in \Multcs{\mfrakA}$. Dann ist $m$ adjungierbar, also insbesondere linear und stetig.
Somit ist $m$ insbesondere abgeschlossen.
Nach \cite{RudinFuncAna:73}, 13.13, hat $1+m^*m$ dichtes Bild.
\dremark{Alternativer Beweis:
  Es ist $D(L_a) = \mfrakA$ ein dichtes Rechtsideal von $\mfrakA$,
  somit gilt offensichtlich: $L_a \in \Multwor{\mfrakA}$.
  Z. 6.2.'08/2}%
\bewitemph{(ii):} Da die Menge $G(\mfrakA)$ der invertierbaren Elemente von $\mfrakA$ offen und nicht-leer ist,
  hat $\mfrakA$ keine echten dichten Rechtsideale.
\dremark{\bewitemph{(iii):}
Setze $\mfrakA := C_0(\Omega)$.
Zeige: $\Phi : C(\Omega) \arrowbij \Multwor{\mfrakA}, f \mapsto M_f$ ist bijektiv
auf $\Multwor{\mfrakA}$ mit $\Phi(f)^* = \Phi(\overline{f})$, wobei
$M_f : D(M_f) \subseteq \mfrakA \to \mfrakA$ f"ur alle $f \in C(\Omega)$.

Zeige: $\Bild(\Phi) \subseteq \Multwor{\mfrakA}$.

Nach \cite{EngelNagelSemigroups}, Prop. I.4.2, ist $M_f$ dicht definiert und abgeschlossen.
Wegen $(M_f)^* = M_{\overline{f}}$ ist $(M_f)^*$ dicht definiert und $\Phi(f)^* = \Phi(\overline{f})$.
\dremark{Es gilt:
\[ D((M_f)^*)
= \{ a \in \mfrakA \setfdg
     \exists b_a \in \mfrakA\,\, \forall c \in D(M_f) :
     \underbrace{M_f(c)^* a = c^* b_a}_{\text{also } b_a = \overline{f} \cdot a} \}. \]
Man erh"alt $(M_f)^*(a) = b_a = \overline{f} \cdot a = M_{\overline{f}}(a)$
f"ur alle $a \in D((M_f)^*)$, also $(M_f)^* = M_{\overline{f}}$.}%
Zeige nun: $(\Id_{\mfrakA} + (M_f)^* M_f) \mfrakA$ ist dicht in $\mfrakA$.
Es gilt
\[ (\Id_\mfrakA + (M_f)^* M_f)g = g + f^*fg = (\mathbbm{1} + \abs{f}^2)g \]
f"ur alle $g \in \mfrakA$.
Setze $h := \mathbbm{1} + \abs{f}^2$.
Wegen $h \geq \mathbbm{1}$ hat $M_h$ nach \cite{EngelNagelSemigroups}, Prop. I.4.2, ein
beschr"anktes Inverses, und es gilt $M_h^{-1} = M_{1/h}$.
Es folgt: $\Bild(\Id_\mfrakA + (M_f)^* M_f) = \Bild(M_h) = D(M_h^{-1}) = D(M_{1/h}) = \mfrakA$.
\smallskip

Zeige: $\Phi$ ist injektiv.
xxx

Zeige: $\Phi$ ist surjektiv.

Sei dazu $T \in \Multwor{C_0(\Omega)}$.
Wegen $\Multcs{C_0(\Omega)} \cong C_b(\Omega)$ via $f \mapsto M_f$ findet man
ein $f \in C_b(\Omega)$ mit $z_T = M_f$.
Sei $a := f (\mathbbm{1} - \abs{f}^2)^{-1/2} \in C(\Omega)$.
Dann gilt: $M_a = T$.
\dremark{Z. 10.9.'07/1, hier ausf"uhrlicherer Beweis}}%

\bewitemph{(iv):} Setze $\mfrakA := \kptOp(H)$.
Man findet ein $y \in \Multcs{\mfrakA} \cong \mLinStet(H)$ mit $z_T = L_y$.
Nach \cite{LanceHmod}, Theorem 10.4, findet man ein $S \in \Multwor{H}$ mit $y = z_S$,
also ist $S$ dicht definiert und abgeschlossen.
\dremark{Z. 6.2.'08/1, dort ausf"uhrlicherer Beweis, basierend auf Beweis bei Lance.
  Eine Richtung fehlt noch.}%
\dmarginpar{zutun}

\dremark{Besser motivieren und formulieren.}\dmarginpar{zutun}%

Sei ein Hom"oomorphismus $z : \mbbC \arrowbij \cmkug_\mbbC(0,1)$ auf $\cmkug_\mbbC(0,1)$ gegeben.
Die Motivation hinter Definition~\ref{defzTransf} ist, \dass man,
indem man eine Funktion aus $C(\Omega)$ mit $z$ komponiert,
$C(\Omega)$ in $C_b(\Omega)$ einbetten kann (allerdings nicht \sterns{}homomorph).
\dremark{$\Omega$: lokalkompakter Hausdorffraum}%
\dremark{Beweis: Sei $\Phi : C(\Omega) \to C_b(\Omega), g \mapsto z \circ g$.
Es gilt: $\norm{\Phi(g)} = \sup_{\omega \in \Omega} \abs{ z \circ g(\omega) } \leq 1$.
Da $z$ invertierbar ist, ist $\Phi$ injektiv.}%
\dremark{\cite{Webster04UnboundedOp}}%

\begin{lemma}\label{zFunctionFromWor}
Sei $z : \mbbC \to \cmkug_\mbbC(0,1), \xi \mapsto \xi(1+\overline{\xi}\xi)^{-1/2}$.
Dann ist $g : \cmkug_\mbbC(0,1) \to \mbbC, \eta \mapsto \eta(1-\overline{\eta}\eta)^{-1/2}$,
die Umkehrabbildung von $z$.
\end{lemma}

\begin{proof}
Sei $\xi \in \mbbC$. Dann gilt
\[    \abs{z(\xi)} < 1
\iff  \frac{\abs{\xi}}{(1+\abs{\xi}^2)^{1/2}} < 1
\iff  \frac{\abs{\xi}^2}{1+\abs{\xi}^2} < 1
\iff  \abs{\xi}^2 < 1 + \abs{\xi}^2.  \]
Weiter gilt:
\[ g \circ z(\xi)
=  z(\xi)(1-\abs{z(\xi)}^2)^{-1/2}
=  z(\xi)(1-\frac{\abs{\xi}^2}{1+\abs{\xi}^2})^{-1/2}
=  z(\xi)(\frac{1}{1+\abs{\xi}^2})^{-1/2}
=  \xi. \]
Analog zeigt man: $z \circ g = \Id_\mbbD$.
\dremark{5.9.'07, S. 1}%
\end{proof}

S. Woronowicz hat in \cite{Woronowicz91} die Funktion $z$ aus Lemma~\ref{zFunctionFromWor} gew"ahlt.
\dremark{Bezug?}
}%

Die Abbildungen $T$ und $z_T$ h"angen wie folgt zusammen:

\begin{bemerkung}[Vgl. \cite{Kustermans97FunctionalCalcRegOp}, Proposition 7.21 und Proposition 7.22]
Seien $E$, $F$ Hilbert-\cstern{}Moduln "uber $\mfrakA$.
Sei $T \in \Multwor{E,F}$.
Dann gilt:
\begin{enumaufz}
\item $1+T^*T$ ist invertierbar mit Inversem $S \in \Adjhm{E}$ und $S \geq 0$.
\item $(1+T^*T)^{-1} = 1 - z_T^* z_T$.
\item $z_T = T\left(\left(1+T^*T\right)^{-1}\right)^{1/2}$ und $T = z_T\left((1-z_T^* z_T)^{1/2}\right)^{-1}$.
\dremark{$1+T^*T$ ist invertierbar mit Inversem in $\Adjhm{E}$,
  siehe \cite{Kustermans97FunctionalCalcRegOp}, Prop. 7.21 und Cor. 5.12.
  Nach \cite{LanceHmod}, Beweis von Lemma 9.1, Zeile 10, ist
  $(1+T^*T)^{-1}$ ein positives Element von $\Adjhm{E}$,
  man kann also die Wurzel ziehen.
  Wegen $(1+T^*T)^{-1} = 1 - z_T^* z_T := z$
  nach \cite{Kustermans97FunctionalCalcRegOp}, Prop. 7.21,
  kann man auch von $z$ die Wurzel ziehen.}%
\dremark{In \cite{Kustermans97FunctionalCalcRegOp}, Prop. 7.22, steht kurz
  $z_T = T(1+T^*T)^{-1/2}$ und $T = z_T(1-z_T^* z_T)^{-1/2}$,
  dies stimmt aber mit dem oben Definierten "uberein, siehe Text nach Prop. 7.22.}%
\end{enumaufz}
\end{bemerkung}

F"ur die Abbildung $T^{**}$ gilt:

\begin{bemerkung}[\cite{LanceHmod}, Corollary 9.4]\label{TssEqT}
Seien $E$, $F$ Hilbert-\cstern{}Moduln "uber $\mfrakA$.
Sei $T \in \Multwor{E,F}$.
Dann gilt: $T^{**} = T$.
\end{bemerkung}

Da nicht-ausgeartete \sterns{}Homomorphismen eine wichtige Rolle
im n"achsten Abschnitt spielen werden, erinnern wir an diesen Begriff:

\begin{definitn}
Sei $\mfrakB$ eine \csalgebra und
$E$ ein Hilbert-\cstern{}Modul "uber $\mfrakA$.
Ein \sterns{}Ho\-mo\-mor\-phis\-mus $\alpha : \mfrakB \to \Adjhm{E}$ wird
\defemph{nicht-ausgeartet} genannt, falls $\erzl{\alpha(\mfrakB) E}$ dicht in $E$ liegt.
\index[B]{nicht-ausgearteter \sterns{}Ho\-mo\-mor\-phis\-mus}%
\dremark{Motivation: \cite{LanceHmod}, S. 18 unten}%
\end{definitn}

F"ur einen beliebigen nicht-ausgearteten \sterns{}Homomorphismus $\alpha : \mfrakB \to \Adjhm{E}$
ergibt sich mit einer Verallgemeinerung des Faktorisierungssatzes von Cohen
(\cite{HewittRossHarmonicAna2}, Theorem 32.22):
$E = \{ \alpha(b)x \setfdg b \in \mfrakB, x \in E \}$.
\dremark{Beachte: $\Adjhm{E}$ ist eine \csalgebra,
  besitzt also eine beschr"ankte approximative Linkseins.}%
\skiptext

Bei der Untersuchung von Operatoren auf \hcsmoduln ist die folgende Topologie wichtig:

\begin{definitn}\label{defStriktStetig}
Seien $E$, $F$ \hcsmoduln "uber $\mfrakA$.
Die \defemphi{strikte Topologie} auf $\Adjhm{E,F}$ ist die Topologie,
die durch die Halbnormen
\[ p_x : \Adjhm{E,F} \to \mbbR, T \mapsto \norm{Tx}, \quad\text{und}\quad
   q_y : \Adjhm{E,F} \to \mbbR, T \mapsto \norm{T^* y}, \]
f"ur alle $x \in E$ und $y \in F$ definiert wird.
\dliter{\cite{LanceHmod}, S. 11 unten}%
\end{definitn}

Identifiziert man $\Multcs{\mfrakA}$ mit $\Adjhm{\mfrakA}$,
so erh"alt man, \dass die strikte Topologie auf $\Multcs{\mfrakA}$
durch die Halbnormen
\[ p_a : \Multcs{\mfrakA} \to \mbbR, x \mapsto \norm{ax}, \quad\text{und}\quad
   q_b : \Multcs{\mfrakA} \to \mbbR, x \mapsto \norm{bx}, \]
f"ur alle $a,b \in \mfrakA$ gegeben ist.
\dremark{Die strikte Topologie auf $\Multcs{\mfrakA}$ wird durch die Halbnormen
  $x \mapsto \norm{ax}$ und $x \mapsto \norm{xa}$ ($x \in \Multcs{\mfrakA}, a \in \mfrakA$)
  definiert (\dliter{\cite{Wegge-OlsenKTheory}, Def. 2.3.1}).
  F"ur alle $x \in \Multcs{\mfrakA}$, $a \in \mfrakA$ gilt
  $\norm{x^* a} = \norm{(a^* x)^*} = \norm{a^* x}$,
  daher erh"alt man durch obige Definition tats"achlich dieselben Topologien.}%
\skiptext

Man kann nicht-ausgeartete \sterns{}Homomorphismen wie folgt charakterisieren:

\begin{bemerkung}[Vgl. \cite{LanceHmod}, Proposition 2.5, und \cite{RaeburnWilliams98MoritaEq}, Proposition 2.50]\label{lanceProp25}
Sei $\mfrakB$ eine \csalgebra{} und $E$ ein Hilbert-\cstern{}Modul "uber~$\mfrakA$.
F"ur jeden \sterns{}Homomorphismus $\alpha : \mfrakB \to \Adjhm{E}$
sind die folgenden Aussagen "aquivalent:
\begin{enumaequiv}
\item $\alpha$ ist nicht-ausgeartet.
\item Es existiert ein eindeutiger unitaler \sterns{}Homomorphismus
  $\overline{\alpha} : \Multcs{\mfrakB} \to \Adjhm{E}$,
  der $\alpha$ fortsetzt und
  strikt stetig auf der abgeschlossenen Einheitskugel
$\overline{\cmkug}_{\Multcs{\mfrakB}}(0,1)$
  ist.
\item F"ur jede  approximative Eins $(e_\lambda)_\lambda$ von $\mfrakB$ gilt:
  $\alpha(e_\lambda) \overset{\lambda}{\longrightarrow} \Id_E$ strikt.
\dremark{Siehe \cite{LanceHmod}, S. 19 unten}%
\item Es existiert eine approximative Eins $(e_\lambda)_\lambda$ von $\mfrakB$ so,
  \dass gilt: $\alpha(e_\lambda) \overset{\lambda}{\longrightarrow} \Id_E$ strikt.
\end{enumaequiv}
\end{bemerkung}

Die Fortsetzung $\overline{\alpha}$ aus (b) wird h"aufig mit $\alpha$ bezeichnet.
\dremark{und wie folgt definiert:
$\overline{\alpha}(b)$ wird definiert durch:
$\sum_i \alpha(a_i) \xi_i \mapsto \sum_i \alpha(b a_i) \xi_i$.

}%
\skiptext

Mit Hilfe von nicht-ausgearteten \sterns{}Homomorphismen erh"alt man regul"are Operatoren:

\begin{bemerkung}[\cite{LanceHmod}, Proposition 10.7]\label{LanceProp10:7}
Sei $\mfrakB$ eine \csalgebra, $E$ ein Hilbert-\cstern{}Modul "uber $\mfrakA$
und $\alpha : \mfrakB \to \Adjhm{E}$ ein nicht-ausgearteter \sterns{}Homomorphismus.
Sei $T \in \Multwor{\mfrakB}$.
Dann findet man einen regul"aren Operator $\alpha(T)$ auf $E$ mit den folgenden Eigenschaften:
\begin{enumaufz}
\item $\Bild(\alpha(a)) \subseteq D(\alpha(T))$ f"ur alle $a \in D(T)$,
\item $\alpha(T)\alpha(a) = \alpha(Ta)$ f"ur alle $a \in D(T)$,
\item $\alpha(z_T) = z_{\alpha(T)}$ und
\item $\alpha(T^*) = \alpha(T)^*$.
\end{enumaufz}
\end{bemerkung}

F"ur selbstadjungierte regul"are Operatoren hat man den folgenden Funktionalkalk"ul
(\cite{LanceHmod}, Theorem 10.9):

\begin{satz}[Funktionalkalk"ul]\label{regOpFkalkuel}
Sei $E$ ein Hilbert-\cstern{}Modul "uber $\mfrakA$ und
$T \in \Multwor{E}$ selbstadjungiert.
Sei $\iota : \mbbR \to \mbbC, \lambda \mapsto \lambda$, und
$f : \mbbR \to \mbbC, \lambda \mapsto \lambda \left(1+\lambda^2\right)^{-1/2}$.
Dann existiert ein \sterns{}Homomorphismus $\varphi_T : C(\mbbR) \to \Multwor{E}$ derart,
\dass $\varphi_T(\iota) = T$ und $\varphi_T(f) = z_T$ gilt.
\dliter{\cite{LanceHmod}, Th. 10.9; \cite{Woronowicz91}, Th. 1.5}%
\index[S]{phiT@$\varphi_T$}%
\end{satz}

\dremww{
Beweisskizze.
Die Abbildung $f$ ist ein Hom"oomorphismus von $\mbbR$ auf $]-1,1[$.
Nach \cite{LanceHmod}, Theorem 10.4, ist $z_T$ selbstadjungiert.
Somit ist
\[ \alpha : C_0(\mbbR) \to \Adjhm{E}, g \mapsto (g \circ f^{-1})(z_T) \]
ein \sterns{}Homomorphismus.
\dremark{$\alpha$ wird mit Hilfe der stetigen Funktionalkalk"uls auf \csalgebren{}
  definiert (Beachte: $\Adjhm{E}$ ist eine \csalgebra{}.).}%
Wegen $\Multwor{C_0(\mbbR)} \cong C(\mbbR)$ erh"alt man mittels $\alpha$ einen
\sterns{}Homomorphismus $\hat{\alpha} : C(\mbbR) \to \Multwor{E}$.
}%

\dremark{Beachte: $\varphi\restring_{C_0(\mbbR)}$ l"a\cms{}t sich nach \cite{LanceHmod}, Prop. 2.1,
  \emph{eindeutig} zu einem \sterns{}Homomorphismus
  $\overline{\varphi_T} : C_b(\mbbR) \to \Adjhm{E}$ fortsetzen.
  Wegen der Eindeutigkeit folgt:
\begin{equation}\label{phiTEqphiRestrCb}
\overline{\varphi_T} = \varphi_T\restring_{C_b(\mbbR)}.
\end{equation}
}%

In \cite{Kustermans97FunctionalCalcRegOp}, Abschnitt 3,
(siehe auch \cite{Woronowicz91}, Abschnitt 1) wird ein Funktionalkalk"ul
f"ur normale regul"are Operatoren eingef"uhrt, wobei ein Operator $T \in \Multwor{E}$
normal hei"st, falls $z_T$ ein normales Element von $\Adjhm{E}$ ist.
\skiptext

Sei $E$ ein \hcsmodul und $T \in \Multwor{E}$ selbstadjungiert.
Es gilt:
\[ \varphi_T(\mathbbm{1})\restring_{\varphi_T(C_0(\mbbR))E} = \Id_{\varphi_T(C_0(\mbbR))E}. \]
\dremark{folgt aus $\eqref{phiTEqphiRestrCb}$ und der Definition der Fortsetzung $\hat{\varphi_T}$, siehe \cmzB{}
  \cite{RaeburnWilliams98MoritaEq}, Beweis von Prop. 2.50}%
Somit hat man:
\dremark{denn $E = \varphi_T(C_0(\mbbR)) E$
  (\cite{HewittRossHarmonicAna2}, Th. 32.22
   (Verallgemeinerung der Faktorisierungssatzes von Cohen))}%
\begin{equation}\label{eqvarphih1EqId}
\varphi_T(\mathbbm{1}) = \Id_E.
\end{equation}
Wie man dem Beweis von \cite{LanceHmod}, Theorem~10.9, entnimmt,
ist $\varphi_T\restring_{C_0(\mbbR)} : C_0(\mbbR) \to \Adjhm{E}$
ein nicht"=ausgearteter \sterns{}Homomorphismus.
Mit \cite{LanceHmod}, Proposition 2.1, folgt:
\dremark{Beachte: $\overline{\varphi_T} : C_b(\mbbR) \to \Adjhm{E}$.}%
\begin{equation}\label{eqvarphihCbSsAdjhm}
\varphi_T(C_b(\mbbR))
\dremarkm{\overset{\eqref{phiTEqphiRestrCb}}{=} \overline{\varphi_T}(C_b(\mbbR)) }
\subseteq \Adjhm{E}.
\end{equation}

Ist $T \in \Adjhm{E}$ selbstadjungiert,
so kann man das Spektrum $\sigma(T)$ von $T$ wie folgt beschreiben:

\begin{lemma}[Vgl. \cite{Kustermans97FunctionalCalcRegOp}, Lemma 3.2]
Sei $E$ ein \hcsmodul und $T \in \Adjhm{E}$ selbstadjungiert.\dremark{normal gen"ugt}
Dann gilt:\dremark{folgt mit \cite{Kustermans97FunctionalCalcRegOp}, Lemma 3.2 und Corollary 3.8}
\[ \sigma(T) = \{ \lambda \in \mbbR \setfdg \forall f \in \Kern(\varphi_T) : f(\lambda) = 0 \}. \]
\end{lemma}

\dremark{
\begin{lemma}
Sei $E$ ein \hcsmodul.
Dann ist $\Adjhm{E}$ eine volle Unteralgebra von $\mLinStet(E)$.
Insbesondere gilt: $\sigma_{\Adjhm{E}}(T) = \sigma_{\mLinStet(E)}(T)$
f"ur alle $T \in \Adjhm{E}$.
\end{lemma}

Beweis
Sei $T \in G(\mLinStet(E)) \cap \Adjhm{E}$.

Zeige: $T^*$ ist injektiv.
Sei $x \in \Kern(T^*)$.
Dann gilt $0 = \skalpr{T^*x}{y} = \skalpr{x}{Ty}$ f"ur alle $y \in E$.
Da $T$ surjektiv auf $E$ ist und $\skalpr{\cdot}{\cdot\cdot}$ nicht-ausgeartet,
folgt: $x=0$.

Zeige: $T^*$ ist surjektiv auf $E$.

Annahme: $T^*$ ist nicht surjektiv auf $E$.

Dann findet man ein $y \in E$ mit: $\forall x \in E : T^*x \neq y$.
Da $T$ injektiv ist, folgt: $\forall x \in E : TT^* x \neq Ty$.
Dies steht im Widerspruch zur Surjektivit"at von $T$ auf $E$.

Seien $x,y \in E$.
Setze $z := (T^*)^{-1} y$.
Dann gilt:
\[ \skalpr{T^{-1}x}{y}
=  \skalpr{T^{-1}x}{T^* z}
=  \skalpr{T T^{-1}x}{z}
=  \skalpr{x}{(T^*)^{-1} y}.  \]
}%

F"ur das Spektrum eines selbstadjungierten regul"aren Operators wird
die folgende Definition eingef"uhrt,
die nach dem vorherigen Lemma konsistent zur Definition des
Spektrums f"ur selbstadjungierte Operatoren aus $\Adjhm{E}$ ist:

\begin{definitn}
Sei $E$ ein \hcsmodul und $T \in \Multwor{E}$ selbstadjungiert.\dremark{normal gen"ugt}
Dann ist das \defemph{Spektrum} von $T$ definiert als die Menge
\index[B]{Spektrum eines Operators}%
\index[S]{sigmaA@$\sigma(A)$ (Spektrum)}%
\[ \sigma(T) := \{ \lambda \in \mbbR \setfdg \forall f \in \Kern(\varphi_T) : f(\lambda)=0 \}. \]
\end{definitn}

Offensichtlich ist $\sigma(T)$ eine abgeschlossene Teilmenge von $\mbbR$.
Au"serdem gilt: $\sigma(T) \neq \emptyset$ (\cite{Kustermans97FunctionalCalcRegOp}, S. 10).

\begin{satz}[Vgl. \cite{Kustermans97FunctionalCalcRegOp}, Theorem 3.4]
Sei $E$ ein \hcsmodul und $T \in \Multwor{E}$ selbstadjungiert.
Dann existiert ein eindeutig bestimmter nicht"=ausgearteter, injektiver \sterns{}Homomorphismus
$\psi_T : C_0(\sigma(T)) \to \Adjhm{E}$ derart, \dass $\psi_T(\iota) = T$ gilt,
wobei $\iota : \sigma(T) \to \mbbC, \lambda \mapsto \lambda$, sei
und $\psi_T(\iota)$ wie in Proposition~\ref{LanceProp10:7} definiert.
\index[S]{psiT@$\psi_T$}%
\end{satz}

Ist $f : D(f) \subseteq \mbbC \to \mbbC$
beschr"ankt auf $\sigma(T)$ mit $\sigma(T) \subseteq D(f)$,
so gilt:
\begin{equation}\label{eqpsiTfinAdj}
\overline{\psi_T}(f\restring_{\sigma(T)}) \in \Adjhm{E}
\end{equation}
(\cite{Kustermans97FunctionalCalcRegOp}, S. 10\dremark{nach Notation 3.5}).

\begin{bemerkung}[Vgl. \cite{Kustermans97FunctionalCalcRegOp}, Lemma 3.6]\label{varphiTeqPsiT}
Sei $E$ ein \hcsmodul und $T \in \Multwor{E}$ selbstadjungiert.\dremark{normal gen"ugt}
Dann gilt: $\varphi_T(f) = \psi_T(f\restring_{\sigma(T)})$ f"ur alle $f \in C(\mbbR)$.
\end{bemerkung}

Man beachte, \dass in der obigen Proposition $\psi_T(f\restring_{\sigma(T)})$
mittels Proposition~\ref{LanceProp10:7} definiert wird,
da nach Beispiel~\ref{bspRegOpC0Omega} gilt:
$f\restring_{\sigma(T)} \in C(\sigma(T)) \cong \Multwor{C_0(\sigma(T))}$.

\dremww{
\begin{bemerkung}
Seien $E$, $F$ Hilbert-\cstern{}Moduln "uber $A$.
Sei $t : D(t) \subseteq E \to F$ ein regul"arer Operator derart,
\dass $t$ und $t^*$ dichtes Bild haben.
Dann sind $t$ und $t^*$ injektiv,
$t^{-1}$ und $(t^*)^{-1}$ sind regul"ar und $(t^{-1})^* = (t^*)^{-1}$.
\dremark{F"ur die Injektivit"at von $t$ braucht man die Vor., \dass $t$ dichtes Bild hat, nicht.}%
\dliter{\cite{LanceHmod}, S. 104 unten, S. 105 oben}
\end{bemerkung}

\begin{proof}
\bewitemph{(i):} Zeige: $t$ und $t^*$ sind injektiv.

Sei $x \in \Kern(t)$.
Da $t^*$ dichtes Bild hat, findet man eine Folge $y$ in $D(t^*)$
mit $\lim t^*(y_n) = x$.
Es folgt
\dremark{$(*)$:
  $\norm{\skalpr{x}{x} - \skalpr{x}{t^*(y_n)}}_A
  =  \norm{\skalpr{x}{x-t^*(y_n)}}
  \leq  \norm{x} \cdot \norm{x-t^*(y_n)}
  \to 0$ f"ur $n \to \infty$}
\[ \skalpr{x}{x}
=  \skalpr{x}{\lim t^*(y_n)}
\overset{\dremarkm{(*)}}{=}  \lim \skalpr{x}{t^*(y_n)}
=  \lim \skalpr{tx}{y_n}
=  0, \]
also $x=0$.
Analog folgt, \dass $t^*$ injektiv ist.

\bewitemph{(ii):} Zeige mit \cite{Lance76TensorProducts}, Prop.~9.5: $t^{-1}$ ist regul"ar.

\bewitemph{(I):} Da $t$ (bzw. $t^*$) dichtes Bild hat,
ist $t^{-1}$ (bzw. $(t^*)^{-1}$) dicht definiert.
Ferner ist $D(t^{-1}) = \Bild(t)$ ein Modul.

\bewitemph{(II):} Zeige: $(t^{-1})^*$ ist dicht definiert.
F"ur alle $a \in D(t)$, $x \in \Bild(t^*)$
(also $(t^*)^{-1} \in D(t^*)$) gilt:
\[ \skalpr{x}{a}
=  \skalpr{t^* (t^*)^{-1}x}{a}
=  \skalpr{(t^*)^{-1}x}{ta}. \]
F"ur alle $a \in D(t)$,
$x \in D((t^{-1})^*)$ gilt:
\[ \skalpr{x}{a}
=  \skalpr{x}{t^{-1}ta}
=  \skalpr{(t^{-1})^*x}{ta}. \]
Nach Definition von $(t^{-1})^*$ folgt somit:
$D((t^*)^{-1}) = \Bild(t^*) \subseteq D((t^{-1})^*)$.
Somit ist $(t^{-1})^*$ dicht definiert.

\bewitemph{(III):} Zeige: $t^{-1}$ ist abgeschlossen.
Definiere
\begin{align*}
\sigma_{E \oplus F} &: E \oplus F \to F \oplus E, (x,y) \mapsto (y,x),\dremarkm{\qquad\text{der Flip}} \\
\rho_{1,E \oplus F} &: E \oplus F \to E \oplus F, (x,y) \mapsto (-x,y), \\
\rho_{2,E \oplus F} &: E \oplus F \to E \oplus F, (x,y) \mapsto (x,-y), \\
v &: E \oplus F \to F \oplus E, (x,y) \mapsto (y,-x), \\
\tilde{v} &: F \oplus E \to E \oplus F, (y,x) \mapsto (x,-y), \\
\end{align*}
Es gilt: $\Graph(t^{-1}) = \Graph(t)^{-1} = \sigma \Graph(t)$.
\dremark{$\Graph(t)^{-1}$ als Relation; 1. \glqq$=$\grqq:
  $(a,b) \in \Graph(t)
   \iff a \in D(t) \wedge b=ta
   \iff b \in D(t^{-1}) \wedge a=t^{-1}b
   \iff (b,a) \in \Graph(t^{-1})$}%
Da $t$ abgeschlossen ist, ist somit $t^{-1}$ abgeschlossen.

\bewitemph{(IV):} Zeige: $\Graph(t^{-1}) \oplus \tilde{v}\Graph((t^{-1})^*) = F \oplus E$.
Es gilt:\dremark{$w$ ist Isometrie}
\begin{equation}\label{eqLan104vtiGraph}
\begin{split}
   \tilde{v}^{-1} \Graph((t^{-1})^*)
&= \tilde{v}^{-1} \tilde{v} \Graph(t^{-1})^\bot
=  \Graph(t^{-1})^\bot
=  (\Graph(t)^{-1})^\bot  \\
&= (w \Graph(t))^\bot
=  w(\Graph(t))^\bot
=  w(v^{-1}v \Graph(t)^\bot)
=  w(v^{-1} \Graph(t^*)).
\end{split}
\end{equation}

\bewitemph{(iii):} Zeige: $(t^*)^{-1}$ ist regul"ar.

\bewitemph{(iv):}

\dremark{1.4.'07, S. 1--3}%
\end{proof}
}%

\section{Der Satz von Stone f"ur Hilbert-\cstern{}Moduln}
\label{subsecSatzVonStoneHCsMod}

In diesem Abschnitt "ubertragen wir Aussagen "uber $C_0$-Gruppen
aus unit"aren Elementen
aus dem Artikel \cite{HollevoetEtAl92StonesTh}, die dort f"ur
\csalgebren formuliert sind, auf \hcsmoduln.
Insbesondere beweisen wir eine Variante des Satzes von Stone
f"ur \hcsmoduln.
\dremark{Wozu: Sp"ater zeigen, was mit den unbeschr"ankten
  Operatoren los ist, falls $X$ ein \hcsmodul ist.}%

\begin{erinnerung}
Operatorgruppen werden in Definition~\ref{defC0Halbgr} definiert.
Eine Operatorgruppe $(T_t)_\indtGr$ auf einem Banachraum $X$
hei"st stark stetig oder $C_0$-Gruppe,
falls gilt: $\lim_{t \to 0} T_t x = x$ f"ur alle $x \in X$.
\end{erinnerung}

\begin{definitn}
Es sei $E$ ein \hcsmodul und $\Omega$ ein topologischer Raum.
Eine Abbildung $f : \Omega \to \Adjhm{E}$ hei"st \defemphi{strikt stetig},
falls die Abbildungen
$\Omega \ni \omega \mapsto f(\omega)x$ und $\Omega \ni \omega \mapsto f(\omega)^*x$
f"ur alle $x \in E$ stetig sind.
\dmarginpar{pr}\dremark{Zusammenhang zu \ref{defStriktStetig}?}%
\dremark{Gilt ebenfalls f"ur einen Operatorraum $X$,
  also f"ur $(X,\Adjlor{X})$ anstelle von $(E,\Adjhm{E})$.}%
\dliter{\cite{HollevoetEtAl92StonesTh}, S. 228 oben}%
\end{definitn}


\begin{bemerkung}
Sei $E$ ein \hcsmodul und $(U_t)_\indtGr$ eine
Operatorgruppe auf $E$ mit $U_t \in \Adjhm{E}$ unit"ar f"ur alle $\indtGr$.
Dann sind die folgenden Aussagen "aquivalent:
\dremark{Gilt ebenfalls f"ur einen Operatorraum $X$,
  also f"ur $(X,\Adjlor{X})$ anstelle von $(E,\Adjhm{E})$.}%
\dliter{F"ur \hcsmoduln siehe \cite{RaeburnWilliams98MoritaEq}, C.8}%
\begin{enumaequiv}
\item $(U_t)_\indtGr$ ist stark stetig.
\item $(U_t)_\indtGr$ ist strikt stetig.
\end{enumaequiv}
\end{bemerkung}

\begin{proof}
\bewitemphqq{(a)$\Rightarrow$(b)} folgt mit der
Charakterisierung der starken Stetigkeit von Halbgruppen (Proposition~\ref{charHGIstC0})
wegen $U_t^* = U_{-t}$ f"ur alle $\indtGr$.
\dremark{Es gelte \bewitemph{(a)}.
Es gilt $U_t U_{-t} = U_0 = \Id_E = U_{-t} U_t$
f"ur alle $\indtGr$, also $U_{-t} = U_t^{-1} = U_t^*$.
Es ist $(U_t^*)_\indtGr$ eine Operatorgruppe.
Weiter gilt f"ur alle $x \in E$: $\norm{ U_t^* x - x } = \norm{ U_{-t} x - x } \to 0$ f"ur $t \to 0$.
Somit folgt \bewitemph{(b)}.}%

\bewitemphqq{(b)$\Rightarrow$(a)} ist offensichtlich.
\end{proof}

In dem Artikel \cite{HollevoetEtAl92StonesTh} werden stets strikt stetige
Gruppen mit uni\-t"a\-ren Elementen aus $\Multcs{\mfrakA}$ betrachtet,
die wegen der obigen Proposition
das gleiche sind wie $C_0$-Gruppen mit unit"aren Elementen aus $\Multcs{\mfrakA}$.
\skiptext

Die folgende Proposition stellt eine Verallgemeinerung einer Aussage
aus \cite{HollevoetEtAl92StonesTh}, S. 229, von \csalgebren{} auf
\hcsmodul{}n dar:

\begin{bemerkung}\label{UtIstSSUnitaereGr}
Sei $E$ ein \hcsmodul{} und $h \in \Multwor{E}$ selbstadjungiert.
Definiere $e_t : \mbbR \to \mbbC, \lambda \mapsto \exp(\mri \lambda t)$,
und $U_t := \varphi_h(e_t) = \exp(\mri th)$ f"ur alle $\indtGr$.
\index[S]{ekleint@$e_t$}%
Dann ist $(U_t)_\indtGr$ eine $C_0$-Gruppe auf $E$
mit $U_t \in \Adjhm{E}$ unit"ar f"ur alle $\indtGr$.
\end{bemerkung}

\dremark{Zusammenhang von \ref{UtIstSSUnitaereGr} und \ref{HollevoetTh21}
zu \cite{Kustermans97FunctionalCalcRegOp}, Prop. 7.13, bzw.
\cite{Kustermans97OneParamRepOnCsAlg}, Th. 5.10:

Bei Kustermans ist $h \in \Multwor{E}$ strikt positiv mit $U_t = h^{\mri t}$.
Hierbei ist $h^{\mri t}$ wie folgt definiert:
F"ur $f : \mbbR \to \mbbC, \lambda \mapsto \lambda^{\mri t}$
definiere $h^{\mri t} := f(h) = \varphi_h(f)$
(\cite{Kustermans97FunctionalCalcRegOp}, Definition 7.10).
In \cite{HollevoetEtAl92StonesTh} wird hingegen $U_t = \varphi_h(e_t)$ betrachtet.
Allerdings soll nach dem Text am Anfang von §5 in \cite{Kustermans97OneParamRepOnCsAlg}
\cite{Kustermans97OneParamRepOnCsAlg}, Th. 5.10, mit \ref{HollevoetTh21} "ubereinstimmen.}%

Zum Beweis formulieren wir das folgende

\begin{lemma}\label{etKonvStrikt}
Es gilt f"ur alle $f \in C_0(\mbbR)$:
$\norm{e_t f - f}_\infty \to 0$ f"ur $t \to 0$.
\end{lemma}

\begin{proof}
Sei $f \in C_0(\mbbR)$ und $\varepsilon \in \mbbR_{>0}$.
Dann findet man ein $a \in \mbbR_{>0}$ derart,
\dass gilt: $\norm{ f\restring_{\mbbR\setminus [-a,a]} } < \frac{\varepsilon}{2}$.
Setze $M := \norm{f}_\infty + 1$.
Man hat:
$\mre^{\mri ta} \to 1$ f"ur $t \to 0$.
Somit findet man ein $t_0 \in \mbbR_{>0}$ so,
\dass gilt:
\[ \forall t \in\, ]0,t_0[\, : \abs{\mre^{\mri t a} - 1} < \frac{\varepsilon}{M}.
   \dremarkm{\qquad (*)}\]

Sei im folgenden $t \in\,]0,t_0[$.
Man hat:
\[ \norm{e_t f - f}_\infty
=  \sup_{\lambda \in \mbbR} \abslr{\mre^{\mri t \lambda} - 1} \, \abs{f(\lambda)}. \]

Sei $\lambda \in \mbbR$.

\bewitemph{1. Fall:} $\lambda \in [-a,a]$.
Dann gilt:\dremark{$\abs{\mre^{\mri t \lambda} - 1} < \frac{\varepsilon}{M}$ nach $(*)$.}
\[ \abslr{\mre^{\mri t \lambda} - 1}\, \abs{f(\lambda)}
\leq  \frac{\varepsilon}{M} M
=  \varepsilon. \]

\bewitemph{2. Fall:} $\lambda \notin [-a,a]$.
Man hat:
\[ \abslr{\mre^{\mri t \lambda} - 1}\, \abs{f(\lambda)}
\leq \dremarkm{2 \cdot \frac{\varepsilon}{2} =} \varepsilon. \]

Somit ergibt sich: $\norm{e_t f - f}_\infty \leq \varepsilon$.
\dremark{Der Beweis geht nicht mittels
  $\norm{e_t f - f}_\infty \leq \norm{e_t - \mathbbm{1}} \, \norm{f}$,
  denn $\norm{e_t - \mathbbm{1}}_\infty \not\to 0$ f"ur $t \to 0$.}%
\end{proof}

\begin{proof}[Beweis von Proposition~\ref{UtIstSSUnitaereGr}]
Mit \eqref{eqvarphihCbSsAdjhm} erh"alt man:
$U_t = \varphi_h(e_t) \in \Adjhm{E}$ f"ur alle $\indtGr$.
Da $\varphi_h$ nach dem Funktionalkalk"ul f"ur regul"are Operatoren (Satz~\ref{regOpFkalkuel})
ein \sterns{}Ho\-mo\-mor\-phis\-mus ist,
folgt: $U_s U_t = U_{s+t}$ f"ur alle $s, \indtGr$.
\dremark{Es gilt:
$  e_{s+t}(\lambda)
=  \mre^{\mri(s+t)\lambda}
=  \mre^{\mri s \lambda} \mre^{\mri t \lambda}
=  (e_s e_t)(\lambda)$.
Damit folgt:
$  U_{s+t}
=  \varphi_h(e_{s+t})
=  \varphi(e_s e_t)
=  \varphi(e_s) \varphi(e_t)
=  U_s U_t$.}%
Nach \eqref{eqvarphih1EqId} gilt $\varphi_h(\mathbbm{1}) = \Id_E$.
Somit folgt
\[ U_t^* U_t
\dremarkm{=  \varphi_h(e_t)^* \varphi_h(e_t)
=  \varphi_h(e_t^* e_t)
=  \varphi_h(\mathbbm{1})}
\overset{\dremarkm{\eqref{eqvarphih1EqId}}}{=}  \Id_E
=  U_t U_t^* \]
f"ur alle $\indtGr$, also ist $U_t$ unit"ar.

Es bleibt zu zeigen: $(U_t)_\indtGr$ ist stark stetig.

Man hat: $e_t \in C_b(\mbbR) \cong \Multcs{C_0(\mbbR)}$ f"ur alle $\indtGr$.
Nach Lemma~\ref{etKonvStrikt} gilt: $e_t \to \mathbbm{1}$ strikt f"ur $t \to 0$.
\dremark{$e_t \cong M_{e_t}$ strikt: $M_{e_t} g \to g$ und $g M_{e_t} g \to g$}%
Da $\varphi_h$ nach Proposition~\ref{lanceProp25}
(Charakterisierung nicht"=ausgearteter \sterns{}Homomorphismen)
auf der abgeschlossenen Einheitskugel $\overline{\cmkug}_{C_b(\mbbR)}(0,1)$
strikt stetig ist, folgt:
\[ U_t = \varphi_h(e_t) \to \varphi_h(\mathbbm{1})
\overset{\dremarkm{\eqref{eqvarphih1EqId}}}{=}  \Id_E = U_0 \quad\text{strikt f"ur } t \to 0. \]

Damit ist $(U_t)_\indtGr$ strikt stetig in $0$.
Mit der Charakterisierung der starken Stetigkeit von Halbgruppen
(Proposition~\ref{charHGIstC0}) folgt,
\dass $(U_t)_\indtGr$ stark stetig ist.
\dremark{Alternativ: Mit der Charakterisierung der starken Stetigkeit von Halbgruppen
  (Proposition~\ref{charHGIstC0}) folgt,
  \dass die Halbgruppe $(U_t)_\indtHG$ stark stetig ist.
  Nach \cite{EngelNagelSemigroups}, Proposition~II.3.11, kann man die
  $C_0$-Halbgruppe $(U_t)_\indtHG$ in eine $C_0$-Gruppe $(V_t)_\indtGr$ einbetten,
  die aber mit $(U_t)_\indtGr$ "ubereinstimmt.}%
\end{proof}

Die beiden folgenden S"atze sind Verallgemeinerungen
von \cite{HollevoetEtAl92StonesTh}, Theorem 2.1 und Proposition 2.2,
von \csalgebren{} auf Hil\-bert-\cstern{}Mo\-duln
und ergeben zusammen eine Variante des Satzes von Stone.
Beide Beweise enthalten Anpassungen an die \hcsmodul-Situation
und werden ausf"uhrlicher als in \cite{HollevoetEtAl92StonesTh} dargestellt.
\dremark{Der Beweis von \ref{HollevoetTh21} wird durch Benutzung der Aussage
  \cite{RaeburnWilliams98MoritaEq}, Proposition C.17, verk"urzt.}%

\begin{satz}[Satz von Stone]\label{HollevoetTh21}
Sei $E$ ein Hilbert-\cstern{}Modul und $(U_t)_\indtGr$ eine $C_0$-Gruppe
auf $E$ mit $U_t \in \Adjhm{E}$ unit"ar f"ur alle $\indtGr$.
Dann gilt:
\begin{enumaufz}
\item Es existiert ein selbstadjungiertes $h \in \Multwor{E}$ derart,
\dass $U_t = \varphi_h(e_t) = \exp(\mri th)$ f"ur alle $\indtGr$ gilt.

\item Ist $(U_t)_\indtGr$ normstetig, so ist $h \in \Adjhm{E}$.
\dmarginpar{DB}\dmarginpar{Bezeichnung}\dremark{Bezeichnung $h$: Sonst $T$.}%
\dliter{Verallgemeinerung von \cite{HollevoetEtAl92StonesTh}, Theorem 2.1,
  auf Hilbert-\cstern{}Moduln}%
\end{enumaufz}
\end{satz}

\dremww{
\begin{bemerkung}\label{striktStetigHM}
Sei $E$ ein Hilbert-\cstern{}Modul "uber $\mfrakA$,
$(U_t)_\indtGr$ eine Operatorgruppe auf $E$ mit $U_t \in \Adjhm{E}$ f"ur alle $\indtGr$.
Dann sind die folgenden Aussagen "aquivalent:
\begin{enumaequiv}
\item $t \mapsto U_t x$ und $t \mapsto U_t^* x$ sind stetig f"ur alle $x \in E$.
\item $t \mapsto U_t T$ und $t \mapsto U_t^* T$ sind stetig f"ur alle $T \in \kptOphm(E)$.
\dremark{3.7.'08, S. 3/4}%
\end{enumaequiv}
\end{bemerkung}

\begin{proof}
\bewitemph{\glqq(a) $\Rightarrow$ (b):\grqq}
Es gelte \bewitemph{(a)}.
Dann ist $t \mapsto U_t \circ \cketbra{x}{y} = \cketbra{U_t x}{y}$ stetig f"ur alle $x,y \in E$.
\dremark{Seien $x,y \in E$, sei $\varepsilon \in \mbbR_{>0}$.
Man findet ein $\delta \in \mbbR_{>0}$ mit:
\[ \forall s,t \in \mbbR :
   \abs{s-t} < \delta \Rightarrow \norm{(U_s - U_t)x} < \frac{\varepsilon}{\norm{y}}. \]
Seien $s,t \in \mbbR$ mit $\abs{s-t} < \delta$.
Es gilt:
\begin{align*}
   \norm{ U_s \circ \cketbra{x}{y} - U_t \circ \cketbra{x}{y} }
=  \norm{ \cketbra{(U_s - U_t) x}{y} }
\leq  \norm{ (U_s - U_t)x } \norm{y}
<  \varepsilon.
\end{align*}
}%
Da $\norm{U_t} = 1$ f"ur alle $\indtGr$, folgt:
$t \mapsto U_t T$ ist stetig f"ur alle $T \in \kptOphm(E)$.
\dremark{\bewitemph{1. Schritt: }
  Zeige: Stetig f"ur alle $T \in \vrerz{ \cketbra{x}{y} \setfdg x,y \in E }$.
Sei $\varepsilon \in \mbbR_{>0}$.
Seien $x,y \in E^n$.
Man findet $\delta \in (\mbbR_{>0})^n$ mit:
\[ \forall s,t \in \mbbR\, \forall i \in \haken{n} :
   \abs{s-t}<\delta_i  \Rightarrow
   \norm{(U_s - U_t)x_i} < \frac{\varepsilon}{n (\norm{y_1} + \dots + \norm{y_n})}. \]
Setze $\delta := \min\{\delta_1, \dots, \delta_n\}$.
Seien $s,t \in \mbbR$ mit $\abs{s-t} < \delta$.
Es gilt:
\begin{align*}
   \norm{ ( U_s - U_t) \circ \sum_{i=1}^n \cketbra{x_i}{y_i} }
&= \norm{ \sum \cketbra{(U_s - U_t)x_i}{y_i} }
\leq  \sum \norm{  \cketbra{(U_s - U_t)x_i}{y_i} }  \\
&\leq  \sum \norm{ (U_s - U_t)x_i } \cdot \norm{y_i}
<  \sum_{i=1}^n \frac{\varepsilon}{n (\norm{y_1} + \dots + \norm{y_n}) } \norm{y_i}
\leq  n \frac{\varepsilon}{n} = \varepsilon.
\end{align*}

\bewitemph{2. Schritt: }
  Zeige: Stetig f"ur alle
  $T \in \kptOphm(E) = \overline{\cmlin}\{ \cketbra{x}{y} \setfdg x,y \in E \}$.
Sei $T \in \kptOphm(E)$.
Dann existiert eine gegen $T$ konvergente Folge $z$
in $\vrerz{ \cketbra{x}{y} \setfdg x,y \in E}$.
F"ur alle $n \in \mbbN$ findet man $x_n^{(1)},\dots,x_n^{(k_n)}, y_n^{(1)},\dots,y_n^{(k_n)} \in E$
mit: $z_n = \sum_{i=1}^{k_n} \cketbra{x_n^{(i)}}{y_n^{(i)}}$.
Man findet ein $m \in \mbbN$ mit $\norm{T-z_m^{(i)}} < \frac{\varepsilon}{3}$.
Es gilt:
\begin{align*}
   \norm{(U_s - U_t) \circ T}
&= \norm{(U_s - U_t) \circ (T - \sum_{i=1}^{k_n} \cketbra{x_m^{(i)}}{y_m^{(i)}} +
                                \sum_{i=1}^{k_n} \cketbra{x_m^{(i)}}{y_m^{(i)}}) }  \\
&\leq  \norm{(U_s - U_t) \circ (T - \sum_{i=1}^{k_n} \cketbra{x_m^{(i)}}{y_m^{(i)}}) } +
      \norm{(U_s - U_t) \circ \sum_{i=1}^{k_n} \cketbra{x_m^{(i)}}{y_m^{(i)}} }  \\
&\leq  \norm{U_s - U_t} \cdot \norm{T - \sum_{i=1}^{k_n} \cketbra{x_m^{(i)}}{y_m^{(i)}}} +
      \norm{ \sum \cketbra{ (U_s - U_t) x_m^{(i)}}{y_m^{(i)}} }  \\
&\leq  2 \cdot \frac{\varepsilon}{3} +
      \sum \norm{ \cketbra{ (U_s - U_t) x_m^{(i)}}{y_m^{(i)}} }
\overset{\dremarkm{\text{wie im 1. Schritt}}}{<}  \varepsilon.
\end{align*}
}%
Analog folgt: $t \mapsto U_t^* T$ ist stetig f"ur alle $T \in \kptOphm(E)$.

\bewitemph{\glqq(b) $\Rightarrow$ (a):\grqq}
Es gelte \bewitemph{(b)}.
Sei $x \in E$.
Nach \cite{RaeburnWilliams98MoritaEq}, 2.31, findet man ein $y \in E$ mit
$x = y \skalpr{y}{y} = \cketbra{y}{y} (y)$.
Setze $T := \cketbra{y}{y} \in \kptOphm(E)$.
Da die Punktauswertung stetig ist, sind auch
$t \mapsto U_t T(y) = U_t x$ und $t \mapsto U_t^* T(y) = U_t^* x$ stetig.
\end{proof}
}

\begin{proof}
\bewitemph{(i):}
Mit $C^*(\mbbR)$ sei die Gruppen-\cstern{}Algebra der lokalkompakten Gruppe $\mbbR$
bezeichnet, also die Vervollst"andigung von $C_c(\mbbR)$
bez"uglich der sogenannten universellen \cstern{}Norm.
\dremark{\cite{RaeburnWilliams98MoritaEq}, S. 280/281}%
Die duale Gruppe der lokalkompakten Gruppe $\mbbR$
wird als $\hat{\mbbR}$ notiert.
Es gilt nach \cite{FollandAbstractHarmonicAna}, Theorem 4.5: $\hat{\mbbR} \cong \mbbR$.
Die Fouriertransformation $\mathscr{F} : L^1(\mbbR) \to C_0(\hat{\mbbR}) \cong C_0(\mbbR)$
l"a\cms{}t sich nach \cite{RaeburnWilliams98MoritaEq}, Example C.20,
zu einem \sterns{}Iso\-mor\-phis\-mus
$\hat{\mathscr{F}} : C^*(\mbbR) \arrowbij C_0(\mbbR)$ auf $C_0(\mbbR)$ fortsetzen.
\dremark{Dies gilt allgemein f"ur eine beliebige abelsche, lokalkompakte Gruppe.
  Zur Fouriertransformation: \cite{FollandAbstractHarmonicAna}, Prop. 4.13.}%
\dremark{3.7.'08/2a}%

Nach \cite{RaeburnWilliams98MoritaEq}, Proposition C.17, wird f"ur alle $f \in C_c(\mbbR)$ durch
\[ \alpha(f) = \int_\mbbR f(t) U_t\, dt \]
ein wohldefinierter Operator in $\Adjhm{E}$ definiert,
und $\alpha$ l"a\cms{}t sich\dremark{\cite{RaeburnWilliams98MoritaEq}, Proposition C.17}
zu einem nicht-ausgearteten \sterns{}Homomorphismus
$\hat{\alpha} : C^*(\mbbR) \to \Adjhm{E}$ fortsetzen.
Nach Proposition~\ref{lanceProp25} (Charakterisierung nicht-ausgearteter \sterns{}Ho\-mo\-mor\-phis\-men)
l"a\cms{}t sich $\hat{\alpha}$ zu einem
\sterns{}Homomorphismus $\check{\alpha} : \Multcs{C^*(\mbbR)} \to \Adjhm{E}$ fortsetzen.

Sei $t \in \mbbR$.
Definiere
\[ \lambda_t : L^1(\mbbR) \to L^1(\mbbR), f \mapsto (s \mapsto f(s-t) ). \]
F"ur alle $f,g \in L^1(\mbbR)$ und $s \in \mbbR$ erh"alt man
wegen der Translationsinvarianz des Lebesgue"=Ma"ses:\dremark{$(*)$: Translation um $s-t$}
\dremark{$\mbbR$ ist eine lokalkompakte Gruppe mit einem Haar-Ma"s.}%
\begin{align*}
   \bigl( (\lambda_t f) * g \bigr)(s)
&= \int_\mbbR (\lambda_t f)(y) g(s-y) \,dy
=  \int_\mbbR f(y-t) g(s-y) \,dy
\dremarkm{=  \int_\mbbR f(t-y) g(y-s) \,dy}  \\
&\dremarkm{\overset{(*)}{=} \int_\mbbR f(s-x) g(x-t) \,dx}
=  \int_\mbbR f(s-x) (\lambda_t g)(x) \,dx
=  \bigl( f * (\lambda_t g) \bigr)(s).
\end{align*}
\dremark{Somit ist $(\lambda_t, \lambda_t)$ ein Bizentralisator von $L^1(\mbbR)$.}%
Daraus folgt: $\lambda_t \in \Multcs{L^1(\mbbR)}$.
Au"serdem gilt: $\norm{\lambda_t} = 1$.

\dremark{Da $L^1(\mbbR)$ eine \sterns{}Banachalgebra mit beschr"ankter approximativer Eins ist,}%
Weiterhin gibt es eine norm-verkleinernde, lineare Injektion
von $\Multcs{L^1(\mbbR)}$ nach
$\Multcs{C^*(\mbbR)} \dremarkm{ \cong \Multcs{C_0(\mbbR)} \cong C_b(\mbbR) }$
(vgl. \cite{FollandAbstractHarmonicAna}, Proposition 4.17).\dremark{Z. 17.2.'09/2}
Indem man ein beliebiges $f \in C^*(\mbbR)$ durch Elemente aus $L^1(\mbbR)$ approximiert,
definiert man $\hat{\lambda}_t \in \Multcs{C^*(\mbbR)}$ mit Hilfe von $\lambda_t$.
\dremark{Geht nach T. Timmermann}\dremark{Z. 12.8.'08/3}%

F"ur alle $f \in C_c(\mbbR)$ und $x \in E$ erh"alt man mit der Translationsinvarianz
des Lebesgue"=Ma"ses:\dremark{$(*)$}
\begin{align*}
   \check{\alpha}(\hat{\lambda}_t) \alpha(f) x
&\overset{\dremarkm{\text{def}}}{=}  \alpha(\lambda_t f) x
\dremarkm{=  \int_\mbbR (\lambda_t f)(s) U_s x\, ds}
=  \int_\mbbR f(s-t) U_s x\, ds  \\
&\overset{\dremarkm{(*)}}{=}  \int_\mbbR f(r) U_{r+t} x\, dr
\dremarkm{=  U_t \int_\mbbR f(r) U_{r} x \, dr}
=  U_t \alpha(f) x.
\end{align*}
Da $C_c(\mbbR)$ dicht in $C^*(\mbbR)$ liegt, folgt:
\begin{equation}\label{eqalphaltEqUt2}
\check{\alpha}(\hat{\lambda}_t) = U_t.
\end{equation}

Da es eine Einbettung von $C_0(\mbbR)$ in $\Multcs{C_0(\mbbR)} \cong C_b(\mbbR)$ gibt,
kann man $\tilde{\mathscr{F}} : C^*(\mbbR) \to \Multcs{C_0(\mbbR)}$ betrachten.
Nach Proposition~\ref{lanceProp25}
(Charakterisierung nicht-ausgearteter \sterns{}Ho\-mo\-mor\-phis\-men)
ist $\tilde{\mathscr{F}}$ nicht-ausgeartet,
\dremark{Zeige: \ref{lanceProp25}.(d).
  Da $\tilde{\mathscr{F}}$ ein Isomorphismus ist, gen"ugt es zu zeigen:
  $e_\lambda f \to f$ f"ur alle $f \in C_0(\mbbR)$.
  Dies gilt aber nach der Definition einer approximativen Eins.}%
kann also zu einem \sterns{}Homomorphismus
$\check{\mathscr{F}} : \Multcs{C^*(\mbbR)} \to \Multcs{C_0(\mbbR)}$ fortgesetzt werden.
\dremark{Evtl. noch ein weiteres Diagramm notieren.}%
\[ \xymatrix{
& & \Adjhm{E}  \\
\\
\Multcs{C^*(\mbbR)} \ar[uurr]^{\check{\alpha}} \ar@(dr,dl)[rrrr]_{\check{\mathscr{F}}} &
\,C^*(\mbbR)\,\, \ar[uur]^{\hat{\alpha}} \ar@{_{(}->}[l] \ar@{>->>}[rr]^{\hat{\mathscr{F}}} &&
C_0(\mbbR) \ar[uul]_\beta \ar@{^{(}->}[r] &
\Multcs{C_0(\mbbR)} \ar[uull]_{\check{\beta} = \varphi_h\restring_{C_b(\mbbR)}}
} \]

Es gilt:
\dremark{Siehe Z. 12.8.'08, S. 3--4.}%
\begin{equation}\label{eqFlambdatEqet}
\check{\mathscr{F}}(\hat{\lambda}_t) = e_t.
\end{equation}

Da $\hat{\alpha}$ nicht-ausgeartet und $\hat{\mathscr{F}}^{-1}$ bijektiv ist,
ist $\beta := \hat{\alpha} \circ \hat{\mathscr{F}}^{-1} : C_0(\mbbR) \to \Adjhm{E}$
ein nicht-ausgearteter \sterns{}Homomorphismus,
den man nach Proposition~\ref{lanceProp25} zu einem \sterns{}Ho\-mo\-mor\-phis\-mus
$\check{\beta} : \Multcs{C_0(\mbbR)} \to \Adjhm{E}$ fortsetzen kann.
Sei $\iota : \mbbR \to \mbbC, \lambda \mapsto \lambda$.
Nach \refb{\cite{Woronowicz91}, Example~2,}{Beispiel~\ref{bspRegOpC0Omega}} gilt
$\Multwor{C_0(\mbbR)} \cong C(\mbbR)$,
also $\iota \in \Multwor{C_0(\mbbR)} \cong C(\mbbR)$.
Somit wird nach \refb{\cite{LanceHmod}, Proposition~10.7,}{Proposition~\ref{LanceProp10:7}}
durch $h := \beta(\iota)$
ein regul"arer Operator auf $E$ definiert,
und es gilt\dremark{nach \ref{LanceProp10:7}} $\beta(\iota^*) = \beta(\iota)^*$,
also ist $h$ selbstadjungiert.
\dremark{$\beta(\iota) = \beta(\iota^*) = \beta(\iota)^*$}%
Weil $\varphi_h\restring_{C_0(\mbbR)} : C_0(\mbbR) \to \Adjhm{E}$
\dremark{siehe Beweis von \cite{LanceHmod}, Th. 10.9}%
und $\beta$ nicht-ausgeartete \sterns{}Homomorphismen sind mit
$\varphi_h(\iota) = h = \beta(\iota)$,
\dremark{$\beta(\iota) \in \Multwor{E}$}%
folgt mit \cite{Kustermans97FunctionalCalcRegOp}, Result 6.3:
$\varphi_h\restring_{C_0(\mbbR)} = \beta$.
Da ein nicht-ausgearteter \sterns{}Homomorphismus eine eindeutige
Fortsetzung auf die Multiplikatoralgebra besitzt
(\refb{\cite{LanceHmod}, Proposition 2.1,}{Proposition~\ref{lanceProp25}}), gilt
$\varphi_h\restring_{C_b(\mbbR)} = \check{\beta}$
und $\check{\beta} \circ \check{\mathscr{F}} = \check{\alpha}$.
Zusammen mit \eqref{eqFlambdatEqet} und \eqref{eqalphaltEqUt2} ergibt sich:
\[ \exp(\mri th)
=  \varphi_h(e_t)
=  \check{\beta}(e_t)
\overset{\dremarkm{\eqref{eqFlambdatEqet}}}{=}
   \check{\beta}\left(\check{\mathscr{F}}(\hat{\lambda}_t)\right)
=  \check{\alpha}(\hat{\lambda}_t)
\overset{\dremarkm{\eqref{eqalphaltEqUt2}}}{=}  U_t. \]

\bewitemph{(ii):}
Mit \cite{RaeburnWilliams98MoritaEq}, Proposition~2.50, folgt,
\dass $\overline{\psi_h} : C_b(\sigma(h)) \to \Adjhm{E}$ ein injektiver \sterns{}Homomorphismus ist,
\dremark{denn nach \cite{MurphyCsalgebras}, Theorem 3.1.8, ist $C_0(\sigma(h))$ ein
  wesentliches Ideal in $\Multcs{C_0(\sigma(h))} = C_b(\sigma(h))$.}%
also eine Isometrie.\dremark{\cite{DoranBelfiCAlg}, Cor. 24.4}
Zusammen mit Proposition~\ref{varphiTeqPsiT}\dremark{und \eqref{eqvarphih1EqId}} erh"alt man:
\dremark{Beachte: $(e_t - \mathbbm{1})\restring_{\sigma(h)} \in C_b(\sigma(h))$}%
\begin{align*}
   \normlr{ U_t - \Id_E}
&\overset{\dremarkm{\eqref{eqvarphih1EqId}}}{=}
   \normlr{\varphi_h(e_t - \mathbbm{1})}
\dremarkm{ \\ &\overset{\dremarkm{\ref{varphiTeqPsiT}}}{=}
   \normlr{\psi_h\left( (e_t - \mathbbm{1})\restring_{\sigma(h)}\right)}}
=  \normlr{\overline{\psi_h}\left( (e_t - \mathbbm{1})\restring_{\sigma(h)}\right)}
=  \normlr{(e_t - \mathbbm{1})\restring_{\sigma(h)}}_\infty.
\end{align*}
Da $(U_t)_\indtGr$ normstetig ist, findet man ein $\delta \in \mbbR_{>0}$
derart, \dass gilt:
\[ \forall t \in\,\, ]-\delta,\delta[\,\, \forall \lambda \in \sigma(h) :
   \abs{\exp(\mri t\lambda) - 1} < \frac{1}{2}. \]
Somit ist $\sigma(h)$ beschr"ankt.
\dremark{Denn: Annahme: Nicht.
  O.B.d.A. sei $\sigma(h)$ nach oben unbeschr"ankt.
  Dann findet man ein $\lambda \in \sigma(h)$ mit: $\lambda \geq \frac{2\pi}{\delta}$.
  Mit $t=\frac{\pi}{\lambda} < \frac{\delta}{2}$ folgt: $\mre^{\mri t \lambda} = -1$.}%
\dremark{also kompakt}%
Daher ist $\iota : \sigma(h) \to \mbbC, \lambda \mapsto \lambda$, beschr"ankt,
also folgt mit \eqref{eqpsiTfinAdj}:
$h = \psi_h(\iota) \in \Adjhm{E}$.
\dremark{Siehe vor \ref{varphiTeqPsiT} oder \cite{Kustermans97FunctionalCalcRegOp}, nach Not. 3.5}%
\end{proof}

Wie im klassischen Satz von Stone f"ur Hilbertr"aume ist im obigen Satz $\mri h$
der Erzeuger der $C_0$-Gruppe $(U_t)_\indtGr$, die aus unit"aren Elementen besteht:

\begin{satz}\label{erzeugerUnitaereGrWor}
Sei $E$ ein Hilbert-\cstern{}Modul und $(U_t)_{t \in \mbbR}$ eine $C_0$-Gruppe
auf $E$ mit $U_t \in \Adjhm{E}$ unit"ar f"ur alle $\indtGr$.
Sei $h \in \Multwor{E}$ selbstadjungiert so,
\dass $U_t = \varphi_h(e_t) = \exp(\mri th)$ f"ur alle $t \in \mbbR$ gilt.
Dann ist $\mri h$ der Erzeuger der $C_0$-Gruppe $(U_t)_\indtGr$,
\dass hei"st, $\mri h = H$, wobei der Operator $H$ auf der Menge
\[ D(H) = \left\{ x \in E \setfdg \lim_{t \to 0} \frac{U_t x - x}{t} \text{ existiert} \right\}
\quad\text{durch}\quad
H x = \lim_{t \to 0} \frac{U_t x - x}{t} \]
f"ur alle $x \in E$ definiert wird.
\dremark{Verallgemeinerung auf $h$ normal m"oglich?}%
\dremark{Siehe auch \cite{WoronowiczNapi92OpThInCsAlg}, Th. 2.3}%
\dliter{Verallgemeinerung von \cite{HollevoetEtAl92StonesTh}, Prop. 2.2,
  auf Hilbert-\cstern{}Moduln}%
\end{satz}

Im Beweis verwenden wir das folgende

\begin{lemma}\label{dtfMinusfTo0}
Sei $t \in \mbbR\setminus\{0\}$.
Definiere \[ d_t : \mbbR \to \mbbC, \lambda \mapsto
\begin{cases}
\frac{\exp(\mri t\lambda) - 1}{\mri t\lambda}, & \text{falls } \lambda \neq 0 \\
1 & \text{sonst}
\end{cases}. \]
Dann gilt:
\begin{enumaufz}
\item $d_t \in C_0(\mbbR)$.
\item F"ur alle $f \in C_0(\mbbR)$ hat man:
\[ \norm{ d_s f - f }_\infty \to 0 \qquad\text{f"ur } s \to 0. \]
\end{enumaufz}
\end{lemma}

\begin{proof}
\bewitemph{(i):}
Es gilt
\[ d_t(\lambda) = \frac{\sum_{n=1}^\infty \frac{1}{n!} (\mri t\lambda)^n}{\mri t\lambda}
= \sum_{n=1}^\infty \frac{1}{n!} (\mri t\lambda)^{n-1}
\dremarkm{= 1 + \frac{1}{2} \mri t\lambda + \frac{1}{6} (\mri t\lambda)^2 + \dots} \]
f"ur alle $\lambda \in \mbbR\setminus\{0\}$,
also folgt: $\lim_{\lambda \to 0} d_t(\lambda) = 1$.
Somit ist $d_t$ stetig.
Da $\exp(\mri t \Id_\mbbR) - \mathbbm{1}$ beschr"ankt ist,
\dremark{und $\lim_{\lambda \to \pm \infty} \abs{\mri t\lambda} = \infty$}%
erh"alt man: $d_t \in C_0(\mbbR)$.
\smallskip

\bewitemph{(ii):}
Sei $f \in C_0(\mbbR) \setminus \{0\}$ und $g := d_1 - \mathbbm{1} \in C_b(\mbbR)$.
\dremark{Also
\[ g : \mbbR \to \mbbC, \mu \mapsto
\begin{cases}
\frac{\exp(\mri\mu) - 1}{\mri\mu} - 1, & \text{falls } \mu \neq 0, \\
0 & \text{sonst}
\end{cases}. \]}%
Es gilt\dremark{$\mu = t\lambda$,
  $(d_t - \mathbbm{1})(\lambda) = (d_1 - \mathbbm{1})(t\lambda) = (d_1 - \mathbbm{1})(\mu)$}%
\[ \norm{ d_t f - f }_\infty
=  \sup_{\lambda \in \mbbR} \abs{ ((d_t - \mathbbm{1}) f)(\lambda) }
=  \sup_{\mu \in \mbbR} \absBig{ g(\mu) f\left(\frac{\mu}{t}\right) }, \]
daher wird im folgenden der letzte Ausdruck untersucht.

\dremark{Beweisidee: Zeige: $g$ ist stetig, beschr"ankt und $\lim_{\mu \to 0} g(\mu) = 0$.
  Sei $\varepsilon > 0$. W"ahle $t$ so klein, \dass $f(\frac{\cdot}{t})$ in die Mulde um $0$
  von $g$ pa\cms{}t.}%

Sei $\varepsilon \in \mbbR_{>0}$.
Dann findet man ein $a \in \mbbR_{>0}$ mit
$\norm{f\restring_{\mbbR\setminus[-a,a]}} < \frac{\varepsilon}{\norm{g}}$.
Da $g$ stetig ist mit $\lim_{\mu \to 0} g(\mu) \dremarkm{= \lim d_1(\mu) - 1} = 0$,
findet man ein $b \in \mbbR_{>0}$ mit
$\norm{g\restring_{[-b,b]}} \leq \frac{\varepsilon}{\norm{f}}$.
Setze $t_0 := \frac{b}{a}$.
Sei $\mu \in \mbbR$ und $t \in\, ]0,t_0[$.

\bewitemph{1. Fall:} $\mu \in [-b,b]$.
Man erh"alt: $\abs{g(\mu) f(\frac{\mu}{t})} \leq \dremarkm{\frac{\varepsilon}{\norm{f}} \norm{f} =} \varepsilon$.

\bewitemph{2. Fall:} $\mu \notin [-b,b]$.
Wegen $\abs{\frac{\mu}{t}} \geq \abs{\frac{\mu}{t_0}} = a \frac{\abs{\mu}}{b} \geq a$ gilt:
$\abs{g(\mu) f(\frac{\mu}{t})}  \leq \dremarkm{\norm{g} \frac{\varepsilon}{\norm{g}} =} \varepsilon$.
\dremark{29.4.'08/2}%
\end{proof}

\begin{proof}[Beweis von Satz~\ref{erzeugerUnitaereGrWor}]
\newcommand{\cmdn}{c_n}%
\bewitemph{\glqq$\mri h \subseteq H$\grqq:}
Sei $x \in D(h)$.
Definiere
\[ g : \mbbR \to \mbbR, \lambda \mapsto (1+\lambda^2)^{-1/2}, \]
also $g \in C_0(\mbbR)$.
Nach \cite{LanceHmod}, S. 107, gilt $\dremarkm{q_h := }\varphi_h(g) = (1 + h^*h)^{-1/2} \in \Adjhm{E}$
und $\Bild(\varphi_h(g)) = D(h)$.
Somit findet man ein $y \in E$ mit $x = \varphi_h(g)y$.
Es gilt:
\begin{equation}\label{eqStonesThUtxx}
\begin{split}
   \frac{U_t x - x}{\mri t}
&=\dremarkm{ \frac{\varphi_h(e_t) x - x}{\mri t}
=}  \frac{\varphi_h(e_t) \varphi_h(g)y - \varphi_h(g)y}{\mri t}  \dremarkm{\\
&= \frac{\varphi_h((e_t - \mathbbm{1})g)y}{\mri t}}
=  \varphi_h\left( \frac{(e_t - \mathbbm{1})g}{\mri t} \right) y
\end{split}
\end{equation}
f"ur alle $t \in \mbbR$.
Mit Lemma~\ref{dtfMinusfTo0} erh"alt man
\begin{align*}
&\phantom{=}\,\,\, \sup_{\lambda \in \mbbR} \abslr{ \left(\frac{\exp(\mri t\lambda) - 1}{\mri t} - \lambda\right)
                               \frac{f(\lambda)}{\sqrt{1+\lambda^2}} }  \\
&=  \sup_{\lambda \in \mbbR\setminus\{0\}} \abslr{ \left(\frac{\exp(\mri t\lambda) - 1}{\mri t\lambda} - 1\right)
                               f(\lambda)} \frac{\abs{\lambda}}{\sqrt{1+\lambda^2}}  \\
&\leq  \dremarkm{\sup_{\lambda \in \mbbR\setminus\{0\}} \absBig{ \Bigl(\frac{\exp(\mri t\lambda) - 1}{\mri t\lambda} - 1\Bigr)
                               f(\lambda)} \\
&=}  \sup_{\lambda \in \mbbR\setminus\{0\}}  \abslr{(d_t f - f)(\lambda)}
\dremarkm{= \norm{d_t f - f}_\infty}
\overset{t \to 0}{\longrightarrow} 0
\end{align*}
f"ur alle $f \in C_0(\mbbR)$, also gilt:
\begin{equation}\label{eqStonesThetKonv}
\frac{(e_t - \mathbbm{1})g}{\mri t} \overset{t \to 0}{\longrightarrow} \Id_\mbbR \cdot g
\qquad\text{strikt}.
\end{equation}
\dremark{denn $(\frac{e_t - \mathbbm{1}}{\mri t} - \Id_\mbbR) f g \to 0$}%
Nach dem Funktionalkalk"ul f"ur regul"are Operatoren (Satz~\ref{regOpFkalkuel})
ist $\varphi_h : C(\mbbR) \to \Multwor{E}$ ein \sterns{}Homomorphismus.
Ferner ist $\varphi_h\restring_{C_0(\mbbR)}$ nicht-ausgeartet
(siehe den Beweis von \cite{LanceHmod}, Theorem 10.9).
Au"serdem ist $\varphi_h$ nach Proposition~\ref{lanceProp25}
(Charakterisierung nicht-ausgearteter \sterns{}Homomorphismen) strikt stetig
auf $\overline{\cmkug}_{C_b(\mbbR)}(0,1)$.
\dremark{Nach \eqref{phiTEqphiRestrCb} gilt:
  $\overline{\varphi_h} = \varphi_h\restring_{C_b(\mbbR)}$.
  Nach \cite{LanceHmod}, Prop. 2.5, ist $\overline{\varphi_h} = \varphi_h\restring_{C_b(\mbbR)}$
  strikt stetig auf $\overline{\cmkug}_{C_b(0,1)}(0,1)$.}%
Weiter gilt f"ur alle $t \in \mbbR\setminus\{0\}$:
\dremark{Die k"urzeste Verbindung zwischen $\mre^{\mri r}$ und $1$ ist die Gerade,
  also $\abs{\mre^{\mri r} - 1} \leq \abs{r}$.}%
\begin{align*}
   \sup_{\lambda \in \mbbR} \abslr{ \left( \frac{(e_t - \mathbbm{1})g}{\mri t} \right) (\lambda) }
&= \sup_{\lambda \in \mbbR} \abslr{ \frac{\exp(\mri t\lambda) - 1}{\mri t \sqrt{1+\lambda^2}} }  \\
&\leq \sup_{\lambda \in \mbbR \setminus \{0\}} \abslr{ \frac{\exp(\mri t\lambda) - 1}{t\lambda}}  \\
&= \sup_{r \in \mbbR \setminus \{0\}} \frac{ \abs{\exp(\mri r) - 1} }{\abs{r}}
\leq 1  .
\end{align*}
Also folgt: $\frac{(e_t - \mathbbm{1})g}{\mri t} \in \overline{\cmkug}_{C_b(0,1)}(0,1)$.
Man erh"alt mit \eqref{eqStonesThUtxx} und \eqref{eqStonesThetKonv}:
\dremark{$(*)$: $\varphi_h$ ist strikt stetig auf $\overline{\cmkug}_{C_b(0,1)}(0,1)$.}%
\begin{align*}
   -\mri Hx
&= \lim_{t \to 0} \frac{U_t x - x}{\mri t}
\dremarkm{\overset{\dremarkm{\eqref{eqStonesThUtxx}}}{=}  \lim_{t \to 0} \varphi_h\left( \frac{(e_t-\mathbbm{1})g}{\mri t} \right)y}
\overset{\dremarkm{(*)}}{=}  \varphi_h\left( \lim_{t \to 0} \frac{(e_t-\mathbbm{1})g}{\mri t} \right)y  \\
&\overset{\dremarkm{\eqref{eqStonesThetKonv}}}{=}  \varphi_h(\Id_\mbbR \cdot g)y
=  \varphi_h(\Id_\mbbR) \varphi_h(g)y
=  hx.
\end{align*}

\bewitemph{\glqq$\mri h \supseteq H$\grqq:}
Sei $x \in D(H)$.
Sei $(f_m)_{m \in \mbbN}$ eine approximative Eins in $C_0(\mbbR)$ so,
\dass $f_m \in C_c(\mbbR)$ f"ur jedes $m \in \mbbN$ ist.
Nach Proposition~\ref{lanceProp25} gilt:
\dremark{da $\varphi_h\restring_{C_0(\mbbR)}$ nicht-ausgeartet}%
\begin{equation}\label{eqStonesThPhifmStrikt}
\varphi_h(f_m) \overset{m \to \infty}{\longrightarrow} \Id_E \qquad\text{strikt}.
\end{equation}
F"ur alle $n \in \mbbN$ definiere $\cmdn := d_{1/n} \in C_0(\mbbR)$,
wobei $d_{1/n}$ in Lemma~\ref{dtfMinusfTo0} definiert wird.
Somit folgt f"ur alle $n,m \in \mbbN$:
\begin{equation}\label{eqStonesThDefxnm}
   x_{n,m}
:= \varphi_h(\cmdn f_m)x
=  \varphi_h(\cmdn) \varphi_h(f_m) x
\overset{m \to \infty}{\longrightarrow}  \dremarkm{\varphi_h(\cmdn) \Id_E x
=} \varphi_h(\cmdn) x
=: x_n.
\end{equation}
Es gilt nach Satz~\ref{regOpFkalkuel}:\dremark{(Funktionalkalk"ul f"ur regul"are Operatoren)}
$\varphi_h(\Id_\mbbR) = h$.
Nach Beispiel~\ref{bspRegOpC0Omega} ist $C(\mbbR) \cong \Multwor{C_0(\mbbR)}$ unter
der Abbildung $f \mapsto M_f : D(M_f) \subseteq C_0(\mbbR) \to C_0(\mbbR), g \mapsto fg$,
also entspricht $\Id_\mbbR$ gerade $M_{\Id_\mbbR}$.

Mit \refb{\cite{LanceHmod}, Proposition~10.7,}{Proposition~\ref{LanceProp10:7}} folgt:
\dremark{\cite{EngelNagelSemigroups}, Beweis von Prop. I.4.2.(i)}%
\dremark{$T=\Id_\mbbR$,$\alpha = \varphi_h$}
\[ \forall \alpha \in C_c(\mbbR) \subseteq D(M_{\Id_\mbbR}):
   \Bild(\varphi_h(\alpha)) \subseteq D(\varphi_h(\Id_\mbbR)) = D(h). \]
Wegen $\cmdn f_m \in C_c(\mbbR)$
\dremark{$\supp(\cmdn f_m)$ ist als abgeschlossene Teilmenge der
  kompakten Menge $\supp(f_m)$ kompakt.}%
ergibt sich somit: $x_{n,m} = \varphi_h(\cmdn f_m) x \in D(h)$
f"ur alle $n,m \in \mbbN$.

Mit \eqref{eqStonesThPhifmStrikt} erh"alt man:
\begin{align*}
   h x_{n,m}
&= \varphi_h(\Id_\mbbR) \varphi_h(\cmdn f_m)x
\dremarkm{=  \varphi_h(\Id_\mbbR \cdot \cmdn f_m)x}  \\
&= \varphi_h(\Id_\mbbR \cdot \cmdn) \varphi_h(f_m) x
\overset{\dremarkm{\eqref{eqStonesThPhifmStrikt}}}{\longrightarrow}  \dremarkm{\varphi_h(\Id_\mbbR \cdot \cmdn) \Id_E x  \\
&=}  \varphi_h(\cmdn \cdot \Id_\mbbR)x \qquad\text{f"ur } m \to \infty.
\end{align*}
Weil $h$ abgeschlossen ist, folgt mit \eqref{eqStonesThDefxnm}:
$x_n \in D(h)$ und $h x_n = \varphi_h(\cmdn \cdot \Id_\mbbR)x$.

Da nach Lemma~\ref{dtfMinusfTo0} gilt $\cmdn \overset{n \to \infty}{\longrightarrow} \mathbbm{1}$ strikt
und da $\varphi_h$ nach Proposition~\ref{lanceProp25}
(Charakterisierung nicht-ausgearteter \sterns{}Homomorphismen)
strikt stetig auf $\overline{\cmkug}_{C_b(\mbbR)}(0,1)$ ist,
\dremark{\cite{LanceHmod}, Prop. 2.5}%
erh"alt man mit \eqref{eqvarphih1EqId}:\dremark{$\norm{\cmdn} \leq 1$, siehe oben}
\[ x_n = \varphi_h(\cmdn)x  \overset{n \to \infty}{\longrightarrow}
   \varphi_h(\mathbbm{1})x
\overset{\dremarkm{\eqref{eqvarphih1EqId}}}{=}  x. \]
Es ergibt sich
\dremark{($*$): F"ur $\lambda \neq 0$ gilt
$ (\Id_\mbbR \cdot \cmdn)(\lambda)
= \frac{\exp(\mri \lambda/n)-1}{\mri \lambda/n} \cdot \lambda
= \frac{\exp(\mri \lambda/n)-1}{\mri/n}$ und
$ (\Id_\mbbR \cdot \cmdn)(0) = 0$.}%
\begin{align*}
   h x_n
&= \varphi_h(\cmdn \cdot \Id_\mbbR)x
\overset{\dremarkm{(*)}}{=}  \varphi_h\left( \frac{e_{1/n} - \mathbbm{1}}{\mri/n} \right)x  \\
&= \frac{\left(\varphi_h(e_{1/n}) \varphi_h(\mathbbm{1}) - \varphi_h(\mathbbm{1})\right)x}{\mri/n}
\dremarkm{=  \frac{\varphi_h(e_{1/n}) \Id_E x - \Id_E x}{\mri/n}}
=  \frac{U_{1/n} x - x}{\mri/n}.
\end{align*}
Wegen $x \in D(H)$ gilt:
\[ \lim_{n \to \infty} h x_n  =  \lim_{n \to \infty} \frac{U_{1/n} x - x}{\mri/n}  =  -\mri Hx. \]
Da $h$ abgeschlossen ist, folgt: $x \in D(h)$ und $\mri hx = Hx$.
\end{proof}

\dremark{
Mit Satz~\ref{erzeugerUnitaereGrWor} erh"alt man:

\begin{satz}
Sei $E$ ein Hilbert-\cstern{}Modul "uber $\mfrakA$.
Sei $T \in \Multwor{E}$ schiefadjungiert.
Dann gilt $T \in \Multlor{E}$.
\dremark{Evtl. als Motivation f"ur die Definition von $\Multlco{X}$ anf"uhren.}%
\dmarginpar{pr}\dremark{Evtl. folgt hier noch mehr.}%
\end{satz}

\begin{proof}
Es ist $h := \mri T$ selbstadjungiert.
Somit wird durch $U_t := \exp(\mri th)$ eine unit"are $C_0$-Gruppe in $\Adjhm{E}$
definiert, deren Erzeuger als $C_0$-Gruppe nach \ref{erzeugerUnitaereGrWor} gleich $T$ ist.
Somit folgt: $T \in \Multlor{E}$.
\end{proof}
}%

\section{Zusammenhang zwischen regul"aren Operatoren auf $E$ und auf $E \oplus E$}

\dremww{
\begin{lemma}\label{dichtInHSModEF}
Seien $E, F$ Hilbert-\cstern{}Moduln "uber $\mfrakA$,
sei $A \subseteq E$ und $B \subseteq F$.
Dann ist $A \oplus B$ genau dann dicht in $E \oplus F$,
wenn $A$ dicht in $E$ und $B$ dicht in $F$ liegt.
\end{lemma}

\begin{proof}
\bewitemph{\glqq$\Rightarrow$\grqq:}
Sei $x_0 \in E$.
Dann findet man eine Folge $(x_n,y_n)_{n \in \mbbN}$ in $A \oplus B$
mit $(x_n,y_n) \to (x_0,0)$.
Es folgt $\norm{x_n - x_0} \leq \norm{ (x_n,y_n) - (x_0,0) } \to 0$,
\dremark{\cite{RaeburnWilliams98MoritaEq}, Example 2.14}%
also ist $A$ dicht in $E$.
Analog folgt: $B$ ist dicht in $F$.

\bewitemph{\glqq$\Leftarrow$\grqq:}
Sei $(x_0, y_0) \in E \oplus F$.
Man findet Folgen $x \in A^\mbbN$ und $y \in B^\mbbN$
mit $x_n \to x_0$ und $y_n \to y_0$.
Es folgt
$  \norm{(x_n, y_n) - (x_0, y_0)}^2
\leq \norm{x_n - x_0}^2 + \norm{y_n - y_0}^2
\to 0$,
\dremark{\cite{RaeburnWilliams98MoritaEq}, Example 2.14}%
also ist $A \oplus B$ dicht in $E \oplus F$.
\end{proof}
}%

Analog zur direkten Summe von Hilbertr"aumen wird eine direkte Summe von \hcsmoduln definiert:

\begin{defBemerkung}[Vgl. \cite{RaeburnWilliams98MoritaEq}, Example 2.14]\label{defDirSumHCSMod}
Seien $E$, $F$ \hcsmoduln "uber $\mfrakA$.
\begin{enumaufz}
\item Auf dem Vektorraum $E \oplus F$ wird durch
\[ \skalpri{(x_1,y_1)}{(x_2,y_2)}{E \oplus F} := \skalpri{x_1}{x_2}{E} + \skalpri{y_1}{y_2}{F} \]
f"ur alle $(x_1,y_1), (x_2,y_2) \in E \oplus F$
ein $\mfrakA$-wertiges inneres Produkt derart definiert,
\dass $E \oplus F$ zu einem \hcsmodul "uber $\mfrakA$ wird,
welchen wir ebenfalls mit $E \oplus F$ bezeichnen.
\index[S]{EgrossplusF@$E \oplus F$}%

\item Es gilt f"ur alle $x \in E$ und $y \in F$:
\begin{equation}\label{eqNormHCsModEPlusE}
     \max\{ \normi{x}{E}, \normi{y}{F} \}
\leq \normi{(x,y)}{E \oplus F}
\leq \sqrt{ \normi{x}{E}^2 + \normi{y}{F}^2 }.
\end{equation}
\end{enumaufz}
\end{defBemerkung}

\begin{notation}\label{defMatABCDHCsMod}
Seien $E$, $F$ \hcsmoduln "uber $\mfrakA$.
Seien $A_1 : D(A_1) \subseteq E \to E$, $A_2 : D(A_2) \subseteq F \to E$,
$A_3 : D(A_3) \subseteq E \to F$ und $A_4 : D(A_4) \subseteq F \to F$ Operatoren.
Definiere
\index[S]{ABCD@$\binom{A_1 \,\,\, A_2}{A_3 \,\,\, A_4}$}%
\begin{align*}
 \cmpmatrix{A_1 & A_2 \\ A_3 & A_4} :
&\,\,(D(A_1) \cap D(A_3)) \cmtimes (D(A_2) \cap D(A_4)) \subseteq E \oplus F \to E \oplus F,  \\
&\,\,(x,y) \mapsto (A_1 x + A_2 y, A_3 x + A_4 y).
\end{align*}
\end{notation}

Seien $E$, $F$ \hcsmoduln "uber $\mfrakA$.
Wir zeigen in diesem Abschnitt,
\dass man zu jedem regul"aren Operator $T \in \Multwor{E,F}$
einen selbstadjungierten regul"aren Operator angeben kann,
n"amlich $\hat{T} := \binom{0 \,\,\, T}{T^* \,\,\, 0}$.
Umgekehrt ist, falls $\hat{T}$ regul"ar ist, auch $T$ regul"ar.
\skiptext

F"ur den Beweis der folgenden beiden Propositionen formulieren wir das
folgende

\begin{lemma}\label{NullTS0sternEq}
Seien $E$, $F$ \hcsmoduln "uber $\mfrakA$.
Seien $T : D(T) \subseteq E \to F$ und $S : D(S) \subseteq F \to E$ dicht definiert
und $\mfrakA$-linear.
Dann gilt: \[ \cmpmatrix{0 & T \\ S & 0}^* = \cmpmatrix{0 & S^* \\ T^* & 0}. \]
\end{lemma}

\begin{proof}
Offensichtlich ist $\hat{T} := \cmpmatrix{0 & T \\ S & 0} \dremarkm{F \oplus E \to F \oplus E}$
dicht definiert\dremark{\ref{dichtInHSModEF}} und $\mfrakA$-linear.

Mit dem Index $1$ (bzw. $2$) bezeichnen wir die erste (bzw. zweite)
Komponente eines Tupels.
F"ur ein beliebiges $x \in F \oplus E$ gilt beispielsweise: $x = (x_1, x_2)$.

\bewitemph{\glqq$\subseteq$\grqq:}
Es gilt f"ur alle $x \in D(\hat{T})$ und $y \in D(\hat{T}^*)$
\begin{align*}
   \skalpri{T x_2}{y_1}{F} + \skalpri{S x_1}{y_2}{E}
&= \skalpri{\hat{T}x}{y}{F \oplus E}
=  \skalpri{x}{\hat{T}^*y}{F \oplus E}  \\
&= \skalpri{x_1}{(\hat{T}^*y)_1}{F} + \skalpri{x_2}{(\hat{T}^*y)_2}{E},
\end{align*}
also $\skalpr{S x_1}{y_2} = \skalpr{x_1}{(\hat{T}^*y)_1}$.
Somit ist $\skalpr{x_1}{(\hat{T}^*y)_1}$ unabh"angig von $y_1$
f"ur alle $x \in D(\hat{T}), y \in D(\hat{T}^*)$,
also gilt $\skalpr{x_1}{(\hat{T}^*y)_1} = \skalpr{x_1}{A y_2}$,
wobei $A : D(A) \subseteq E \to F$ sei.
Da $\skalpri{\cdot}{\cdot\cdot}{E}$ nicht-ausgeartet ist, folgt: $(\hat{T}^*y)_1 = A y_2$
f"ur alle $y \in D(\hat{T}^*)$.
Analog erh"alt man: $(\hat{T}^*y)_2 = B y_1$ f"ur alle $y \in D(\hat{T}^*)$,
wobei $B : D(B) \subseteq F \to E$ sei.
Somit hat man $A \subseteq S^*$ und $B \subseteq T^*$,
\dremark{denn $\skalpr{S x_1}{y} = \skalpr{x_1}{A y_2}$}%
also $\hat{T}^* \subseteq \cmpmatrix{0 & S^* \\ T^* & 0}$.

\bewitemph{\glqq$\supseteq$\grqq:}
Es gilt f"ur alle $x \in D(\hat{T})$ und $y \in D(T^*) \cmtimes D(S^*)$
\begin{align*}
   \skalpr{\hat{T}x}{y}
&= \skalpr{T x_2}{y_1} + \skalpr{S x_1}{y_2}  \\
&= \skalpr{x_2}{T^* y_1} + \skalpr{x_1}{S^* y_2}
=  \skalpr{x}{\cmpmatrix{0 & S^* \\ T^* & 0}y},
\end{align*}
also $\cmpmatrix{0 & S^* \\ T^* & 0} \subseteq \hat{T}^*$.
\dremark{27.6.'08/2}%
\end{proof}

\begin{bemerkung}\label{NullTTs0Selbstadj}
Seien $E$, $F$ \hcsmoduln "uber $\mfrakA$.
Dann ist f"ur alle $T \in \Multwor{E,F}$ der Operator
$\begin{pmatrix} 0 & T \\ T^* & 0 \end{pmatrix}$ regul"ar und selbstadjungiert.
\dliter{Aus Beweis von Lemme 2.1 in \cite{Hilsum89FonctorialiteEnKth}, dort ohne Beweis}%
\end{bemerkung}

\begin{proof}
Setze $\hat{T} := \begin{pmatrix} 0 & T \\ T^* & 0 \end{pmatrix}$.
Da $T$ und $T^*$ dicht definiert sind, ist wegen
$D(\hat{T}) = D(T^*) \cmtimes D(T)$ auch $\hat{T}$ dicht definiert.\dremark{\ref{dichtInHSModEF}}

Es gilt $T^{**} = T$ (\refb{\cite{LanceHmod}, Corollary 9.4,}{Proposition~\ref{TssEqT}}).
Damit und mit Lemma~\ref{NullTS0sternEq} folgt:
\[ \hat{T}^*
\overset{\dremarkm{\ref{NullTS0sternEq}}}{=}  \cmpmatrix{0 & T^{**} \\ T^* & 0}
\dremarkm{=  \cmpmatrix{0 & T \\ T^* & 0}}
=  \hat{T}. \]

Da $T^*$ abgeschlossen ist,\dremark{\cite{LanceHmod}, S. 95 unten}
ist auch $\hat{T}$ abgeschlossen.
\dremark{Sei $(x_n, y_n)_{n \in \mbbN} \in D(\hat{T})^\mbbN$
  konvergent gegen $(x_0, y_0) \in E \oplus E$ und
  $\hat{T}(x_n, y_n)$ konvergent gegen $(\xi_0, \eta_0) \in E \oplus E$.
  Zu zeigen: $(x_0,y_0) \in D(\hat{T})$ und $\hat{T}(x_0,y_0) = (\xi_0,\eta_0)$.
  Es gilt $\max \{ \norm{x_n - x_0}, \norm{y_n - y_0} \} \leq \norm{ (x_n - x_0, y_n - y_0) }$,
  also folgt $x_n \to x_0$ und $y_n \to y_0$.
  Wegen $\hat{T}(x_n,y_n) = (T y_n, T^* x_n)$ erh"alt man
  $T y_n \to \xi_0$ und $T^* x_n \to \eta_0$.
  Da $T$ und $T^*$ abgeschlossen sind, folgt $x_0 \in D(T^*)$, $y_0 \in D(T)$ und
  $T^* x_0 = \eta_0$, $T y_0 = \xi_0$.
  Man erh"alt $(x_0,y_0) \in D(T^*) \oplus D(T) = D(\hat{T})$ und
  $\hat{T}(x_0,y_0) = (T y_0, T^* x_0) = (\xi_0, \eta_0)$.}%

Es bleibt zu zeigen: $1 + \hat{T}^* \hat{T}$ hat dichtes Bild.
Es gilt:
\begin{align*}
   (1 + \hat{T}^* \hat{T})(x,y)
\dremarkm{=  (x,y) + \hat{T} (T y, T^* x)}
=  \bigl((1+TT^*)x, (1+T^*T)y\bigr)
\end{align*}
f"ur alle $(x,y) \in D(1 + \hat{T}^* \hat{T})$.
Nach \cite{LanceHmod}, Corollary 9.6, ist $T^*$ regul"ar.
Somit folgt, \dass
\[ \Bild(1+\hat{T}^*\hat{T}) = \Bild\bigl(1+(T^*)^* T^*\bigr) \cmtimes \Bild(1+T^*T) \]
dicht in $F \oplus E$ liegt.\dremark{\ref{dichtInHSModEF}}%
\dremark{Z. 28.4.'08/2}%
\end{proof}

\begin{bemerkung}\label{ThatRegFolgtTReg}
Seien $E$, $F$ \hcsmoduln "uber $\mfrakA$.
Sei $T : D(T) \subseteq E \to F$ dicht definiert und $\mfrakA$-linear.
Ist $\cmpmatrix{0 & T \\ T^* & 0}$ regul"ar,
so ist $T$ regul"ar.
\dremark{Man mu\cms{} dicht definiert und $\mfrakA$-linear voraussetzen,
  damit $T^*$ existiert.}%
\end{bemerkung}

\begin{proof}
Offensichtlich ist $T^*$ dicht definiert,
\dremark{$D(\hat{T}) = D(T^*) \cmtimes D(T)$, \ref{dichtInHSModEF}}%
und $T$ ist abgeschlossen.
\dremark{abgeschlossen: Analog zum Beweis von \ref{zshgRegOpUnbeschrMult}, $(\ddagger)$.}%
Nach Lemma~\ref{NullTS0sternEq} gilt
$\hat{T} := \cmpmatrix{0 & T \\ T^* & 0} = \cmpmatrix{0 & T^{**} \\ T^* & 0}$,
somit erh"alt man:
\begin{align*}
   (1 + \hat{T}^*\hat{T})(x,y)
&= \dremarkm{(x,y) + \hat{T}^* \hat{T}(x,y)
=}  (x,y) + \hat{T}^*(Ty,T^*x)  \\
&\overset{\dremarkm{\ref{NullTS0sternEq}}}{=} (x,y) + (T^{**}T^*x,T^*Ty)
=  \bigl((1 + T^{**}T^*)x, (1 + T^*T)y\bigr)
\end{align*}
f"ur alle $(x,y) \in D(1 + \hat{T}^*\hat{T})$.
Wegen $\norm{y} \leq \norm{(x,y)}$ f"ur alle $(x,y) \in F \oplus E$ folgt\dremark{mit \ref{dichtInHSModEF}}
aus der Dichtheit von $\Bild(1 + \hat{T}^*\hat{T})$ in $F \oplus E$
die Dichtheit von $\Bild(1+T^*T)$ in $E$.
\dremark{4.6.'08/1}%
\end{proof}

\dremww{
\begin{definitn}
Sei $\mfrakA$ eine \csalgebra, $T$ ein dicht definierter, linearer Operator auf $\mfrakA$.
Dann hei"st $T$ \defemph{(\cstern{}algebraisch) affiliiert} zu $\mfrakA$, \cmiZ $T \eta \mfrakA$,
falls ein $z_T \in M(\mfrakA)$ mit $\norm{z_T} \leq 1$ so existiert, \dass gilt:
\[ x \in D(T)
\iff  \exists a \in \mfrakA : x = (1-z_T^* z_T)^{1/2} a  \cap  Tx = z_T a. \]
Man nennt $z_T$ die $z$-Transformierte von $T$.
\dliter{\cite{Woronowicz91}, s. auch \cite{WoronowiczNapi92OpThInCsAlg},
  \cite{Woronowicz95CsalgGenByUnbEl};
  die Theorie stammt von S. Baaj
  (s. \cite{Baaj80Multiplicateurs}, \cite{BaajJulg83TheorieBivariante})
  eine "aquivalente Definition steht in \cite{LanceHmod}.}%
\end{definitn}

\begin{anmerkung}
\begin{enumaufz}
\item Nach \cite{Webster04UnboundedOp} gibt es einen Zusammenhang zwischen den unbeschr"ankten Multiplikatoren
nach E. C. Lance (bzw. nach Woronowicz) und den Multiplikatoren des Pedersen Ideals einer \csalgebra.
Dort wird au"serdem eine weitere M"oglichkeit vorgestellt, mit Hilfe der $q$-Stetigkeit
unbeschr"ankte Multiplikatoren zu definieren.
\dremark{Angucken!}

\item Literatur:
\cite{LanceHmod} (auf S. 107 steht, \dass in §10 bei die Hauptergebnisse von der Artikel
  \cite{Woronowicz91} und \cite{WoronowiczNapi92OpThInCsAlg} vorgestellt werden),
  \cite{Woronowicz91}, \cite{WoronowiczNapi92OpThInCsAlg}, \cite{Woronowicz95CsalgGenByUnbEl}.
  Evtl. n"utzlich f"ur Beispiele (unklar): \cite{Turowska02ComplexityOfDescription}.

\item In \cite{Trout00GradedKTh}, §3, wird eine Art umgekehrter Funktionalkalk"ul f"ur
selbstadjungierte, regul"are Operatoren vorgestellt.

\item In \cite{Pal99RegOpOnHilbertCsMod} wird eine Verallgemeinerung von regul"aren Operatoren
auf Hilbert-\cstern{}Moduln vorgestellt, sogenannte semiregul"are Operatoren.
Es wird gezeigt, \dass nicht jeder abgeschlossene semiregul"are Operator regul"ar ist (Prop. 2.2).
Dies ist allerdings in kommutativen \csalgebren (Prop. 4.1) und in der
kompakten Operatoren "uber einem komplexen, separablen Hilbertraum (Prop. 5.1) der Fall.
Au"serdem wird am Anfang von Abschnitt 3 gezeigt,
\dass es gen"ugt, regul"are Operatoren auf \csalgebren zu untersuchen.

\item Literatur zu Quantengruppen:
  Nach Thomas Timmermann: \cite{Woronowicz91}, \cite{WoronowiczNapi92OpThInCsAlg}.
  Au"serdem n"utzlich: \cite{Woronowicz95CsalgGenByUnbEl}, mit vielen Beispielen.
  Etwas findet man bei \cite{LanceHmod}.
  Einen groben "Uberblick "uber Quantengruppen liefert Thomas Timmermann in seinem
  Vortrag im \cstern{}Oberseminar am 29.1.'08.

\item In \cite{Kustermans97FunctionalCalcRegOp}, 12.18 und 12.19, wird die
  Summe und das Produkt von zwei kommutierenden normalen regul"aren Operatoren definiert.
\end{enumaufz}
\end{anmerkung}
}%

\section{Multiplikatoren des Pedersen-Ideals}

\dremark{Aus \cite{LazarTaylorMultPedId}}

Im vorherigen Abschnitt haben wir mit den regul"aren Operatoren eine M"oglichkeit
kennengelernt, auf einer \csalgebra unbeschr"ankte Multiplikatoren zu definieren.
In diesem Abschnitt erinnern wir an einen anderen Zugang,
n"amlich an Multiplikatoren des Pedersen-Ideals.
Es wird sich herausstellen, \dass die regul"aren Operatoren die Multiplikatoren
des Pedersen-Ideals umfassen.
\skiptext

Eine grundlegende Philosophie beim Studium von \csalgebren{} ist, \dass
diese nicht"=kommutative Analoga von $C_0(\Omega)$ sind,
wobei $\Omega$ ein lokalkompakter Hausdorffraum $\Omega$ sei.
Die Multiplikatoralgebra $\Multcs{\mfrakA}$ einer \csalgebra{} $\mfrakA$ entspricht bei dieser Sichtweise
der Algebra $C_b(\Omega)$ der stetigen, beschr"ankten Funktionen auf $\Omega$.
Eine nat"urliche Frage ist dann: Was ist das nicht"=kommutative Analogon
zu der Algebra $C(\Omega)$ der stetigen Funktionen auf~$\Omega$?

Falls $\mfrakA$ unital ist, dann entspricht dieses der Situation, \dass $\Omega$ kompakt ist,
und daher erwartet man, \dass das Analogon von $C(\Omega)$ die Algebra $\mfrakA$ selbst ist.
Daher liegt das eigentliche Interesse in nicht-unitalen \csalgebren.

Im nicht-unitalen Fall betrachtet man die zugeh"orige Algebra von Multiplikatoren.
Die Multiplikatoralgebra von $C_c(\Omega)$,
den stetigen Funktionen auf $\Omega$ mit kompaktem Tr"ager,
ist in der klassischen Situation $C(\Omega)$.
Ferner ist $C_c(\Omega)$ ein minimales dichtes Ideal von $C(\Omega)$,
und in der Situation von \csalgebren{} existiert dieses minimale dichte beidseitige Ideal
und hei"st Pedersen-Ideal von~$\mfrakA$ (\cite{Pedersen66MeasureTh}, \cite{PedersenCAlg}).

A. J. Lazar und D. C. Taylor studieren in \cite{LazarTaylorMultPedId}
die Multiplikatoralgebra (Algebra der Bizentralisatoren) dieses Ideals,
also unbeschr"ankte Multiplikatoren.
N.~C.~Phillips (\cite{Phillips88MultPedId}) vereinfacht ihre
Arbeit, indem er feststellt, \dass die Multiplikatoralgebra
des Pedersen-Ideals eine Pro-\csalgebra{} ist, also ein inverser Limes von \csalgebren.

\begin{definitn}
\begin{enumaufz}
\item Ein \defemphi{Ordnungsideal}\dremark{(engl.\ \glqq order ideal\grqq{} oder \glqq face\grqq)}
  von $\mfrakA$ ist ein Unterkegel $J$ von $\mPosEl{\mfrakA}$ derart,
  \dass gilt: $\forall x \in \mPosEl{\mfrakA} \, \forall y \in J: x \leq y  \Rightarrow  x \in J$.
\dremark{Evtl. Unterkegel definieren.}%
\item Ein Ordnungsideal $J$ von $\mfrakA$ hei"st \defemph{invariant}, falls $a^*Ja \subseteq J$
  f"ur alle $a \in \mfrakA$ gilt.
\index[B]{invariantes Ordnungsideal}%
\item Eine \sterns{}Unteralgebra $\mfrakB$ von $\mfrakA$ hei"st \defemph{hereditär},
\index[B]{hereditäre \sterns{}Unteralgebra}%
\dremark{(engl.\ \glqq hereditary\grqq{}, \glqq order related\grqq{} oder \glqq facial\grqq)}%
  falls $\mfrakB^+ := \mfrakB \cap \mPosEl{\mfrakA}$ ein Ordnungsideal und
  $\mfrakB$ der Vektorraumaufspann von~$\mfrakB^+$ ist.
\end{enumaufz}
\end{definitn}

\begin{satzDefinition}[\cite{LazarTaylorMultPedId}, Theorem 2.2]
Sei $\mathfrak{J}$ die Menge aller dichten, invarianten Ordnungsideale von $\mPosEl{\mfrakA}$
und $K_\mfrakA^+ := \bigcap \mathfrak{J}$.
\begin{enumaufz}
\item $K_\mfrakA^+$ ist ein dichtes, invariantes Ordnungsideal von $\mPosEl{\mfrakA}$.
\item Der Vektorraumaufspann $K_\mfrakA$ von $K_\mfrakA^+$ ist ein dichtes,
  heredit"ares beidseitiges Ideal von $\mfrakA$,
  welches minimal unter allen dichten beidseitigen Idealen von $\mfrakA$ ist
  und das \defemph{Pedersen-Ideal} von $\mfrakA$ genannt wird.
\index[S]{KA@$K_\mfrakA$ (Pedersen-Ideal)}%
\dremark{Da $K_\mfrakA$ der Vektorraumaufspann von $K_\mfrakA^+$ ist,
  ist $K_\mfrakA$ ein \sterns{}Ideal.}%
\dremark{(Die Dichtheit von $K_\mfrakA$ wird sp"ater bewiesen.)}%
\dremark{D.h.: F"ur jedes dichte, beidseitige Ideal $I$ von $\mfrakA$ gilt:
  $K_\mfrakA \subseteq I$}%
\dliter{vgl. \cite{PedersenCAlg}, Th. 5.6.1}%
\end{enumaufz}
\end{satzDefinition}

In \cite{Ara01MoritaEqAndPedersenId}, Definition~1.2,
wird eine Verallgemeinerung des Pedersen"=Ideals
auf Hilbert-\cstern{}Moduln $E$ "uber $\mfrakA$ definiert.
Genauer wird $P_E = E K_J$ Pedersen"=Untermodul von $E$ genannt,
wobei $J := \overline{\skalpr{E}{E}}$ ein Ideal in $\mfrakA$ ist.
Es ist aber nicht offensichtlich, wie man diese Definition auf Operatorr"aume verallgemeinern kann,
da man beispielsweise in Operatorr"aumen den Begriff des Ideals nicht zur Verf"ugung hat.
\dmarginpar{zutun}\dremark{Genauer angeben, was nicht geht.}%
\skiptext

\dremark{
\begin{satz}\label{lazarTaylor2.3}
Sei $A$ eine \csalgebra.
Setze
\begin{align*}
K^+_{00}(A) &:= \{ x \in \mPosEl{A} \setfdg \exists y \in \mPosEl{A} : xy = x \} \quad\text{und} \\
J_0      &:= \{ x \in \mPosEl{A} \setfdg \exists n \in \mbbN \, \exists y_1, \dots y_n \in K^+_{00}(A) :
  x \leq \sum_{i=1}^n y_i \}.
\end{align*}
Dann gilt: $J_0 = K_A^+$.
\end{satz}

\begin{lemma}
Sei $\mbbR_{\geq 0}$ mit der von $\mbbR$ induzierten Topologie versehen.
Sei $K \subseteq \mbbR_{\geq 0}$ kompakt.
Dann gilt: $\exists a,b \in \mbbR_{>0}: K \subseteq\, ]a,b[$.
\end{lemma}

\begin{proof}
Beachte: Die Menge $U := \{ ]\frac{1}{n},n[ \setfdg n \in \mbbN \}$ ist eine offene "Uberdeckung von $K$.
Somit findet man eine endliche Teil"uberdeckung $T \subseteq U$ von $K$.
\end{proof}

\begin{satz}\label{PedIdCharMitFunkionalkalkuel}
Sei $A$ eine \csalgebra.
Dann gilt: $K^+_{00}(A) = \{ f(x) \setfdg x \in \mPosEl{A}, f \in \mPosEl{C_c(]0,\infty[, \mbbC)} \}$.
\dliter{\cite{PedersenCAlg}, 5.6.1}%
\end{satz}

\begin{proof}
Setze $M := \{ f(x) \setfdg x \in \mPosEl{A}, f \in \mPosEl{C_c(]0,\infty[, \mbbC)} \}$.

\bewitemph{\glqq$\subseteq$\grqq:} ?

\bewitemph{\glqq$\supseteq$\grqq:}
Sei $z \in M$.
Dann findet man ein $x \in \mPosEl{A}$ mit $z = f(x)$.
Es ist $\supp(f)$ kompakt, also beschr"ankt, denn
$\{ ]\frac{1}{n},n[ \setfdg n \in \mbbN \}$ ist eine offene "Uberdeckung von $\supp(f)$.
Somit findet man $a,b \in \mbbR_{>0}$ mit $\supp(f) \subseteq \,]a,b[$.
Definiere $\hat{f} := f \cup \{(0,0)\}$.
Dann ist $\hat{f}$ stetig mit $\supp(\hat{f}) = \supp(f)$.
Weil f"ur alle $r \in \mbbR_{>0}$ die Menge $\{ [r,r+\frac{1}{n}] \setfdg n \in \mbbN \}$
eine kompakte Umgebungsbasis von $r$ ist, ist $\mbbR_{>0}$ ein lokalkompakter Hausdorffraum.
Nach \cite{PedersenAnalysisNow}, 1.7.5, findet man ein $g \in C_c(\mbbR_{>0}, [0,1])$ so, \dass
$g\restring_{\supp(f)} \equiv 1$ und $\supp(g) \subseteq \,]a,b[$.
Insbesondere ist $g \geq 0$.
Setze $\hat{g} := g \cup \{(0,0)\}$.
Es gilt $\hat{f} \cdot \hat{g} = \hat{f}$.
Setze $\tilde{f} := \hat{f}\restring_{\sigma'(x)}$, $\tilde{g} := \hat{g}\restring_{\sigma'(x)}$,
Sei $B$ die von $x$ erzeugte \cstern{}Unteralgebra von $A$.
Sei $C_0(\sigma'(x)) := \{f \in C(\sigma'(x)) \setfdg f(0) = 0 \}$.
Dann gibt es einen \sterns{}Isomorphismus $\Phi$ von $C_0(\sigma'(x))$ auf $B$ mit $\Phi(\Id_{\sigma'(x)}) = x$.
\dremark{\cite{SunderFuncAna}, 3.3.10}%
Ferner gilt: $z = f(x) = \Phi(\tilde{f})$.
Setze $y := \Phi(\tilde{g}) \in B$.
Dann gilt:
\[ zy = \Phi(\tilde{f}) \Phi(\tilde{g}) = \Phi(\tilde{f}\tilde{g}) = \Phi(\tilde{f}) = z. \]
Da $\Phi$ ein \sterns{}Isomorphismus ist, gilt
$\sigma'(y) = \sigma'_{C_0(\sigma'(x))}(\tilde{g}) \subseteq \mbbR_{\geq 0}$ und
$y = \Phi(\tilde{g}) = \Phi(\tilde{g}^*) = \Phi(\tilde{g})^* = y^*$.
Also ist $y \in \mPosEl{A}$.

Zeige: $z \in \mPosEl{A}$.
Da $\Phi$ ein \sterns{}Homomorphismus ist, folgt:
$z = \Phi(\tilde{f}) = \Phi(\tilde{f}^*) = \Phi(\tilde{f})^* = z^*$.
Da $\Phi$ ein Algebrenisomorphismus ist, gilt ferner:
$\sigma'_A(z) = \sigma'_{C_0(\sigma(x))}(f) \subseteq \mbbR_{\geq 0}$.
Somit gilt: $z \in \mPosEl{A}$.\dremark{3.5.'05}
\end{proof}

\begin{satz}
Sei $A$ eine \csalgebra.
Das Pedersen-Ideal $K_A$ von $A$ liegt dicht in $A$.
\dliter{G. K. Pedersen, Measure theory for \cstern{}algebras, Math. Scand. 19 (1966) 131--145;
        \cite{PedersenCAlg}, 5.6.1}%
\end{satz}

\begin{proof}
(Unvollst"andig!)
F"ur alle $n \in \mbbN$ definiere
\[ \varphi_n : \mbbR \rightarrow \mbbR, t \mapsto
\begin{cases}
0,                 & \text{falls } t \leq \frac{1}{n} \\
2(t-\frac{1}{n}),  & \text{falls } \frac{1}{n} < t \leq 2\frac{1}{n} \\
t                  & \text{sonst}
\end{cases}.
\]
Sei $n \in \mbbN$.
Es ist $\varphi_n$ stetig.
Sei $a \in \mPosEl{A}$ und $B$ die von $a$ erzeugte \cstern{}Unteralgebra von $A$.
Dann findet man nach \cite{SunderFuncAna}, 3.3.10, einen isometrischen \sterns{}Isomorphismus $\Phi$
von $C_0(\sigma'(a))$ auf $B$ so, \dass $\Phi(\Id_{\sigma'(a)}) = a$.
Es folgt:
\[
   \norm{a - \varphi_n(a)}
=  \norm{(\Id_{\sigma'(a)} - \varphi_n)(a)}
=  \norm{ \Phi(\Id_{\sigma'(a)} - \varphi_n) }
=  \norm{ \Id_{\sigma'(a)} - \varphi_n\restring_{\sigma'(a)} }_\infty
\leq \frac{1}{n}.
\]
Es ist noch z.z.: $\varphi_n(a) \in K^+_{00}(A)$.
\end{proof}
}%

Im folgenden werden einige Beispiele f"ur das Pedersen-Ideal angegeben:

\pagebreak
\begin{beispiel}
\begin{enumaufz}
\item Ist $\mfrakA$ eine unitale \csalgebra, so gilt: $K_\mfrakA = \mfrakA$.
\item Sei $\Omega$ ein lokalkompakter Hausdorffraum.
  Dann gilt: $K_{C_0(\Omega)} = C_c(\Omega)$.
\dremark{$C_c(\Omega) := \{ f \in C(\Omega) \setfdg \supp(f) \text{ kompakt} \}$}%
\item Sei $H$ ein Hilbertraum.
  Dann ist $K_{\kptOp(H)}$
  gleich der Menge\dremark{$\mLinStet_f(H)$} der Operatoren auf $H$
  mit endlich-dimensionalem Bild.
\dremark{(endlichem Rang).
  $\mLinStet_f(H)$ liegt dicht in $\kptOp(H)$ bzgl. der Operatornorm, genauer:
  $\kptOp(H)$ ist der Abschlu\cms{} von $\mLinStet_f(H)$ in der Operatornorm
  (\cite{WernerFunkana3}, VI.3.7).}%
\dremark{Die Operatoren auf $H$ mit endlichem Rang, also die Operatoren, die endliche
  Summen von Rang-1-Operatoren sind,
  sind gleich den Operatoren auf $H$ mit endlich-dimensionalem Bild.}%
\end{enumaufz}
\end{beispiel}

\begin{proof}
\bewitemph{(i):} Sei $\mfrakA$ eine unitale \csalgebra.
Da $K_\mfrakA$ ein dichtes beidseitiges Ideal von $\mfrakA$ ist,
hat $K_\mfrakA$ nicht-leeren Schnitt mit den invertierbaren Elementen von $\mfrakA$.
\dremark{denn $\exists \varepsilon \in \mbbR_{>0} :
               \cmkug(e_\mfrakA,\varepsilon) \subseteq G(\mfrakA)$}%
Somit folgt (i).

\bewitemph{(ii):} \cite{LazarTaylorMultPedId}, Example 2.4.

\bewitemph{(iii):} \cite{LazarTaylorMultPedId}, Example 2.5.
\end{proof}

\dremark{
Beweis.
\bewitemph{(i):} Die Bezeichnungen $K^+_{00}$ und $J_0$ stammen aus \ref{lazarTaylor2.3}.
Da $e_A \in \mPosEl{A}$, gilt: $K^+_{00} = \mPosEl{A}$.
Es folgt mit $n=1$ und $y_1=x$: $J_0 = K^+_{00} = \mPosEl{A}$.
Da $K_A$ der Vektorraumaufspann von $K_A^+$ ist und man jedes Element von $A$ als Summe
von vier positiven Elementen schreiben kann\dremark{s. D/B oder Pedersen, C*-Algebras}, folgt: $K_A = A$.

\bewitemph{(ii):} \bewitemph{\glqq$\subseteq$\grqq:}
  Bekanntlich ist $B:=C_c(X)$ ein dichtes, beidseitiges Ideal von $A$.
  Da $K_A$ minimal ist, folgt: $K_A \subseteq C_c(X)$.

\bewitemph{\glqq$\supseteq$\grqq:}
Falls $X$ kompakt ist, ist $A$ unital, also folgt mit (i): $C_c(X) = A = K_A$
Sei daher $X$ nicht kompakt.
Sei $f \in B$.
Sei $\hat{X} = X \cup \{\infty\}$ die Ein-Punkt-Kompaktifizierung von $X$.
Dann existiert $\hat{f} \in C(\hat{X})$ mit $f=\hat{f}\restring_X$ und $\hat{f}(\infty)=0$.
$\hat{X}$ ist als kompakter Hausdorffraum normal.\dremark{Pedersen, Analysis now, 1. ed., 1.6.6}
Somit existiert zu den disjunkten, abgeschlossenen Mengen $\supp(\hat{f})$ ($=\supp(f)$) und $\{\infty\}$
ein $\hat{g} \in C(\hat{X},[0,1])$ mit $\hat{g}\restring_{\supp(\hat{f})} \equiv 1$ und
$\hat{g}(\infty) = 0$.\dremark{Bou., Gen. Top., IX.4.1 Def. 1}
Somit ist $g := \hat{g}\restring_X \in C_0(X,[0,1])$ mit $g \cdot f = f$.\dremark{6.4.'05}

\begin{bemerkung}
Seien $A,B$ \csalgebren.
Sei $\varphi$ ein surjektiver \sterns{}Homomorphismus von $A$ auf $B$.
Dann gilt: $\varphi(K_A) = K_B$.
\end{bemerkung}
}%

\dremww{
\begin{bemerkung}
Sei $J$ eine normierte \sterns{}Algebra mit approximativer Eins.
Sei $(S,T) \in \Gamma(J)$.
\begin{enumaufz}
\item Falls $S$ beschr"ankt ist, dann ist auch $T$ beschr"ankt, und
es gilt: $\norm{S} = \norm{T}$.
\item Ein Multiplikator $(S,T)$ hei"st \defemph{beschr"ankt}, falls $S$ beschr"ankt ist.
\item Die Menge aller beschr"ankten Multiplikatoren von $J$ sei mit $\mathcal{DC}(J)$ bezeichnet
  (von englisch \glqq double centralizer\grqq).\index[S]{DCJ@$\mathcal{DC}(J)$}
\item $\mathcal{DC}(J)$ ist, versehen mit der Norm $\norm{(S,T)} := \norm{S}$,
  eine Banach-\sterns{}Unteralgebra von $\Gamma(J)$
\end{enumaufz}
\end{bemerkung}

Sei $\mfrakA$ eine \csalgebra.
Dann ist $K_\mfrakA$ ein dichtes, beidseitiges Ideal von $\mfrakA$ und besitzt somit eine aufsteigende
approximierende Eins, die beschr"ankt ist durch 1 und aus hermiteschen Elementen von $K_\mfrakA$ besteht.
\dremark{\cite{DoranBelfiCAlg}, 13.1}%

\begin{beispiel}
F"ur Beispiele von $\mathcal{DC}(A)$, wobei $A$ eine \csalgebra{} ist, siehe \ref{bspMultiplikatoralg}.
\end{beispiel}

\begin{bemerkung}
Sei $J$ eine normierte \sterns{}Algebra mit approximativer Eins.
Sei $(S,T) \in \Gamma(J)$.
\begin{enumaufz}
\item $S$ und $T$ sind lineare Abbildungen mit abgeschlossem Graphen.
\item $S(ab) = S(a)b$ und $T(ab) = aT(b)$ f"ur alle $a,b \in J$.
\item Falls $J$ vollst"andig ist, so gilt: $\mathcal{DC}(J) = \Gamma(J)$.
\end{enumaufz}
\end{bemerkung}

\begin{definitn}
Sei $A$ eine Algebra und $a \in A$.
Definiere $L_a : A \rightarrow A, b \mapsto ab$ und $R_a : A \rightarrow A, b \mapsto ba$.
\end{definitn}

\begin{bemerkung}
Sei $J$ eine normierte \sterns{}Algebra mit approximativer Eins.
Sei $\mu_0 : J \rightarrow \mathcal{DC}(J), a \mapsto (L_a, R_a)$.
\begin{enumaufz}
\item $\mu_0$ ist ein \sterns{}Isomorphismus von $J$ nach $\mathcal{DC}(J)$, und
  $\mu_0(J)$ ist ein beidseitiges Ideal von $\Gamma(J)$.
\item $\mu_0$ ist surjektiv auf $\mathcal{DC}(J)$ genau dann, wenn $J$ ein Einselement hat.
\item Ist $J$ kommutativ, so auch $\Gamma(J)$.
\end{enumaufz}
\end{bemerkung}

\begin{satz}
Sei $A$ eine \csalgebra{}.
Dann ist $\mathcal{DC}(A)$ eine \csalgebra.
\end{satz}

\begin{bemerkung}
Seien $A,B$ normierte \sterns{}Algebren mit approximativer Eins.
Sei $\varphi$ ein surjektiver \sterns{}Homomorphismus von $A$ auf $B$.
Dann induziert $\varphi$ einen \sterns{}Homomorphismus $\overline{\varphi}$
von $\Gamma(A)$ nach $\Gamma(B)$.
F"ur alle $(S,T) \in \Gamma(A)$ ist $(U,V) := \overline{\varphi}(S,T)$ definiert durch
\[ U(\varphi(x)) = \varphi(S(x)) \quad\text{und}\quad V(\varphi(x)) = \varphi(T(x)) \]
f"ur alle $x \in A$.
Ferner gilt: $\overline{\varphi}(\mathcal{DC}(A)) \subseteq \mathcal{DC}(B)$.
\end{bemerkung}

Sei $A$ eine \csalgebra.
Da $K_A$ dicht in $A$ liegt, kann man jeden beschr"ankten Multiplikator von $K_A$ eindeutig zu einem
Multiplikator von $A$ fortsetzen. Die Umkehrung gilt ebenfalls:

\begin{bemerkung}
Sei $A$ eine \csalgebra.
Sei $(S,T) \in \mathcal{DC}(A)$ und $K := K_\mfrakA$.
\begin{enumaufz}
\item $(S\restring_K, T\restring_K) \in \Gamma(K)$.
\item Durch $\rho : \mathcal{DC}(A) \rightarrow \mathcal{DC}(K), (S,T) \mapsto (S\restring_K, T\restring_K)$ wird ein
isometrischer \sterns{}Isomorphismus auf $\mathcal{DC}(K)$ definiert.
\end{enumaufz}
\end{bemerkung}
}%

Die Elemente der Multiplikatoralgebra $\cmGamma(K_\mfrakA)$ des Pedersen-Ideals
bezeichnet man auch als unbeschr"ankte Multiplikatoren.
Beispiele hierf"ur sind:

\begin{beispiel}[\cite{LazarTaylorMultPedId}, 4.1 und 4.2]\label{bspMultPedId}
\begin{enumaufz}
\item Sei $\Omega$ ein lokalkompakter Hausdorffraum.
  Dann gilt:  $\cmGamma(K_{C_0(\Omega)}) \cong C(\Omega)$.
\item Sei $H$ ein Hilbertraum.
  Dann kann man $\cmGamma(K_{\kptOp(H)})$ mit $\mLinStet(H)$ identifizieren.
\end{enumaufz}
\end{beispiel}

\dremark{
Beweis.
\bewitemph{(i):} Setze $K := C_c(X)$.
Sei $f \in C(X)$.
Definiere $L_f : K \rightarrow C(X), g \mapsto f \cdot g$.
Sei $g \in K$.
Es gilt $\supp(f \cdot g) = \overline{ \{ x \in X \setfdg (fg)(x) \neq 0 \} } \subseteq \supp(g)$.
Da $\supp(g)$ kompakt ist, ist $\supp(fg)$ als abgeschlossene Teilmenge einer kompakten
Menge kompakt. Es folgt: $\Bild(L_f) \subseteq K$.

Definiere $R_f : K \rightarrow K, g \mapsto g \cdot f$.
Zeige:
\[ T : C(X) \rightarrow \Gamma(K), f \mapsto (L_f, R_f) \]
ist ein Algebrenisomorphismus auf $\Gamma(K)$.

\bewitemph{linear:} $f \mapsto L_f$ und $f \mapsto R_f$ sind linear, also auch $T$.

\bewitemph{injektiv:} Seien $f,g \in C(X)$ mit $f \neq g$.
Somit existiert ein $x \in X$ mit $f(x) \neq g(x)$.
Sei $y \in X \setminus \{x\}$.
Da $X$ hausdorffsch ist, findet man disjunkte, offene Mengen $U,V$ in $X$ mit $x \in U$ und
$y \in V$.
Da $X$ lokalkompakt ist, findet man kompakte Menge $A,B$ in $X$ mit $x \in A \subseteq U$ und
$y \in B \subseteq V$.\dremark{Bou., Gen. Top., bei Def. lokalkompakt}
Nach Pedersen, Analysis Now, 2. ed., 1.7.5 findet man eine Funktion $h \in C(X,[0,1])$ mit
$h\restring_A \equiv 1$ und $\supp(h)$ kompakt und $\supp(h) \subseteq U$.
Also ist $h \in K$.
Es folgt
\[ L_f(h)(x) = f(x)h(x) = f(x) \neq g(x) = L_g(h)(x), \]
also $L_f \neq L_g$, also $T(f) \neq T(g)$.

\bewitemph{surjektiv:} Siehe \cite{Johnson64Centralisers}, S. 604, Theorem 2.

\bewitemph{Algebrenhomomorphismus:}
Es gilt:\dremark{13.4.'05}
\[ T(fg) = (L_{fg}, R_{fg}) = (L_f \circ L_g, R_g \circ R_f)
=  (L_f,R_f)(L_g,R_g) = T(f)T(g).  \qedhere \]
}%

\begin{bemerkung}[Vgl. \cite{LazarTaylorMultPedId}, S. 8]
Man kann $K_\mfrakA$, $\mfrakA$ und $\Multcs{\mfrakA}$ als Unteralgebren von $\cmGamma(K_\mfrakA)$
auffassen.
\dmarginpar{zutun}\dremark{Die Literaturangabe stimmt nicht ganz.}%
\dliter{vgl. \cite{LazarTaylorMultPedId}, S. 8}%
\end{bemerkung}

Wie man in Beispiel~\ref{bspRegOp}.(iv) und \ref{bspMultPedId}.(ii) sieht,
stimmen $\Multwor{\mfrakA}$ und $\cmGamma(K_\mfrakA)$ \iallg{} nicht "uberein.
Vielmehr sind die Multiplikatoren des Pedersen-Ideals
in den regul"aren Operatoren enthalten:

\begin{satz}[\cite{Webster04UnboundedOp}, Theorem 3.1]
\begin{enumaufz}
\item F"ur jedes $a \in \cmGamma(K_\mfrakA)$ existiert ein $T \in \Multwor{\mfrakA}$ mit
  $K_\mfrakA \subseteq D(T)$ und $Tx = ax$ f"ur alle $x \in D(T)$.
\dremark{Gemeint ist: Nach \cite{LazarTaylorMultPedId}, S. 8 unten, existiert eine
  Einbettung $\iota : \mfrakA \hookrightarrow \Gamma(K_\mfrakA)$.
  Also gilt: $\iota(Tx) = a \iota(x)$.}%
\dremark{Ist $T$ eindeutig?}%
\dremark{Ist $D(T)$ ein Ordnungsideal?}%
\item F"ur jedes $T \in \Multwor{\mfrakA}$ mit $K_\mfrakA \subseteq D(T)$
findet man ein $a \in \cmGamma(K_\mfrakA)$ so, \dass gilt: $Tx = ax$ f"ur alle $x \in D(T)$.
\end{enumaufz}
\end{satz}

Auch wenn die Multiplikatoralgebra $\cmGamma(K_\mfrakA)$ des Pedersen-Ideals
\iallg{} nicht alle regul"aren Operatoren enth"alt,
hat sie gegen"uber $\Multwor{\mfrakA}$ den Vorteil, eine \sterns{}Algebra zu sein.

\dremww{
\begin{anmerkung}
\begin{enumaufz}
\item Eine alternative Charakterisierung des Pedersen-Ideals findet man in \cite{Phillips88MultPedId}.
\item In \cite{Ara01MoritaEqAndPedersenId} wird ein Pedersen-Untermodul $P_E$
  f"ur Hilbert-\cstern{}Moduln $E$ definiert,
  welcher eine Verallgemeinerung des Pedersen-Ideals ist.
  Setze $\mfrakB := \kptOphm_\mfrakA(E)$.
  Dann gilt: $\Gamma(K_\mfrakB) = \Adjhmi{\mfrakA}{P_E}$.
  Kann man dies nicht verwenden, um unbeschr"ankte Operatoren in
  Operatorr"aumen zu definieren?
\end{enumaufz}
\end{anmerkung}
}%

\chapter{Operatorr"aume}

In diesem Kapitel wird an verschiedene Definitionen und Aussagen erinnert,
die wir in dieser Arbeit verwenden.
Im ersten Abschnitt werden elementare Begriffe in Operatorr"aumen wiederholt.
Im Anschlu\cms werden Beispiele f"ur Operatorraumstrukturen auf verschiedenen R"aumen aufgef"uhrt.
Unitale Operatorsysteme werden im dritten Abschnitt behandelt.
Im vierten Abschnitt wird auf den Begriff der Injektivit"at f"ur Operatorr"aume eingegangen.
Au"serdem wird an Eigenschaften der injektiven H"ulle eines Operatorraumes erinnert.
Tern"are Ringe von Operatoren und Tripelsysteme werden im darauf\/folgenden Abschnitt eingef"uhrt,
ebenso der Begriff der tern"aren H"ulle eines Operatorraumes.
Im sechsten Abschnitt werden Grundlagen der Theorie der selbstadjungierten Operatorr"aume notiert.
Anschlie"send beweisen wir f"ur solche R"aume die Existenz und Eindeutigkeit
einer injektiven H"ulle.

\dremark{Einige Aussagen k"onnen gestrichen werden. Evtl. straffen}%

\section{Grundlagen}

In diesem Abschnitt wird an Grundlagen der Theorie der Operatorr"aume erinnert,
also insbesondere an elementare Definitionen.
\skiptext

Sei $X$ ein normierter Raum.
Nach dem Satz von Hahn-Banach kann man $X$ isometrisch in $C(\Omega)$ einbetten,
wobei $\Omega = \overline{\cmkug}_{X^*}(0,1)$ ist.
Also wird $X$ als Funktionenraum,
\dheisst als Untervektorraum der stetigen Funktionen auf~$\Omega$, realisiert.
Analog wird ein (konkreter) Operatorraum
als abgeschlossener Untervektorraum von Operatoren auf einem Hilbertraum definiert.

\dremark{Weitere Motivation: \cmvb{e} Abbildungen werden bei Strukturuntersuchungen von
Operatoralgebren verwendet und sind die Morphismen zwischen
Operatorr"aumen.}%

\begin{definitn}
\begin{enumaufz}
\item Als einen (\defemph{konkreten}) \defemph{matrixnormierten Raum} bezeichnet man einen
Untervektorraum $X$ von $\mLinStet(H)$, wobei $H$ ein Hilbertraum sei.
\item Ist $X$ abgeschlossen, so spricht man von einem (\defemph{konkreten}) \defemphi{Operatorraum}.
\dremark{Norm auf $M_n(E)$: Da $M_n(\mLinStet(H)) \cong \mLinStet(H^n)$ gilt,
  kann man $M_n(E)$ mit einer Norm versehen.}%
\dremark{Motivation, warum man Operatorr"aume betrachtet:
  \cite{EffrosRuan00OperatorSpaces}, S. viii, S. 19/20, \cite{WittstockWasSindOpr},
  \cite{Pisier03OpSpaceTheory}, S. 1ff}%
\dremark{Beispiele f"ur Operatorr"aume: Banachr"aume, Hilbertr"aume, Hilbert-\cstern{}Moduln}%
\end{enumaufz}
\end{definitn}

\begin{definitn}
Seien $m,n \in \mbbN$, sei $V$ ein Vektorraum.
Mit $M_{m,n}(V)$ wird der Raum der $m \times n$-Matrizen "uber $V$ notiert.
Setze $M_n(V) := M_{n,n}(V)$, $M_{m,n} := M_{m,n}(\mbbC)$ und $M_n := M_{n,n}$.
\index[S]{MmnV@$M_{m,n}(V)$ ($m \times n$-Matrizen über V)}%
\index[S]{MnV@$M_n(V)$ ($:= M_{n,n}(V)$)}%
\index[S]{Mmn@$M_{m,n}$ ($:= M_{m,n}(\mbbC)$)}%
\index[S]{Mn@$M_n$ ($:= M_{n,n}(\mbbC)$)}%
\end{definitn}

Sei $X \subseteq \mLinStet(H)$ ein Operatorraum.
Seien $m,n \in \mbbN$.
Durch
\[ M_n(X) \subseteq M_n(\mLinStet(H)) \cong \mLinStet(H^n) \]
wird $M_n(X)$ mit einer Norm $\norm{\cdot}_n$ versehen.
Mit Hilfe solcher Normen wird ein abstrakter Operatorraum definiert.

Setze $p := \max\{m,n\}$.
Wir fassen den Raum $M_{m,n}$ als Untervektorraum von $M_p$ auf
und erhalten so eine Norm auf $M_{m,n}$,
die als $\norm{\cdot}_{M_{m,n}}$ oder kurz als $\norm{\cdot}$ notiert wird.

\begin{definitn}\label{defOpraumAbstrakt}
\begin{enumaufz}
\item Unter einer Matrixnorm $\norm{\cdot}$ auf einem $\mbbC$"=Vektorraum $X$
versteht man eine Folge $(\norm{\cdot}_n)_{n \in \mbbN}$ von Normen derart, \dass $\norm{\cdot}_n$
eine Norm auf $M_n(X)$ f"ur jedes $n \in \mbbN$ ist.

\item Ein (\defemph{abstrakter}) \defemphi{matrixnormierter Raum} ist ein $\mbbC$"=Vektorraum $X$,
der so mit einer Matrixnorm $\norm{\cdot}$ versehen ist, \dass gilt:
\begin{enumaufzB}
\item[($N_1$)] $\norm{x \oplus y}_{m+n} = \max\bigl\{\norm{x}_m, \norm{y}_n\bigr\}$,
wobei mit $x \oplus y$ die $(m+n) \times (m+n)$-Matrix
 $\begin{smallpmatrix} x & 0 \\ 0 & y \end{smallpmatrix}$ bezeichnet wird,
 \index[S]{xplusy@$x \oplus y$}%
\item[($N_2$)]
$ \norm{ \alpha x \beta}_n
\leq  \norm{\alpha}_{M_{n,m}} \cdot \norm{x}_m \cdot \norm{\beta}_{M_{m,n}}$
\end{enumaufzB}
f"ur alle $m, n \in \mbbN$, $x \in M_m(X), y \in M_n(X)$,
$\alpha \in M_{n,m}$ und $\beta \in M_{m,n}$.

\item Sei $X$ ein matrixnormierter Raum.
Ist f"ur ein $n \in \mbbN$ (und damit f"ur alle)
$\left( M_n(X), \norm{\cdot}_n \right)$ vollst"andig,
so hei"st $X$ ein (\defemph{abstrakter}) \defemphi{Operatorraum}.
\dremark{In $(N_1)$ gen"ugt es, \glqq$\leq$\grqq{} zu fordern
  (\cite{EffrosRuan00OperatorSpaces}, 2.3.6).
  Anstelle von $(N_1)$ und $(N_2)$ kann man auch verlangen:
  \cite{Pisier03OpSpaceTheory}, S. 35, \cite{PaulsenComplBoundedMapsPS}, Exerc. 13.2}%
\dliter{Begriffe aus \cite{WittstockWasSindOpr}, 2.1}%
\dremark{In Analogie zu den Begriffen normierter Raum und Banachraum}%
\end{enumaufz}
\end{definitn}

Offensichtlich ist jede \csalgebra{} ein Operatorraum.
Wir werden in Abschnitt \ref{BeispieleOR} sehen,
\dass man beispielsweise Banachr"aume, Hilbertr"aume und Hilbert-\cstern{}Moduln
mit Operatorraumstrukturen versehen kann.

\begin{definitn}
Seien $V$, $W$ Vektorr"aume.
Sei $\alpha : V \rightarrow W$ und $n \in \mbbN$.
Definiere
\[ \alpha_n : M_n(V) \rightarrow M_n(W), v \mapsto (\alpha(v_{ij}))_{i,j \in \haken{n}}, \]
genannt $n$-te \defemphi{Amplifikation} von $\alpha$.
Anstelle von $\alpha_n$ schreiben wir auch $\alpha^{(n)}$.
\index[S]{alphan@$\alpha_n$ (Amplifikation)}%
\index[S]{alphan@$\alpha^{(n)}$ (Amplifikation)}%
\end{definitn}

\begin{definitn}
Seien $X$, $Y$ matrixnormierte R"aume, sei $\alpha : X \rightarrow Y$ linear.
\begin{enumaufz}
\item $\alpha$ hei"st \defemph{vollständig kontraktiv}, 
wenn $\alpha_n$ f"ur alle $n \in \mbbN$ kontraktiv ist, \dheisst $\norm{\alpha_n} \leq 1$.
\item $\alpha$ bezeichnet man als \defemph{vollständig isometrisch},
falls $\alpha_n$ f"ur alle $n \in \mbbN$ isometrisch ist.
\index[B]{vollständig kontraktive Abbildung}%
\index[B]{vollständig isometrische Abbildung}%
\item Man nennt $\alpha$ \defemph{vollständig beschränkt},\dremark{(kurz \defemph{v.\,b.})}
  falls gilt:
\index[B]{vollständig beschränkte Abbildung}%
  \[ \norm{\alpha}_\text{cb} := \sup\bigl\{ \norm{\alpha_n} \setfdg n \in \mbbN \bigr\} < \infty \]
  (\glqq cb\grqq{} steht f"ur \glqq completely bounded\grqq).
\index[S]{normcb@$\normcb{\cdot}$}%
\dremark{Dies sind die Morphismen zwischen Operatorr"aumen.}%
\item Die Menge der vollst"andig beschr"ankten Abbildungen von $X$ nach $Y$
  wird mit $\mCB(X,Y)$ bezeichnet.
  Setze $\mCB(X) := \mCB(X,X)$.
\index[S]{CBX@$\mCB(X)$}%
\index[S]{CBXY@$\mCB(X,Y)$}%
\end{enumaufz}
\end{definitn}

Die Operatorr"aume bilden zusammen mit den vollst"andig beschr"ankten
(beziehungsweise vollst"andig kontraktiven) Abbildungen eine Kategorie.
\skiptext

Wegen des folgenden Resultats wird h"aufig nicht zwischen abstrakten und konkreten
matrixnormierten R"aumen unterschieden (vgl. \cite{EffrosRuan00OperatorSpaces}, Theorem 2.3.5):

\begin{satz}[Darstellungssatz f"ur matrixnormierte R"aume]
Sei $X$ ein abstrakter matrixnormierter Raum.
Dann findet man einen Hilbertraum $H$,
einen konkreten matrixnormierten Raum $Y \subseteq \mLinStet(H)$ und
eine surjektive vollst"andige Isometrie $\Phi : X \arrowbij Y$ auf $Y$.
\dremark{Falls $X$ als normierter Raum separabel ist, kann man $H = \ell^2$ w"ahlen.}%
\dliter{\cite{PaulsenComplBoundedMapsPS}, 13.4}%
\dremark{von Ruan, 1988}%
\end{satz}

Wir erinnern an die folgende bekannte Absch"atzung:

\begin{lemma}[\cite{EffrosRuan00OperatorSpaces}, S. 22]
\label{normaLeqNormSumAij}\label{OpraumNormabschaetzung}\label{normxijLeqNormx}
Sei $X$ ein Operatorraum und $x \in M_n(X)$.
F"ur alle $i,j \in \haken{n}$ gilt:
\dremark{Diese Ungl. kann man verbessern, siehe z.B. \cite{MurphyCsalgebras}, 3.4.1,
  \cite{PaulsenComplBoundedMapsPS}, Ex. 3.10.(i)}%
\begin{equation}\label{eqNormAijZuNormA}
\norm{x_{ij}} \leq \norm{x}_n \leq \sum_{k,\ell = 1}^n \norm{x_{k \ell}}.
\end{equation}
\end{lemma}

Aus der Absch"atzung \eqref{eqNormAijZuNormA} folgt,
\dass eine Folge $(x^{(m)})_{m \in \mbbN}$ in $M_n(X)$ genau dann konvergiert,
wenn ihre Eintr"age $(x_{ij}^{(m)})_{m \in \mbbN}$ f"ur alle $i,j \in \haken{n}$ konvergieren.

\dremww{
\begin{proof}
\bewitemph{(i):} Es gilt:
$\norm{a_{ij}} = \norm{e_i^T a e_j} \leq \norm{e_i^T} \cdot \norm{a} \cdot \norm{e_j} = \norm{a}$.

\bewitemph{(ii):}
F"ur alle $i,j \in \haken{n}$ sei $E_{ij} = ( \delta_{ik} \delta_{j\ell} )_{k,\ell \in \haken{n}}$.
Es gilt:
$  \norm{a}
=  \norm{ \sum_{i,j = 1}^n a_{ij} \cdot E_{ij} }
\leq  \sum_{i,j = 1}^n \norm{a_{ij} \cdot E_{ij}}
=  \sum_{i,j = 1}^n \norm{a_{ij}}$.
\end{proof}

\begin{lemma}\label{TnInvEqTInvn}
Seien $X$, $Y$ Operatorr"aume, sei $T : X \to Y$ invertierbar.
Dann gilt: $(T_n)^{-1} = (T^{-1})_n$.
\end{lemma}

\begin{proof}
Es gilt:
$  T_n (T^{-1})_n x
=  (T T^{-1}x_{ij})_{i,j}
=  (x_{ij})_{i,j}
=  (T^{-1} T x_{ij})_{i,j}
=  (T^{-1})_n T_n x$.
\end{proof}

\begin{lemma}
Seien $U$, $V$, $W$ Vektorr"aume,
seien $\alpha : U \rightarrow V$ und $\beta : V \rightarrow W$ linear.
Dann gilt $(\beta \circ \alpha)_n = \beta_n \circ \alpha_n$ f"ur alle $n \in \mbbN$.
\end{lemma}

\begin{proof}
F"ur alle $a \in M_n(U)$ gilt:
\begin{align*}
   (\alpha \circ \beta)_n(a)
&= ((\alpha \circ \beta)(a_{ij}))_{ij}
=  \begin{pmatrix}
  (\alpha \circ \beta)(a_{11}) & \dots & (\alpha \circ \beta)(a_{1n}) \\
  \vdots & \ddots & \vdots \\
  (\alpha \circ \beta)(a_{n1}) & \dots & (\alpha \circ \beta)(a_{nn})
  \end{pmatrix}  \\
&= \alpha_n \Big( \begin{pmatrix}
  \beta(a_{11}) & \dots & \beta(a_{1n}) \\
  \vdots & \ddots & \vdots \\
  \beta(a_{n1}) & \dots & \beta(a_{nn})
  \end{pmatrix}  \Big)
=  (\alpha_n \circ \beta_n)(a).  \qedhere
\end{align*}
\end{proof}

\begin{lemma}
Sei $A$ eine \csalgebra oder gleich $\mLinStet(H,K)$, wobei $H$, $K$ Hilbertr"aume seien.
Sei $X$ ein Unterraum von $A$.
Dann gilt f"ur alle $x_1, \dots, x_n \in X$:
\dliter{\cite{BlecherLeMerdy04OpAlg}, 1.2.5, \cite{Zarikian01Thesis}, Lemma 1.3.1}%
\[ \normlr{\begin{pmatrix} x_1 \\ \vdots \\ x_n \end{pmatrix}}_{C_n(X)}
= \normBig{\sum_{k=1}^n x^*_k x_k}^{\frac{1}{2}}  \quad\text{und}\quad
  \norm{(x_1 \dots  x_n)}_{R_n(X)}
= \normBig{\sum_{k=1}^n x_k x^*_k}^{\frac{1}{2}}. \]
\end{lemma}

\begin{proof}
F"ur alle $x_1, \dots, x_n \in X$ gilt:\dremark{Benutzt wird: $\norm{a}^2 = \norm{a^*a}$}
\begin{align*}
   \normlr{\begin{pmatrix} x_1 \\ \vdots \\ x_n \end{pmatrix}}_{C_n(X)}^2
&= \normlr{\begin{pmatrix}
       x_1 & 0 & \cdots & 0 \\ \vdots & \vdots & \ddots & \vdots \\ x_n & 0 & \cdots & 0
     \end{pmatrix}}_{M_n(X)}^2
=  \normlr{\begin{pmatrix}
       x_1^* & \dots & x_n^* \\ 0 & \dots & 0 \\ \vdots & \ddots & \vdots \\ 0 & \dots & 0
     \end{pmatrix}
     \begin{pmatrix}
       x_1 & 0 & \cdots & 0 \\ \vdots & \vdots & \ddots & \vdots \\ x_n & 0 & \cdots & 0
     \end{pmatrix}} \\
&= \normBig{\sum_{k=1}^n x^*_k x_k}_X.  \qedhere
\end{align*}
\end{proof}

\begin{lemma}\label{unitaeresElErhNorm}
Sei $E$ ein Operatorraum und $x \in M_n(E)$.
Dann gilt f"ur alle unit"aren $\alpha \in M_n$: $\norm{\alpha v} = \norm{v} = \norm{v \alpha}$.
Insbesondere wird durch Vertauschen von Zeilen oder Spalten von $v$ dessen Norm nicht ver"andert.
\dliter{\cite{EffrosRuan00OperatorSpaces}, (2.1.5), S. 21 unten}%
\end{lemma}

\begin{lemma}\label{DdichtFolgtMnDdicht}
Sei $E$ ein Operatorraum, $D$ eine dichte Teilmenge von $E$ und $m,n \in \mbbN$.
Dann ist $M_{m,n}(D)$ dicht in $M_{m,n}(E)$.
\end{lemma}

\begin{proof}
Sei $a \in M_{m,n}(E)$ und $\varepsilon \in \mbbR_{>0}$.
F"ur alle $i \in \haken{m}$, $j \in \haken{n}$ findet man ein $d_{i,j} \in D$ mit
$\norm{a_{i,j} - d_{i,j}} < \frac{\varepsilon}{mn}$.
Mit $b := (d_{i,j})_{i \in \haken{m},j \in \haken{n}} \in M_{m,n}(D)$
und \ref{normaLeqNormSumAij} erh"alt man:\dremark{6.10.'05}%
\[ \norm{a - b}
\overset{\ref{normaLeqNormSumAij}}{\leq}  \sum_{i=1}^m \sum_{j=1}^n \norm{a_{i,j} - d_{i,j}}
<  mn \frac{\varepsilon}{mn}
=  \varepsilon.  \qedhere \]
\end{proof}

\begin{lemma}\label{konvergenzInMnX}
Sei $X$ ein Operatorraum, $n \in \mbbN$.
Sei $x^{(i,j)} \in X^\mbbN$ eine Folge und $x^{(i,j)}_0 \in X$ f"ur alle $i,j \in \haken{n}$.
Dann konvergieren die Folge $x^{(i,j)}$ gegen $x^{(i,j)}_0$ f"ur alle $i,j \in \haken{n}$
genau dann, wenn die Folge $(x^{(i,j)})_{i,j}$ in $M_n(X)$ gegen $(x^{(i,j)}_0)_{i,j}$ konvergiert.
\end{lemma}

\begin{proof}
\bewitemph{\glqq$\Leftarrow$\grqq:}
Die Folge $(x^{(i,j)})_{i,j}$ in $M_n(X)$ konvergiere gegen $(x^{(i,j)}_0)_{i,j}$.
Dann gilt f"ur alle $k, \ell \in \haken{n}$:
\[ \lim_m x^{(k,\ell)}_m
=  \lim \pr_{k,\ell}((x^{(i,j)}_m)_{i,j})
\overset{\ref{prVollstKontraktiv}}{=}
   \pr_{k,\ell}(\lim (x^{(i,j)}_m)_{i,j})
=  \pr_{k,\ell} ((x^{(i,j)}_0)_{i,j})
=  x^{(k,\ell)}_0. \]

\bewitemph{\glqq$\Rightarrow$\grqq:}
Es gelte: Die Folgen $x^{(i,j)}$ konvergieren gegen $x^{(i,j)}_0$ f"ur alle $i,j \in \haken{n}$.
Dann gilt:
\[ \norm{ (x^{(i,j)}_m)_{i,j} - (x^{(i,j)}_0)_{i,j} }
\leq  \sum_{i,j=1}^n  \underbrace{\norm{x^{(i,j)}_m - x^{(i,j)}_0 }}_{\overset{m \rightarrow \infty}{\longrightarrow} 0}
\overset{m \rightarrow \infty}{\longrightarrow} 0. \]
\end{proof}

\begin{lemma}\label{AabgFolgtAnabg}
Seien $X$, $Y$ Operatorr"aume, sei $A : D(A) \subseteq X \to Y$ linear und abgeschlossen.
Dann ist $A_n$ abgeschlossen.
\end{lemma}

\begin{proof}
Sei $(x^{(k)})_{k \in \mbbN} \in D(A_n)^\mbbN$ konvergent gegen ein $x^{(0)} \in M_n(X)$ und
$(A_n x^{(k)})_{k \in \mbbN}$ konvergent gegen $y^{(0)} \in M_n(Y)$.
\dremark{Zu zeigen: $x^{(0)} \in D(A_n)$ und $A_n x^{(0)} = y^{(0)}$.}%
Nach \ref{konvergenzInMnX} konvergiert $(x^{(k)}_{ij})_{k \in \mbbN}$ gegen $x^{(0)}_{ij}$
und $((A_n x^{(k)})_{ij})_{k \in \mbbN} = (A x^{(k)}_{ij})_{k \in \mbbN}$
gegen $(y^{(0)})_{ij}$ f"ur alle $i,j \in \haken{n}$.
Da $A$ abgeschlossen ist, folgt $x^{(0)}_{ij} \in D(A)$ und
$(A_n x^{(0)})_{ij} = (y^{(0)})_{ij}$ f"ur alle $i,j \in \haken{n}$,
also $x^{(0)} \in M_n(D(A)) = D(A_n)$ und $A_n x^{(0)} = y^{(0)}$.
\dremark{3.12.'08/1}%
\end{proof}

\begin{lemma}\label{TstetigFolgtTnSt}
Seien $X$, $Y$ Operatorr"aume, sei $T \in L(X,Y)$.
Dann gilt: $T_n \in L(M_n(X),M_n(Y))$.
\end{lemma}

\begin{proof}
Zeige: $T_n$ ist stetig in $0$.
Sei $\varepsilon \in \mbbR_{>0}$.
Da $T$ stetig in $0$ ist, findet man f"ur alle $i,j \in \haken{n}$
ein $\delta_{i,j} \in \mbbR_{>0}$ mit:
\[ \forall x \in X: \norm{x}<\delta_{i,j} \Rightarrow \norm{Tx}<\frac{\varepsilon}{n^2}. \]
Setze $\delta_0 := \min_{i,j \in \haken{n}} \delta_{i,j}$.
Sei $x \in M_n(X)$ mit $\norm{x} < \delta_0$.
Mit \ref{normxijLeqNormx} gilt: $\norm{x_{ij}} \leq \norm{x} < \delta_0$.
Somit folgt:
\[ \norm{T_n x}
\leq  \sum_{i,j=1}^n \norm{Tx_{ij}}
<  \sum_{i,j=1}^n \frac{\varepsilon}{n^2}
= \varepsilon. \]
\end{proof}

\begin{lemma}\label{TnVertauscht}\label{lem191104}
Seien $X,Y$ Vektorr"aume "uber $\mbbC$.
Sei $T: X \rightarrow Y$ linear.
Sei $\alpha \in M_{n,m}$, $\beta \in M_{m,n}$ und $v \in M_m(X)$.
Dann gilt: $T_n(\alpha v \beta) = \alpha T_m(v) \beta$.
\end{lemma}

\begin{proof}
Es gilt:\dremark{19.11.'04, S.2}
\begin{align*}
   (T_n(\alpha v \beta))_{i,\ell}
=  T(\sum_{j,k=1}^m \alpha_{ij} v_{jk} \beta_{k\ell})
=  \sum_{j,k} T(\alpha_{ij} v_{jk} \beta_{k\ell})
=  \sum_{j,k} \alpha_{ij} T(v_{jk}) \beta_{k\ell}
=  (\alpha \cdot T_m(v) \cdot \beta)_{i,l}.
\end{align*}
\end{proof}

\begin{bemerkung}\label{prVollstKontraktiv}
Sei $X$ Operatorraum, $n \in \mbbN$, seien $i,j \in \haken{n}$.
Dann ist $\pr_{ij} : M_n(X) \rightarrow X, x \mapsto x_{ij}$ vollst"andig kontraktiv.
\end{bemerkung}

\begin{proof}
Folgt aus \ref{OpraumNormabschaetzung}.
\dremark{Alternativer Beweis:
\bewitemph{1. Fall:} $n=1$. \bewitemph{Unterfall (a):} $i=j=1$.
Es gilt:
\begin{gather*}
   \pr_{11}(x)
=  x_{11}
=  \begin{pmatrix} x_{11} & x_{12} \end{pmatrix}
   \begin{pmatrix} 1 \\ 0 \end{pmatrix}
=  \begin{pmatrix} 1 & 0 \end{pmatrix} x
   \begin{pmatrix} 1 \\ 0 \end{pmatrix}
=  e_1^T x \, e_1, \qquad\text{also}  \\
   \norm{\pr_{11}(x)}
\leq \underbrace{\norm{\begin{pmatrix} 1 & 0 \end{pmatrix}}}_{=
        \normlr{\begin{pmatrix} 1 & 0 \\ 0 & 0 \end{pmatrix}}_{M_2} = 1}
   \cdot \norm{x} \cdot
   \normlr{\begin{pmatrix} 1 \\ 0 \end{pmatrix}}
=  \norm{x}.
\end{gather*}

\bewitemph{Unterfall (b):} $i=j=2$. Analog zum Unterfall 1.(a).

\bewitemph{Unterfall (c):} $i=1$, $j=2$. Es gilt:
\[ \pr_{12}(x)
=  x_{12}
=  \begin{pmatrix} x_{11} & x_{12} \end{pmatrix}
   \begin{pmatrix} 0 \\ 1 \end{pmatrix}
=  \begin{pmatrix} 1 & 0 \end{pmatrix} x
   \begin{pmatrix} 0 \\ 1 \end{pmatrix}
=  e_1^T x \, e_2, \]
also $\norm{\pr_{12}(x)} \leq \norm{x}$.

\bewitemph{Unterfall (d):} $i=2$, $j=1$. Analog zum Unterfall 1.(c).

\bewitemph{1. Fall:} $n=2$. \bewitemph{Unterfall (a):} $i=j=1$. Es gilt:
\dremark{(a): $M_2(x)$ ist Operatorraum, benutze $n=1$, 1. Fall;
         (b): \cite{EffrosRuan00OperatorSpaces}, S. 21 unten}%
\begin{align*}
   (\pr_{11})_2(x)
&= \normlr{ \begin{pmatrix} x_{11} & x_{13} \\ x_{31} & x_{33} \end{pmatrix} }  \\
&\overset{\dremarkm{(a)}}{\leq}
   \normlr{ \begin{pmatrix}
   x_{11} & x_{13} & x_{12} & x_{14} \\
   x_{31} & x_{33} & x_{32} & x_{34} \\
   x_{21} & x_{23} & x_{22} & x_{24} \\
   x_{41} & x_{43} & x_{42} & x_{44} \end{pmatrix} }
\overset{\dremarkm{(b)}}{=}
   \normlr{ \begin{pmatrix}
   x_{11} & x_{12} & x_{13} & x_{14} \\
   x_{31} & x_{32} & x_{33} & x_{34} \\
   x_{21} & x_{22} & x_{23} & x_{24} \\
   x_{41} & x_{42} & x_{43} & x_{44} \end{pmatrix} }
\overset{\dremarkm{(b)}}{=}  \norm{x}.
\end{align*}

\bewitemph{Unterfall (b) bis (d):} Analog.

\bewitemph{3. Fall:} $n>2$. Analog.\dremark{6.9.'05, S. 3}}%
\end{proof}

\begin{bemerkung}
Sei $X$ ein Operatorraum, $n \in \mbbN$, seien $i,j \in \haken{n}$.
Die Abbildung $\iota_{ij} : X \to M_n(X), x \mapsto (\delta_{ik}\delta_{j\ell} x)_{k,\ell \in \haken{n}}$
\dremark{Also ist $\iota_{ij}(x)$ die Matrix, die "uberall $0$ ist bis auf den Eintrag $(i,j)$,
  dort steht $x$.}%
ist eine lineare Isometrie.\dremark{Evtl. vollst. isometrisch}%
\end{bemerkung}

\begin{bemerkung}
Sei $E$ ein Operatorraum, $P : E \to E$ eine Projektion.
Dann ist $P_n$ f"ur alle $n \in \mbbN$ eine Projektion,
insbesondere ist $P$ vollst"andig kontraktiv.
\dremark{$P_n^2(x) = P_n((P(x_{ij}))_{ij}) = (P^2(x_{ij}))_{ij} = (P(x_{in}))_{ij} = P_n(x)$}%
\end{bemerkung}
}%

\dremww{
\begin{beispiel}
Sei $X := \mbbC \oplus \mbbC$ mit der Struktur der direkten Summe von \csalgebren versehen,
\dremark{$\norm{\cdot}_\text{max}$}%
$u : X \to X, (a,b) \mapsto (b,a)$.
Dann ist $u$ ein \sterns{}Isomorphismus, also insbesondere vollst"andig kontraktiv.
Weiter gilt:
$\begin{smallpmatrix} u & \Id_X \\ 0 & 0 \end{smallpmatrix} : M_2(X) \to M_2(X)$
ist \emph{nicht} kontraktiv.
\dremark{20.4.'07, S. 1/2}%
\end{beispiel}

\begin{proof}
Setze $y := \begin{smallpmatrix} (1,0) & (0,1) \\ 0 & 0 \end{smallpmatrix}$.
Es gilt:
\begin{align*}
   \norm{y}^2
&= \norm{y^*y}
=  \normlr{ \begin{pmatrix} (1,0) & 0 \\ (0,1) & 0 \end{pmatrix}
            \begin{pmatrix} (1,0) & (0,1) \\ 0 & 0 \end{pmatrix} }  \\
&= \normlr{ \begin{pmatrix} (1,0) & 0 \\ 0 & (0,1) \end{pmatrix} }
=  \max\{ \norm{(1,0)}, \norm{(0,1)} \}
=  1.
\end{align*}
Weiter gilt
\begin{align*}
   \normlr{ \begin{pmatrix} u & \Id_X \\ 0 & 0 \end{pmatrix} y }^2
&= \normlr{ \begin{pmatrix} (0,1) & (0,1) \\ 0 & 0 \end{pmatrix} }^2  \\
&= \normlr{ \begin{pmatrix} (0,1) & 0 \\ (0,1) & 0 \end{pmatrix}
            \begin{pmatrix} (0,1) & (0,1) \\ 0 & 0 \end{pmatrix} }
=  \normlr{ \begin{pmatrix} (0,1) & (0,1) \\ (0,1) & (0,1) \end{pmatrix} }
\geq  \sqrt{2},
\end{align*}
denn mit $\xi := (0,1)$, $\eta:=0$ folgt
\begin{align*}
   \normlr{ \begin{pmatrix} (0,1) & (0,1) \\ (0,1) & (0,1) \end{pmatrix} }^2
&\geq  \normlr{ \begin{pmatrix} (0,1) & (0,1) \\ (0,1) & (0,1) \end{pmatrix} \binom{\xi}{\eta} }^2  \\
&= \normlr{ \binom{(0,\xi_2 + \eta_2)}{(0,\xi_2 + \eta_2)} }^2
=  \normlr{ \binom{(0,1)}{(0,1)} }^2
=  2\norm{(0,1)}^2
=  2.
\end{align*}
\end{proof}
}%

\section{Beispiele}\label{BeispieleOR}

Beispiele f"ur Operatorraumstrukturen auf verschiedenen R"aumen wie
$C_0(\Omega)$, Banachr"aumen, \hcsmoduln und Hilbertr"aumen werden
in diesem Abschnitt angegeben.

\begin{beispiel}\label{EoplusFIstOpraum}
Seien $X,Y$ Operatorr"aume.
F"ur alle $(x,y) \in X \oplus Y$ sei
$\norm{(x,y)}_{X \oplus Y,n} := \max\{ \norm{x}_{X,n}, \norm{y}_{Y,n} \}$.
\dremark{$\norm{(x,y)}_{E \oplus F,n} := \sqrt{ \norm{x}_n^2, \norm{y}_n^2 }$ scheint nicht
  zu funktionieren}%
\dremark{siehe auch \cite{BlecherLeMerdy04OpAlg}, 1.2.17}%
Dann ist $X \oplus Y$, versehen mit der Matrixnorm $(\norm{\cdot}_{X \oplus Y,n})_{n \in \mbbN}$,
  ein Operatorraum.
\end{beispiel}

\dremww{
\begin{beispiel}\label{xxEoplusFIstOpraum}
Seien $X,Y$ matrixnormierte R"aume.
F"ur alle $(x,y) \in X \oplus Y$ sei
$\norm{(x,y)}_{Y \oplus X,n} := \max\{ \norm{x}_{X,n}, \norm{y}_{Y,n} \}$.
\dremark{$\norm{(x,y)}_{E \oplus F,n} := \sqrt{ \norm{x}_n^2, \norm{y}_n^2 }$ scheint nicht
  zu funktionieren}%
\dremark{siehe auch \cite{BlecherLeMerdy04OpAlg}, 1.2.17}%
\begin{enumaufz}
\item Versehen mit der Matrixnorm $(\norm{\cdot}_{X \oplus Y,n})_{n \in \mbbN}$
  ist $X \oplus Y$ ein matrixnormierter Raum.
\item Sind $X$ und $Y$ vollst"andig, so ist auch $Y \oplus Y$ vollst"andig,
  also ein Operatorraum.
\end{enumaufz}
\end{beispiel}

\begin{proof}
Setze $G := E \oplus F$.
Da die anderen Normeigenschaften offensichtlich sind, wird nur die Dreiecksungleichung gezeigt.
Seien $(v,w), (x,y) \in M_n(G)$.
Es gilt:
\begin{align*}
   \norm{ (v,w) + (x,y) }_n
&= \max\{ \norm{v+x}_n, \norm{w+y}_n \}  \\
&\leq \max\{ \norm{v}_n, \norm{w}_n \} + \max\{ \norm{x}_n, \norm{y}_n \}
=  \norm{(v,w)}_n + \norm{(x,y)}_n.
\end{align*}

\bewitemph{($N_1$):}
F"ur alle $(v,w) \in M_m(G)$, $(x,y) \in M_n(G)$ gilt:
\begin{align*}
   \norm{ (v,w) \oplus (x,y) }_{m+n}
&= \max \{ \norm{v \oplus x}_{m+n}, \norm{w \oplus y}_{m+n} \}  \\
&= \max \{ \max \{ \norm{v}_m, \norm{x}_m \}, \max \{ \norm{w}_m, \norm{y}_m \} \}  \\
&\overset{\ref{maxVon2Max}}{=}
   \max \{ \max \{ \norm{v}, \norm{w} \}, \max \{ \norm{x}, \norm{y} \} \}
=  \max \{ \norm{(v,w)}, \norm{(x,y)} \}.
\end{align*}

\bewitemph{($N_2$):}
F"ur alle $\alpha \in M_{n,m}$, $(x,y) \in M_m(G)$ und $\beta \in M_{m,n}$ gilt:
\dremark{22.7.'05, S. 3}%
\begin{align*}
   \norm{\alpha (x,y) \beta}_n
&= \norm{ (\alpha x \beta, \alpha y \beta) }_n
=  \max \{ \norm{\alpha x \beta}, \norm{\alpha y \beta} \}  \\
&\leq \max \{ \norm{\alpha} \norm{x} \norm{\beta}, \norm{\alpha} \norm{y} \norm{\beta} \}
=  \norm{\alpha} \cdot \norm{(x,y)} \cdot \norm{\beta}.  \qedhere
\end{align*}
\end{proof}
}%

\begin{beispiel}\label{MnEIstOpraum}
Sei $X$ ein Operatorraum mit Matrixnorm $(\norm{\cdot})_{n \in \mbbN}$,
und sei $k \in \mbbN$.
Da man $M_m(M_k(X))$ mit $M_{mk}(X)$ f"ur alle $m \in \mbbN$ identifizieren kann,
erh"alt man so eine Norm auf $M_m(M_k((X))$.
Versehen mit dieser Norm ist
$(M_k(X), (\norm{\cdot}_{mk})_{m \in \mbbN})$ ein Operatorraum.
\end{beispiel}

\dremww{
Beweis.
Da $E$ ein Operatorraum ist, ist $E$ isomorph zu einem konkreten Operatorraum
$\hat{E} \leq \mLinStet(H)$.
Dann gilt: $M_k(E) \cong M_k(\hat{E}) \leq M_k(\mLinStet(H)) \cong \mLinStet(H^k)$.
}%

\begin{defBemerkung}
Sei $X$ ein Operatorraum, seien $m,n \in \mbbN$.
Setze $p := \max\{m,n\}$.
\begin{enumaufz}
\item Die Norm auf den Matrixr"aumen $M_{m,n}(X)$ erh"alt man, indem man
$M_{m,n}(X)$ als Unterraum von $M_p(X)$ auf\/fa\cms{}t.
Durch den linearen Isomorphismus $M_p(M_{m,n}(X)) \cong M_{pm,pn}(X)$ erh"alt man
eine Operatorraumstruktur auf $M_{m,n}(X)$.

\item Mit $C_n(X) := M_{n,1}(X)$ wird der \defemphi{Spaltenoperatorraum} "uber $X$ bezeichnet
(\glqq{}C\grqq{} f"ur englisch \glqq{}column\grqq).
\index[S]{CnX@$C_n(X)$ ($=M_{n,1}(X)$)}%
\dliter{\cite{EffrosRuan00OperatorSpaces}, S. 22, Mitte}%
\end{enumaufz}
\end{defBemerkung}


\dremww{
\begin{bemerkung}
Sei $X$ ein Operatorraum.
Dann ist $C_2(C_2(X))$ vollst"andig isometrisch isomorph zu $C_4(X)$.
\dremark{Allgemein m"u\cms{}te gelten: $M_{m,n}(M_{p,q}(X)) \cong M_{mp,nq}(X)$.}%
\end{bemerkung}

Beweis.
Es gelte: $X \subseteq L(H)$.
Es gilt: $C_2(X) \hookrightarrow M_2(X) \subseteq M_2(L(H)) \cong L(H^2)$.
Sei $x \in C_2(C_2(X))$ mit $x = \binom{x_1}{x_2}$ und $x_1 = \binom{a}{b}$, $x_2 = \binom{c}{d}$.
Es gilt:
\begin{align*}
   \norm{x}^2_{C_2(C_2(X))}
xxx
\end{align*}
}%

\begin{beispiel}[\cite{EffrosRuan00OperatorSpaces}, S. 50; Operatorraumstruktur auf $C_0(\Omega)$]
Sei $\Omega$ ein lokalkompakter Hausdorffraum und $\mfrakA := C_0(\Omega)$.
Sei $\delta : \Omega \rightarrow \mLinStet(\mfrakA,\mbbC),
              \omega \mapsto ( \nu \mapsto \nu(\omega) )$.
Versehen mit der Norm
\begin{align*}
   \norm{a}
&= \sup\bigr\{ \norm{ (f(a_{ij}))_{i,j \in \haken{n}} } \setfdg f \in \Bild(\delta) \bigl\}  \\
&= \sup\bigr\{ \norm{ (a_{ij}(\omega))_{i,j \in \haken{n}} }_{M_n} \setfdg \omega \in \Omega \bigr\}
\end{align*}
ist $M_n(\mfrakA)$ eine \csalgebra.
Ferner gilt: $M_n(\mfrakA) \cong C_0(\Omega,M_n)$.\dremark{genauer: isometrisch \sterns{}iso\-morph}
\end{beispiel}

\dremww{
Beweis.
Zur Isomorphie:
\[ \varphi : M_n(A) \rightarrow C_0(\Omega,M_n),
   a \mapsto ( \omega \mapsto (a_{ij}(\omega))_{i,j \in \haken{n}}) \]
ist ein \sterns{}Isomorphismus auf $C_0(\Omega, M_n)$, also nach \cite{DoranBelfiCAlg}, 24.4,
eine Isometrie.\dremark{13.10.'05, S. 1}%
}%

Sei $X$ ein Banachraum.
Jede lineare Einbettung von $X$ in $\mLinStet(H)$, wobei $H$ ein Hilbertraum sei,
definiert eine Operatorraumstruktur auf $X$.
Wenn man von einer solchen Einbettung zus"atzlich verlangt, \dass sie isometrisch ist,
findet man unter allen solchen Einbettungen eine,
die $X$ so mit einer Operatorraumstruktur versieht, \dass die entstehende Norm minimal ist
(siehe die folgende Definitions-Proposition, Punkt (iv)).

\begin{defBemerkung}[Vgl. \cite{EffrosRuan00OperatorSpaces}, S. 47/48;
  minimale Operatorraumstruktur $\min(X)$]\label{defminX}
Seien $X, Y$ Banachr"aume. Sei $Z$ ein Operatorraum.
\begin{enumaufz}
\item Es sei $X^*_1$ die abgeschlossene Einheitskugel der Dualraumes $X^*$ von $X$,
  versehen mit der schwach*-Topologie.
  Dann ist $X^*_1$ ein kompakter Hausdorffraum.
  Die Abbildung $j : X \to C(X^*_1), x \mapsto (f \mapsto f(x))$, ist eine lineare Isometrie.
  Da $j(X)$ als abgeschlossener Unterraum einer \csalgebra{}
  die Struktur eines Operatorraumes tr"agt,
  erh"alt man so eine Operatorraumstruktur auf $X$, die man mit $\min(X)$ bezeichnet.
\index[S]{minX@$\min(X)$}%
\item Es gilt f"ur alle $x \in M_n(X)$:
\[ \norm{x}_{\min(X),n} = \sup\bigr\{ \norm{f_n(x)}_{M_n} \setfdg f \in X^*_1 \bigl\}. \]
\item F"ur jede lineare Abbildung $\varphi : Z \to X$ gilt:
\dremark{Dies folgt, da $\min(E)$ vollst"andig isometrisch in einer kommutativen \csalgebra{} liegt.}%
\[ \norm{\varphi}_{\text{CB}(Z,\min(X))} = \norm{\varphi}_{\mLinStet(Z,X)}. \]
\item Es gilt: $\norm{x}_{\min(Z),n} \leq \norm{x}_n$ f"ur alle $x \in M_n(Z)$.
\dremark{Wende (iii) auf $\Id_Z$ an.}%
\item Sei $\varphi : X \to Y$ eine Kontraktion.
  Dann ist $\varphi : \min(X) \to Y$ eine Kontraktion und
\[ \varphi : \min(X) \to \min(Y) \]
vollst"andig kontraktiv.
\dremark{Man hat also einen Funktor.}%
\dremark{Insbes. ist $\min(X)$ ein homogener Operatorraum:
  $\norm{\varphi : \min(X) \to \min(X)}_\text{cb} = \norm{\varphi : \min(X) \to X}
   = \norm{\varphi : \min(X) \to \min(X)}$}%
\dliter{\cite{Pisier03OpSpaceTheory}, Chap. 3, S. 71,
  \cite{PaulsenComplBoundedMapsPS}, vor 14.1, S. 195,
  \cite{EffrosRuan00OperatorSpaces}, §3.3, S. 47/48}%
\end{enumaufz}
\end{defBemerkung}

\hcsmoduln sind mit einer nat"urliche Operatorraumstruktur ausgestattet
(vgl. \cite{BlecherLeMerdy04OpAlg}, 8.2.1):

\begin{satz}[Operatorraumstruktur auf Hilbert-\cstern{}Moduln]\label{HCsModORStruktur}
Sei $E$ ein \dremark{Rechts-}\hcsmodul "uber $\mfrakA$ und $n \in \mbbN$.
Dann ist f"ur alle $x \in M_n(E)$ die Matrix
$\bigl(\sum_{k=1}^n \skalpr{x_{ki}}{x_{kj}}\bigr)_{i,j \in \haken{n}}$
ein positives Element in $M_n(\mfrakA)$,
und $E$, versehen mit der durch
\[ \norm{x}_n = \normlr{\left(\sum_{k=1}^n \skalpr{x_{ki}}{x_{kj}}\right)_{i,j}}^{1/2} \]
gegebenen Matrixnorm, ein Operatorraum.
\dliter{\cite{BlecherLeMerdy04OpAlg}, 8.2.1, (8.6); \cite{PaulsenComplBoundedMapsPS}, Th. 14.10}%
\end{satz}

F"ur ein \hcsmodul $E$ h"angt die Norm von $E \oplus E$
(siehe Definitions"=Proposition~\ref{defDirSumHCSMod})
und von $C_2(E)$, wobei $E$ bei $C_2(E)$ als Operatorraum aufgefa\cms{}t wird,
wie folgt zusammen:

\begin{bemerkung}\label{normHCsModEoplusEEqC2E}
Sei $E$ ein Hilbert-\cstern{}Modul "uber $\mfrakA$.
Dann gilt: $E \oplus E \cong C_2(E)$.
Genauer ist $\varphi : E \oplus E \arrowbij C_2(E), x \mapsto \binom{x_1}{x_2}$
eine vollst"andige Isometrie.
\dremark{
  Es gilt $\norm{x}_{M_n(E \oplus E)} = \norm{\varphi_n(x)}_{M_n(C_2(E))}$
  f"ur alle $n \in \mbbN$ und $x \in M_n(E \oplus E)$.}%
\end{bemerkung}

\begin{proof}
\bewitemph{1. Fall: $n=1$.}
Sei $x \in E \oplus E$.
Mit Satz~\ref{HCsModORStruktur} (Operatorraumstruktur auf $E$) erh"alt man:
\begin{align*}
   \norm{ \varphi(x) }^2_{C_2(E)}
&= \normBigg{ \underbrace{\cmpmatrix{ x_1 & 0 \\ x_2 & 0 }}_{:=y} }^2_{M_2(E)}
\overset{\dremarkm{\ref{HCsModORStruktur}}}{=}
   \normlr{ \left(\sum_{k=1}^2 \skalpri{y_{ki}}{y_{kj}}{E}\right)_{i,j} }_{M_2(\mfrakA)}  \\
&= \normlr{ \cmpmatrix{ \skalpr{x_1}{x_1} + \skalpr{x_2}{x_2} & 0 \\ 0 & 0 } }_{M_2(\mfrakA)}  \\
\dremarkm{&= \skalpr{x_1}{x_1} + \skalpr{x_2}{x_2}  \\}
&=  \skalpr{x}{x}
=  \norm{x}^2_{E \oplus E}.
\end{align*}

\bewitemph{2. Fall: $n=2$.}
Sei $x \in M_2(E \oplus E)$.
Seien $a, b \in M_2(E)$ mit $x_{ij} = (a_{ij},b_{ij})$ f"ur alle $i,j \in \haken{2}$.
Setze $y := \cmpmatrix{ a_{11} & a_{12} & 0 & 0 \\
                        b_{11} & b_{12} & 0 & 0 \\
                        a_{21} & a_{22} & 0 & 0 \\
                        b_{21} & b_{22} & 0 & 0 }$.
Mit Satz~\ref{HCsModORStruktur} folgt:
\begin{align*}
&\,\,  \norm{\varphi_2(x)}^2_{M_2(C_2(E))}
\dremarkm{=  \norm{\varphi_2(x)}^2_{M_{4,2(E)}} }
=  \normlr{ \cmpmatrix{ a_{11} & a_{12} \\
                        b_{11} & b_{12} \\
                        a_{21} & a_{22} \\
                        b_{21} & b_{22} } }_{M_{4,2}(E)}  \\
=&\, \normlr{ \cmpmatrix{ a_{11} & a_{12} & 0 & 0 \\
                        b_{11} & b_{12} & 0 & 0 \\
                        a_{21} & a_{22} & 0 & 0 \\
                        b_{21} & b_{22} & 0 & 0 } }_{M_4(E)}
\overset{\dremarkm{\ref{HCsModORStruktur}}}{=}
   \normlr{ \left(\sum_{k=1}^4 \skalpri{y_{ki}}{y_{kj}}{E}\right)_{i,j} }_{M_4(\mfrakA)}  \\
=&\, \normlr{ \cmpmatrix{
   \sum_{k=1}^4 \skalpr{y_{k,1}}{y_{k,1}} & \sum_{k=1}^4 \skalpr{y_{k,1}}{y_{k,2}} & 0 & 0 \\
   \sum_{k=1}^4 \skalpr{y_{k,2}}{y_{k,1}} & \sum_{k=1}^4 \skalpr{y_{k,2}}{y_{k,2}} & 0 & 0 \\
   0 & 0 & 0 & 0 \\
   0 & 0 & 0 & 0 } }_{M_4(\mfrakA)}  \\
\dremarkm{=&\, \normlr{ \cmpmatrix{
   {\scriptstyle \skalprb{a_{11}} + \skalprb{b_{11}} + \skalprb{a_{21}} + \skalprb{b_{21}} } &
     {\scriptstyle \skalpr{a_{11}}{a_{12}} + \skalpr{b_{11}}{b_{12}} + \skalpr{a_{21}}{a_{22}} + \skalpr{b_{21}}{b_{22}} } &
     0 & 0 \\
   {\scriptstyle \skalpr{a_{12}}{a_{11}} + \skalpr{b_{12}}{b_{11}} + \skalpr{a_{22}}{a_{21}} + \skalpr{b_{22}}{b_{21}} } &
     {\scriptstyle \skalprb{a_{12}} + \skalprb{b_{12}} + \skalprb{a_{22}} + \skalprb{b_{22}} } &
     0 & 0  \\
   0 & 0 & 0 & 0 \\
   0 & 0 & 0 & 0 } }  \\
=&\, \normlr{ \cmpmatrix{
   \skalpr{x_{11}}{x_{11}} + \skalpr{x_{21}}{x_{21}} & \skalpr{x_{11}}{x_{12}} + \skalpr{x_{21}}{x_{22}} \\
   \skalpr{x_{12}}{x_{11}} + \skalpr{x_{22}}{x_{21}} & \skalpr{x_{12}}{x_{12}} + \skalpr{x_{22}}{x_{22}} } }_{M_2(\mfrakA)}  \\ }
=&\, \normlr{ \left(\sum_{k=1}^2 \skalpri{x_{ki}}{x_{kj}}{E \oplus E}\right)_{i,j} }_{M_2(\mfrakA)}
=  \norm{x}^2_{M_2(E \oplus E)}.
\end{align*}

\bewitemph{3. Fall: $n > 2$.} Der Beweis verl"auft analog zum 2. Fall.
\end{proof}

\dremww{
\begin{bemerkung}\label{HCSModEoplusEEqC2E}
Sei $X$ ein Operatorraum.
Da $Z := \ternh{X}$ als TRO auch ein Hilbert-\cstern{}Modul ist, tr"agt der Hilbert-\cstern{}Modul
$Z \oplus Z$ zugleich die Struktur eines Operatorraumes,
wodurch $X \oplus X$ eine Operatorraumstruktur erh"alt,
die mit der von $C_2(X)$ "ubereinstimmt.
\dmarginpar{Wof"ur?}\dremark{Wof"ur?}%
\end{bemerkung}

\begin{proof}
\newcommand{\cmjX}{\tilde{X}}%
\dremark{Zeige: $j(X) \oplus j(X)$ ist abgeschlossen in $Z \oplus Z$.
Setze $\cmjX := j(X)$.
Sei $(x,y) \in (X \oplus X)^\mbbN$,
sei $(j(x_n),j(y_n))_n$ konvergent gegen $(x_0,y_0) \in Z \oplus Z$.
Es gilt
$  \norm{j(x_n) - x_0}_Z
\leq  \norm{(j(x_n) - x_0, j(y_n) - y_0)}_{Z \oplus Z} \to 0$ f"ur $n \to \infty$.
Da $\cmjX$ abgeschlossen ist, folgt $x_0 \in \cmjX$.
Analog erh"alt man $y_0 \in \cmjX$.}%
Es ist $\cmjX := j(X)$ abgeschlossen in $Z$.
Da $\cmjX \oplus \cmjX$ im Hilbert-\cstern{}Modul und Operatorraum $Z \oplus Z$ abgeschlossen ist,
erh"alt $\cmjX \oplus \cmjX$ und somit auch $X \oplus X$ eine Operatorraumstruktur,
die mit $X \oplushm X$ bezeichnet wird.
\smallskip

Sei $\varphi : X \oplus X \to C_2(X), x \mapsto \binom{x_1}{x_2}$ und
$\rho : Z \oplus Z \to C_2(Z), x \mapsto \binom{x_1}{x_2}$.
Da $Z$ ein Hilbert-\cstern{}Modul ist,
folgt mit \ref{normHCsModEoplusEEqC2E} f"ur alle $x \in M_n(X \oplushm X)$:
\begin{align*}
   \norm{x}_{n, X \oplushm X}
&= \norm{(j \oplus j)_n(x)}_{n, Z \oplus Z}
\overset{\ref{normHCsModEoplusEEqC2E}}{=}
   \norm{\rho_n(j \oplus j)_n(x)}_{M_n(C_2(Z))}  \\
&= \norm{\varphi_n(x)}_{M_n(C_2(X))}.
\end{align*}
\dremark{
Nach dem oben Gezeigten ist die Operatorraumstruktur von $X \oplus X$
unabh"angig von der Einbettung des Operatorraumes $X$.
Ein alternativer Beweis hierf"ur:

Sei $u : X \to \mLinStet(K)$ vollst"andig isometrisch.
Setze $Y := u(X)$.
Nach \cite{BlecherLeMerdy04OpAlg}, 4.4.5.(2), l"a\cms{}t sich $u$ zu einem \sterns{}Isomorphismus
$\theta : I(\mcalS(X)) \arrowbij I(\mcalS(Y))$ auf $I(\mcalS(Y))$ fortsetzen,
der ecken\-er\-hal\-tend ist.
Mit $\cdot_\Psi$ sei die Multiplikation auf $I(\mcalS(Y))$ bezeichnet.
Es sind $\ternh{X}$ und $\ternh{Y}$ tripelisomorph und somit insbesondere
vollst"andig isometrisch isomorph.
Sei $j : X \to Z$ die kanonische Einbettung und
$(Z_2,\kappa)$ eine tern"are H"ulle von $Y$.
Man hat:
\begin{equation}\label{eqDirSumORtheta}
\theta \circ j = \kappa \circ u.
\end{equation}
\dremark{$\theta \circ j(x) = \theta(\cmpmatrix{ 0 & x \\ 0 & 0 })
=  \cmpmatrix{0 & u(x) \\ 0 & 0} = \kappa \circ u (x)$}%
Sei $y = ((c_{ij},d_{ij}))_{i,j} \in M_n(Y \oplus Y)$.
Dann findet man ein $x = ((a_{ij},b_{ij}))_{i,j} \in M_n(X \oplus X)$ mit
$(u(a_{ij}), u(b_{ij}))_{i,j} = y$.
Es gilt:
\begin{align*}
   \norm{y}^2_{n, Y \oplushm Y}
&= \norm{(\kappa \oplus \kappa)_n(y)}^2_{n, Z_2 \oplushm Z_2}
=  \norm{ ( \sum_{k=1}^n
      \skalpr{(\kappa \oplus \kappa)(y_{ki})}{(\kappa \oplus \kappa)(y_{kj})} )_{i,j} }_{M_n(Z_2^* Z_2)}  \\
&= \norm{ ( \sum ( \skalpr{\kappa(c_{ki})}{\kappa(c_{kj})} +
                   \skalpr{\kappa(d_{ki})}{\kappa(d_{kj})} )_{i,j} }  \\
&= \norm{ ( \sum ( \kappa(u(a_{ki}))^* \cdot_\Psi \kappa(u(a_{kj})) +
                   \kappa(u(b_{ki}))^* \cdot_\Psi \kappa(u(b_{kj})) )_{i,j} }  \\
&\overset{\eqref{eqDirSumORtheta}}{=}
   \norm{ ( \sum ( \theta(j(a_{ki}))^* \cdot_\Psi \theta(j(a_{kj})) +
                   \theta(j(b_{ki}))^* \cdot_\Psi \theta(j(b_{kj})) )_{i,j} }  \\
&= \norm{ ( \sum ( \theta(j(a_{ki})^* \cdot_\Phi j(a_{kj}) +
                          j(b_{ki})^* \cdot_\Phi j(b_{kj})) )_{i,j} }  \\
&= \norm{ \theta_n(( \sum \skalpr{j(a_{ki})}{j(a_{kj})} + \skalpr{j(b_{ki})}{j(b_{kj})} )_{i,j} ) }  \\
&= \norm{ (\sum \skalpr{(j \oplus j)(x_{ki})}{(j \oplus j)(x_{kj})})_{i,j} }_{M_n(Z^*Z)}  \\
&= \norm{(j \oplus j)_n(x)}^2_{n, Z \oplushm Z}
=  \norm{x}^2_{n, X \oplushm X}.
\end{align*}
Somit ist die Operatorraumstruktur von $X \oplus X$ unabh"angig von der Einbettung von $X$.
}%
\dremark{Geht so wohl nicht: Mit $\cdot_\Phi$ sei die Multiplikation auf $I(\mcalS(X))$ bezeichnet.
Verm"oge $\skalpr{x}{y} := x^* \cdot_\Phi y$ ist $Z$ und somit auch $Z \oplus Z$
ein Rechts-Hilbert-\cstern{}Modul "uber $Z^*Z \subseteq I_{22}(X)$.
Es gilt f"ur alle $x \in X \oplus X$:
\dremark{18.6.'08, S. 1/2}%
\[ \norm{x}^2
=  \skalpr{x}{x}
=  \skalpr{x_1}{x_1} + \skalpr{x_2}{x_2}
=  x_1^* \cdot_\Phi x_1 + x_2^* \cdot_\Phi x_2
\overset{?,\dremarkm{\cite{BlecherLeMerdy04OpAlg}, 1.2.5}}{=}  \norm{x}^2_{C_2(X)}.   \]
Andere Beweisidee: Es ex. unitale, isometrische \sterns{}Darst. $\pi : I(S(X)) \to L(K)$,
  also $\pi(x \cdot_\Phi y) = \pi(x) \pi(y)$.
}%
\end{proof}
}%

\begin{definitn}
Sei $H$ ein Hilbertraum.
Mit $\overline{H}$ sei der Anti-Hilbertraum von $H$ bezeichnet,
also der Hilbertraum auf der Menge $H$ mit der Addition von $H$,
der Skalarmultiplikation
\[ (\lambda, \xi) \mapsto \overline{\lambda} \xi \]
und dem Skalarprodukt
\dliter{\cite{EffrosRuan00OperatorSpaces}, S. 336/337}%
\[ \skalpr{\overline{\xi}}{\overline{\eta}} = \skalpr{\eta}{\xi}. \]
\end{definitn}

Hilbertr"aume kann man mit verschiedenen Operatorraumstrukturen versehen,
wie die folgenden beiden Propositionen zeigen
(\cite{EffrosRuan00OperatorSpaces}, S. 55 bzw. S.~54):

\begin{bemerkung}[Zeilen-Hilbertoperatorraum]\label{defHr}
Sei $H$ ein Hilbertraum.
\begin{enumaufz}
\item Sei $\Phi : \overline{H} \to \dualr{H}$ der nat"urliche Isomorphismus.
Die nat"urliche Isometrie
\[ R : H \to H^{**} = \mLinStet(\dualr{H},\mbbC) \cong \mLinStet(\overline{H},\mbbC) \]
ist gegeben durch
\[ R(\zeta)(\overline{\xi}) := \Phi(\overline{\xi})(\zeta) = \skalpr{\zeta}{\xi} \]
f"ur alle $\zeta, \xi \in H$.
Da $\mLinStet(\overline{H}, \mbbC)$ ein Operatorraum ist, wird hierdurch $H$ mit
der Struktur eines Operatorraumes versehen, der \defemphi{Zeilen-Hilbert\-opera\-tor\-raum}
genannt und mit $H^r$ bezeichnet wird (\glqq r\grqq{} f"ur englisch \glqq row\grqq).
\index[S]{Hr@$H^r$}%

\item F"ur alle $\xi \in M_{m,n}(H^r)$ gilt:
\[ \norm{\xi}
=  \normlr{ \left(\sum_{k=1}^n \skalpr{\xi_{ik}}{\xi_{jk}} \right)_{i,j} }^{1/2}_{M_{m,n}(\mbbC)}. \]
\dliter{\cite{EffrosRuan00OperatorSpaces}, S. 55}%
\dremark{Ferner ist $H^r$ injektiv.}%
\dremark{Sei $q$ Rang-1-Projektion auf $H$.
  Dann gilt $H^r \cong \mLinStet(H)q$ (vollst"andig isometrisch).
  Genauer: \cite{Robertson91InjectiveMatHilbertSp}}%
\dremark{Interessante Literatur zum Zeilen- und Spalten-Hilbert\-opera\-tor\-raum:
  \cite{Mathes94CharOfRowHilbertSp} (stellt verschiedene Charakterisierungen bereit),
  \cite{Mathes94CharOfRowHilbertSp} ("uber injektive Hilbert-Operatorr"aume)}%
\end{enumaufz}
\end{bemerkung}

\begin{bemerkung}[Spalten-Hilbertoperatorraum]
Sei $H$ ein Hilbertraum.
\begin{enumaufz}
\item Definiere f"ur alle $\xi \in H$
\[ C^\xi : \mbbC \to H, a \mapsto a\xi, \quad\text{und}\quad
   C : H \to \mLinStet(\mbbC,H), \xi \mapsto C^\xi. \]
Dann ist $C$ ein isometrischer Vektorraumisomorphismus auf $\mLinStet(\mbbC,H)$.
F"ur alle $\xi \in M_n(H)$ wird durch
\[ \norm{\xi}_n = \norm{C_n(\xi)}_n \]
eine Matrixnorm auf $H$ definiert, durch die $H$ mit der Struktur eines
Operatorraumes versehen wird, der \defemphi{Spalten-Hilbertoperatorraum} genannt und
mit $H^c$ bezeichnet wird (\glqq c\grqq{} f"ur englisch \glqq column\grqq).
\index[S]{Hc@$H^c$}%

\item F"ur alle $\xi \in M_{m,n}(H^c)$ gilt:
\[ \norm{\xi}
=  \normlr{ \left(\sum_{k=1}^m \skalpr{\xi_{kj}}{\xi_{ki}} \right)_{i,j} }^{1/2}_{M_{m,n}(\mbbC)}. \]
\dliter{\cite{EffrosRuan00OperatorSpaces}, S. 54}%
\dremark{Ferner ist $H^c$ injektiv.}%
\dremark{Man kann sich $H$ als eine Spalte in $\mLinStet(H)$ vorstellen,
  siehe \cite{BlecherLeMerdy04OpAlg}, 1.2.23, Beschreibung "uber
  Rang-1-Operatoren und Z. 22.8.'07, S. 1/2.
  F"ur $H=\ell^2_n:=\ell_2-\bigoplus_{i=1}^n \mbbC$ ist $H^c$
  vollst"andig isometrisch zu $M_{n,1}$, man kann sich somit $H^c$ als
  eine Spalte in $M_n$ vorstellen. Siehe Z. 15.6.'07, S. 1/2}%
\dliter{\cite{BlecherLeMerdy04OpAlg}, 1.2.23}%
\dremark{Sei $p := \Id_h$, $q$ Rang-1-Projektion auf $H$.
  Es gilt $H^c \cong \mLinStet(H)q$ (vollst"andig isometrisch).
  Da $\mLinStet(H)$ eine injektive \csalgebra ist und $p$ und $q$ \cmvk{e} Projektionen sind,
  ist $p\mLinStet(H)q$ ein injektiver Operatorraum.}%
\dremark{Sei $O_n$ die von $s_1,\dots,s_n$ erzeugte Cuntz-Algebra.
  Dann ist $\vrerz{s_1,\dots,s_n}$ vollst"andig isometrisch zu
  $C_n = (\ell^2_n)^c$ und
  $\vrerz{s_1^*,\dots,s_n^*}$ vollst"andig isometrisch zu
  $R_n = (\ell^2_n)^r$.}%
\dliter{\cite{Robertson91InjectiveMatHilbertSp}, S. $184^6$.}%
\end{enumaufz}
\end{bemerkung}

\dremww{
\begin{bemerkung}
Sei $H$ ein Hilbertraum.
\begin{enumaufz}
\item Dann ist $H$ zugleich ein \dremark{Rechts-}\hcsmodul{} "uber $\mbbC$,
welcher nach \ref{HCsModORStruktur} eine Operatorraumstruktur tr"agt.
Es sei $H$, versehen mit dieser Operatorraumstruktur, als $X$ bezeichnet.
Dann gilt: $H^c \cong X$ vollst"andig isometrisch.
\dremark{Beweis: \cite{EffrosRuan00OperatorSpaces}, (3.4.2), S. 54}%

\item Definiere $\skalprd{x}{y} := \skalpr{x}{y}_H$ f"ur alle $x,y \in H$.
Versehen mit dem inneren Produkt $\skalprd{\cdot}{\cdot\cdot}$
kann man $H$ als Links-\hcsmodul $\tilde{H}$ "uber $\mbbC$ auffassen.
Es gilt: $H^r \cong \tilde{H}$ vollst"andig isometrisch.
\end{enumaufz}
\end{bemerkung}

\begin{proof}
\bewitemph{(i):}
Sei $x \in M_n(H)$.
Nach \ref{HCsModORStruktur} ist
$( \sum_{k=1}^n \skalpr{x_{ki}}{x_{kj}} )_{i,j} \in M_n(\mbbC)_+$,
also insbesondere normal.
Es folgt:
\[ \norm{x}^2_{M_n(H^c)}
\overset{\cite{EffrosRuan00OperatorSpaces}, (3.4.2), S. 54}{=}
   \norm{ ( \sum_{k=1}^n \skalpr{x_{kj}}{x_{ki}} )_{i,j} }
\overset{\ref{CsalgNormatEqnorma}}{=}
   \norm{ ( \sum_{k=1}^n \skalpr{x_{kj}}{x_{ki}} )_{j,i} }
=  \norm{ ( \sum_{k=1}^n \skalpr{x_{ki}}{x_{kj}} )_{i,j} }
=  \norm{x}^2_{M_n(X)}. \]

\bewitemph{(ii):}
Es gilt:
\dremark{15.7.'08/1}%
\[ \norm{x}^2_{M_n(H^r)}
\overset{\cite{EffrosRuan00OperatorSpaces}, S. 56}{=}
   \norm{ ( \sum_{k=1}^n \skalprd{x_{ik}}{x_{jk}} )_{j,i} }
\overset{\cite{PaulsenComplBoundedMapsPS}, Th. 14.10}{=}
   \norm{x}_{M_n(\tilde{H})}. \]
\end{proof}

\begin{beispiel}[Operatorraumstruktur auf $L^p$]
Die Operatorraumstruktur auf $L^p$ wird in Pisier, Operator Space Theory, S. 178 (9.5) geschildert.
Dort (1. Absatz) findet man einen Literaturverweis auf den weiterf"uhrenden Artikel
\cite{PisierNoncomLpSpaces}, von dem eine Kurzfassung bei arxiv erschienen ist.
\end{beispiel}
}%

\dremww{
\section{Homogene Operatorr"aume}

\begin{definitn}
Ein Operatorraum $E$ hei"st \defemph{homogen},
falls f"ur jede beschr"ankte, lineare Abbildung $T : E \to E$ gilt:
$\norm{T}_\text{cb} = \norm{T}$.
\index[B]{homogener Operatorraum}%
\dliter{\cite{Pisier03OpSpaceTheory}, S. 172 oben, §9.2}%
\end{definitn}

\begin{definitn}
Sei $\lambda \in \mbbR_{\geq 1}$.
Ein Operatorraum $E$ hei"st \defemphi{homogen} (bzw. \defemphi{$\lambda$-homogen}),
falls f"ur jede beschr"ankte, lineare Abbildung $T : E \to E$ gilt:
$\norm{T}_\text{cb} = \norm{T}$ (bzw. $\norm{T}_\text{cb} \leq \lambda\norm{T}$).
\dliter{\cite{Pisier03OpSpaceTheory}, S. 172 oben, §9.2}%
\end{definitn}

\begin{beispiel}\label{HrHomogen}
\begin{enumaufz}
\item Sei $H$ ein Hilbertraum.
Dann sind der Spalten"=Hilbertoperatorraum $H^c$ und der Zeilen-Hilbertoperatorraum $H^r$
\dremark{ und Pisiers Operatorhilbertraum $H^\text{oh}$}%
homogen.
\dliter{\cite{Paulsen02ComplBoundedMaps}, S. 200}%
\dremark{Es gilt sogar mehr: F"ur jeden Hilbertraum $K$, f"ur jede Abbildung
$u : H \to K$ gilt: $\norm{u} = \norm{u}_\text{cb}$.}%

\item Sei $X$ ein Banachraum.
Dann ist $\min(X)$\dremark{und $\max(E)$} homogen.
\dliter{\cite{Pisier03OpSpaceTheory}, S. 172}%
Es gilt sogar f"ur jeden Operatorraum $F$:
$\mCB(F, \min(E)) \cong \mLinStet(F,E)$ (isometrisch) und
$\mCB(max(E), F) \cong \mLinStet(E,F)$ (isometrisch).
\dremark{Weitere Beispiele in \cite{Pisier03OpSpaceTheory}, §9}%
\dliter{\cite{EffrosRuan00OperatorSpaces}, (3.3.7), S. 48 und (3.3.9), S. 49}%

\item Der Schnitt und die Summe zweier homogener hilbertscher Operatorr"aume
sind wieder homogen und hilbertsch.
\dliter{\cite{Pisier96OpHilbertSpInterpol}}%

\dremww{\item Aus \cite{Pisier03OpSpaceTheory}: Der Raum $\Phi(I)$ (§9.3) ist
(isometrisch) Hilbertsch und homogen.\dliter{Th. 9.3.1}%
Die R"aume $\mcalG_p$ (9.8.1) und $\mcalV_I$ (9.9.5) sind hilbertsch und homogen.}%
\end{enumaufz}
\end{beispiel}

\begin{proof}
\bewitemph{(i):} \cite{Paulsen02ComplBoundedMaps}, S. 200.

\bewitemph{(ii):} \cite{Pisier03OpSpaceTheory}, S. 172.
\end{proof}

\begin{proof}
\bewitemph{(i):}
\bewitemph{$H_c$:} Nach \cite{EffrosRuan00OperatorSpaces}, 3.4.1, gilt
$\mLinStet(H,K) \overset{\text{v.i.}}{\cong} \mCB(H_c,K_c)$, also ist insbesondere
$H_c$ homogen.

\bewitemph{$H_r$}: Nach \cite{BlecherLeMerdy04OpAlg}, S. $11_3$, gilt
$\mLinStet(H,K) = \mCB(H_r,K_r)$ isometrisch.

\bewitemph{$H_\text{oh}$}: Man findet eine Menge $I$ so, \dass $H$ isomorph zu $\ell^2(I)$ ist.
Ferner ist $\text{OH}(I)$ als Banachraum isometrisch zu $\ell^2(I)$
(\cite{Pisier03OpSpaceTheory}, Th. 7.1), und es gilt
f"ur jede Abbildung $u : \text{OH}(I) \to \text{OH}(I)$: $\norm{u} = \norm{u}_\text{cb}$.
(\cite{Pisier03OpSpaceTheory}, Prop. 7.2).
\end{proof}

\begin{definitn}
Zwei Operatorr"aume $X$ und $Y$ hei"sen \defemphi{vollständig isomorph},
falls es einen linearen Isomorphismus $u : X \arrowbij Y$ auf $Y$ so gibt,
\dass $u$ und $u^{-1}$ \cmvb{} sind.
\dliter{\cite{Pisier03OpSpaceTheory}, Def. 1.4}%
\end{definitn}

\begin{satz}
Jeder $\lambda$-homogene Operatorraum $X$ ist vollst"andig isomorph zu einem
homogenen Operatorraum $\tilde{X}$ mit $d_\text{cb}(X,\tilde{X}) \leq \lambda$.
\dliter{\cite{Pisier03OpSpaceTheory}, Prop. 9.2.2, S. 172.
  Ein ausf"uhrlicherer Beweis: Z. 15.7.'07, S. 1}%
\end{satz}

\begin{proof}
\dremark{Nur die Konstruktion:}
Sei $\beta$ die Familie aller Operatoren $u : X \to X$ mit $\norm{u} \leq 1$.
Man findet einen Hilbertraum $H$ mit $X \subseteq \mLinStet(H)$.
Betrachte die direkte Summe $\bigoplus_{u \in \beta} \mLinStet(H)$,
versehen mit der Operatorraumstruktur der direkten Summe.
Die Abbildung $J : X \to \bigoplus_{u \in \beta}, x \mapsto (u(x))_{u \in \beta}$
ist eine Isometrie.
Setze $\tilde{X} := \Bild(J)$.
\end{proof}

\begin{satz}\label{homogOpraumVollstIsometrie}
Sei $E$ ein homogener Operatorraum und $f : E \rightarrow E$ eine lineare, surjektive
Isometrie auf $E$.
Dann ist $f$ eine vollst"andige Isometrie.\dliter{Pisier, 9.2, S. 172 oben}
\end{satz}

\begin{proof}
\dmarginpar{Bew. pr}
Da $E$ homogen ist, gilt:
$\norm{f} = \norm{f}_\text{cb} = \sup \{ \norm{f_n} \setfdg n \in \mbbN \}$.
Mit $(N_1)$ erh"alt man:
\[ \norm{f} \leq \norm{f_2} \leq \norm{f_3} \leq \dots \leq \norm{f_n} \leq \dots
   \leq \norm{f}_\text{cb}. \]
Also gilt f"ur alle $n \in \mbbN$: $\norm{f_n} = \norm{f}$.
F"ur alle $x \in E$ gilt $\norm{f(x)} = \norm{x}$, also $1 = \norm{f} = \norm{f}_\text{cb}$
und $1 = \norm{f} = \norm{f_n}$.
F"ur alle $a \in M_n(E) \setminus \{0\}$ erh"alt man:
\begin{equation}\label{eqPisier9:2a}
\norm{f_n(a)} = \norm{a} \cdot \norm{ f_n\bigl(\frac{a}{\norm{a}}\bigr) }
\leq \norm{a} \cdot \norm{\frac{a}{\norm{a}}} = \norm{a}.
\end{equation}
Da $f$ eine surjektive Isometrie ist, gilt
$\norm{f^{-1}(b)} = \norm{f(f^{-1}(b))} = \norm{b}$ f"ur alle $b \in E$, also ist $f^{-1}$ eine Isometrie.
Ferner ist $f^{-1}$ linear und surjektiv auf $E$, also gilt:
$\norm{f^{-1}}_\text{cb} = \norm{f^{-1}} = \norm{(f^{-1})_n} = 1$.
Analog folgt:
\begin{equation}\label{eqPisier9:2b}
\norm{(f^{-1})_n(b)} \leq \norm{b} \qquad\text{f"ur alle } b \in M_n(E).
\end{equation}
Man erh"alt f"ur alle $a \in M_n(E)$:
\[ \norm{a}_n
=  \norm{(f^{-1}_n \circ f_n)(a)}
\overset{\eqref{eqPisier9:2b}}{\leq}  \norm{f_n(a)}
\overset{\eqref{eqPisier9:2a}}{\leq}  \norm{a}_n,  \]
also $\norm{f_n(a)} = \norm{a}$.
Somit ist $f$ eine vollst"andige Isometrie.
\dremark{8.6.'05, S. 2}
\end{proof}
}%

\section{Unitale Operatorsysteme}

Wir erinnern in diesem Abschnitt an einige grundlegende Definitionen und Aussagen
zu unitalen Operatorsystemen.

\begin{definitn}
\begin{enumaufz}
\item Unter einem \defemphi{Operatorsystem} versteht man
  einen abgeschlossenen, selbstadjungierten Untervektorraum
  der linearen, stetigen Operatoren auf einem Hilbertraum $H$.
\item Ein Operatorsystem $X \subseteq \mLinStet(H)$ hei"st \defemph{unital}, falls $\Id_H \in X$ gilt.
  In diesem Fall bezeichnet man $\Id_H$ als Einselement von $X$ und
  notiert es als $e_X$.
\index[S]{ekleinX@$e_X$ (Einselement)}%
\dremark{Achtung: In der Literatur werden Operatorsysteme h"aufig als unital vorausgesetzt.}%
\dremark{Welche guten Eigenschaften haben unitale Operatorsysteme?
  Z. B.: Seien $E$, $F$ unitale Operatorsysteme,
  sei $\alpha : E \rightarrow F$ linear und unital.
  Dann ist $\alpha$ genau dann vollst"andig positiv,
  wenn $\alpha$ eine vollst"andige Kontraktion ist.}%
\end{enumaufz}
\end{definitn}

\begin{definitn}
Sei $X \subseteq \mLinStet(H)$ ein unitales Operatorsystem.
Dann wird ein beliebiges $x \in X$ \defemph{positiv} genannt,
falls $x$ in $\mLinStet(H)$ positiv ist.
\index[B]{positives Element eines unitalen Operatorsystems}%
\end{definitn}

\begin{definitn}
Seien $X$, $Y$ unitale Operatorsysteme.
Dann hei"st eine lineare Abbildung $\varphi : X \to Y$
\begin{enumaufz}
\item \defemph{unital}, falls $\varphi(e_X) = e_Y$ gilt,
\index[B]{unitale Abbildung}%
\item \defemph{positiv}, falls die positiven Elemente von $X$ in
  die positiven Elemente von $Y$ abgebildet werden,
\index[B]{positive Abbildung}%
\item \defemph{vollständig positiv},\dremark{(kurz \defemph{v.\,p.})}
  falls $\varphi_n$ f"ur alle $n \in \mbbN$ positiv ist,
\index[B]{vollständig positive Abbildung}%
\item \defemphi{Ordnungsmonomorphismus},
  falls $\varphi$ injektiv ist und $\varphi$ und $\varphi^{-1}$ positiv sind, und
\item \defemphi{Ordnungsisomorphismus},
  falls $\varphi$ ein surjektiver Ordnungsmonomorphismus ist.
\end{enumaufz}
\end{definitn}

\begin{bemerkung}[\cite{EffrosRuan00OperatorSpaces}, Corollary 5.1.2]
Seien $X \subseteq \mLinStet(H), Y \subseteq \mLinStet(K)$ unitale Operatorsysteme.
Sei $\varphi : X \rightarrow Y$ eine unitale, lineare Abbildung.
Dann ist $\varphi$ genau dann vollst"andig positiv, wenn $\varphi$ vollst"andig kontraktiv ist.
\end{bemerkung}

Auf die folgende Art und Weise kann man f"ur jeden Operatorraum
ein zugeh"origes unitales Operatorsystem definieren:

\begin{defBemerkung}[Vgl. \cite{BlecherLeMerdy04OpAlg}, 1.3.14 und Lemma 1.3.15]\label{defPaulsenSys}
Seien $X \subseteq \mLinStet(H)$ und $Y \subseteq \mLinStet(K)$ Operatorr"aume.
\begin{enumaufz}
\item Es ist
\begin{align*}
   \mcalS(X)
&:= \begin{pmatrix} \mbbC \Id_H & X \\ X^* & \mbbC \Id_H \end{pmatrix}
=  \left\{ \begin{pmatrix} \lambda \Id_H & x \\ y^* & \mu \Id_H \end{pmatrix}
           \setfdg x,y \in X, \lambda, \mu \in \mbbC \right\}  \\
&\subseteq  M_2(\mLinStet(H)),
\end{align*}
versehen mit der von $\mLinStet(H \oplus H) \cong M_2(\mLinStet(H))$ induzierten Struktur,
ein unitales Operatorsystem, genannt \defemphi{Paulsen-System} von $X$.
\index[B]{Paulsen-System}%
\index[S]{SX@$\mcalS(X)$ (Paulsen-System)}%

\item Weiter ist
\[ \varphi : X \to \mcalS(X), x \mapsto \cmpmatrix{0 & x \\ 0 & 0}, \]
eine vollst"andig isometrische Einbettung.

\item Sei $u : X \to Y$. Definiere
\[ \mcalS(u) : \mcalS(X) \to \mcalS(Y),
   \begin{pmatrix} \lambda\Id_H & x_1 \\ x_2^* & \mu\Id_H \end{pmatrix} \mapsto
   \begin{pmatrix} \lambda\Id_K & u(x_1) \\ u(x_2)^* & \mu\Id_K \end{pmatrix}. \]
Ist $u$ \cmvk{} (bzw. vollst"andig isometrisch), so ist $\mcalS(u)$ unital und \cmvp{}
(bzw. ein unitaler, vollst"andiger Ordnungsmonomorphismus).
\dremark{Zur Motivation: \ref{normxLeq1IffMatGeq0}. Hiermit beweist man auch (ii).}%
\dliter{\cite{PaulsenComplBoundedMapsPS}, 8.1}%
\dremark{Abstraktes Paulsen-System in \cite{PaulsenComplBoundedMapsPS}, 13.3}%
\end{enumaufz}
\end{defBemerkung}

Als unitales Operatorsystem, also bis auf unitale, vollst"andige Ordnungsisomorphismen,
h"angt $\mcalS(X)$ nach Definitions"=Proposition~\ref{defPaulsenSys}.(iii) nur von der Operatorraumstruktur von $X$ ab
und nicht von der Einbettung von $X$ in~$\mLinStet(H)$.
\dremark{$X \subseteq \mLinStet(H)$ isomorph zu $u(X) = Y \subseteq \mLinStet(K)$
  $\Rightarrow$ $\mcalS(X) \cong \mcalS(Y)$}%

\dremww{
\begin{satz}[\cite{Hamana79}, S. 779]
Seien $X$, $Y$ unitale Operatorsysteme.
Sei $\varphi : X \to Y$ linear und unital.
Dann ist $\varphi$ genau dann ein vollst"andiger Ordnungsmonomorphismus,
wenn $\varphi$ eine vollst"andige Isometrie ist.
\end{satz}

\begin{bemerkung}
Sei $X$ ein Operatorraum
Im allgemeinen gilt: $M_2(\mcalS(X)) \ncong \mcalS(M_2(X))$.
\end{bemerkung}

Beweis.
Setze $X := \mbbC$.
Es gilt $\dim(\mcalS(X)) = 4$,
also $\dim(M_2(\mcalS(X))) = 16 \neq 10 = \dim(\mcalS(M_2(X)))$.
}%

\section{Injektive Operatorr"aume}

Eine wichtige Rolle in der Theorie der Operatorr"aume spielen die
injektiven Operatorr"aume.
Da wir uns mit Multiplikatoren auf Operatorr"aumen besch"aftigen,
ist f"ur uns die injektive H"ulle eines Operatorraumes von besonderem Interesse.

\begin{definitn}\label{defInjOpraum}
Ein Operatorraum $Z$ hei"st \defemph{injektiv}, falls
f"ur jeden Operatorraum $X$, jede \cmvk{e} Abbildung $\varphi : X \to Z$
und f"ur jeden Operatorraum~$Y$, der $X$ als abgeschlossenen Unterraum enth"alt,
eine \cmvk{e} Fortsetzung $\Phi : Y \to Z$ von $\varphi$ existiert.
\index[B]{injektiver Operatorraum}%
\[
\xymatrix{
Y \ar@{-->}[rd]^\Phi & \\
X \ar[u]^\subseteq \ar[r]_{\varphi} & Z}
\]
\dremark{Man betrachtet die Kategorie $\mcalO_1$ der Operatorr"aume mit \cmvk{en}
  Abbildungen als Morphismen. Die Einbettungen sind vollst"andige Isometrien.
  Nach \cite{PaulsenComplBoundedMapsPS}, S. 218 kann man stattdessen
  auch \cmvb{e} Abbildungen als Morphismen betrachten,
  wobei die Fortsetzung dieselbe cb-Norm haben mu\cms{}.}%
\dliter{\cite{BlecherLeMerdy04OpAlg}, 1.2.9, \cite{PaulsenComplBoundedMapsPS}, S. 218,
  siehe auch \cite{EffrosRuan00OperatorSpaces}, S. 70 unten}%
\end{definitn}

\dremark{
Ein Operatorraum ist somit genau dann injektiv,
falls er ein \glqq{}injektives Objekt\grqq{} in der Kategorie
der Operatorr"aume und vollst"andig kontraktiven Abbildungen ist.
\dliter{\cite{BlecherLeMerdy04OpAlg}, 1.2.9}%
}%

\begin{beispiel}[\cite{BlecherLeMerdy04OpAlg}, Theorem 1.2.10]\label{LHIstInjOR}
Es seien $H$, $K$ Hilbertr"aume.
Dann ist $\mLinStet(H,K)$ ein injektiver Operatorraum.
\end{beispiel}

\dremww{
\begin{bemerkung}
In der Kategorie $\mcalO_1$ sind die Monomorphismen injektiv.
\end{bemerkung}

\begin{proof}
Sei $m : F \to G$ ein Monomorphismus in $\mcalO_1$.

Annahme: $m$ ist nicht injektiv.

Dann findet man ein $x \in \Kern(m) \setminus \{0\}$.
Sei $\mfrakB$ eine Basis von $F$ mit $x \in \mfrakB$.
Setze $f := \Id_F$ und $g := \Id_{\cmlin{\mfrakB \setminus \{x\}}} \cup\, 0_{\vrerz{x}}$.
Dann gilt $m \circ f = m \circ g$, aber $f \neq g$.
Dies steht im Widerspruch dazu, \dass $m$ ein Monomorphismus ist.
\end{proof}

\begin{definitn}[injektiv in \csalgebren]
Eine unitale \csalgebra{} $\mfrakB$ hei"st injektiv, falls, gegeben ein beliebiger
selbstadjungierten unitalen Untervektorraum $S$ einer unitalen \csalgebra{} $\mfrakC$,
jede unitale, vollst"andig positive lineare Abbildung von $S$ nach $\mfrakB$
fortgesetzt werden kann
zu einer vollst"andig positiven, linearen Abbildung von $\mfrakC$ nach $\mfrakB$.
\dremark{Kategorie, deren Objekte unitale \csalgebren{} sind und
deren Morphismen unitale, vollst"andig positive, lineare Abbildungen}%
\end{definitn}

\begin{definitn}[injektiv in Banachr"aumen]
Ein Banachraum $Y$ hei"st injektiv, falls jede stetige lineare Abbildung von
einen Untervektorraum eines Banachraumes $Z$ nach $Y$ zu einer stetigen linearen
Abbildung der selben Norm auf ganz $Z$ fortgesetzt werden kann.
\dremark{Kategorie, deren Objekte Banachr"aume sind und deren Morphismen
kontraktive lineare Abbildungen. Einbettung: lineare Isometrie.}%
\dliter{\cite{CohenInjEnvBanachSp}}%
\end{definitn}
}%

\begin{definitn}
Sei $X$ ein Operatorraum.
Eine \defemph{Erweiterung} von $X$ ist ein Paar, bestehend aus einem Operatorraum
$Y$ und einer vollst"andigen Isometrie $\kappa : X \rightarrow Y$.
\dremark{Idee: $X \subseteq Y$ soll ja bedeuten,
  \dass die Struktur von $Y$ auf $X$ vererbt wird.}%
\index[B]{Erweiterung eines Operatorraumes}%
\end{definitn}

\begin{definitn}
Sei $X$ ein Operatorraum.
Dann hei"st $(Y,\kappa)$ \defemphi{injektive Hülle} von $X$, falls
$(Y,\kappa)$ eine Erweiterung von $X$ ist,
$Y$ injektiv ist und es keinen echten injektiven Unterraum von $Y$ gibt, der $\kappa(X)$ enth"alt.
\dliter{\cite{BlecherLeMerdy04OpAlg}, 4.2.3; \cite{PaulsenComplBoundedMapsPS}, S. 220}%
\end{definitn}

Die Konstruktion einer injektiven H"ulle f"ur einen Operatorraum $X$
wird nicht auf dem Wege durchgef"uhrt, \dass man eine absteigende Kette
$(Y_\lambda)_\lambda$ von injektiven Operatorr"aumen mit $X \subseteq Y_\lambda$ betrachtet
und $\bigcap_\lambda Y_\lambda$ bildet, denn es ist nicht klar, \dass der Schnitt
wieder injektiv ist.
Stattdessen benutzt man sogenannte $X$-Halbnormen.

\begin{definitn}\label{defFAbbOR}
Sei $X \subseteq \mLinStet(H)$.\dremark{Sp"ater wird $X$ stets ein Operatorraum sein}%
\begin{enumaufz}
\item Eine Abbildung $\varphi : \mLinStet(H) \rightarrow \mLinStet(H)$ hei"st
\defemph{$X$-Abbildung}, falls $\varphi$ vollst"andig kontraktiv ist und
$\varphi\restring_X = \Id_X$ erf"ullt.
\index[B]{XAbbildung@$X$-Abbildung}%

\item Unter einer \defemph{$X$-Projektion} versteht man eine $X$-Abbildung $\varphi$,
f"ur die $\varphi \circ \varphi = \varphi$ gilt.
\index[B]{XProjektion@$X$-Projektion}%
\dliter{analog zu \cite{PaulsenComplBoundedMapsPS}, S. 221 oben;
  \cite{Paulsen02ComplBoundedMaps}, S. 209 unten}%
\end{enumaufz}
\end{definitn}

Somit ist eine $X$-Projektion $\varphi$ eine Projektion auf $Y := \varphi(\mLinStet(H))$,
und es gilt: $X \subseteq Y$.

\begin{defBemerkung}
Sei $X \subseteq \mLinStet(H)$.
Auf der Menge der $X$"=Projektionen wird durch
\[ \varphi \preceq \psi
\overset{\text{def}}{\iff} \varphi \circ \psi = \varphi = \psi \circ \varphi \]
eine Ordnungsrelation definiert.\dremark{$\preceq$ ist i.a. nicht konnex}
\end{defBemerkung}

\dremark{
\begin{proof}
Seien $\varphi, \psi, \rho$ $F$-Projektionen.

\bewitemph{Reflexivit"at} ist klar.

\bewitemph{Antisymmetrie:}
Es gelte $\varphi \preceq \psi$ und $\psi \preceq \varphi$.
Also gilt $\varphi = \varphi \circ \psi = \psi$.

\bewitemph{Transitivit"at:}
Es gelte $\varphi \preceq \psi$ und $\psi \preceq \rho$.
Also gilt
\begin{align}
\varphi \circ \psi &= \varphi = \psi \circ \varphi \qquad\text{und} \label{eqOrdnungAufFProj:a}\\
\psi \circ \rho    &= \psi    = \rho \circ \psi \label{eqOrdnungAufFProj:b}.
\end{align}
Es folgt:
\begin{align*}
   \varphi \circ \rho
\overset{\eqref{eqOrdnungAufFProj:a}}{=}  \varphi \circ \psi \circ \rho
\overset{\eqref{eqOrdnungAufFProj:b}}{=}  \varphi \circ \psi
\overset{\eqref{eqOrdnungAufFProj:a}}{=}  \varphi  \quad\text{und}\quad
    \rho \circ \varphi
\overset{\eqref{eqOrdnungAufFProj:a}}{=}  \rho \circ \psi \circ \varphi
\overset{\eqref{eqOrdnungAufFProj:b}}{=}  \psi \circ \varphi
\overset{\eqref{eqOrdnungAufFProj:a}}{=}  \varphi,
\end{align*}
also $\varphi \preceq \rho$.\dremark{16.2.'06, S. 2}
\end{proof}
}%

\begin{definitn}
Sei $X \subseteq \mLinStet(H)$.
\begin{enumaufz}
\item Eine Abbildung $p : \mLinStet(H) \rightarrow \mbbR$ hei"st \defemph{$X$-Halbnorm},
  falls eine $X$"=Abbildung $\varphi : \mLinStet(H) \rightarrow \mLinStet(H)$ so existiert,
  \dass gilt: $p(x) = \norm{\varphi(x)}$ f"ur alle $x \in \mLinStet(H)$.
\item Sei $\varphi : \mLinStet(H) \to \mLinStet(H)$ eine $X$-Abbildung.
  Man nennt die Halbnorm
\[ p_\varphi : \mLinStet(H) \rightarrow \mbbR, x \mapsto \norm{\varphi(x)} \]
die zu $\varphi$ geh"orige $X$-Halbnorm.
\index[B]{XHalbnorm@$X$-Halbnorm}%
\end{enumaufz}
\end{definitn}

Wie man leicht sieht, gilt:

\begin{bemerkung}
Sei $X \subseteq \mLinStet(H)$.
Auf der Menge der $X$-Halbnormen wird durch
\[ p \leq q  \overset{\text{def}}{\iff}  p(x) \leq q(x) \text{ f"ur alle } x \in \mLinStet(H) \]
eine Ordnungsrelation definiert.\dremark{$\leq$ ist i.a. nicht konnex}
\end{bemerkung}

\begin{satz}[\cite{Paulsen02ComplBoundedMaps}, Proposition 15.3.]\label{Paulsen15.3:Opraum}
Sei $X \subseteq \mLinStet(H)$ ein Operatorraum.
Dann existiert eine bez"uglich \glqq$\leq$\grqq{} minimale $X$-Halbnorm auf $\mLinStet(H)$.
\dliter{\cite{Paulsen02ComplBoundedMaps}, Proposition 15.3}%
\end{satz}

\dremark{
\begin{proof}
Sei $\Lambda$ eine Kette von $F$-Halbnormen.
F"ur alle $\lambda \in \Lambda$ findet man eine $F$-Abbildung
$\varphi_\lambda : \mLinStet(H) \rightarrow \mLinStet(H)$ mit $\lambda = p_{\varphi_\lambda}$.
Auf $\Lambda$ wird wie folgt eine Ordnung $\leq_1$ definiert:
\[ \lambda_1 \leq_1 \lambda_2
\overset{\text{def}}{\iff}  \lambda_1 \geq \lambda_2. \]
Dann gilt f"ur alle $\lambda, \mu \in \Lambda$:
$\lambda \leq_1 \mu \Rightarrow p_{\varphi_\lambda} = \lambda \geq \mu = p_{\varphi_\mu}$.
Somit ist $(p_{\varphi_\lambda})_{\lambda \in \Lambda}$ ein absteigendes Netz.

Definiere\dremark{Siehe \cite{PaulsenComplBoundedMapsPS}, S. 92}
\[
CB_1(\mLinStet(H),H)
= \{ T \in \mLinStet(\mLinStet(H),\mLinStet(H)) \setfdg \norm{T}_\text{cb} \leq 1 \}.
\]
Nach \cite{PaulsenComplBoundedMapsPS}, 7.4, ist $CB_1(\mLinStet(H),H)$
kompakt in der BW-Topologie.
\dremark{Def. in \cite{PaulsenComplBoundedMapsPS}, S. 92}%
Somit findet man ein Teilnetz $(\varphi_{\lambda_\mu})_{\mu \in M}$ von $(\varphi_\lambda)_\lambda$,
welches gegen $\varphi_0 \in CB_1(\mLinStet(H),H)$
in der BW-Topologie konvergiert.
Offensichtlich ist $\varphi_0$ eine $F$-Abbildung.
\dremark{Da $(\varphi_{\lambda_\mu})_{\mu \in M}$ konvergiert, ist dieses Netz beschr"ankt.
Sei $x \in F$, seien $h,k \in H$. Es gilt:
$\skalpr{\varphi_0(x)h}{k}
= \lim \skalpr{\varphi_{\lambda_\mu}(x)h}{k}
= \lim \skalpr{x(h)}{k}
= \skalpr{x(h)}{k}
\,\Rightarrow\,  \skalpr{(\varphi_0(x)-x)h}{k}=0
\,\overset{\cite{PaulsenComplBoundedMapsPS}, 7.3}{\Rightarrow} \varphi_0(x)=x$}%
F"ur alle $\lambda \in \Lambda$ setze $p_\lambda := p_{\varphi_\lambda}$.
Sei $x \in \mLinStet(H)$, $h \in H$.
F"ur alle $k \in H$ folgt mit der Cauchy-Schwarzschen Ungleichung:
\[ \abs{\skalpr{\varphi_0(x)h}{k}}
=  \lim_\mu \abs{\skalpr{\varphi_{\lambda_\mu}(x)h}{k}}
=  \liminf_\mu \abs{\skalpr{\varphi_{\lambda_\mu}(x)h}{k}}
\leq \liminf_\mu \norm{\varphi_{\lambda_\mu}(x)h} \cdot \norm{k}. \]

Mit $k = \varphi_0(x)h$ erh"alt man:
$\norm{\varphi_0(x)h}^2
=  \skalpr{\varphi_0(x)h}{\varphi_0(x)h}
\leq  \liminf_\mu \norm{\varphi_{\lambda_\mu}(x)h} \cdot \norm{\varphi_0(x)h}$,
also $\norm{\varphi_0(x)} \leq \liminf_\mu \norm{\varphi_{\lambda_\mu}(x)}$.
Da $(\varphi_{\lambda_\mu})_{\mu \in M}$ ein Teilnetz ist,
findet man ein $\tilde{\mu} \in M$ mit $\lambda \leq_1 \lambda_{\tilde{\mu}}$.
\dremark{Sei $g : M \rightarrow \Lambda, \mu \mapsto \lambda_\mu$.
Dann ist $g(M)$ kofinal in $\Lambda$}%
Also gilt: $p_\lambda \geq p_{\lambda_{\tilde{\mu}}}$.
Da $(p_\lambda)_{\lambda \in \Lambda}$ eine absteigende Kette ist,
gilt insbesondere f"ur alle $\hat{\mu} \in M$, $\hat{\mu} \geq \tilde{\mu}$
\begin{gather*}
   \norm{\varphi_{\lambda_{\tilde{\mu}}}(x)}
=  p_{\lambda_{\tilde{\mu}}}(x)
\geq  p_{\lambda_{\hat{\mu}}}(x)
=  \norm{\varphi_{\lambda_{\hat{\mu}}}(x)}, \quad\text{also} \\
     \norm{\varphi_{\lambda_{\tilde{\mu}}}(x)}
\geq \liminf_{\mu \in M} \norm{\varphi_{\lambda_\mu}(x)}.
\dremarkm{= \lim_{\mu \in M} \inf\{ \norm{\varphi_{\lambda_{\check{\mu}}}(x)} \setfdg
                                    \check{\mu} \in M, \check{\mu} \geq \mu \}}%
\end{gather*}
Es folgt
\[ p_\lambda(x)
\geq  p_{\lambda_{\tilde{\mu}}}(x)
=  \norm{\varphi_{\lambda_{\tilde{\mu}}}(x)}
\geq  \liminf_{\mu \in M} \norm{\varphi_{\lambda_\mu}(x)}
\geq  \norm{\varphi_0(x)}
=  p_{\varphi_0}(x), \]
also $p_\lambda \geq p_{\varphi_0}$.
Somit besitzt jede absteigende Kette von $F$-Halbnormen eine untere Schranke.
\dremark{n"amlich $p_{\varphi_0}$}%
Mit den Lemma von Zorn erh"alt man, \dass eine minimale $F$-Halbnorm existiert.
\dremark{Der Beweis ist deutlich ausf"uhrlicher als der in \cite{PaulsenComplBoundedMapsPS}}%
\dremark{30.11.'06, S. 1}%
\end{proof}

\dremark{Zum Beweis: Frage von W. Werner: Kann man den Beweis nicht auch mit $F$-Projektionen
  anstelle von $F$-Halbnormen und $F$-Abbildungen f"uhren?
  Nach \cite{PaulsenComplBoundedMapsPS}, nach Th. 15.4, ist nicht klar, \dass gilt:
  $\varphi$ minimale $F$-Projektion $\Rightarrow$ $p_\varphi$ minimale $F$-Halbnorm.
  Die Eigenschaft, \dass $p_\varphi$ minimale $F$-Halbnorm ist, wird aber im Beweis
  der Existenz der injektiven H"ulle (\cite{PaulsenComplBoundedMapsPS}, Th. 15.4)
  verwendet.}%
}%

\dremark{
\begin{anmerkung}[zur BW-Topologie]
\dmarginpar{in Div}%
Seien $X$, $Y$ Banachr"aume.
Sei $x \in X$, $y \in Y$.
Definiere ein lineares Funktional auf $\mLinStet(X,Y^*)$ durch
$x \otimes y(L) = L(x)(y)$.
Sei $Z$ der abgeschlossene Aufspann von Elementartensoren in $\mLinStet(X,Y^*)^*$.

Dann gilt: $\mLinStet(X,Y^*)$ ist isometrisch isomorph zu $Z^*$ verm"oge der Dualit"at
\[ \skalpr{L}{x \otimes y} = x \otimes y(L). \]

Die schwach*-Topologie, die durch diese Identifizierung auf $\mLinStet(X,Y^*)$ induziert
wird, wird BW-Topologie genannt (f"ur \glqq bounded weak\grqq).
\dliter{\cite{Paulsen02ComplBoundedMaps}, S. 85 oben}%
\end{anmerkung}
}%


\begin{satz}[\cite{Paulsen02ComplBoundedMaps}, Theorem~15.4; Existenz einer injektiven H"ulle]\label{exInjHuelle}
Sei $X \subseteq \mLinStet(H)$ ein Operatorraum und
$\varphi : \mLinStet(H) \to \mLinStet(H)$ eine $X$-Abbildung derart,
\dass $p_\varphi$ eine bez"uglich \glqq$\leq$\grqq{} minimale $X$-Halbnorm ist.
Sei $\iota : X \to \varphi(\mLinStet(H)), x \mapsto x$.
Dann ist $\varphi$ eine bez"uglich \glqq$\preceq$\grqq{} minimale $X$-Projektion und
$(\varphi(\mLinStet(H)), \iota)$ eine injektive H"ulle von $X$.
\dliter{\cite{PaulsenComplBoundedMapsPS}, 15.4}%
\end{satz}

\dremww{
\begin{lemma}\label{RigideKonkretOpraum}
Sei $F \subseteq \mLinStet(H)$ ein Operatorraum und
$\varphi : \mLinStet(H) \rightarrow \mLinStet(H)$ eine $F$-Abbildung derart,
\dass $p_\varphi$ eine minimale $F$-Halbnorm ist.
Setze $Z := \varphi(\mLinStet(H))$.
Sei $\alpha : Z \rightarrow Z$ vollst"andig kontraktiv
mit $\alpha\restring_F = \Id_F$.
Dann gilt: $\alpha = \Id_Z$.
\dremark{$\varphi$ existiert nach \ref{Paulsen15.3:Opraum}}%
\dliter{\cite{PaulsenComplBoundedMapsPS}, 15.5}%
\dremark{Beachte zum Beweis: Da $\varphi$ eine Projektion ist,
also $\varphi\restring_{\Bild{\varphi}} = \Id_{\Bild{\varphi}}$ gilt, erh"alt man:
$ \gamma \circ \varphi
= \gamma \circ \varphi \circ \gamma \circ \varphi
= \gamma \circ \gamma \circ \varphi$.}%
\end{lemma}
}%

\begin{satz}[\cite{Paulsen02ComplBoundedMaps}, Theorem 15.6; Eindeutigkeit der injektiven H"ulle]\label{EindInjHuelleOpraum}\label{EindInjHuelleKonkreterOpraum}
Sei $X$ ein Operatorraum.
Seien $(Y_1,\kappa_1)$, $(Y_2,\kappa_2)$ injektive H"ullen
von $X$.
Dann l"a\cms{}t sich die Abbildung
$i : \kappa_1(X) \to \kappa_2(X)$, gegeben durch
$i(\kappa_1(x)) = \kappa_2(x)$ f"ur alle $x \in X$,
eindeutig zu einer vollst"andigen Isometrie
von $Y_1$ auf $Y_2$ fortsetzen.
\[ \xymatrix{
\,\,Y_1\,\, \ar@{>-->>}[rr]^{\cong} & & Y_2 \\
\kappa_1(X) \ar[rr]^{i} \ar[u]^\subseteq & & \kappa_2(X) \ar[u]_\subseteq \\
& X \ar[ul]^{\kappa_1}  \ar[ur]_{\kappa_2} &
} \]
\end{satz}

\dremark{
\begin{satz}[Eindeutigkeit der injektiven H"ulle]
Seien $(Y_1,\kappa_1)$, $(Y_2,\kappa_2)$ injektive H"ullen eines Operatorraumes $X$.
Dann findet man eine surjektive, vollst"andige Isometrie $\varphi : Y_1 \to Y_2$ derart,
\dass $\varphi \circ \kappa_1 = \kappa_2$ gilt
(\cite{BlecherLeMerdy04OpAlg}, Lemma 4.2.5).
\dliter{\cite{PaulsenComplBoundedMapsPS}, 15.6.}%
\[ \xymatrix{
Y_1 \ar[rr]^{\varphi} & & Y_2 \\
& X \ar[ul]^{\kappa_1}  \ar[ur]_{\kappa_2} &
} \]
\end{satz}

\begin{proof}
Nach \ref{Paulsen15.3:Opraum} findet man eine $F$-Abbildung
$\varphi : \mLinStet(H) \rightarrow \mLinStet(H)$ so, \dass $p_\varphi$ eine
minimale $F$-Halbnorm ist.
Setze $Z := \varphi(\mLinStet(H))$.
Falls wir zeigen k"onnen, \dass $\kappa_1$ eine eindeutige Fortsetzung zu einem
vollst"andig isometrischen Isomorphismus
$\gamma_1 : Z \arrowbij I_1$ auf $I_1$ hat,
dann folgt die Behauptung, da $\gamma_2 \circ \gamma_1^{-1}$ die Abbildung liefert.
\[ \xymatrix{
I_1 & \,\,Z\,\, \ar@{>->>}[l]_{\gamma_1} \ar@{>->>}[r]^{\gamma_2} & I_2 \\
& F \ar[ul]^{\kappa_1} \ar[u]_\subseteq \ar[ur]_{\kappa_2} &
} \]

\bewitemph{Existenz:}
Da $I_1$ injektiv ist, findet man eine \cmvk{e} Abbildung
$\gamma_1 : Z \rightarrow I_1$, die $\kappa_1$ fortsetzt.
Weil $Z$ injektiv ist, findet man eine \cmvk{e} Abbildung
$\beta : I_1 \rightarrow Z$ mit $\beta \circ \kappa_1 = \Id_F$.
\[ \xymatrix{
Z \ar[r]^{\gamma_1} & I_1 \ar[r]^\beta & Z \\
& F \ar[ul]^{\Id_F} \ar[u]_{\kappa_1} \ar[ur]_{\Id_F} &
} \]

Dann ist $\beta \circ \gamma_1 : Z \rightarrow Z$
\cmvk{} mit $(\beta \circ \gamma_1)\restring_F = \Id_F$.
\dremark{$\beta \circ \gamma_1(x) = \beta \circ \kappa_1(x) = x$ f. a. $x \in F$}%
Mit \ref{RigideKonkretOpraum} erh"alt man:
$\beta \circ \gamma_1 = \Id_Z$.


Da $\beta$ und $\gamma_1$ beide \cmvk{} sind, ist $\gamma_1$ eine vollst"andige Isometrie.
\dremark{$\norm{x} = \norm{\beta \circ \gamma_1(x)} \leq \norm{\gamma_1(x)} \leq \norm{x}$}%
Dann ist $\Bild(\gamma_1)$ ein injektiver Operatorunterraum von $I_1$,
wegen der Minimalit"at von $I_1$ also gleich $I_1$.
Somit ist $\gamma_1$ surjektiv auf $I_1$.

\bewitemph{Eindeutigkeit:}
Seien $\gamma_1, \delta_1 : Z \rightarrow I_1$
\cmvk{e} Fortsetzungen von $\kappa_1$.
Dann gilt: $\delta_1^{-1} \circ \gamma_1\restring_F = \Id_F$.
Mit \ref{RigideKonkretOpraum} folgt:
$\delta_1^{-1} \circ \gamma_1 = \Id_Z$.
Analog erh"alt man: $\gamma_1^{-1} \circ \delta_1 = \Id_Z$.
Also folgt: $\gamma_1 = \delta_1$.
\[ \xymatrix{
Z \ar[r]^{\gamma_1} & I_1 & Z \ar[l]_{\delta_1} \\
& F \ar[ul]^{\Id_F} \ar[u]_{\kappa_1} \ar[ur]_{\Id_F} &
} \]
\end{proof}
}%

\begin{notation}\label{notInjHuelleOR}
Nach Satz~\ref{exInjHuelle} besitzt jeder Operatorraum $X$ eine injektive H"ulle,
die nach Satz~\ref{EindInjHuelleKonkreterOpraum} bis auf vollst"andige
Isometrie eindeutig ist und mit $(I(X),j)$ oder kurz mit $I(X)$ bezeichnet wird.
\index[S]{IX@$I(X)$ (injektive Hülle)}%
\index[S]{j@$j$}%
\end{notation}

Mit Hilfe vollst"andig kontraktiver Projektionen
erh"alt man injektive Operatorr"aume:

\begin{bemerkung}[\cite{EffrosRuan00OperatorSpaces}, Proposition 4.1.6]\label{BildProjIstInjOpraum}
Sei $X \subseteq \mLinStet(H)$ ein injektiver Operatorraum.
\begin{enumaufz}
\item Ist $P : X \rightarrow X$ eine vollst"andig kontraktive
Projektion,\dremark{\dheisst{} $P \circ P = P$}
dann ist $P(X)$ ein injektiver Operatorraum.
\item Umgekehrt gilt: Man findet eine \cmvk{e} Projektion von $\mLinStet(H)$ auf $X$.
\end{enumaufz}
\end{bemerkung}

\dremww{
\begin{proof}
\bewitemph{(i):}
\dremark{Nach \cite{WernerFunkana3}, Lemma IV.6.1, ist $P(E)$ abgeschlossen,
  also ein Operatorraum.}%
Seien $F \subseteq G \subseteq \mLinStet(K)$ Operatorr"aume,
sei $\varphi : F \rightarrow P(E)$ vollst"andig kontraktiv.
Dann ist auch $\varphi : F \rightarrow E$ vollst"andig kontraktiv.
Da $E$ injektiv ist, findet man eine \cmvk{e} Fortsetzung $\Phi : G \to E$
von $\varphi$.
Dann ist $\tilde{\Phi} := P \circ \Phi : G \to P(E)$ die gesuchte \cmvk{e} Fortsetzung von $\varphi$.
\[
\xymatrix{
G \ar@{-->}[rr]^\Phi && E \ar[d]^P \\
F \ar[rr]^\varphi \ar[u]^\subseteq && P(E)
}
\]

\bewitemph{(ii):} Sei $E \subseteq \mLinStet(H)$ ein injektiver Operatorraum.
Dann findet man eine \cmvk{e} Fortsetzung  $P : \mLinStet(H) \to E$ von $\Id_E$.
Somit ist $P$ insbesondere eine Projektion.
\[
\xymatrix{
\mLinStet(H) \ar@{-->}[rd]^P & \\
E \ar[u]^\subseteq \ar[r]_{\Id_E} & E}
\]
\end{proof}
}%

Da $\mLinStet(H)$ injektiv ist
(\refb{\cite{BlecherLeMerdy04OpAlg}, Theorem 1.2.10,}{Beispiel~\ref{LHIstInjOR}}),
erh"alt man mit der obigen Proposition:\dremark{\ref{BildProjIstInjOpraum}}

\begin{beispiel}
Sei $H$ Hilbertraum.
Dann sind $H^c$ und $H^r$ injektive Operatorr"aume.
\end{beispiel}

\dremww{
\begin{folgerung}[Rigidit"at]
Sei $F$ ein Operatorraum und $(E, \kappa)$ eine injektive H"ulle von $F$.
Sei $\psi : E \to E$ \cmvk{} mit $\psi\restring_{\kappa(F)} = \Id\restring_{\kappa(F)}$.
Dann gilt: $\psi = \Id_E$.
\dliter{\cite{PaulsenComplBoundedMapsPS}, 15.7}%
\end{folgerung}

\dremark{Man erh"alt alle Aussagen auch f"ur unitale Operatorsysteme,
  indem man \cite{PaulsenComplBoundedMapsPS}, Prop. 15.1 und S. 224 unten verwendet.}%

\begin{bemerkung}
Sei $E$ ein endlich-dimensionaler injektiver Operatorraum.
Dann findet man $p \in \mbbN$ und $m_k, n_k \in \haken{p}$ so, \dass $E$
vollst"andig isometrisch ist zu
\dliter{\cite{EffrosRuan00OperatorSpaces}, 6.1.8;
  siehe auch \cite{Pisier03OpSpaceTheory}, S. 357, Remark 1}%
\[ M_{m_1,n_1} \oplus \dots \oplus M_{m_p,n_p}. \]
\end{bemerkung}

\begin{bemerkung}
Sei $\mfrakA$ eine \csalgebra{} ohne Einselement und
$\tilde{\mfrakA}$ eine Unitalisierung von $\mfrakA$.
Dann ist die injektive H"ulle im Operatorraum-Sinne $I(\mfrakA)$ von $\mfrakA$ eine
unitale \csalgebra.
Ferner gilt: $I(\mfrakA) = I(\tilde{\mfrakA})$.
\dliter{\cite{BlecherPaulsen00MultOfOpSpaces}, Prop. 2.8}%
\end{bemerkung}
}%

\begin{bemerkung}[\cite{BlecherLeMerdy04OpAlg}, 4.2.10]\label{injHuelleMnX}
Seien $m,n \in \mbbN$, und sei $X$ ein Operatorraum.
Dann gilt:
\[ M_{m,n}(I(X)) \cong I(M_{m,n}(X)), \]
\dheisst $M_{m,n}(I(X))$ ist eine injektive H"ulle von $M_{m,n}(X)$.
\dremark{Gilt dies auch f"ur unitale Operatorsysteme?}%
\dliter{\cite{BlecherEffrosZarikian02}, Lemma 4.3}%
\end{bemerkung}

\begin{definitn}\label{defUnitalerOpraum}
Ein Operatorraum $X$ hei"st \defemph{unital}, falls er ein ausgezeichnetes Element
$e_X$, genannt Einselement von $X$, so besitzt,
\dass eine unitale \csalgebra $\mfrakA$ und eine vollst"andige Isometrie $u : X \to \mfrakA$
existieren mit $u(e_X) = e_\mfrakA$.
\dremark{Die Definition von $e_X$ ist konsistent mit der Definition von $e_Y$ f"ur ein
  unitales Operatorsystem $Y \subseteq \mLinStet(H)$,
  denn $u : Y \to \mLinStet(H)$ ist eine vollst"andige Isometrie mit $u(e_Y) = \Id_H$.}%
\index[B]{unitaler Operatorraum}%
\index[S]{ekleinX@$e_X$ (Einselement)}%
\dliter{\cite{BlecherLeMerdy04OpAlg}, 1.3.1}%
\end{definitn}

\begin{anmerkung}\label{AnmUnitalerORunabh}
Die obige Definition h"angt nicht von einer speziellen unitalen \csalgebra $\mfrakA$ ab:
Sei $\mfrakB$ eine weitere unitale \csalgebra, $v : X \to \mfrakB$
eine vollst"andige Isometrie mit $v(e_X) = e_\mfrakB$.
Dann existiert ein eindeutiger vollst"andiger Ordnungsisomorphismus
von dem unitalen Operatorsystem $u(X) + u(X)^*$ auf das unitale Operatorsystem $v(X) + v(X)^*$,
der die vollst"andige Isometrie $v \circ u^{-1} : u(X) \to v(X)$ fortsetzt.
\dliter{\cite{Blecher01ShilovBoundary}, S. $4^{17}$}%
\end{anmerkung}

Eine abstrakte Charakterisierung unitaler Operatorr"aume wird
in \cite{BlecherNeal08MetricChar}, Theorem 1.1,
und in \cite{HuangNg08AbstractCharUnitalOS} aufgef"uhrt.
\dremark{Theorem 2.9}%
\skiptext

\dremark{Unitalisierung (unitization) eines Operatorraumes:
  Siehe \cite{Pisier03OpSpaceTheory}, S. 163 und Prop. 8.19.}%

F"ur die injektive H"ulle eines unitalen Operatorraumes gilt:

\begin{bemerkung}[\cite{BlecherLeMerdy04OpAlg}, Corollary 4.2.8]\label{unitalerORFolgtjUnital}
Sei $X$ ein unitaler Operatorraum.
Dann gibt es eine injektive H"ulle $(I(X),j)$ von $X$ so,
\dass $I(X)$ eine unitale \csalgebra{} ist und $j$ eine unitale Abbildung.
\end{bemerkung}

In der obigen Proposition kann man nicht auf die Voraussetzung verzichten,
\dass $X$ unital ist, wie das folgende Beispiel zeigt:

\begin{beispiel}[\cite{EffrosRuan00OperatorSpaces}, S. 108]
Setze $E_{11} := \begin{smallpmatrix} 1 & 0 \\ 0 & 0 \end{smallpmatrix}$.
Die Abbildung
\[ \Phi : M_2 \to M_2, a \mapsto E_{11}a, \]
ist eine \cmvk{e} Projektion von $M_2$ auf den Operatorraum
\[ E := E_{11} M_2 = \left\{ \begin{pmatrix} \alpha & \beta \\ 0 & 0 \end{pmatrix}
  \setfdg \alpha, \beta \in \mbbC \right\}, \]
somit ist $E$ injektiv.
Ferner ist $E$ nicht vollst"andig isometrisch zu einer \csalgebra.
\dremark{Im Gegensatz zu einem Operatorsystem, vgl. \cite{EffrosRuan00OperatorSpaces}, 6.1.3}%
\end{beispiel}

\dremark{Beweisskizze. Annahme: $E$ ist vollst"andig isomorph zu der \csalgebra $\mfrakA$.
  Wegen $\dim \mfrakA = 2$ folgt: $\mfrakA = \ell^\infty_2$.
  Somit ist jede kontraktive Abbildung nach $E$ \cmvk{}.
  Die Abbildung $\mathbf{t} : M_2 \to M_2, a \mapsto a^t$ ist eine Isometrie.
  $\varphi := \Phi \circ \mathbf{t} : M_2 \to E$ ist kontraktiv, aber nicht \cmvk{}.}%

\dremww{
\begin{bemerkung}
Sei $\mfrakA$ eine \csalgebra ohne Einselement und
$\tilde{\mfrakA}$ eine Unitalisierung von $\mfrakA$.
Dann gilt $I(\mfrakA) = I(\tilde{\mfrakA})$, wobei jeweils die
injektive H"ulle im Operatorraumsinne gebildet wird.
\dliter{\cite{BlecherPaulsen00MultOfOpSpaces}, Prop. 2.8}%
\end{bemerkung}

\begin{beispiel}
Sei $\mcalT(\mbbT)$ der Raum aller Toeplitz-Operatoren auf dem zur Einheitskreisscheibe
assoziierten Hardy-Raum $H^2(\mbbT)$.
Dann ist $\mcalT(\mbbT)$ ein injektiver Operatorraum.
\dliter{\cite{BeltitaPrunaru07AmenabilitycbProj}, Cor. 3.7}%
\end{beispiel}
}%

\section{Tern"are Ringe von Operatoren}

Wir wiederholen in diesem Abschnitt Grundlagen der Theorie der tern"aren Ringe
von Operatoren, kurz TROs.
Jeder TRO ist ein \hcsmodul, und umgekehrt kann
jeder \hcsmodul treu als TRO dargestellt werden.
Des weiteren kann man jeden TRO $Z$ nicht-ausgeartet einbetten,
\dheisst man findet Hilbertr"aume $H$ und $K$ mit $\erzl{ZH}$ dicht in $K$
und $\erzl{Z^* K}$ dicht in $H$.
Au"serdem ist jeder Operatorraum in einem kleinsten TRO enthalten,
der sogenannten tern"aren H"ulle.

\begin{definitn}
Seien $H, K$ Hilbertr"aume.
\begin{enumaufz}
\item Unter einem \defemph{ternären Ring von Operatoren} (kurz \defemphi{TRO})
versteht man einen abgeschlossenen Untervektorraum $Z$ von $\mLinStet(H,K)$
(oder einer \csalgebra) mit der Eigenschaft: $ZZ^*Z \subseteq Z$.
F"ur alle $x,y,z \in Z$ schreibt man anstelle von $xy^*z$ auch $[x,y,z]$,
dieses Produkt wird \defemphi{Tripelprodukt} auf $Z$ genannt.
\dremark{Wird bei \cite{Zet83} ohne Abschlu\cms{} definiert.}%
\item Ein \defemphi{Untertripel} eines TRO $Z$ ist ein abgeschlossener
Untervektorraum $Y$ von $Z$, der $YY^*Y \subseteq Y$ erf"ullt.
\item Ein \defemphi{Tripelmorphismus} zwischen TROs $Y$ und $Z$ ist
eine lineare Abbildung $T : Y \to Z$ mit der Eigenschaft:
$T([x,y,z]) = [Tx,Ty,Tz]$ f"ur alle $x,y,z \in Y$.
\index[B]{ternärer Ring von Operatoren}%
\dliter{\cite{BlecherLeMerdy04OpAlg}, 4.4.1}%
\end{enumaufz}
\end{definitn}

In einem TRO $Z$ gilt bereits $ZZ^*Z = Z$ (\cite{Exel97TwistedPartialActions}, Corollary 4.10).
\skiptext

In \cite{NeR03} werden TROs mit Hilfe von Operatorr"aumen und beschr"ankten
symmetrischen Gebieten charakterisiert (siehe Theorem 5.3).
\dremark{also intrinsisch charakterisiert}%

\dmarginpar{DB}\dremark{
Dies "uberarbeiten: Alle TROs sind von dieser Gestalt.
\begin{beispiel}\label{bspTROCsAlg}
Ein Hauptbeispiel f"ur einen TRO ist $p \mfrakA (1-p)$,
wobei $\mfrakA$ eine \csalgebra und $p$ eine Projektion in $\mfrakA$ oder $M(\mfrakA)$ sei
(vgl. \cite{BlecherLeMerdy04OpAlg}, 4.4.1).
\end{beispiel}
}%
\dremark{Beispiel: Obere Dreiecksmatrizen; $H^c$: $H^c$ ist injektiv, also TRO.}%
\dremark{Weitere Beispiele? Siehe \cmzB \cite{BlecherLeMerdy04OpAlg}, 4.4.7, am Ende ($C^*_e(X)$)}%

\dremark{Einen Teil der Tripelmorphismen erh"alt man auf nat"urliche Art und Weise:
Sei $Z$ ein Untertripel einer beliebigen \csalgebra $\mfrakA$
und $\varphi : \mfrakA \to \mfrakA$ ein \sterns{}Homomorphismus mit $\varphi(Z) \subseteq Z$.
Schr"ankt man $\varphi$ auf $Z$ ein, so erh"alt man einen Tripelmorphismus.}%

\begin{bemerkung}[Vgl. \cite{BlecherLeMerdy04OpAlg}, Corollary 4.4.6 und Lemma 8.3.2]\label{TripelisomIffSurjVIso}
Seien $Z_1$, $Z_2$ TROs.
Dann ist $\varphi : Z_1 \to Z_2$ genau dann ein Tripelisomorphismus,
wenn $\varphi$ eine surjektive vollst"andige Isometrie ist.
\end{bemerkung}

\begin{bemerkung}[Vgl. \cite{Zet83}, S. 123]
Sei $Z$ ein TRO.
Dann sind $ZZ^*$ und $Z^*Z$ \csalgebren.
Ferner ist $Z$ ein voller \hcsmodul "uber $Z^*Z$ verm"oge
$\skalpri{x}{y}{Z^*Z} := x^*y$.
\dremark{Ferner ist $Z$ ein Hilbert-$ZZ^*$-$Z^*Z$-Bimodul verm"oge
  $\skalpril{x}{y}{ZZ^*} := xy^*$ und $\skalpri{x}{y}{Z^*Z} := x^*y$.}%
\dremark{Die Norm des Rechts-Hilbert-\cstern{}Moduls $Z$
  ist gleich der Operatornorm von $L(H)$.}%
\dremark{$Z^*Z$ ist gleich der von $\vrerz{ x^*y \setfdg x,y \in Z }$
  in $L(H)$ erzeugten \csalgebra{}.
  Es ist $xy = x \cdot_{L(H)} y$.
  Alternativ kann man anstelle von $x^*y$ auch die
  Abbildung $\alpha_{z,y} : Z \to Z, x \mapsto [x,y,z]$ betrachten
  und die von $\vrerz{ \alpha_{z,y} \setfdg z,y \in Z }$ in $L(Z)$
  erzeugte \csalgebra, die man mit $Z^*Z$ identifizieren kann.
  Dies wird \cmzB{} in tern"aren \cstern{}Ringen gemacht,
  siehe beispielsweise \cite{Zet83}, S. 118 und Prop. 3.2.}%
\dliter{vgl. \cite{BlecherLeMerdy04OpAlg}, 8.1.2, S. 298 oben}%
\dremark{Z. 21.5.'08/1}%
\end{bemerkung}

Umgekehrt kann jeder Hilbert-\cstern{}Modul treu als TRO dargestellt werden
(\cite{Zet83}, Theorem 2.6).\dliter{\cite{BlecherLeMerdy04OpAlg}, 8.2.8}
Alternativ kann man die Linking-Algebra eines Hil\-bert-\cstern{}Mo\-duls $E$ bilden.
Dann ist $E$ als rechte obere Ecke einer \csalgebra{} ein TRO.
\dliter{\cite{BlecherLeMerdy04OpAlg}, 8.1.19, 2. Absatz}%

Die Morphismen zwischen TROs sind die Tripelmorphismen.
Hingegen sind die Morphismen zwischen \hcsmoduln die adjungierbaren Abbildungen.
\dliter{\cite{WWerner08MorphForHCsManifolds},1.4; \cite{Wegge-OlsenKTheory}, Def. 15.1.5}%
\skiptext

"Ahnlich wie bei der nicht-ausgearteten Darstellung einer \csalgebra
kann man auch jeden TRO nicht-ausgeartet einbetten:

\begin{bemerkung}[Vgl. \cite{Harris81GeneralizationOfCsAlg}, S. 341]\label{TROnichtAusgeartet}
Sei $Z \subseteq \mLinStet(H,K)$ ein TRO,
$H_1 := \{ \xi \in H \setfdg Z \cdot \xi = \{0\} \}$
und $K_1 := \{ \eta \in K \setfdg Z^* \cdot \eta = \{0\} \}$.
\begin{enumaufz}
\item Es ist $\erzl{Z H_1^\perp}$ dicht in $K_1^\perp$
und $\erzl{Z^* K_1^\perp}$ dicht in $H_1^\perp$.
\item Es ist
$W := \left\{ z\restring_{H_1^\perp} \setfdg z \in Z \right\} \subseteq \mLinStet(H_1^\perp, K_1^\perp)$
ein TRO mit $\erzl{W H_1^\perp}$ dicht in $K_1^\perp$ und
$\erzl{W^* K_1^\perp}$ dicht in $H_1^\perp$.
\item $\varphi : Z \arrowbij W, z \mapsto z\restring_{H_1^\perp}$, ist
ein tern"arer Isomorphismus auf $W$.
\dremark{Ist $W$ eindeutig? Unklar.}%
\dremark{Resultat ist allgemein bekannt,
  siehe \cmzB \cite{Harris81GeneralizationOfCsAlg}, S. 341 unten,
  \cite{KaurRuan02LocalPropOfTROs}, S. 269.}%
\dremark{Hierzu etwas schreiben, \cmzB Vergleich mit nicht-ausgearteter Darstellung
  einer \csalgebra.}%
\dremark{Z. 2.6.'08/2}%
\end{enumaufz}
\end{bemerkung}

\dremww{
\begin{proof} \bewitemph{(iii):}
F"ur alle $z \in Z$ gilt:
\[ \norm{z}_Z
=  \sup\{ \norm{z(\xi)} \setfdg \xi \in H, \norm{\xi} \leq 1 \}
=  \sup\{ \norm{z(\xi)} \setfdg \xi \in H_1^\perp, \norm{\xi} \leq 1 \}
=  \norm{z\restring_{H_1^\perp}}_W
=  \norm{\varphi(z)}_W. \]
Analog folgt, \dass $\varphi$ vollst"andig isometrisch ist.
Nach \ref{TripelepiIffVollstIsometr} ist $\varphi$ ein tern"arer Morphismus.
\end{proof}

\begin{bemerkung}[Harris/Kaup]\label{TripelepiIffVollstIsometr}
Seien $Y,W$ TROs, sei $\Phi : Y \to W$ ein Tripelmorphismus.
\begin{enumaufz}
\item $\Phi$ ist vollst"andig kontraktiv.
\item $\Phi$ ist genau dann eine vollst"andige Isometrie, wenn $\Phi$ surjektiv auf $W$ ist.
\dliter{\cite{BlecherLeMerdy04OpAlg}, 8.3.2}%
\end{enumaufz}
\end{bemerkung}
}%

\begin{bemerkung}[\cite{BlecherLeMerdy04OpAlg}, 4.4.5.(1)]\label{XTRO:TProdGleich}
Sei $X$ ein TRO.
Dann ist die Einbettung $j : X \to I(X)$ ein Tripelmorphismus.
\dremark{Ohne $j$: $[\cdot,\cdot\cdot,\cdots]_X$ und
$[\cdot,\cdot\cdot,\cdots]_{I(X)}$ stimmen "uberein.}%
\end{bemerkung}

\dremark{
Definiere $\sigma : X \to \begin{pmatrix} 0 & X \\ 0 & 0 \end{pmatrix},
                    x \mapsto \begin{pmatrix} 0 & x \\ 0 & 0 \end{pmatrix}$.
Dann gilt f"ur alle $x,y,z \in X$:
\begin{align*}
   [\sigma(x),\sigma(y),\sigma(z)]
&= \sigma(x) \sigma(y)^* \sigma(z)
=  \sigma(xy^*z)
=  \sigma([x,y,z]).
\begin{pmatrix} \end{pmatrix}
\end{align*}
Wegen $\begin{pmatrix} 0 & X \\ 0 & 0 \end{pmatrix} \cong  X$ folgt die Behauptung.
}%

\begin{definitn}
\begin{enumaufz}
\item Ein \defemphi{Tripelsystem} ist ein Operatorraum $Y$,
der eine Abbildung $[\cdot,\cdot\cdot,\cdots] : Y \times Y \times Y \to Y$,
genannt \defemphi{Tripelprodukt}, derart besitzt,
\dass ein TRO $Z$ und eine vollst"andige Isometrie $\Phi : Y \to Z$ existieren mit der Eigenschaft,
\dass $\Phi$ ein \defemphi{Tripelmorphismus} ist,
\dheisst linear ist und folgendes erf"ullt:
\[ \Phi([x,y,z]) = [\Phi(x),\Phi(y),\Phi(z)] \]
f"ur alle $x,y,z \in Y$.
\dremark{Man kann o.B.d.A. (siehe \ref{BildTripelmorIstAbg}) verlangen,
  \dass $\Phi$ surjektiv auf $Z$ ist.}%

\item Unter einem \defemphi{Untertripel} eines Tripelsystems $Z$ versteht man einen
abgeschlossenen Unterraum von $Z$, der unter dem Tripelprodukt abgeschlossen ist.
\dliter{\cite{BlecherLeMerdy04OpAlg}, 8.3.1}%
\end{enumaufz}
\end{definitn}

Insbesondere sind TROs Tripelsysteme.
Ferner ist ein Untertripel eines Tripelsystems wieder ein Tripelsystem.

\begin{bemerkung}[Vgl. \cite{Hamana99TripleEnv}, Proposition~2.1.(iii)]\label{BildTripelmorIstAbg}
Sei $Y$ ein Tripelsystem, $Z$ ein TRO und $\Phi : Y \to Z$ ein Tripelmorphismus.
Dann ist $\Phi(Y)$ ein TRO.
\end{bemerkung}

Mit der obigen Proposition\dremark{\ref{BildTripelmorIstAbg}} und
Proposition~\ref{TripelisomIffSurjVIso} folgt,
\dass es auf einem Operatorraum $Y$ h"ochstens ein Tripelprodukt derart geben kann,
\dass $Y$ ein Tripelsystem ist.
\dremark{Seien $Y_1$, $Y_2$ Tripelsysteme derart, \dass die Operatorr"aume $Y_1$, $Y$ und $Y_2$
  "ubereinstimmen.
  Nach \ref{BildTripelmorIstAbg} findet man vollst"andig isometrische Tripelisomorphismen
  $\Phi_1 : Y_1 \arrowbij Z_1$, $\Phi_2 : Y_2 \arrowbij Z_2$.
  Insbesondere sind $Y$, $Z_1$ und $Z_2$ vollst"andig isomorph.
  Mit \ref{TripelisomIffSurjVIso} folgt, \dass $Z_1$ tripelisomorph zu $Z_2$ ist.
  Also ist $Y_1$ tripelisomorph zu $Y_2$.}%

\begin{definitn}
Sei $X$ ein Operatorraum.
Eine \defemphi{Tripelerweiterung} von $X$ ist ein Paar,
bestehend aus einem Tripelsystem $Y$ und
einer vollst"andigen Isometrie $\kappa : X \to Y$ derart,
\dass $Y$ das kleinste Untertripel von $Y$ ist, welches $\kappa(X)$ enth"alt.
\dremark{$Y$ ist hierdurch noch keine tern"are H"ulle,
  siehe \cite{BlecherLeMerdy04OpAlg}, 8.3.8}%
\dliter{\cite{BlecherLeMerdy04OpAlg}, 8.3.8}%
\end{definitn}

\begin{defSatz}[Vgl. \cite{BlecherLeMerdy04OpAlg}, 8.3.9]\label{defTernHuelle}
Sei $X$ ein Operatorraum.
\begin{enumaufz}
\item Es existiert eine Tripelerweiterung $(Y,\kappa)$ von $X$ mit der folgenden
universellen Eigenschaft:
F"ur jede Tripelerweiterung $(Z,\rho)$ von $X$ existiert ein (notwendig eindeutiger und surjektiver)
Tripelmorphismus $\Phi : Z \to Y$ derart, \dass $\Phi \circ \rho = \kappa$ gilt.
\[ \xymatrix{
Z \ar@{-->}[rr]^\Phi && Y \\
& X \ar[ul]^{\rho} \ar[ur]_\kappa}
\]

\item Eine Tripelerweiterung $(Y,\kappa)$ wie oben ist bis auf Tripelisomorphismen eindeutig,
hei"st \defemphi{ternäre Hülle} von $X$
und wird mit $(\ternh{X},j)$ bezeichnet.
\index[S]{j@$j$}%
\index[S]{TX@$\ternh{X}$ (ternäre Hülle)}%
\dremark{\item $\mcalT(X)$ ist gleich dem von $j(X)$ in $I(X)$ erzeugten Untertripel.}%
\dremark{Wegen \ref{TXeqSpan} stimmt dieses $j$ mit dem aus \ref{notInjHuelleOR} "uberein.}%
\end{enumaufz}
\end{defSatz}

Man kann die tern"are H"ulle wie folgt beschreiben:

\begin{bemerkung}[Vgl. \cite{BlecherLeMerdy04OpAlg}, 8.3.8]\label{TXeqSpan}
Sei $X$ ein Operatorraum.
Dann gilt:
\begin{equation}\label{eqTXeqSpan}
   \ternh{X}
\cong  \overline{\cmlin}\left\{x_1 \multis x_2^* \multis x_3 \multis x_4^* \cdots x_{2n+1} \setfdg
     n \in \mbbN_0, x \in j(X)^{2n+1} \right\} \subseteq I(X),
\end{equation}
wobei $j : X \to I(X)$ wie in Notation~\ref{notInjHuelleOR} sei.
Insbesondere folgt:
\begin{equation}\label{ITXeqIX}
I(\ternh{X}) \cong I(X).
\end{equation}
\end{bemerkung}

\dremark{Beweis.
Offensichtlich gilt \glqq$\supseteq$\grqq.
Au"serdem ist $\overline{\cmlin}\{ \dots \}$ ein Untertripel, welches $j(X)$ enth"alt.}%

Wir betrachten $j : X \to \ternh{X}$ wie in Definitions-Satz \ref{defTernHuelle}
und $j : X \to I(X)$ wie in Notation~\ref{notInjHuelleOR}.
Wegen der obigen Proposition\dremark{\ref{TXeqSpan}} ist es gerechtfertigt,
f"ur beide Abbildungen dieselbe Bezeichnung zu benutzen.
\skiptext

Die tern"are H"ulle besitzt die folgenden Eigenschaften:

\begin{bemerkung}[Vgl. \cite{BlecherLeMerdy04OpAlg}, 4.4.7]\label{CsHuelleEqTernHuelle}
\begin{enumaufz}
\item Sei $X$ ein TRO. Dann gilt: $\mcalT(X) \cong X$.
\item Sei $X$ ein unitaler Operatorraum.
Dann ist $\mcalT(X)$ eine unitale \csalgebra, die $j(X)$ als unitalen Untervektorraum enth"alt.
Ferner ist $\mcalT(X)$ isomorph zu der von $j(X)$ erzeugten \cstern{}Unteralgebra in $I(X)$.
\dremark{Au"serdem ist $C^*_\text{env}(X)$ eine tern"are H"ulle von $X$.}%
\dliter{\cite{BlecherLeMerdy04OpAlg}, 4.4.7 unten}%
\end{enumaufz}
\end{bemerkung}

\dremark{Beweis.
(i) folgt mit \ref{XTRO:TProdGleich} und \eqref{eqTXeqSpan}.

(ii): (I): $X+X^*$ ist ein Operatorsystem, also ist $I(X+X^*)$ eine unitale \csalgebra.
Ferner gilt: $I(X) = I(X+X^*)$.

(II): Setze in \eqref{eqTXeqSpan} $x_{2i} = e$ f"ur alle $i \in \haken{n}$.
}%

\dremww{
\begin{beispiel}
Seien $\Omega$ und $A$ wie in \ref{bspPaulsenS223}.
Nach \ref{bspPaulsenS223} ist $A$ ein unitales Operatorsystem mit \cstern{}H"ulle $C(\Omega)$.
Nach \ref{CsHuelleEqTernHuelle} gilt $\ternh{X} = C(\Omega)$.
Da $\Omega$ nicht extrem unzusammenh"angend ist, ist $C(\Omega)$ nicht injektiv.
Insbesondere gilt: $\ternh{X} \neq \injh{X}$.
\dremark{Beispiel daf"ur, \dass $\ternh{X} \neq \injh{X}$.}%
\end{beispiel}

\begin{bemerkung}
Sei $X$ ein Operatorraum, $T : X \to X$ linear und $B := \mcalT(X)^*\mcalT(X)$.
Dann sind die folgenden Aussagen "aquivalent:
\begin{enumaequiv}
\item $T \in M_\ell(X)$ mit $\norm{T}_{M_\ell(X)} \leq 1$.
\item $T$ ist die Einschr"ankung auf $X$ einer (notwendig eindeutigen) kontraktiven
  $B$-Rechtsmodulabbildung $u$ auf $\mcalT(X)$ mit $u(X) \subseteq X$.
\dliter{\cite{BlecherLeMerdy04OpAlg}, 4.5.2}%
\end{enumaequiv}
\end{bemerkung}
}%

\dremark{Beweis.
\bewitemph{\glqq$\Rightarrow$\grqq:}
Man findet $a \in I_{11}$ mit $Tx = ax$ f"ur alle $x \in X$.
Definiere $u : \mcalT(X) \to \mcalT(X), z \mapsto az$.

Rechtsmodulabbildung: F"ur alle $x,y,z \in \mcalT(X)$ gilt:
$u(x(y^*z)) = a(x(y^*z)) = (ax)(y^*z) = u(x)(y^*z)$.

Kontraktiv: $\norm{u(z)} = \norm{az} \leq \norm{a} \cdot \norm{z} \leq \norm{z}$.

\bewitemph{\glqq$\Leftarrow$\grqq:}
Benutzt wird: $u(z)^*u(z) \leq z^*z$ f"ur alle $z \in \mcalT(X)$.
}%

\dremark{Literatur zu TROs:
  \cite{Hamana99TripleEnv}, \cite{Exel97TwistedPartialActions}, §§4--6,
  \cite{Zet83}}%

\section{Selbstadjungierte Operatorr"aume}

In diesem Abschnitt notieren wir die Definition eines selbstadjungierten
Operatorraumes.
Au"serdem erinnern wir an \sterns{}TROs und an die tern"are \sterns{}H"ulle.

\begin{definitn}
Seien $V$, $W$ Vektorr"aume mit Involution.
Eine Abbildung $\alpha : V \to W$ hei"st \defemph{\sterns{}linear},
falls $\alpha$ linear ist und $\alpha(v^*) = \alpha(v)^*$ f"ur alle $v \in V$ gilt.
\dremark{\cite{BlecherKNW07OrderedInvOpSpaces}, S. 2}%
\index[B]{sternlinear@\sterns{}linear}%
\end{definitn}

\begin{definitn}
Unter einem \defemphi{selbstadjungierter Operatorraum} (kurz \defemph{\cmsa Operatorraum})
versteht man einen Operatorraum $X$, versehen mit einer Involution $\cdot^* : X \to X$ derart,
\dass eine \sterns{}lineare, vollst"andige Isometrie $i : X \to L(H)$ existiert.
\dremark{\cite{BlecherKNW07OrderedInvOpSpaces}, S. 2}%
\index[B]{sa Operatorraum@\cmsa Operatorraum}%
\end{definitn}

\begin{bemerkung}[\cite{BlecherKNW07OrderedInvOpSpaces}, S. 2]
Sei $X$ ein Operatorraum.
Dann ist $X$ genau dann ein selbstadjungierter Operatorraum,
falls eine Involution $\cdot^* : X \to X$ derart existiert,
\dass gilt:
\begin{equation}\label{eqNormxjiEqNormx}
   \norm{ (x_{ji}^*)_{i,j} }_n = \norm{ x }_n
\end{equation}
f"ur alle $n \in \mbbN$ und $x \in M_n(X)$.
\end{bemerkung}

\begin{definitn}
\begin{enumaufz}
\item Einen selbstadjungierten TRO bezeichnet man auch als \defemph{\sterns{}TRO}.
\index[B]{sternTRO@\sterns{}TRO}%
\item Ein \defemph{\sterns{}Untertripel} eines \sterns{}TRO $Z$ ist ein
  abgeschlossener, selbstadjungierter Unterraum $Y$ von $Z$,
  der $YY^*Y \subseteq Y$ erf"ullt.
\index[B]{sternUntertripel@\sterns{}Untertripel}%
\item Ein \defemph{ternärer \sterns{}Morphismus} zwischen \sterns{}TROs ist
ein \sterns{}linearer Tripelmorphismus.
\index[B]{ternärer sternMorphismus@ternärer \sterns{}Morphismus}%
\dliter{\cite{BlecherWerner06OrderedCsModules}, S. 3}%
\end{enumaufz}
\end{definitn}

\dremark{
Einige Aussagen f"ur tern"are \sterns{}Morphismen findet man in
\cite{BlecherWerner06OrderedCsModules}, S. 3.}%

\begin{definitn}
Sei $X$ ein selbstadjungierter Operatorraum.
\begin{enumaufz}
\item Eine \defemph{tern"are \sterns{}Erweiterung} von $X$
ist ein Paar $(Z,\kappa)$,
bestehend aus einem \sterns{}TRO $Z$ und einer \sterns{}linearen, vollst"andigen Isometrie
$\kappa : X \to Z$ derart, \dass es keine nicht-trivialen \sterns{}Untertripel von $Z$ gibt,
die $\kappa(X)$ enthalten.
\index[B]{ternäre sternErweiterung@ternäre \sterns{}Erweiterung}%

\item Eine \defemph{ternäre \sterns{}Hülle} von $X$ ist
eine beliebige tern"are \sterns{}Erweiterung $(Y,j)$ von $X$ mit der universellen Eigenschaft des
folgenden Satzes.
\dremark{Hier w"are es besser, die tern"are \sterns{}H"ulle analog zur tern"aren H"ulle
  "uber \sterns{}Tripelsysteme zu definieren.
  Wird aber in \cite{BlecherWerner06OrderedCsModules} nicht gemacht,
  w"are au"serdem mehr Arbeit.}%
\index[B]{ternäre sternHülle@ternäre \sterns{}Hülle}%
\dliter{\cite{BlecherWerner06OrderedCsModules}, S. 3}%
\end{enumaufz}
\end{definitn}

\begin{satz}[\cite{BlecherWerner06OrderedCsModules}, S. 3]
Es sei $X$ ein selbstadjungierter Operatorraum.
Dann existiert eine tern"are \sterns{}Er\-wei\-terung $(Y,j)$ von $X$
mit der folgenden universellen Eigenschaft:
F"ur jede tern"are \sterns{}Erweiterung $(Z,\kappa)$ von $X$ existiert ein
(notwendigerweise eindeutiger und surjektiver) tern"arer \sterns{}Mor\-phis\-mus
$\Phi : Z \to Y$ derart, \dass $\Phi \circ \kappa = j$ gilt.
\[ \xymatrix{
Z \ar@{-->}[rr]^{\Phi} & & Y \\
& X \ar[ul]^{\kappa}  \ar[ur]_{j} &
} \]
\end{satz}

\begin{notation}
Sei $X$ ein \cmsa Operatorraum.
Dann wird mit $(\ternsh{X},j)$ die (bis auf tern"are \sterns{}Isomorphismen)
eindeutige tern"are \sterns{}H"ulle von $X$ bezeichnet.
\end{notation}

\begin{bemerkung}[\cite{BlecherKNW07OrderedInvOpSpaces}, nach Theorem 2.1]\label{ternsHuelleIstTernHuelle}
Sei $X$ ein \cmsa Operatorraum.
Dann ist die tern"are \sterns{}H"ulle von $X$ isomorph
zu der tern"aren H"ulle des Operatorraumes $X$, kurz: $\ternsh{X} \cong \ternh{X}$.
\end{bemerkung}

\section{Die injektive H"ulle eines \cmsa Operatorraumes}

In diesem Abschnitt wird zun"achst die Existenz einer injektiven
\sterns{}H"ulle eines selbstadjungierten Operatorraumes $X$ gezeigt.
Unser Vorgehen lehnt sich dem in \cite{Paulsen02ComplBoundedMaps}, Kapitel 15, an,
in dem die Existenz einer injektiven H"ulle eines Operatorraumes gezeigt wird.
Anschlie"send wird bewiesen, \dass die injektive H"ulle eindeutig ist.
\dremark{und \dass f"ur jeden \cmsa Operatorraum gilt: $I^*(X) = I(X)$.
\dremark{Erw"ahnen, welche \cmsa Operatorraumstruktur $I(X)$ tr"agt.}%
}%
\skiptext

Zun"achst definieren wir die Begriffe injektiver selbstadjungierter Operatorraum
und injektive \sterns{}H"ulle.

\dremark{Idee: Eine injektive H"ulle konstruieren f"ur \cmsa Operatorr"aume.
  Die Objekte sind \cmsa Operatorr"aume,
  die Morphismen \sterns{}lineare, \cmvk{e} Abbildungen.
  Beweis analog zu \cite{PaulsenComplBoundedMapsPS}, 15.4,
  allerdings ben"otigt man angepa\cms{}te Aussagen.}%

\begin{definitn}
Ein \cmsa Operatorraum $Z$ hei"st \defemph{injektiv}, falls
f"ur jeden \cmsa Operatorraum $X$ und jede \sterns{}lineare, \cmvk{e} Abbildung $\varphi : X \to Z$
und f"ur jeden \cmsa Operatorraum~$Y$, der $X$ als abgeschlossenen Unterraum enth"alt,
eine \sterns{}lineare, \cmvk{e} Fortsetzung $\Phi : Y \to Z$ von $\varphi$ existiert.
\index[B]{injektiver \cmsa Operatorraum}%
\[
\xymatrix{
Y \ar@{-->}[rd]^\Phi & \\
X \ar[u]^\subseteq \ar[r]_{\varphi} & Z}
\]
\end{definitn}

\begin{definitn}
Sei $X$ ein \cmsa Operatorraum.
Eine \defemph{\sterns{}Erweiterung} von $X$ ist ein Paar, bestehend aus einem \cmsa Operatorraum
$Y$ und einer \sterns{}linea\-ren, vollst"andigen Isometrie $\kappa : X \rightarrow Y$.
\dremark{Idee: $X \subseteq Y$ soll ja bedeuten,
  \dass die Struktur von $Y$ auf $X$ vererbt wird.}%
\index[B]{sternErweiterung eines sa Operatorraumes@\sterns{}Erweiterung eines \cmsa Operatorraumes}%
\end{definitn}

\begin{definitn}
Sei $X$ ein \cmsa Operatorraum.
Dann hei"st $(Y,\kappa)$ \defemphi{injektive \sterns{}Hülle} von $X$, falls
$(Y,\kappa)$ eine \sterns{}Erweiterung von $X$ ist,
$Y$ injektiv ist und es keinen echten injektiven selbstadjungierten Unterraum von $Y$ gibt,
der $\kappa(X)$ enth"alt.
\end{definitn}

Wir sprechen hier von einer injektiven \sterns{}H"ulle,
um diese besser von der injektiven H"ulle eines Operatorraumes unterscheiden zu k"onnen.
\skiptext

Im folgenden wird die Existenz einer injektiven \sterns{}H"ulle
f"ur einen beliebigen \cmsa Operatorraum $X$ gezeigt.
Hierf"ur werden zun"achst die Begriffe \sterns$X$-Abbildung, \sterns$X$-Projektion
und \sterns$X$-Halbnorm entsprechend den Begriffen $Y$-Abbildung, $Y$-Projektion
und $Y$-Halbnorm f"ur einen beliebigen Operatorraum $Y$ eingef"uhrt.
\dremark{Wir gehen dabei "ahnlich vor wie Paulsen beim Beweis der Existenz einer
  injektiven H"ulle f"ur Operatorr"aume in \cite{Paulsen02ComplBoundedMaps}, Kapitel 15.}%
\dremark{Analog zu \cite{Paulsen02ComplBoundedMaps}, S. 209 unten
\dremark{\cite{PaulsenComplBoundedMapsPS}, S. 221 oben; \ref{defFAbbOR}}%
werden $X$-Abbildung und $X$-Projektion definiert,
allerdings wird zus"atzlich verlangt, \dass diese Abbildungen \sterns{}linear sind.}

\begin{definitn}\label{defXAbbsaOR}
Sei $X \subseteq \mLinStet(H)$ ein \cmsa Operatorraum.
\begin{enumaufz}
\item Eine Abbildung $\varphi : \mLinStet(H) \rightarrow \mLinStet(H)$ hei"st
\defemph{\sterns{}$X$-Abbildung},
falls $\varphi$ eine \cmvk{e}, \sterns{}lineare Abbildung mit
$\varphi\restring_X = \Id_X$ ist.
\index[B]{sternXAbbildung@\sterns{}$X$-Abbildung}%

\item Unter einer \defemph{\sterns{}$X$-Projektion} versteht
man eine \sterns{}$X$-Abbildung $\varphi$,
f"ur die $\varphi \circ \varphi = \varphi$ gilt.
\index[B]{sternXProjektion@\sterns{}$X$-Projektion}%
\dliter{analog zu \cite{PaulsenComplBoundedMapsPS}, S. 221 oben;
  \cite{Paulsen02ComplBoundedMaps}, S. 209 unten}%
\end{enumaufz}
\end{definitn}

Wie man leicht sieht, gilt:

\begin{defBemerkung}
Sei $X \subseteq \mLinStet(H)$ ein \cmsa Operatorraum.
Auf der Menge der \sterns{}$X$-Projektionen wird durch
\[ \varphi \preceq \psi
\overset{\text{def}}{\iff} \varphi \circ \psi = \varphi = \psi \circ \varphi \]
eine Ordnungsrelation definiert.
\dremark{$\preceq$ ist i.a. nicht konnex}%
\end{defBemerkung}

\begin{definitn}
Sei $X \subseteq \mLinStet(H)$ ein \cmsa Operatorraum.
\begin{enumaufz}
\item Eine Abbildung $p : \mLinStet(H) \rightarrow \mbbR$ hei"st \defemph{\sterns{}$X$-Halbnorm},
  falls eine \sterns{}$X$-Abbildung $\varphi : \mLinStet(H) \rightarrow \mLinStet(H)$ so existiert,
  \dass gilt: $p(x) = \norm{\varphi(x)}$ f"ur alle $x \in \mLinStet(H)$.
\item Sei $\varphi : \mLinStet(H) \to \mLinStet(H)$ eine \sterns{}$X$-Abbildung.
  Man nennt die Halbnorm
\[ p_\varphi : \mLinStet(H) \rightarrow \mbbR, x \mapsto \norm{\varphi(x)} \]
\index[B]{sternXHalbnomr@\sterns{}$X$-Halbnorm}%
  die zu $\varphi$ geh"orige \sterns{}$X$-Halbnorm.
\end{enumaufz}
\end{definitn}

\begin{defBemerkung}
Sei $X \subseteq \mLinStet(H)$ ein \cmsa Operatorraum.
Auf der Menge der \sterns$X$-Halbnormen wird durch
\[ p \leq q  \overset{\text{def}}{\iff}  p(x) \leq q(x) \text{ f"ur alle } x \in \mLinStet(H) \]
eine Ordnungsrelation definiert.\dremark{$\leq$ ist i.a. nicht konnex}
\end{defBemerkung}

\begin{definitn}
Seien $X$, $Y$ \cmsa Operatorr"aume, sei $\varphi : X \to Y$.
Definiere
\index[S]{phis@$\abbsa{\varphi}$}%
\[ \abbsa{\varphi} : X \to Y, x \mapsto \left(\varphi(x^*)\right)^*. \]
\end{definitn}

\begin{lemma}\label{normPhiEqNormPhis}
Seien $X$, $Y$ \cmsa Operatorr"aume, sei $\varphi : X \to Y$ linear.
Dann gilt: $\norm{\abbsa{\varphi}} = \norm{\varphi}$.
\end{lemma}

\begin{proof}
F"ur alle $x \in X$ erh"alt man wegen $\norm{x^*} = \norm{x}$ (siehe \eqref{eqNormxjiEqNormx}):
\[ \norm{\abbsa{\varphi}(x)}
=  \norm{ \varphi(x^*)^* }
=  \norm{ \varphi(x^*) }
\leq \norm{\varphi} \cdot \norm{ x^* }
=  \norm{\varphi} \cdot \norm{x}. \]
Es folgt: $\norm{\abbsa{\varphi}} \leq \norm{\varphi}$.
Wegen
$\norm{\varphi}
=  \norm{\abbsa{(\abbsa{\varphi})}}
\leq \norm{\abbsa{\varphi}}
\leq \norm{\varphi}$ folgt die Behauptung.
\end{proof}

\begin{bemerkung}\label{LHinjAlsSAOpRaum}
Sei $H$ ein Hilbertraum.
Dann ist $L(H)$ injektiv als \cmsa Operatorraum.
\end{bemerkung}

\begin{proof}
Seien $X$, $Y$ \cmsa Operatorr"aume mit $X \subseteq Y$.
Sei $\varphi : X \to L(H)$ vollst"andig kontraktiv und \sterns{}linear.
Da $L(H)$ ein injektiver Operatorraum ist
(\refb{\cite{BlecherLeMerdy04OpAlg}, Theorem 1.2.10,}{Beispiel~\ref{LHIstInjOR}}),
findet man eine \cmvk{e} Abbildung $\Phi : Y \to \mLinStet(H)$ mit $\Phi\restring_X = \varphi$.
\[
\xymatrix{
Y \ar@{-->}[r]^\Phi & L(H) \\
X \ar[u]^\subseteq \ar[ru]_{\varphi} & }
\]
Setze $\Psi := \frac{\Phi + \abbsa{\Phi}}{2}$.
Dann ist $\Psi$ \sterns{}linear und eine Fortsetzung von $\varphi$.
\dremark{$2\Psi(y^*) = \Phi(y^*) + \abbsa{\Phi}(y^*) = (\Phi(y^*)^*)^* + \Phi(y)^*
   = (\abbsa{\Phi}(y))^* + \Phi(y)^* = 2\Psi(y)^*$,
   $2\Psi(x) = \Phi(x) + \abbsa{\Phi}(x) = \varphi(x) + \Phi(x^*)^*
   =  \varphi(x) + \varphi(x^*)^* = 2 \varphi(x)$ f"ur alle $x \in X$.}%
Mit Lemma~\ref{normPhiEqNormPhis} erh"alt man,
\dass $\abbsa{\Phi}$ vollst"andig kontraktiv ist.
\dremark{Beweis ohne Lemma~\ref{normPhiEqNormPhis}:
Weiter gilt mit \eqref{eqNormxjiEqNormx} f"ur alle $n \in \mbbN$ und $y \in M_n(Y)$:
\begin{align*}
   \normlr{ (\abbsa{\Phi})_n(y) }_n
&= \dremarkm{\normlr{ \left(\abbsa{\Phi}(y_{ij})\right)_{i,j} }
=}  \normlr{ \left(\Phi(y^*_{ij})^*\right)_{i,j} }
\overset{\dremarkm{\eqref{eqNormxjiEqNormx}}}{=}  \normlr{ \left(\Phi(y^*_{ij})\right)_{j,i} }
\dremarkm{= \normlr{ \left(\Phi(y^*_{ji})\right)_{i,j} } }  \\
&= \normlr{ \Phi_n\left((y^*_{ji})_{i,j}\right) }
\leq  \normlr{ (y^*_{ji})_{i,j} }
\overset{\dremarkm{\eqref{eqNormxjiEqNormx}}}{=}  \norm{y}.
\end{align*}}%
Hiermit folgt:
$\normcb{\Psi} \leq \frac{1}{2}\bigl(\normcb{\Phi} + \normcb{\abbsa{\Phi}}\bigr)  \leq  1$.
\end{proof}

Analog zu Proposition~\ref{BildProjIstInjOpraum} im Operatorraum gilt in
\cmsa Operatorr"aumen:

\begin{bemerkung}\label{ProjInjSAOpRaum}
Sei $X$ ein injektiver \cmsa Operatorraum und
$P : X \rightarrow X$ eine \sterns{}lineare, \cmvk{e} Projektion.
Dann ist $P(X)$ ein injektiver \cmsa Operatorraum.
\end{bemerkung}

\begin{proof}
\dremark{Nach \cite{WernerFunkana3}, Lemma IV.6.1, ist $P(X)$ abgeschlossen,
  also ein Operatorraum.}%
F"ur alle $x \in X$ gilt $P(x)^* = P(x^*) \in P(X)$, also ist $\cdot^*\restring_{P(X)}$
eine Involution auf $P(X)$.
\dremark{Es existiert eine \sterns{}lineare, vollst"andige Isometrie $i : X \to L(H)$.
Dann ist auch $i\restring_{P(X)}$ eine \sterns{}lineare, vollst"andige Isometrie.}%

Seien $Y \subseteq Z \subseteq \mLinStet(K)$ \cmsa Operatorr"aume,
sei $\varphi : Y \rightarrow P(X)$ vollst"andig kontraktiv und \sterns{}linear.
Dann ist auch $\varphi : Y \rightarrow X$ vollst"andig kontraktiv und \sterns{}linear.
Da $X$ injektiv ist, findet man eine \sterns{}lineare, \cmvk{e} Fortsetzung $\Phi : Z \to X$
von $\varphi$.
Somit ist $\Psi := P \circ \Phi : Z \to P(X)$ die gesuchte \sterns{}lineare,
\cmvk{e} Fortsetzung von $\varphi$.
\dremark{Es gilt f"ur alle $y \in Y$:
  $\Psi(y) = P \circ \Phi(y) = P \circ \varphi(y) = \varphi(y)$
  (Beachte: $P\restring_{P(X)} = \Id_{P(X)}$).}%
\[
\xymatrix{
Z \ar@{-->}[rr]^\Phi && X \ar[d]^P \\
Y \ar[rr]^\varphi \ar[u]^\subseteq && P(X) }
\]
\end{proof}

\begin{satz}\label{SAOpRaumExMinHalbnorm}
Sei $X \subseteq \mLinStet(H)$ ein \cmsa Operatorraum.
Dann existiert eine bez"uglich \glqq$\leq$\grqq{} minimale \sterns{}$X$-Halbnorm auf $\mLinStet(H)$.
\end{satz}

\begin{proof}
\newcommand{\cmvarphi}{\varphi}%
Dies beweist man analog zum Beweis f"ur Operatorr"aume,
siehe \cmzB \cite{Paulsen02ComplBoundedMaps}, Beweis von Proposition 15.3.
\dremark{\ref{Paulsen15.3:Opraum}}%
Wir skizzieren kurz den dortigen Beweis, angepa\cms{}t an unsere Situation:
Sei $(p_\lambda)_{\lambda \in \Lambda}$ eine absteigende Kette von \sterns{}$X$-Halbnormen,
wobei $\varphi_\lambda : \mLinStet(H) \to \mLinStet(H)$ eine \sterns{}$X$-Abbildung
f"ur alle $\lambda \in \Lambda$ sei.
Dann besitzt $(\varphi_\lambda)_{\lambda \in \Lambda}$ ein Teilnetz
$(\varphi_{\lambda_\mu})_\mu$,
welches gegen die $X$-Abbildung $\cmvarphi$ konvergiert.
Man erh"alt: $p_\cmvarphi \leq p_{\varphi_\lambda}$ f"ur alle $\lambda \in \Lambda$.

Es bleibt zu zeigen, \dass $\cmvarphi$ \sterns{}linear ist.
Sei $x \in X$.
Man hat f"ur alle $h,k \in H$:
\begin{align*}
   \skalpr{\cmvarphi(x^*)h}{k}
&= \lim_\mu \skalpr{\varphi_{\lambda_\mu}(x^*)h}{k}
=  \lim_\mu \skalpr{h}{\varphi_{\lambda_\mu}(x)k}
=  \lim_\mu \overline{\skalpr{\varphi_{\lambda_\mu}(x)k}{h}}  \\
&= \overline{\lim_\mu \skalpr{\varphi_{\lambda_\mu}(x)k}{h}}
=  \overline{\skalpr{\cmvarphi(x)k}{h}}
=  \skalpr{h}{\cmvarphi(x)k}
=  \skalpr{\cmvarphi(x)^*h}{k}.
\end{align*}
Es folgt:
$\cmvarphi(x^*) = \cmvarphi(x)^*$.
\dremark{2.10.'08/2}%
\end{proof}

Unter Verwendung von Satz~\ref{SAOpRaumExMinHalbnorm}, Proposition~\ref{LHinjAlsSAOpRaum} und
Proposition~\ref{ProjInjSAOpRaum} beweist man
analog zu \cite{Paulsen02ComplBoundedMaps}, Theorem~15.4, den folgenden Satz:

\begin{satz}[Existenz einer injektiven \sterns{}H"ulle]\label{exInjsHuelle}
Sei $X \subseteq \mLinStet(H)$ ein \cmsa Operatorraum und
$\varphi : \mLinStet(H) \to \mLinStet(H)$ eine \sterns{}$X$-Abbildung derart,
\dass $p_\varphi$ eine minimale \sterns{}$X$-Halbnorm ist.
Sei $\iota : X \to \varphi(\mLinStet(H)), x \mapsto x$.
Dann ist $\varphi$ eine minimale \sterns{}$X$-Projektion und
$(\varphi(\mLinStet(H)), \iota)$ eine injektive \sterns{}H"ulle von $X$,
die ein \sterns{}TRO ist.
\end{satz}

\dremww{Beweis.
Zeige: $I^*(X)$ ist \sterns{}TRO.
Da $\varphi$ \sterns{}linear ist, folgt:
$I^*(X)^* = \varphi(L(H))^* = \varphi(L(H)^*) = \varphi(L(H)) = I^*(X)$.
}%

Analog zu \cite{Paulsen02ComplBoundedMaps}, Lemma 15.5, zeigt man
f"ur \cmsa Operatorr"aume mit Hilfe von Satz~\ref{exInjsHuelle}:

\begin{lemma}\label{lemmaRigidesaOR}
Sei $X \subseteq L(H)$ ein \cmsa Operatorraum.
Sei $\varphi : L(H) \to L(H)$ eine \sterns{}$X$-Abbildung derart,
\dass $p_\varphi$ eine minimale \sterns{}$X$-Halbnorm ist.
Falls $\gamma : \varphi(L(H)) \to \varphi(L(H))$ \cmvk{} und \sterns{}linear ist
mit $\gamma\restring_X = \Id_X$, dann gilt:
$\gamma(\varphi(x)) = \varphi(x)$ f"ur alle $x \in X$.
\end{lemma}

Mit diesem Lemma beweist man,
analog zu \cite{Paulsen02ComplBoundedMaps}, Theorem 15.6,
den folgenden

\begin{satz}[Eindeutigkeit der injektiven \sterns{}H"ulle]\label{EindInjsHuelleOpraum}
Sei $X$ ein \cmsa Operatorraum.
Seien $(Y_1,\kappa_1)$, $(Y_2,\kappa_2)$ injektive \sterns{}H"ullen von $X$.
Dann l"a\cms{}t sich die Abbildung
$i : \kappa_1(X) \to \kappa_2(X)$, gegeben durch
$i(\kappa_1(x)) = \kappa_2(x)$ f"ur alle $x \in X$,
eindeutig zu einer \sterns{}linearen, surjektiven vollst"andigen Isometrie
von $Y_1$ auf $Y_2$ fortsetzen.
\[ \xymatrix{
\,\,Y_1\,\, \ar@{>-->>}[rr]^{\cong} & & Y_2 \\
\kappa_1(X) \ar[rr]^{i} \ar[u]^\subseteq & & \kappa_2(X) \ar[u]_\subseteq \\
& X \ar[ul]^{\kappa_1}  \ar[ur]_{\kappa_2} &
} \]
\end{satz}

\dremark{
\begin{notation}
Jeder \cmsa Operatorraum $X$ besitzt nach Satz~\ref{exInjsHuelle} eine injektive
\sterns{}H"ulle,
die nach Satz~\ref{EindInjsHuelleOpraum} bis auf vollst"andige
Isometrie eindeutig ist und mit $(I^*(X),j)$ oder kurz mit $I^*(X)$ bezeichnet wird.
\index[S]{IsX@$I^*(X)$ (injektive \sterns{}Hülle)}%
\index[S]{j@$j$}%
\end{notation}
}%

Mit Lemma~\ref{lemmaRigidesaOR} und Satz~\ref{EindInjsHuelleOpraum} beweist man:

\begin{folgerung}[Rigidit"at]
Sei $X$ ein \cmsa Operatorraum.
Sei $(Y,\kappa)$ eine injektive \sterns{}H"ulle von $X$
und $\psi : Y \to Y$ \cmvk{} und \sterns{}linear mit
$\psi(\kappa(x)) = \kappa(x)$ f"ur alle $x \in X$.
Dann gilt: $\psi = \Id_Y$.
\dremark{Beweis analog zu \cite{Paulsen02ComplBoundedMaps}, Corollary 15.7}%
\end{folgerung}

\dremark{
Analog zum Resultat f"ur tern"are \sterns{}H"ullen
($\ternsh{X} = \ternh{X}$, Proposition~\ref{ternsHuelleIstTernHuelle})
gilt f"ur die injektive \sterns{}H"ulle:

\begin{bemerkung}\label{XsaFolgtIsXEqIX}
Sei $X$ ein \cmsa Operatorraum.
Dann gilt: $I^*(X) = I(X)$.
\dremark{In Worte fassen, was genau gemeint ist. (WW)}%
\end{bemerkung}

Zum Beweis dieser Proposition erinnern wir zun"achst an eine Definition und formulieren dann drei Lemmata:

\begin{definitn}
Sei $X \subseteq \mLinStet(H)$ ein Operatorraum.
Dann wird
\[ \adjor{X} := \{ x^* \setfdg x \in X \} \]
der \defemph{adjungierte Operatorraum} von $X$ genannt.
\index[S]{Xs@$\adjor{X}$ (adjungierter Operatorraum)}%
\index[B]{adjungierter Operatorraum}%
\end{definitn}

Als ein abstrakter Operatorraum ist $\adjor{X}$ unabh"angig von
der jeweiligen Einbettung von $X$ in $\mLinStet(H)$ (\cite{BlecherLeMerdy04OpAlg}, 1.2.25).

\begin{lemma}\label{XProjIffXsProj}
Sei $X \subseteq L(H)$ ein Operatorraum und $\varphi : L(H) \to L(H)$.
Dann ist $\varphi$ genau dann eine $X$-Abbildung (bzw. $X$-Projektion),
wenn $\abbsa{\varphi}$ eine $\adjor{X}$-Abbildung (bzw. $\adjor{X}$-Projektion) ist.
\end{lemma}

\begin{proof}
\bewitemphq{$\Rightarrow$}
Sei $\varphi$ eine $X$-Abbildung.
\dremark{Es gilt:
  $\abbsa{\varphi}(x^*) = \varphi(x^{**})^* = \varphi(x)^* \in \adjor{X}$
  f"ur alle $x \in X$.}%
Offensichtlich ist $\abbsa{\varphi}$ linear.
\dremark{Die Skalarmultiplikation auf $\adjor{X}$ ist definiert durch
  $\lambda \diamond x^* = (\overline{\lambda} x)^*$.
  Es gilt f"ur alle $x \in X$ und $\lambda \in \mbbC$:
\begin{align*}
   \abbsa{\varphi}(\lambda \diamond x^*)
&= \varphi(((\overline{\lambda}x)^*)^*)^*
=  \varphi(\overline{\lambda}x)^*
=  \overline{\overline{\lambda}} \varphi(x)^*
=  \left(\overline{\lambda} \varphi(x)\right)^*
=  \lambda \diamond \varphi(x)^*
=  \lambda \diamond \abbsa{\varphi}(x^*).
\end{align*}}%
Da $\varphi$ vollst"andig kontraktiv ist,\dremark{$(\ddagger)$}
folgt f"ur alle $n \in \mbbN$ und $x \in M_n(\adjor{X})$ mit \eqref{eqNormxjiEqNormx}
\begin{align*}
   \norm{ (\abbsa{\varphi})_n(x) }_n
&= \dremarkm{\norm{ (\abbsa{\varphi}(x_{ij}) )_{i,j} }_n
=}   \norm{ (\varphi(x^*_{ij})^* )_{i,j} }_n  \\
&\overset{\dremarkm{\eqref{eqNormxjiEqNormx}}}{=}  \norm{ (\varphi(x^*_{ji}) )_{i,j} }_n
\overset{\dremarkm{(\ddagger)}}{\leq}  \norm{ (x^*_{ji})_{i,j} }_n
=  \norm{x}_n,
\end{align*}
also ist $\abbsa{\varphi}$ vollst"andig kontraktiv.

Sei $x \in \adjor{X}$.
Wegen $\varphi\restring_X = \Id_X$\dremark{$(*)$} erh"alt man:
\[ \abbsa{\varphi}(x) \dremarkm{= \varphi(x^*)^* \overset{\dremarkm{(*)}}{=} x^{**}} = x. \]
Somit ist $\abbsa{\varphi}$ eine $X$-Abbildung.

Ist $\varphi$ eine $X$-Projektion,\dremark{$(\dagger)$} so gilt
\[ \abbsa{\varphi} \circ \abbsa{\varphi}(x)
\dremarkm{=  \abbsa{\varphi}(\varphi(x^*)^*)
=  \varphi(\varphi(x^*)^{**})^*}
\overset{\dremarkm{(\dagger)}}{=}  \varphi(\varphi(x^*))^*
\dremarkm{=  \varphi(x^*)^*}
=  \abbsa{\varphi}(x), \]
daher ist $\abbsa{\varphi}$ eine $\adjor{X}$-Projektion.

\bewitemphq{$\Leftarrow$}
Ist $\abbsa{\varphi}$ eine $\adjor{X}$-Abbildung (bzw. $\adjor{X}$-Projektion),
dann folgt wegen $\adjor{(\adjor{X})} = X$ mit \glqq$\Rightarrow$\grqq,
\dass $\abbsa{(\abbsa{\varphi})} = \varphi$ eine $X$-Abbildung (bzw. $X$-Projektion) ist.
\dremark{25.11.'08/1}%
\end{proof}




\begin{lemma}\label{leqXProjIffleqXsProj}
Sei $X \subseteq \mLinStet(H)$ ein Operatorraum.
Seien $\varphi$, $\psi$ $X$-Pro\-jek\-tio\-nen.
Dann gilt:
\[ \varphi \preceq \psi \quad\text{(als $X$-Projektionen)}
\iff  \abbsa{\varphi} \preceq \abbsa{\psi} \quad\text{(als $\adjor{X}$-Projektionen)}. \]
\end{lemma}

\begin{proof}
F"ur alle $x \in \adjor{X}$ gilt:
\[ \abbsa{\varphi} \circ \abbsa{\psi}(x)
\dremarkm{=  \abbsa{\varphi}(\psi(x^*)^*)
=  \varphi( \psi(x^*)^{**} )^*
=  \varphi( \psi(x^*) )^*}
=  \abbsa{(\varphi \circ \psi)}(x). \dremarkm{\quad(*)} \]
Somit folgt:\dremark{25.11.'08/1}%
\begin{align*}
   \varphi \preceq \psi
&\iff \varphi \circ \psi = \varphi = \psi \circ \varphi
\overset{\dremarkm{(*)}}{\iff}
   \abbsa{(\varphi \circ \psi)} = \abbsa{\varphi} = \abbsa{(\psi \circ \varphi)}  \\
&\iff \abbsa{\varphi} \circ \abbsa{\psi} = \abbsa{\varphi} = \abbsa{\psi} \circ \abbsa{\varphi}
\iff \abbsa{\varphi} \preceq \abbsa{\psi}.  \qedhere
\end{align*}
\end{proof}

\dremark{Im obigen Beweis wurde gezeigt:
\begin{equation}\label{eqphispsisEqphipsis}
   \abbsa{\varphi} \circ \abbsa{\psi}
=  \abbsa{(\varphi \circ \psi)}.
\end{equation}}%

Durch unmittelbares Nachrechnen erh"alt man das folgende

\begin{lemma}\label{leqXHalbnIffleqXsHalbn}
Sei $X \subseteq \mLinStet(H)$ ein Operatorraum.
Seien $\varphi$, $\psi$ $X$-Abbildungen.
Dann gilt:
\dremark{Beachte: $\abbsa{\varphi}$ ist nach \ref{XProjIffXsProj} eine $\adjor{X}$-Abbildung,
  also ist $p_{\abbsa{\varphi}}$ eine $\adjor{X}$-Halbnorm.}%
\[ p_\varphi \leq p_\psi \quad\text{(als $X$-Halbnormen)}
\iff  p_{\abbsa{\varphi}} \leq p_{\abbsa{\psi}} \quad\text{(als $\adjor{X}$-Halbnormen)}. \]
\end{lemma}

\dremww{Beweis.
Es gilt:
\dremark{Beachte: $\norm{\varphi(x^*)} = \norm{\varphi(x^*)^*} = \norm{\abbsa{\varphi}(x)}$}%
\begin{align*}
   p_\varphi \leq p_\psi
&\iff  \forall x \in \mLinStet(H) : \norm{\varphi(x)} = \norm{\psi(x)}  \\
&\iff  \forall x \in \mLinStet(H) : \norm{\varphi(x^*)} = \norm{\psi(x^*)}  \\
&\iff  \forall x \in \mLinStet(H) :
  \norm{\abbsa{\varphi}(x)} = \norm{\abbsa{\psi}(x)}  \\
&\iff  p_{\abbsa{\varphi}} \leq p_{\abbsa{\psi}}
\end{align*}
}%

\begin{satz}
Sei $X \subseteq \mLinStet(H)$ ein \cmsa Operatorraum und
$\varphi : \mLinStet(H) \to \mLinStet(H)$ eine $X$-Abbildung derart,
\dass $p_\varphi$ eine minimale $X$-Halbnorm ist.
Dann ist $\tilde{\varphi} := \frac{\varphi + \abbsa{\varphi}}{2}$
eine minimale \sterns{}$X$-Projektion und
$p_{\tilde{\varphi}}$ eine minimale \sterns{}$X$-Halbnorm und
eine minimale $X$-Halbnorm.
Insbesondere folgt: $I^*(X) = I(X)$.
\end{satz}

\begin{proof}
\dremark{Zeige: $\tilde{\varphi}$ ist $X$-Projektion.}%
Offensichtlich ist $\tilde{\varphi}$ eine $X$-Abbildung.
\dremark{$\abbsa{(\varphi + \abbsa{\varphi})}  =  \abbsa{\varphi} + \varphi$,
  also \sterns{}linear.
  Da $\varphi$ und $\abbsa{\varphi}$ $X$-Projektionen sind,
  gilt: $\tilde{\varphi}$ ist eine $X$-Abbildung.}%
Es gilt  wegen \dremark{$(*)$:}$\varphi\restring_{\Bild(\varphi)} = \Id_{\Bild(\varphi)}$
f"ur alle $x \in \mLinStet(H)$:
\begin{equation}\label{eqphiphisEqphis}
   (\varphi \circ \abbsa{\varphi})(x)
=  \varphi(\varphi(x^*)^*)
=  \varphi(x^*)^*
=  \abbsa{\varphi}(x).
\end{equation}
Hiermit folgt:
\begin{equation}\label{eqphisphiEqphi}
   \abbsa{\varphi} \circ \varphi
\overset{\dremarkm{\eqref{eqphispsisEqphipsis}}}{=}  \abbsa{(\varphi \circ \abbsa{\varphi})}
\overset{\dremarkm{\eqref{eqphiphisEqphis}}}{=}  \abbsa{(\abbsa{\varphi})}
=  \varphi.
\end{equation}
Zusammen mit \eqref{eqphiphisEqphis} erh"alt man:
\[ \tilde{\varphi} \circ \tilde{\varphi}
=  \frac{1}{4}\left( \varphi \circ \varphi + \varphi \circ \abbsa{\varphi} +
                     \abbsa{\varphi} \circ \varphi + \abbsa{\varphi} \circ \abbsa{\varphi} \right)
\overset{\dremarkm{\eqref{eqphiphisEqphis},\eqref{eqphisphiEqphi}}}{=}
   \tilde{\varphi}. \]
Also ist $\tilde{\varphi}$ eine $X$-Projektion.
\smallskip

\dremark{Zeige: $p_{\tilde{\varphi}}$ ist minimale \sterns{}$X$-Halbnorm.}%
Sei $\psi$ eine \sterns{}$X$-Abbildung.
Dann ist $\psi$ insbesondere eine $X$-Abbildung.
Da $p_\varphi$ eine minimale $X$-Halbnorm ist,
xxx
\dmarginfz\dremark{Problem: Vergleichbarkeit!}%

Sei $\psi$ eine \sterns{}$X$-Projektion.
Da $\varphi$ eine minimale $X$-Projektion ist und $\psi$ insbesondere eine $X$-Projektion,
\dremark{$(*)$: $\varphi \circ \psi = \varphi = \psi \circ \varphi$}%
erh"alt man f"ur alle $x \in \mLinStet(H)$:\dremark{$\varphi$ ist \sterns{}linear}
\[ (\abbsa{\varphi} \circ \psi)(x)
=  \varphi( \psi(x)^* )^*
=  \varphi( \psi(x^*) )^*
=  \left( (\psi \circ \varphi)(x^*) \right)^*
\overset{\dremarkm{(*)}}{=}  \varphi(x^*)^*
=  \abbsa{\varphi}(x) \]
und
\[ (\psi \circ \abbsa{\varphi})(x)
=  \psi(\varphi(x^*)^*))
=  \psi\left(\varphi(x^*)\right)^*
=  \left((\psi \circ \varphi)(x^*)\right)^*
\overset{\dremarkm{(*)}}{=}  \varphi(x^*)^*
=  \abbsa{\varphi}(x). \]
\dmarginpar{pr}\dremark{Wof"ur? Weg?}%

\end{proof}

\begin{lemma}\label{IXsEqIXs}
Sei $X \subseteq \mLinStet(H)$ ein Operatorraum.
Dann gilt: $I(\adjor{X}) = \adjor{I(X)}$.
\end{lemma}

\begin{proof}
Nach Satz~\ref{Paulsen15.3:Opraum} (Existenz einer minimalen $X$-Halbnorm) und
Satz~\ref{exInjHuelle} (Existenz einer injektiven H"ulle) findet man
eine minimale $X$-Projektion $\varphi$ mit $I(X) = \varphi(\mLinStet(H))$.

Mit Lemma~\ref{XProjIffXsProj} und Lemma~\ref{leqXProjIffleqXsProj} folgt,
\dass $\varphi \mapsto \abbsa{\varphi}$ eine ordnungserhaltende Bijektion
von der Menge der $X$-Projektionen auf die Menge der $\adjor{X}$-Projektionen ist.
Somit ist $\abbsa{\varphi}$ eine minimale $\adjor{X}$-Projektion.
Es folgt:
\dremark{25.11.'08/1}%
\[ I(\adjor{X})
=  \abbsa{\varphi}(\mLinStet(H))
\dremarkm{=  \adjor{\adjor{\varphi(\mLinStet(H)})} }
=  \adjor{\varphi(\mLinStet(H))}
=  \adjor{I(X)}.  \qedhere \]
\end{proof}

\begin{proof}[Beweis von Proposition~\ref{XsaFolgtIsXEqIX}]
Da $X$ selbstadjungiert ist, hat man: $X = \adjor{X}$.
Mit Lemma~\ref{IXsEqIXs} folgt: $I(X) = I(\adjor{X}) = \adjor{I(X)}$.
\dmarginfz\dremark{Weiter?}%
\dremark{25.11.'08/1}%
\end{proof}
}%

\dmarginpar{Frage}\dremark{Frage: Wie h"angen die tern"are \sterns{}H"ulle und die
  injektive \sterns{}H"ulle zusammen?}%

\dremark{
\begin{anmerkung}[Operatorr"aume mit Involution]
Operatorr"aume mit Involution werden teilweise behandelt in
\cite{BlecherWerner06OrderedCsModules}, \cite{BlecherKNW07OrderedInvOpSpaces}.
\end{anmerkung}
}%

\dremww{
\section{Pisiers selbstdualer Hilbert-Operatorraum}

\begin{definitn}
Sei $H$ ein Hilbertraum.
Dann wird durch
\[ \underbrace{(\sum_{m,\ell=1}^m \eta_{k\ell} \otimes e_{k\ell})}_{=: \eta} \otimes
   \underbrace{(\sum_{i,j=1}^n \xi_{i,j} \otimes e_{ij})}_{=: \xi}
   \mapsto \sum_{k,\ell} \sum_{i,j} \skalpr{\eta_{k\ell}}{\xi_{ij}} e_{k\ell} \otimes e_{ij}
   =: \skalprd{\eta}{\xi} \]
eine sesquilineare Abbildung
$\skalprd{\cdot}{\cdot\cdot} : M_m(H) \otimes M_n(H) \to M_m \otimes M_n$ definiert.
\dliter{\cite{EffrosRuan00OperatorSpaces}, S. 60 unten}%

Die OH-Matrixnorm auf $H$ wir definiert durch
\[ \norm{\xi}_o := \norm{\skalprd{\xi}{\xi}}^{1/2} \]
f"ur alle $\xi \in M_n(H)$.\dliter{\cite{EffrosRuan00OperatorSpaces}, vor 3.5.2}
\dremark{Zu Pisiers Operatorhilbertraum siehe auch \cite{BlecherLeMerdy04OpAlg}, 5.3.4 und 5.3.5.
  Es gilt: $H_\text{oh}$ ist vollst"andig isometrisch zu einer Operatoralgebra.}%
\end{definitn}

\begin{bemerkung}
Sei $H$ ein Hilbertraum.
Seien $\xi, \eta \in M_n(H) = M_n \otimes H$ so, \dass
\[ \xi = \sum_{h=1}^p \alpha^{(h)} \otimes e_h \quad\text{und}\quad
   \eta = \sum_{h=1}^p \beta^{(k)} \otimes e_k, \]
wobei $e_h \in H$ orthonormal sei und
$\alpha^{(h)}, \beta^{(h)} \in M_n$ f"ur alle $h \in \haken{p}$.
Dann gilt:\dliter{\cite{EffrosRuan00OperatorSpaces}, vor 3.5.2}
\[ \skalprd{\eta}{\xi}
=  \sum_{k,\ell,i,j=1}^n  \sum_{h=1}^p  \beta^{(h)}_{k\ell} \alpha^{(h)}_{ij}
   e_{k\ell} \otimes e_{ij}
=  \sum_{h=1}^p \beta^{(h)} \otimes \overline{\alpha^{(h)}}
\in  M_n \otimes \overline{M_n}. \]
\end{bemerkung}

\begin{bemerkung}\label{matrixeinheitenFuerMmn}
Es gilt: $M_m \otimes M_n \cong M_{mn}$.
Konkret: Seien $(e_{ij})_{ij}$ Matrixeinheiten von $M_m$ und
$(f_{st})_{st}$ Matrixeinheiten von $M_n$.
Dann werden durch $g_{(i,s),(j,t)} := e_{ij} \otimes e_{st}$
f"ur alle $i,j \in \haken{m}$, $s,t \in \haken{n}$
Matrixeinheiten von $M_{mn}$ definiert.
\end{bemerkung}

\dremww{
\begin{beispiel}
Seien $\eta, \xi \in M_2$.
Dann gilt mit den Bezeichnungen aus \ref{matrixeinheitenFuerMmn}:
\dremark{$(*)$: \cite{BourbakiAlgebra1to3:1989}, II.10.10, S. 357, vor (43)}%
\begin{align*}
    \eta \otimes \xi
&= (\eta_{ij} e_{ij})_{ij} \otimes (\xi_{st} e_{st})_{st}
=  (\eta_{ij} \xi_{st} g_{(i,s),(j,t)})_{(i,s),(j,t) \in \haken{2} \times \haken{2}}  \\
&\overset{\dremarkm{(*)}}{=}  \begin{pmatrix}
      \eta_{11}\xi_{11} & \eta_{11}\xi_{12} & \eta_{12}\xi_{11} & \eta_{12}\xi_{12} \\
      \eta_{11}\xi_{21} & \eta_{11}\xi_{22} & \eta_{12}\xi_{21} & \eta_{12}\xi_{22} \\
      \eta_{21}\xi_{11} & \eta_{21}\xi_{12} & \eta_{22}\xi_{11} & \eta_{22}\xi_{12} \\
      \eta_{21}\xi_{21} & \eta_{21}\xi_{22} & \eta_{22}\xi_{21} & \eta_{22}\xi_{22}
    \end{pmatrix}
=  \begin{pmatrix} \eta_{11}\xi & \eta_{12}\xi \\ \eta_{21}\xi & \eta_{22}\xi \end{pmatrix}.
\end{align*}
\end{beispiel}
}%
}%

\chapter{Unbeschr"ankte Multiplikatoren auf Operatorr"aumen}

Dieses Kapitel ist der zentrale Teil dieser Arbeit.
Hier werden drei verschiedene Begriffe
von unbeschr"ankten Multiplikatoren auf einem Operatorraum untersucht.

Da ein Operatorraum weniger Struktur als ein \hcsmodul besitzt,
ist nicht ohne weiteres ersichtlich, wie man die Definition
eines regul"aren Operators auf Operatorr"aume "ubertragen kann.
Bekanntlich ist die Menge $\Adjhm{E}$ der adjungierbaren Abbildungen auf $E$ isomorph zu
der Menge $\Adjlor{E}$ der von links adjungierbaren Multiplikatoren
auf $E$, wobei $E$ als Operatorraum aufgefa\cms{}t wird.
Da au"serdem schiefadjungierte regul"are Operatoren mit Hilfe von $C_0$-Gruppen
charakterisiert werden (siehe Satz~\ref{HollevoetTh21} (Satz von Stone)
und Satz~\ref{erzeugerUnitaereGrWor}),
wird hierdurch die Definition eines unbeschr"ankten schiefadjungierten Multiplikators
auf einem Operatorraum $X$ motiviert:

\begin{definitn}\label{defAdjlcox}
\begin{enumaufz}
\item Ein Operator $A : D(A) \subseteq X \to X$ hei"st
\defemph{unbeschr"ankter schiefadjungierter Multiplikator} auf $X$,
falls $A$ Erzeuger einer $C_0$-Gruppe $(U_t)_\indtHG$ auf $X$
mit $U_t \in \Adjlor{X}$ unit"ar f"ur alle $\indtGr$ ist.

\item Die Menge aller unbeschr"ankten schiefadjungierten Multiplikatoren auf $X$ wird mit
$\Adjlco{X}$ bezeichnet.\dremark{\glqq{}s\grqq{} f"ur schiefadjungiert}
\end{enumaufz}
\end{definitn}

Im ersten Abschnitt des Kapitels wird bewiesen,
\dass die Menge der schiefadjungierten Elemente von $\Adjlor{X}$ in $\Adjlco{X}$ enthalten ist
und \dass in unitalen Operatorr"aumen die umgekehrte Inklusion wahr ist.
Au"serdem wird der Begriff des $C_0$-Linksmultiplikators eingef"uhrt, eines Operators,
der eine $C_0$-Halb\-grup\-pe von Elementen aus $\Multlor{X}$,
der Menge der Linksmultiplikatoren auf $X$, erzeugt.
Wiederum gilt, \dass $\Multlor{X}$ in der Menge der $C_0$"=Linksmultiplikatoren enthalten ist
und \dass in unitalen Operatorr"aumen Gleichheit gilt.\\

Im darauf\/folgenden Unterabschnitt werden Beispiele f"ur $C_0$"=Linksmultiplikatoren
und f"ur unbeschr"ankte schiefadjungierte Multiplikatoren festgehalten.
Insbesondere wird gezeigt, \dass in einem \hcsmodul die schiefadjungierten regul"aren Operatoren
mit den unbeschr"ankten schiefadjungierten Multiplikatoren "ubereinstimmen.\\

Um den Begriff eines beliebigen unbeschr"ankten Multiplikators einzuf"uhren,
wird benutzt, \dass f"ur jeden regul"aren Operator $T$ auf $E$
durch $\mri \cmsmallpmatrix{ 0 & T \\ T^* & 0 }$ ein schiefadjungierter regul"arer
Operator auf $E \oplus E$ gegeben ist.
Dies f"uhrt zu der folgenden Definition im zweiten Abschnitt des Kapitels:

\begin{definitn}\label{defAdjloruEinl}
\begin{enumaufz}
\item Eine Abbildung $T : D(T) \subseteq X \to X$ hei"st
\defemph{unbeschränkter Multiplikator} auf $X$,
falls ein $S : D(S) \subseteq X \to X$ existiert mit der Eigenschaft:
$\mri \begin{pmatrix} 0 & T \\ S & 0 \end{pmatrix} \in \Adjlco{C_2(X)}$,
wobei mit $C_2(X)$ der Spaltenoperatorraum bezeichnet wird.

\item Mit $\Adjloru{X}$ wird die Menge der unbeschr"ankten Multiplikatoren auf $X$
bezeichnet.
\end{enumaufz}
\end{definitn}

Ein zentrales Resultat ist, \dass $\Adjloru{E}$ mit der Menge der regul"aren Operatoren
auf $E$ "ubereinstimmt.
Da ferner $S$ aus der obigen Definition gleich $T^*$ ist,
sind die unbeschr"ankten Multiplikatoren eine Verallgemeinerung
der regul"aren Operatoren auf Operatorr"aume.
Weiterhin gilt: $\Adjlor{X} \subseteq \Adjloru{X}$.\\

Im dritten Abschnitt wird untersucht, unter welchen Voraussetzungen
die Operatoren $A$ und $\matnull{A} := \cmsmallpmatrix{ A & 0 \\ 0 & 0 }$
auf $X$ (bzw. $C_2(X)$) eine $C_0$"=Halbgruppe erzeugen und
welche Zusammenh"ange zwischen diesen beiden Operatoren bestehen.
Mit Hilfe von $\matnull{A}$ erh"alt man intrinsische Charakterisierungen
der drei Begriffe von unbeschr"ankten Multiplikatoren.
Beispielsweise ist $A$ genau dann ein unbeschr"ankter schiefadjungierter Multiplikator,
wenn $\matnull{A}$ eine vollst"andig kontraktive $C_0$-Gruppe auf $C_2(X)$ erzeugt,
\dheisst eine $C_0$-Gruppe $(T_t)_\indtGr$, f"ur die gilt:
$T_t$ ist vollst"andig kontraktiv f"ur alle $\indtGr$.\\

Daher ist von Interesse, wann ein Operator eine vollst"andig kontraktive
$C_0$"=Halbgruppe (bzw. $C_0$-Gruppe) erzeugt.
Charakterisierungen hierzu werden im n"achsten Abschnitt bewiesen,
in dem die S"atze von Hille-Yosida und von Lumer-Phillips
auf Operatorr"aume "ubertragen werden.
Hiermit ergeben sich weitere Charakterisierungen f"ur die
drei Begriffe von unbeschr"ankten Multiplikatoren.\\

Im f"unften Abschnitt wird gezeigt,
\dass man jeden $C_0$"=Linksmultiplikator auf $X$
als Einschr"ankung eines $C_0$"=Linksmultiplikators auf der tern"aren H"ulle $\ternh{X}$
von $X$ auf\/fassen kann.
Entsprechende Aussagen werden auch f"ur die unbeschr"ankten schiefadjungierten
Multiplikatoren und die unbeschr"ankten Multiplikatoren bewiesen.
Hiermit folgt, \dass der Operator $S$ aus Definition~\ref{defAdjloruEinl}
eindeutig ist und somit als Adjungierte von $T$ angesehen werden kann.\\

Im anschlie"senden Abschnitt wird
festgehalten, wie man $C_0$"=Linksmultiplikatoren (bzw. unbeschr"ankte
schiefadjungierte Multiplikatoren) und die hiervon erzeugten $C_0$"=Halbgruppen (bzw. $C_0$-Gruppen)
von $X$ auf einen Hilbertraum $H$ "uberf"uhren kann.
Hiermit erh"alt man eine weitere Charakterisierung der Elemente
von $\Multlco{X}$ und $\Adjlco{X}$.\\

Im siebten Abschnitt wird gezeigt, unter welchen Voraussetzungen
man eine $C_0$-Halbgruppe von einem Hilbertraum $H$
auf einen Operatorraum $X \subseteq \mLinStet(H)$ "uberf"uhren kann.
Hierf"ur wird die sogenannte strikte $X$-Topologie eingef"uhrt und untersucht.\\

Im letzten Abschnitt des Kapitels wird ein St"orungsresultat
von S. Damaville (\cite{Damaville04RegulariteDesOp}) f"ur regul"are Operatoren
auf unbeschr"ankte Multiplikatoren verallgemeinert.
Au"serdem werden einige Resultate aus der St"orungstheorie f"ur Erzeuger
von $C_0$"=Halbgruppen auf Operatorr"aumen formuliert und
auf unbeschr"ankte Multiplikatoren angewendet.\\

\section[$C_0$-Linksmultiplikatoren und unbeschr"ankte schiefadjungierte Multiplikatoren]{$C_0$-Linksmultiplikatoren und unbeschr"ankte\newlinef{}schiefadjungierte Multiplikatoren}

Wir f"uhren in diesem Abschnitt zwei Arten von unbeschr"ankten Multiplikatoren auf
Operatorr"aumen ein, n"amlich die $C_0$-Linksmultiplikatoren
(als Menge: $\Multlco{X}$) und die unbeschr"ankten schiefadjungierten Multiplikatoren
(als Menge: $\Adjlco{X}$).
Letztere sind Erzeuger einer $C_0$-Gruppe von unit"aren Elementen
aus $\Adjlor{X}$.
Wir zeigen, \dass $\Adjlco{X}$  alle schiefadjungierten Elemente von $\Adjlor{X}$ enth"alt
und \dass in unitalen Operatorr"aumen Gleichheit gilt.
Des weiteren formulieren wir entsprechende Aussagen f"ur $\Multlco{X}$.
Ferner beleuchten wir den Zusammenhang zu den regul"aren Operatoren.
\skiptext

Wir erinnern kurz an einige Definitionen aus dem Anhang:

\begin{erinnerung}
Sei $X$ ein Operatorraum (siehe Definition~\ref{defOpraumAbstrakt}).
Wir bezeichnen mit $\Multlor{X}$ die Menge der Linksmultiplikatoren
auf $X$ (siehe Definition~\ref{defLinksmultOR})
und mit $\Adjlor{X}$ die Menge der von links adjungierbaren Multiplikatoren auf $X$
(siehe Definition~\ref{defAdjAbbOR}).

Es ist
\[ \mcalS(X) := \begin{pmatrix} \mbbC \Id_H & X \\ X^* & \mbbC \Id_H \end{pmatrix} \]
ein unitales Operatorsystem, welches Paulsen-System von $X$ genannt wird
(siehe Definitions"=Proposition~\ref{defPaulsenSys}).
F"ur die Ecke $(k,\ell)$ der injektiven H"ulle $I(\mcalS(X))$ von $\mcalS(X)$ verwenden wir
die Schreibweise $I_{k\ell}(X)$ f"ur alle $k,\ell \in \haken{2}$
(siehe Proposition~\ref{BleLeM4.4.2}).
Es gilt somit:
\[ I(\mcalS(X)) \cong \begin{pmatrix} I_{11}(X) & I_{12}(X) \\ I_{21}(X) & I_{22}(X) \end{pmatrix}. \]
Mit $j$ sei die kanonische Abbildung von $X$ nach $I_{12}(X)$ bezeichnet.

Es ist
\[ \mIMl(X) := \bigl\{ a \in I_{11}(X) \setfdg a \multis j(X) \subseteq j(X) \bigr\} \]
eine unitale Banachalgebra und
\[ \mIMl^*(X) := \mIMl(X) \cap \mIMl(X)^* \subseteq I(\mcalS(X)) \]
eine unitale \csalgebra.
F"ur alle $a \in \mIMl(X)$ wird durch
\[ \cmlmult_a : X \to X, x \mapsto j^{-1}(a \multis j(x)), \]
ein Linksmultiplikator definiert.
Dann ist
\[ \isoIMl : \mIMl(X) \arrowbij \Multlor{X}, a \mapsto \cmlmult_a, \]
ein unitaler, isometrischer Isomorphismus von $\mIMl(X)$ auf
die Banachalgebra $(\Multlor{X}, \normMl{\cdot}{X})$
(siehe Satz~\ref{IsomIMlAufMl}) und
$\isoIMl\restring_{\mIMls(X)} : \mIMls(X) \arrowbij \Adjlor{X}$
ein unitaler \sterns{}Iso\-mor\-phis\-mus auf $\Adjlor{X}$ (Proposition~\ref{AdjlorIsomIMls}).
\end{erinnerung}

\pagebreak
\begin{definitn}\label{defAdjlco}
Sei $X$ ein Operatorraum.
\begin{enumaufz}
\item Ein Operator $A : D(A) \subseteq X \to X$ hei"st
\defemph{$C_0$-Linksmultiplikator} auf $X$,
falls $A$ Erzeuger einer $C_0$-Halbgruppe $(T_t)_\indtHG$ auf $X$ mit $T_t \in \Multlor{X}$
f"ur alle $\indtHG$ ist.
\index[B]{C0Linksmultiplikator@$C_0$-Linksmultiplikator}%

\item Die Menge aller $C_0$-Linksmultiplikatoren auf $X$ wird mit
$\Multlco{X}$ bezeichnet.
\index[S]{Mscrlc0@$\Multlco{X}$}%

\item Ein Operator $A : D(A) \subseteq X \to X$ hei"st
\defemph{unbeschr"ankter schiefadjungierter Multiplikator} auf $X$,
falls $A$ Erzeuger einer $C_0$-Gruppe $(U_t)_\indtHG$ auf $X$
mit $U_t \in \Adjlor{X}$ unit"ar f"ur alle $\indtGr$ ist.
\index[B]{unbeschränkter schiefadjungierter Multiplikator}%
\dremark{Bezeichnung "uberdenken:
  Kann man so leicht mit unbeschr"ankten Multiplikatoren verwechseln.
  Etwas zur Bezeichnung notieren.}%
\dremark{Es gilt $\norm{U_t} = \normcb{U_t} = 1$,
  denn $U_t$ ist unit"ar, also $\normMl{U_t}{X} = 1$,
  und mit \ref{inAlXNormTEqNormcb} folgt: $\norm{U_t} = \normcb{U_t} = \normMl{U_t}{X} = 1$.}%

\item Die Menge aller unbeschr"ankten schiefadjungierten Multiplikatoren auf $X$ wird mit
$\Adjlco{X}$ bezeichnet.\dremark{\glqq{}s\grqq{} f"ur schiefadjungiert}
\index[S]{Adjlc0X@$\Adjlco{X}$}%
\end{enumaufz}
\end{definitn}

Analog zu Beispiel~\ref{bspRegOp}.(i) gilt:

\begin{bemerkung}\label{MultlorSubseteqMultlco}\label{AdjlorSubseteqAdjlco}
Sei $X$ ein Operatorraum.
Es gilt:
\dremark{Gilt eine entsprechende Aussage auch f"ur $\Adjloru{X}$?}%
\begin{enumaufz}
\item $\Adjlco{X} \subseteq \Multlco{X}$.
\item $\Multlor{X} = \{ A \in \Multlco{X} \setfdg D(A) = X \}$.
\dremark{Insbesondere gilt: $\Multlor{X} \subseteq \Multlco{X}$.}%
\item $\{ A \in \Adjlor{X} \setfdg A \text{ schiefadjungiert} \}
  = \{ A \in \Adjlco{X} \setfdg D(A) = X \}$.
\end{enumaufz}
\end{bemerkung}

\begin{proof}
\bewitemph{(i)} folgt unmittelbar.

\bewitemph{(ii):}
\bewitemphq{$\subseteq$}
Sei $A \in \Multlor{X}$.
Da $\Multlor{X}$ eine unitale Banachalgebra ist
(\refb{\cite{Zarikian01Thesis}, Theorem 1.6.2,}{Satz~\ref{MlXunitaleBA}}),
wird durch $T_t := \mre^{tA} \in \Multlor{X}$ f"ur alle $\indtHG$ eine
normstetige Halbgruppe auf $X$ mit Erzeuger $A$ definiert.
\dremark{Alternativer, l"angerer Beweis f"ur $\Multlor{X}$:
Dann findet man ein $a \in \mIMl(X)$ mit:
$j(Ax) = a \,j(x)$ f"ur alle $x \in X$.
Setze $T_t := \mre^{tA}$ f"ur alle $\indtHG$.
Dann ist $(T_t)_\indtHG$ eine $C_0$-Halbgruppe auf $X$.
Es gilt
\begin{align*}
   j(T_t x)
=  j(\mre^{tA}(x))
=  j\Bigl( \sum_{n=0}^\infty \frac{t^n}{n!} A^n(x) \Bigr)
=  \sum \frac{t^n}{n!} j(A^n x)
=  \sum \frac{t^n}{n!} a^n \multis j(x)
=  \mre^{ta} \multis j(x)
\end{align*}
f"ur alle $x \in X$, $\indtHG$, also $T_t \in \Multlor{X}$.
\dremark{Z. 21.8.'08/1, evtl. auch 8.1.'08}%
}%

\bewitemphq{$\supseteq$}
Sei $A \in \Multlco{X}$ mit $D(A) = X$.
Sei $(T_t)_\indtHG$ die von $A$ erzeugte $C_0$-Halbgruppe.
Nach \cite{WernerFunkana6}, Satz VII.4.9, gilt: $T_t = \mre^{tA}$ f"ur alle $\indtHG$.
Mit \cite{EngelNagelSemigroups}, Proposition I.3.5, erh"alt man:
$A = (t \mapsto T_t)'(0)$.
\dremark{$\Multlor{X}$ ist abgeschlossen,
  $\frac{T_t - T_0}{t} \in \Multlor{X}$,
  also $\lim_{t \downarrow 0} \frac{T_t - T_0}{t} \in \Multlor{X}$}%
Also folgt: $A \in \Multlor{X}$.

\bewitemph{(iii):}
\bewitemphq{$\subseteq$}
Weil $\Adjlor{X}$ eine unitale \csalgebra{} ist
(\refb{\cite{Zarikian01Thesis}, Proposition 1.7.4,}{Proposition~\ref{AlXunitaleCsAlg}}),
folgt dies analog zu (ii).
\dremark{unit"ar: $T_t^* = (\mre^{tA})^*
  = (\sum_{n=0}^\infty \frac{1}{n!}(tA)^n)^*
  = \sum \frac{1}{n!} t^n (A^*)^n
  = \sum \frac{1}{n!} t^n (-A)^n
  = \mre^{-tA}
  = T_{-t}$,
  also $T_t^* T_t = T_{-t} T_t = \mre^{-tA} \mre^{tA} = \mre^{-tA + tA} = \mre^0 = \Id_X$.}%

\bewitemphq{$\supseteq$}
Sei $A \in \Adjlco{X}$ mit $D(A) = X$.
Sei $(T_t)_\indtGr$ die von $A$ erzeugte $C_0$"=Gruppe.
Analog zu (ii) erh"alt man: $A \in \Adjlor{X}$.
Da $T_t = \mre^{tA}$ unit"ar f"ur alle $\indtGr$ ist,
folgt mit \cite{Upmeier85SymBanachManifolds}, Lemma~8.12: $A^* = -A$.
\end{proof}

\dfrage{Bei Woronowicz gilt:
  $\Multcs{\mfrakA} = \{ T \in \Multwor{\mfrakA} \setfdg \norm{T} < \infty \}$.
Gilt etwas "Ahnliches f"ur $\Multlco{X}$?\dremark{6.1.'08/1}}%

Das folgende Lemma ben"otigen wir f"ur den Beweis von Proposition~\ref{AdjlorsaEqAdjlco}:

\pagebreak
\begin{lemma}\label{atIstOpHalbgr}
Sei $X$ ein Operatorraum und $(T_t)_\indtHG$ eine $C_0$-Halb\-grup\-pe auf $X$
mit $T_t \in \Multlor{X}$ f"ur alle $\indtHG$.
F"ur jedes $\indtHG$ setze $a_t := \isoIMl^{-1}(T_t) \in \mIMl(X)$,
also hat man $j(T_t x) = a_t \multis j(x)$ f"ur alle $x \in X$.
Es gilt:
\begin{enumaufz}
\item $a_0 = e_{I_{11}(X)}$,
  wobei mit $e_{I_{11}(X)}$ das Einselement von $I_{11}(X)$ notiert wird.
\item $a_{s+t} = a_s \multis a_t$ f"ur alle $s,\indtHG$.
\dremark{$a_s \multis a_t = \Phi(a_s a_t)$}%
\item $\lim_{t \to 0} a_t \multis y = y$ f"ur alle $y \in \overline{\alg\left(j(X)\right)}$,
  wobei $\alg(j(X))$ die von $j(X)$ in $I(\mcalS(X))$ erzeugte Unteralgebra bezeichnet.
\item Definiere $R_t : j(X) \to j(X), z \mapsto a_t \multis z$, f"ur alle $\indtHG$.
Dann ist $(R_t)_\indtHG$ eine $C_0$-Halbgruppe auf $j(X)$.
\end{enumaufz}
\end{lemma}

\begin{proof}
\bewitemph{(i):} F"ur alle $x \in X$ gilt nach \eqref{eqeI11xEqx} die Gleichung
$e_{I_{11}(X)} \multis j(x) = j(x)$, also ergibt sich:
$0 = j(T_0 x) - j(x) \dremarkm{= a_0 \multis j(x) - e_{I_{11}(X)} \multis j(x)} = (a_0 - e_{I_{11}(X)}) \multis j(x)$.
Mit \refb{\cite{BlecherLeMerdy04OpAlg}, Proposition 4.4.12,}{Proposition~\ref{Lcinj}} folgt:
$a_0 = e_{I_{11}(X)}$.

\bewitemph{(ii)} erh"alt man analog zu (i).
\dremark{$a_{s+t} \multis j(x) = j(T_{s+t} x) = j(T_s T_t x)
  = a_s \multis j(T_t x) = a_s \multis a_t \multis j(x)$}%

\bewitemph{(iii):} F"ur alle $x,y \in X$ gilt:
\[ \lim_{t \downarrow 0} a_t \multis j(x) \multis j(y)
\dremarkm{=  \lim_{t \downarrow 0} j(T_t x) \multis j(y)}
=  j\left(\lim_{t \downarrow 0} T_t x \right) \multis j(y)
=  j(x) j(y). \]
Wenn man dies analog auf Elemente aus $\alg(j(X))$ anwendet,
erh"alt man:
 $\lim_{t \to 0} a_t \multis y = y$ f"ur alle $y \in \alg(j(X))$.
Indem man die Stetigkeit von $j$ benutzt, folgt (iii).
\dremark{(iii): Genauer: Siehe Z. 26.4.'07, S. 1.
  Beweis ist dort recht technisch, daher nicht eingetippt.}%

\bewitemph{(iv)} folgt mit (i), (ii) und (iii).
\dremark{Mit (i) und (ii) erh"alt man,
\dass $(R_t)_\indtHG$ eine Operatorhalbgruppe ist.
\dremark{$R_0(z) = a_0 \multis z = e_{I_{11}} \multis z \overset{\eqref{eqeI11xEqx}}{=} z$,
  $R_s R_t z = a_s \multis a_t \multis z = a_{s+t} \multis z = R_{s+t} z$}%
Da $(T_t)_\indtHG$ stark stetig ist, ist auch $(R_t)_\indtHG$ stark stetig.}%
\dremark{$\lim R_t(z) = \lim a_t \multis z \overset{(iii)}{=} z$}%
\dremark{Weiter gilt f"ur alle $x \in X$:
\[ \normlr{ R_t(j(x)) - j(x) }
=  \norm{ a_t \multis j(x) - j(x) }
=  \normlr{ j(T_t x) - j(x) }
=  \norm{ T_t x - x }
\to 0 \]
f"ur $t \downarrow 0$.}%
\end{proof}

\dfrage{Unter welchen Voraussetzungen ist $(a_t)_\indtHG$ eine $C_0$-Halbgruppe?}%

Eine entsprechende Fassung von Lemma~\ref{atIstOpHalbgr} f"ur $C_0$-Gruppen lautet:

\begin{lemma}\label{atIstOpGr}
Sei $X$ ein Operatorraum und $(T_t)_\indtGr$ eine $C_0$-Gruppe auf $X$
mit $T_t \in \Adjlor{X}$ f"ur alle $\indtGr$.
F"ur jedes $\indtGr$ setze $a_t := \isoIMl^{-1}(T_t) \in \mIMls(X)$.
Es gilt:
\begin{enumaufz}
\item $a_0 = e_{I_{11}(X)}$.
\item $a_{s+t} = a_s \multis a_t$ f"ur alle $s,\indtGr$.
\dremark{$a_s \multis a_t = \Phi(a_s a_t)$}%
\item $\lim_{t \to 0} a_t \multis y = y$ f"ur alle $y \in \overline{\alg\left(j(X)\right)}$.
\item Definiere $R_t : j(X) \to j(X), z \mapsto a_t \multis z$, f"ur alle $\indtGr$.
Dann ist $(R_t)_\indtGr$ eine $C_0$-Gruppe auf $j(X)$.
\item F"ur alle $\indtGr$ gilt: Ist $T_t$ unit"ar, so auch $a_t$.
\end{enumaufz}
\end{lemma}

\begin{proof}
\bewitemph{(i)--(iv)} erh"alt man analog zum Beweis von Lemma~\ref{atIstOpHalbgr}.

\bewitemph{(v)} folgt mit
\refb{\cite{BlecherLeMerdy04OpAlg}, Proposition 4.4.12,}{Proposition~\ref{Lcinj}}.
\dremark{$a_t^* \multis a_t \multis j(x) = j(T_t^* T_t x) = j(x)$
  f"ur alle $x \in X$, also folgt mit \ref{Lcinj}: $a_t^* a_t = e_{I_{11}}$.}%
\end{proof}

Es sei daran erinnert, \dass ein Operatorraum $X$ unital hei"st,
falls es eine unitale vollst"andige Isometrie von $X$ in eine unitale \csalgebra gibt.
\skiptext

In einer unitalen \csalgebra $\mfrakA$ gilt:
$\Multwor{\mfrakA} \cong \Multcs{\mfrakA}$
(\refb{vgl. \cite{LanceHmod}, S. 117,}{Beispiel~\ref{bspRegOp}}).
Nach der folgenden Proposition gilt ein entsprechendes Resultat in unitalen Operatorr"aumen.
Somit erh"alt man insbesondere eine notwendige Bedingung daf"ur, \dass ein Operatorraum unital ist.

\begin{bemerkung}\label{AdjlorsaEqAdjlco}
Sei $X$ ein unitaler Operatorraum.
Dann gilt:
\begin{enumaufz}
\item $\Multlor{X} = \Multlco{X}$.
\item $\bigl\{ A \in \Adjlor{X} \setfdg A \text{ schiefadjungiert} \bigr\} = \Adjlco{X}$.
\end{enumaufz}
\end{bemerkung}

\begin{proof}
\bewitemph{(i):}
\bewitemph{\glqq$\subseteq$\grqq} folgt mit Proposition~\ref{MultlorSubseteqMultlco}.

\bewitemph{\glqq$\supseteq$\grqq:}
Mit $e_X$ sei das Einselement von $X$ bezeichnet.
Da die Definition eines unitalen Operatorraumes
unabh"angig von der Einbettung in eine spezielle \csalgebra{} ist
(Anmerkung~\ref{AnmUnitalerORunabh}),
findet man o.B.d.A. einen Hilbertraum $H$ mit $X \subseteq L(H)$ und $e_X = \Id_H$.
\dremark{Sei $\varphi : X \to \mfrakA$ eine unitale vollst"andige Isometrie,
  $\psi : \mfrakA \to \mLinStet(H)$ eine unitale, injektive \sterns{}Darstellung,
  also insbesondere eine vollst"andige Isometrie.
  Dann ist $\psi \circ \phi$ eine unitale vollst"andige Isometrie.}%
Sei $A \in \Multlco{X}$.
Dann ist $A$ Erzeuger einer $C_0$-Halbgruppe
$(T_t)_\indtHG$ auf $X$ mit $T_t \in \Multlor{X}$ f"ur alle $\indtHG$.
F"ur jedes $\indtHG$ gibt es ein $a_t \in \mIMl(X)$ mit
$j(T_t x) = a_t \multis j(x)$.
Da nach Proposition~\ref{BleLeM4.4.2} gilt $I_{11}(X) \subseteq M_2(\mLinStet(H))$,
findet man f"ur alle $\indtHG$ ein $\tilde{a}_t \in \mLinStet(H)$ mit
$a_t = \begin{smallpmatrix} \tilde{a}_t & 0 \\ 0 & 0 \end{smallpmatrix}$.
Da $X$ unital ist, gilt nach \cite{Blecher01ShilovBoundary}, S. 20,
f"ur alle $S \in \Multlor{X}$:\dremark{$(\ddagger)$} $\normcb{S} = \normMl{S}{X}$.
Au"serdem gilt:\dremark{nach \eqref{eqeIXEq0IdH00}}
$e_{I(X)} = \begin{smallpmatrix} 0 & \Id_H \\ 0 & 0 \end{smallpmatrix}$.
Da $j : X \to I(X)$ unital
(\refb{\cite{BlecherLeMerdy04OpAlg}, Corollary 4.2.8,}{Proposition~\ref{unitalerORFolgtjUnital}})
und $\isoIMl^{-1}$ nach Satz~\ref{IsomIMlAufMl} eine unitale Isometrie ist,
erh"alt man:
\begin{align*}
   \norm{T_t - \Id_X}
&\leq  \normcb{T_t - \Id_X}
\overset{\dremarkm{(\ddagger)}}{=}  \normMl{T_t - \Id_X}{X}  \nonumber \\
&= \norm{\isoIMl^{-1}(T_t - \Id_X)}
=  \norm{a_t - e_{I_{11}(X)}}  \nonumber \\
&= \normlr{ \cmpmatrix{\tilde{a}_t & 0 \\ 0 & 0} - \cmpmatrix{\Id_H & 0 \\ 0 & 0} }
=  \norm{\tilde{a}_t - \Id_H}    \\
&\overset{\dremark{(*)}}{=} \norm{ a_t \multis e_{I(X)} - e_{I(X)} }
=  \norm{ a_t \multis j(e_X) - j(e_X) }  \nonumber \\
&= \norm{ T_t(e_X) - e_X }
\to 0  \qquad\text{f"ur } t \to 0.  \nonumber
\end{align*}
\dmarginpar{pr}\dremark{Beweis ist an der Stelle, wo $(*)$ verwendet wird, arg kurz (DB).}%
\dremark{$(\dagger)$: $
   \cmpmatrix{ 0 & \tilde{a}_t \\ 0 & 0 }
=  \Phi(\cmpmatrix{ 0 & \tilde{a}_t \Id_H \\ 0 & 0 })
=  \Phi(\begin{pmatrix} \tilde{a}_t & 0 \\ 0 & 0 \end{pmatrix}
             \begin{pmatrix} 0 & \Id_H \\ 0 & 0 \end{pmatrix})
=  \begin{pmatrix} \tilde{a}_t & 0 \\ 0 & 0 \end{pmatrix} \multis
             \begin{pmatrix} 0 & \Id_H \\ 0 & 0 \end{pmatrix}
=  a_t \multis e_{I(X)}$,
$(*)$:
$  \norm{\tilde{a}_t - \Id_H}
=  \normBig{ \begin{pmatrix} 0 & \tilde{a}_t - \Id_H \\ 0 & 0 \end{pmatrix} }
=  \normBig{ \begin{pmatrix} 0 & \tilde{a}_t \\ 0 & 0 \end{pmatrix} -
            \begin{pmatrix} 0 & \Id_H \\ 0 & 0 \end{pmatrix} }
\overset{(\dagger)}{=}  \norm{ a_t \multis e_{I(X)} - e_{I(X)} }$}%
Somit ist $(T_t)_\indtHG$ eine normstetige Halbgruppe.
Nach \cite{WernerFunkana6}, Satz VII.4.9, ist $A \in L(X)$
und $T_t = \mre^{tA}$ f"ur alle $\indtHG$.
Mit \cite{EngelNagelSemigroups}, Proposition I.3.5, erh"alt man $A = (t \mapsto T_t)'(0)$,
also folgt: $A \in \Multlor{X}$.
\dremark{$\Multlor{X}$ ist unitale Banachalgebra,
  $\frac{T_t - T_0}{t} \in \Multlor{X}$,
  also $\lim_{t \downarrow 0} \frac{T_t - T_0}{t} \in \Multlor{X}$.}%
\smallskip

\bewitemph{(ii):}
\bewitemph{\glqq$\subseteq$\grqq} folgt mit Proposition~\ref{MultlorSubseteqMultlco}.

\bewitemph{\glqq$\supseteq$\grqq:}
Sei $A \in \Adjlco{X}$.
Wir verwenden die Bezeichnungen aus dem Beweisteil~(i).
Analog zu (i) folgt: $A \in \Adjlor{X}$.
Nach \cite{WernerFunkana6}, Satz VII.4.9, gilt: $T_t = \mre^{tA}$ f"ur alle $\indtGr$.
Da $T_t$ f"ur alle $\indtGr$ unit"ar ist,
folgt mit \cite{Upmeier85SymBanachManifolds}, Lemma~8.12: $A^* = -A$.
%
\dremark{Sei $\mfrakA$ unitale Banachalgebra.
Es ist $U(\mfrakA) = \{ a \in G(\mfrakA) \setfdg \norm{a} = \norm{a^{-1}} \leq 1 \}
\overset{(*)}{=}  \{ a \in G(\mfrakA) \setfdg a^* a = e_\mfrakA \}$
($(*)$ gilt nach Dennis) eine Banach-Lie-Untergruppe von $G(\mfrakA)$
(siehe auch \cite{Upmeier85SymBanachManifolds}, S. 116).
Weiter gilt:
$\{ X \in \mfrakA \setfdg X^* + X = 0 \}
=  u\ell(\mfrakA)
=  \operatorname{Lie}(U(\mfrakA))
\overset{\text{Lemma 8.12}}{=}
  \{ X \in \mfrakA \setfdg \forall \indtGr : \exp(tX) \in U(\mfrakA) \}$,
wobei $\operatorname{Lie}(U(\mfrakA))$ die Lie-Algebra von $U(\mfrakA)$ ist.
\dremark{Beweisidee f"ur \cite{Upmeier85SymBanachManifolds}, Lemma 8.12: Bei $0$ ableiten.}%
\dremark{$\rho(A^* x) = B^* \circ \rho(x) = -B \circ \rho(x) = \rho(-Ax)$.}%
}%
\dremark{Z. 18.1.'07/1, 19.8.'08/1}%
\end{proof}

\subsection{Beispiele}

In diesem Abschnitt f"uhren wir Beispiele f"ur $C_0$-Linksmultiplikatoren
und unbeschr"ankte schiefadjungierte Multiplikatoren auf.
\skiptext

In \hcsmoduln gilt:

\begin{satz}\label{MultworInAdjlco}
Sei $E$ ein Hilbert-\cstern{}Modul.
Dann gilt:
\dremark{Man h"atte gern \glqq$\supseteq$\grqq, Problem:
  Auf $\Adjlco{E}$ hat man keine Adjungierte.}%
\[ \left\{ T \in \Multwor{E} \setfdg T \text{ schiefadjungiert} \right\} = \Adjlco{E}. \]
\end{satz}

In den Beweis gehen entscheidend der Satz von Stone f"ur \hcsmoduln (Satz~\ref{HollevoetTh21})
und Satz~\ref{erzeugerUnitaereGrWor} ein.

\begin{proof}
\bewitemphq{$\subseteq$}
Sei $T \in \Multwor{E}$ schiefadjungiert.
Dann ist $h := -\mri T$ selbstadjungiert.
Es gilt $\Adjhm{E} \cong \Adjlor{E}$
(\refb{vgl. \cite{BlecherLeMerdy04OpAlg}, Corollary 8.4.2,}{Beispiel~\ref{HCsModAdjMultEqAdj}}).
Nach Proposition~\ref{UtIstSSUnitaereGr} wird durch $U_t := \varphi_h(e_t) = \exp(\mri th)$
eine $C_0$-Gruppe auf $E$ mit
$U_t \in \Adjhm{E} \cong \Adjlor{E}$ unit"ar f"ur alle $\indtGr$ definiert,
deren Erzeuger als $C_0$-Gruppe nach Satz~\ref{erzeugerUnitaereGrWor} gleich $\mri h = T$ ist.

\dmarginpar{Bew. pr}%
\bewitemphq{$\supseteq$}
Sei $T \in \Adjlco{E}$ und $T$ erzeuge die $C_0$-Gruppe $(U_t)_\indtGr$
auf $E$ mit $U_t \in \Adjlor{E} \cong \Adjhm{E}$ unit"ar f"ur alle $\indtGr$.
Dann findet man nach Satz~\ref{HollevoetTh21} (Satz von Stone)
ein selbstadjungiertes $h \in \Multwor{E}$ so,
\dass $U_t = \varphi_h(e_t) = \exp(\mri th)$ f"ur alle $\indtGr$ gilt.
Da $\mri h$ nach Satz~\ref{erzeugerUnitaereGrWor} Erzeuger der $C_0$-Gruppe $(U_t)_\indtGr$ ist,
folgt wegen der Eindeutigkeit des Erzeugers: $T = \mri h$.
Also ist $T$ schiefadjungiert.
\end{proof}

\begin{beispiel}\label{bspMultlcoC0Omega}
Sei $\Omega$ ein lokalkompakter Hausdorffraum.
F"ur alle $q \in C(\Omega)$ definiere
$M_q : D(M_q) \subseteq C_0(\Omega) \to C_0(\Omega), f \mapsto qf$.
Dann gilt:
\dmarginpar{zutun}\dremark{$\Adjlco{X}$ bestimmen.}%
\dremark{Beachte: $\Multlor{X} = \Adjlor{X}$}%
\dremark{Im Text mit $\Multwor{X}$ vergleichen.}%
\begin{align*}
   \Multlco{C_0(\Omega)}
&= \left\{ M_q \setfdg
   q \in C(\Omega) \text{ mit } \sup_{\omega \in \Omega} \Re q(\omega) < \infty \right\}  \dremarkm{\\
&\cong  \left\{ q \in C(\Omega) \setfdg \sup_{\omega \in \Omega} \Re q(\omega) < \infty \right\}} \text{ und} \\
   \Adjlco{C_0(\Omega)}
&= \Bigl\{ M_q \setfdg
   q \in C(\Omega) \text{ mit } \overline{q} = -q \Bigr\}  \dremarkm{\\
&\cong  \left\{ q \in C(\Omega) \setfdg  \overline{q} = -q \right\}}.
\end{align*}
\end{beispiel}

Setze $X := C_0(\Omega)$.
Dann ist die Menge $\Multlor{X}$ der Linksmultiplikatoren auf dem Operatorraum $X$
isomorph zu der Menge $\Multlcs{X}$ der Linksmultiplikatoren auf der \csalgebra $X$
\dremark{Kurz: $\Multlor{X} \cong \Multlcs{X}$.}%
und $\Adjlor{X}$ isomorph zu der Menge $\Multcs{X}$
der Multiplikatoren auf der \csalgebra $X$
(\refb{vgl. \cite{BlecherLeMerdy04OpAlg}, Corollary 8.4.2,}{Beispiel~\ref{bspMlCsAlgEqMlOR}}).
Da $X$ eine kommutative \csalgebra ist, folgt somit:
\begin{align}
\Multlor{X} &\cong \Multlcs{X} = \Multcs{X} \cong C_b(\Omega)\quad \text{und} \label{eqMultlorC0MEqCbM} \\
\Adjlor{X} &\cong \Multcs{X} \cong C_b(\Omega). \nonumber
\end{align}
Nach dem obigen Beispiel\dremark{\ref{bspMultlcoC0Omega}} gilt daher im allgemeinen:
$\Multlor{X} \neq \Multlco{X}$ und
$\{ T \in \Adjlor{X} \setfdg T \text{ schiefadjungiert} \} \neq \Adjlco{X}$.

\begin{proof}[Beweis von Beispiel~\ref{bspMultlcoC0Omega}]
\bewitemph{(i):}
\bewitemph{\glqq$\subseteq$\grqq:}
Sei $A$ Erzeuger einer $C_0$-Halbgruppe $(T_t)_\indtHG$ auf $X := C_0(\Omega)$ mit
$T_t \in \Multlor{X}$ f"ur alle $\indtHG$.
Nach \eqref{eqMultlorC0MEqCbM} gilt: $\Multlor{X} \cong C_b(\Omega)$.
Somit findet man f"ur jedes $\indtHG$ ein $m_t \in C_b(\Omega)$ mit
$T_t = M_{m_t}$.
Nach \cite{EngelNagelSemigroups}, Proposition~I.4.6, findet man ein $q \in C(\Omega)$ mit
$\sup_{\omega \in \Omega} \Re q(\omega) < \infty$ derart, \dass gilt:
$m_t(\omega) = \mre^{t q(\omega)}$ f"ur alle $\indtHG$ und $\omega \in \Omega$.
Nach \cite{EngelNagelSemigroups}, II.2.9, ist somit $M_q$ Erzeuger von $(T_t)_\indtHG$.
\smallskip

\bewitemph{\glqq$\supseteq$\grqq:}
Sei $q \in C(\Omega)$ mit $\sup_{\omega \in \Omega} \Re q(\omega) < \infty$.
Nach \cite{EngelNagelSemigroups}, S. 27 und Proposition~I.4.5, erzeugt $M_q$ eine $C_0$-Halbgruppe
$(T_t)_\indtHG$ auf $X$, wobei $T_t = M_{\mre^{tq}} \in \Multlor{X}$ f"ur alle $\indtHG$ ist.
\smallskip

\bewitemph{(ii)} folgt durch Nachrechnen oder
mit Satz~\ref{MultworInAdjlco} und Beispiel~\ref{bspRegOpC0Omega}.
\dremark{Beweis ohne \ref{MultworInAdjlco}:
  Man erh"alt, \dass $\mre^{t(q+\overline{q})} = \mathbbm{1}$,
  also $\overline{q} = -q$ gelten mu\cms.}%
\dremark{11.1.'08}%
\end{proof}

F"ur die kompakten Operatoren $\kptOp(H)$ auf einem Hilbertraum $H$ gilt:

\begin{beispiel}\label{bspMultlcoKH}
Sei $H$ ein Hilbertraum.
Dann gilt:
\dmarginpar{zutun}\dremark{$\Adjlco{\kptOp(H)}$ bestimmen.}%
\dremark{Folgt dies auch mit einem Charakterisierungssatz?
  Wenn ja, erw"ahnen. (WW)}%
\begin{align*}
\Multlco{\kptOp(H)}
\cong
\bigl\{ & A : D(A) \subseteq H \to H \setfdg \\
& A \text{ ist Erzeuger einer $C_0$-Halbgruppe auf $H$} \bigr\} \text{ und} \\
\Adjlco{\kptOp(H)}
\cong
\bigl\{ & A : D(A) \subseteq H \to H \setfdg A \text{ ist linear,} \\
& \text{dicht definiert und abgeschlossen mit } A^* = -A \bigr\}.
\end{align*}
\end{beispiel}

Setze $X := \kptOp(H)$.
Wie man sieht, h"angen unbeschr"ankte Multiplikatoren auf dem
Operatorraum $X \subseteq \mLinStet(H)$ mit Erzeugern von \cohgn (bzw. $C_0$-Gruppen)
auf $H$ zusammen.
Im Abschnitt~\ref{secStrikteXTop} wird untersucht,
wie man f"ur einen beliebigen Operatorraum $Y \subseteq \mLinStet(K)$
\cohgn auf $K$ in \cohgn auf $Y$ "uberf"uhren kann.
\skiptext

Zum Beweis des ersten Teils von Beispiel~\ref{bspMultlcoKH} verwenden wir das folgende

\begin{lemma}[\cite{BlackadarOpAlg}, Proposition I.8.1.4]\label{KHStarkKonvFolgtKonv}
Seien $X$, $Y$, $Z$ Banachr"aume.
Sei $T_0 \in \mLinStet(Y,Z)$ und
$(T_\lambda)_\lambda$ ein gleichm"a"sig beschr"anktes Netz in $\mLinStet(Y,Z)$
mit $T_\lambda \overset{\lambda}{\to} T_0$ stark.
Sei $K \in \kptOp(X,Y)$.
Dann gilt: $\norm{T_\lambda K - T_0 K} \overset{\lambda}{\to} 0$.
\dremark{Kann man auf Beschr"anktheit verzichten? Z.B. falls X,Y,Z Hilbertr"aume?}%
\dremark{Evtl. \glqq{}glm. beschr"ankt\grqq{} definieren.}%
\dliter{\cite{BlackadarOpAlg}, I.8.1.4, dort ohne Beweis; in \cite{Weidmann00LinOperatoren}, "Ubung 3.5.(a)
  wird Beschr"anktheit nicht verlangt.}%
\dremark{Bew.: 25.1.'08}%
\end{lemma}

\begin{proof}[Beweis von Beispiel~\ref{bspMultlcoKH}]
\bewitemph{(i):}
\bewitemphq{$\subseteq$}
Sei $B \in \Multlco{\kptOp(H)}$.
Sei $(T_t)_\indtHG$ die von $B$ erzeugte $C_0$-Halbgruppe auf $\kptOp(H)$ mit
$T_t \in \Multlor{\kptOp(H)}$ f"ur alle $\indtHG$.
Da $\kptOp(H)$ eine \csalgebra ist, gilt:
$\Multlor{\kptOp(H)} \cong  \Multlcs{\kptOp(H)}$
(\refb{vgl. \cite{BlecherLeMerdy04OpAlg}, Corollary 8.4.2,}{Beispiel~\ref{bspMlCsAlgEqMlOR}}).
Wegen $\Multlcs{\kptOp(H)} \cong L(H)$ findet man f"ur alle $\indtHG$ ein $S_t \in L(H)$ mit
$T_t = L_{S_t}$, wobei
\[ L_{S_t} : \kptOp(H) \to \kptOp(H), K \mapsto S_t \circ K. \]
Dann ist $(S_t)_\indtHG$ eine Operatorhalbgruppe auf $H$.
\dremark{$L_{S_0} = T_0 = \Id_{\kptOp(H)}$, also $S_0 = \Id_H$.
  $S_s S_t = L_{S_s} L_{S_t} (\Id_H) = L_{S_{s+t}}(\Id_H) = S_{s+t}$.}%
Da $(T_t)_\indtHG$ stark stetig ist, gilt f"ur alle $K \in \kptOp(H)$:
\[ \norm{S_t \circ K - K}
=  \norm{T_t(K) - K}
\to 0  \qquad\text{f"ur } t \to 0. \]
Insbesondere folgt: $S_t \to \Id_H$ f"ur $t \to 0$ stark.
\dremark{Sei $\xi \in H$.
  Wegen $P := P_{\vrerz{\xi}} \in \kptOp(H)$ folgt:
$\norm{S_t(\xi) - \xi}
=  \norm{S_t(P(\xi)) - P(\xi)}
=  \norm{T_t(P)(\xi) - P(\xi)}
\to 0$ f"ur $t \to 0$.}%
Somit ist $(S_t)_\indtHG$ stark stetig.
\smallskip

\bewitemphq{$\supseteq$}
Sei $A$ Erzeuger einer $C_0$-Halbgruppe $(S_t)_\indtHG$ auf $H$.
Nach der Charakterisierung der starken Stetigkeit von Halbgruppen (Proposition~\ref{charHGIstC0})
findet man ein $\delta \in \mbbR_{>0}$
und ein $M \in \mbbR_{\geq 1}$ so, \dass gilt:
\[ \forall t \in [0,\delta] : \norm{S_t} \leq M. \]
Definiere $T_t := L_{S_t}$
f"ur alle $\indtHG$.
Somit ist $(T_t)_\indtHG$ eine Operatorhalbgruppe auf $\kptOp(H)$ mit
$T_t \in \Multlor{\kptOp(H)}$ f"ur alle $\indtHG$.
Da $(S_t)_{t \in [0,\delta]}$ ein gleichm"a"sig beschr"anktes Netz ist,
erh"alt man mit Lemma~\ref{KHStarkKonvFolgtKonv} f"ur alle $K \in \kptOp(H)$:
\[ \norm{T_t(K) - K}
=  \norm{S_t \circ K - K}
\to 0  \qquad\text{f"ur } t \to 0. \]
Also ist $(T_t)_\indtHG$ eine $C_0$-Halbgruppe auf $\kptOp(H)$.
Offensichtlich gilt: $T_t \in \Multlor{\kptOp(H)}$ f"ur alle $\indtHG$.
\dremark{Linksmultiplikator mit $\sigma = \Id_H$}%
\dremark{Bew.: Z. 28.1.'08/1. Benutzt wird \ref{KHStarkKonvFolgtKonv}.}
\smallskip

\bewitemph{(ii)} folgt mit Satz~\ref{MultworInAdjlco} und Beispiel~\ref{bspRegOpC0Omega}.
\dremark{Folgt ohne \ref{MultworInAdjlco} mit dem Satz von Stone f"ur Hilbertr"aume.}%
\end{proof}

\dfrage{Auf $C_0(\Omega)$ ist jede Multiplikationshalbgruppe quasikontraktiv.
Gilt etwas "Ahnliches f"ur $C_0$-Halbgruppen in $\Multlor{X}$?}%

\dremark{
\begin{definitn}
Sei $X$ ein Operatorraum.
\begin{enumaufz}
\item Ein Operator $A : D(A) \subseteq X \to X$ hei"st \defemphi{$C_0$-Linksmultiplikator},
falls ein $B : D(B) \subseteq X \to X$ so existiert, \dass
$\begin{pmatrix} 0 & \mri A & 0 & 0 \\ \mri B & 0 & 0 & 0 \\ 0 & 0 & 0 & 0 \\ 0 & 0 & 0 & 0 \end{pmatrix}$
Erzeuger einer \cmvk{en} $C_0$-Halbgruppe auf $C_4(X)$ ist.
\dremark{Kann man evtl. stets $B=0$ w"ahlen?}%

\item Die Menge aller $C_0$-Linksmultiplikatoren auf $X$ wird mit
$\Multlco{X}_4$ bezeichnet.
\end{enumaufz}
\end{definitn}

\begin{satz}\label{REsubseteqMlco}
Sei $E$ ein Hilbert-\cstern{}Modul "uber $\mfrakA$.
Dann gilt $\Multwor{E} \subseteq \Multlco{E}_4$.
\dremark{Evtl. als Motivation f"ur die Definition von $\Multlco{X}$ anf"uhren.}%
\dremark{Evtl. folgt hier noch mehr.}%
\end{satz}

\begin{proof}
Es ist $h := \begin{pmatrix} 0 & T \\ T^* & 0 \end{pmatrix}$ selbstadjungiert.
Somit wird durch $U_t := \exp(\mri th)$ eine unit"are $C_0$-Gruppe
in $\Adjhm{E \oplus E} \subseteq \Multlor{E \oplus E}$ definiert,
deren Erzeuger als $C_0$-Gruppe nach \ref{erzeugerUnitaereGrWor} gleich $\mri h$ ist.
Es gilt
\[ \normlr{\binom{U_t}{\Id_{E \oplus E}}\binom{x}{y}}
=  \normlr{\binom{U_t x}{y}}
=  \normlr{ \begin{pmatrix} u_t & 0 \\ 0 & 1 \end{pmatrix} \begin{pmatrix} x \\ y \end{pmatrix} }
=  \max\{\norm{u_t},1\} \cdot \normlr{ \binom{x}{y} }
\leq  \normlr{ \binom{x}{y} }, \]
also erzeugt $\begin{pmatrix} h & 0 \\ 0 & 0 \end{pmatrix}$
eine \cmvk{e} $C_0$-Halbgruppe.
\end{proof}
}%

\section{Unbeschr"ankte Multiplikatoren auf Operatorr"aumen}

In diesem Abschnitt wird eine weitere Art von Multiplikatoren
auf einem Operatorraum $X$ definiert, n"amlich die unbeschr"ankten
Multiplikatoren (als Menge: $\Adjloru{X}$), die wir
mit Hilfe der unbeschr"ankten schiefadjungierten Multiplikatoren $\Adjlco{X}$
aus Definition~\ref{defAdjlco} definieren.
Wir verwenden verschiedene Aussagen "uber $\Adjlco{X}$,
um Eigenschaften von $\Adjloru{X}$ zu zeigen.

F"ur die unbeschr"ankten Multiplikatoren auf einem beliebigen \hcsmodul $E$
gilt das wichtige Resultat: $\Multwor{E} = \Adjloru{E}$.
Also sind die unbeschr"ankten Multiplikatoren
eine Verallgemeinerung der regul"aren Operatoren auf Operatorr"aume.

Weiterhin beweisen wir, \dass $\Adjlor{X} \subseteq \Adjloru{X}$ gilt
und \dass in unitalen Operatorr"aumen Gleichheit vorliegt.
\skiptext

F"ur einen beliebigen \hcsmodul $E$ gilt $E \oplus E \cong C_2(E)$
(Proposition~\ref{normHCsModEoplusEEqC2E}),
wobei mit $C_2(X) = M_{2,1}(X)$ f"ur einen beliebigen Operatorraum $X$
der Spaltenoperatorraum bezeichnet wird.
Daher benutzen wir in Operatorr"aumen analog zu Definition~\ref{defMatABCDHCsMod}
die folgende Notation:

\begin{notation}
Sei $X$ ein Operatorraum.
Sei $A_i : D(A_i) \subseteq X \to X$ ein Operator f"ur alle $i \in \haken{4}$.
Definiere
\index[S]{ABCD@$\binom{A_1 \,\,\, A_2}{A_3 \,\,\, A_4}$}%
\begin{align*}
 \cmpmatrix{A_1 & A_2 \\ A_3 & A_4} :
&\,\,(D(A_1) \cap D(A_3)) \cmtimes (D(A_2) \cap D(A_4)) \subseteq C_2(X)  \to  C_2(X),  \\
&\,\,(x,y) \mapsto (A_1 x + A_2 y, A_3 x + A_4 y).
\end{align*}
\end{notation}

\begin{definitn}\label{defAdjloru}
Sei $X$ ein Operatorraum.
\begin{enumaufz}
\item Eine Abbildung $T : D(T) \subseteq X \to X$ hei"st
\defemphi{unbeschränkter Multiplikator} auf $X$,
falls ein $S : D(S) \subseteq X \to X$ existiert mit der Eigenschaft:
$\tilde{T} := \mri \begin{pmatrix} 0 & T \\ S & 0 \end{pmatrix} \in \Adjlco{C_2(X)}$,
\dheisst $\tilde{T}$ erzeugt eine
$C_0$-Gruppe $(U_t)_\indtGr$ auf $C_2(X)$ mit $U_t \in \Adjlor{C_2(X)}$ unit"ar
f"ur alle $\indtGr$.
\dremark{Etwas zum Namen sagen.}
\dremark{Notwendige/hinreichende Voraussetzung f"ur die Existenz von $S$?}

\item Mit $\Adjloru{X}$ wird die Menge der unbeschr"ankten Multiplikatoren auf $X$
bezeichnet.
\index[S]{AdjluX@$\Adjloru{X}$}%
\end{enumaufz}
\end{definitn}

Der Operator $S$ aus (i) ist, wie in Proposition~\ref{AdjloruSIstEind} gezeigt wird, eindeutig.
Ist $X$ ein \hcsmodul, so erh"alt man: $S = T^*$ (Satz~\ref{zshgRegOpUnbeschrMult}).
Somit kann man $S$ als die Adjungierte von $T$ ansehen.
\skiptext

\dremww{
Mit \ref{normxijLeqNormx} erh"alt man:

\begin{bemerkung}\label{dichtInC2}
Sei $X$ ein Operatorraum,
seien $U$, $V$ Untervektorr"aume von $X$ mit $U \oplus V$ dicht in $C_2(X)$.
Dann sind $U$ und $V$ dicht in $X$.
\end{bemerkung}

Beweis.
Sei $x_0 \in X$.
Dann findet man eine gegen $\binom{x_0}{0}$ konvergente Folge
$x \in (U \oplus V)^\mbbN$.
Es gilt:\dremark{17.10.'08/1}
\[ \norm{x_n^{(1)} - x_0}
\overset{\eqref{normxijLeqNormx}}{\leq}
   \normBig{ x_n - \binom{x_0}{0} }_{C_2(X)}
\to 0 \qquad\text{f"ur } n \to \infty. \]
}%

Ein beliebiges $T \in \Adjloru{X}$ besitzt die folgenden Eigenschaften:

\begin{bemerkung}\label{TinAdjloruIstabg}
Sei $X$ ein Operatorraum und
$T \in \Adjloru{X}$.
Dann gilt:
\begin{enumaufz}
\item $T$ ist dicht definiert und abgeschlossen.
\item Nach Definition findet man ein $S : D(S) \subseteq X \to X$ derart,
\dass $\mri \cmpmatrix{0 & T \\ S & 0} \in \Adjlco{X}$ ist.
Dann ist $S$ ebenfalls dicht definiert und abgeschlossen.
\end{enumaufz}
\end{bemerkung}

Diese Proposition\dremark{\ref{TinAdjloruIstabg}} beweist man
mit dem folgenden Lemma,\dremark{und \ref{dichtInC2}}
welches man leicht durch Nachrechnen zeigt:

\begin{lemma}\label{aequivAbg0AB0}
Sei $X$ ein Operatorraum,
seien $A : D(A) \subseteq X \to X$ und $B : D(B) \subseteq X \to X$ linear.
Dann sind die folgenden Aussagen "aquivalent:
\begin{enumaequiv}
\item $A$ und $B$ sind abgeschlossen.
\item $\cmpmatrix{0 & A \\ B & 0} : D(B) \cmtimes D(A) \subseteq C_2(X) \to C_2(X)$
  ist abgeschlossen, wobei $D(B) \cmtimes D(A)$ mit der von $C_2(X)$ induzierten
  Norm versehen wird.
\end{enumaequiv}
\end{lemma}

\dremww{Beweis.
\bewitemphq{(a)$\Rightarrow$(b)}
Setze $\hat{A}  := \cmpmatrix{0 & A \\ B & 0}$.
Sei $\binom{x_n}{y_n}_n \in D(\hat{A})^\mbbN$ konvergent gegen
$\binom{x_0}{y_0} \in C_2(X)$ und
$(\hat{A}\binom{x_n}{y_n})_n$ konvergent gegen $z \in C_2(X)$.
Dann gilt:
\begin{align*}
   \norm{x_n - x_0}
&\leq  \normBig{ \binom{x_n}{y_n} - \binom{x_0}{y_0} }
\to 0 \quad\text{f"ur } n \to \infty \quad\text{und}  \\
   \norm{B x_n - z_2}
&\leq  \normBig{ \binom{A y_n}{B x_n} - z }
=  \normBig{ \hat{A}\binom{x_n}{y_n} - z }
\to 0 \quad\text{f"ur } n \to \infty.
\end{align*}
Es folgt: $x_0 \in D(B)$ und $B x_0 = z_2$.
Analog erh"alt man $y_0 \in D(A)$ und $A y_0 = z_1$,
also $\binom{x_0}{y_0} \in D(\hat{A})$ und $\hat{A}\binom{x_0}{y_0} = z$.

\bewitemphq{(b)$\Rightarrow$(a)}
Sei $\hat{A}$ abgeschlossen.
Sei $x \in D(B)^\mbbN$ konvergent gegen $x_0 \in X$ und
$(B x_n)_{n \in \mbbN}$ konvergent gegen $y_0 \in X$.
Dann gilt:
\begin{align*}
   \normBig{ \binom{x_n}{0} - \binom{x_0}{0} }
&= \norm{x_n - x_0}  \to 0 \quad\text{f"ur } n \to \infty \quad\text{und}  \\
   \normBig{ \hat{A}\binom{x_n}{0} - \binom{0}{y_0} }
&= \normBig{ \binom{0}{B x_n} - \binom{0}{y_0} }
=  \norm{B x_n - y_0}  \to 0 \quad\text{f"ur } n \to \infty.
\end{align*}
Es folgt: $\binom{0}{y_0} \in D(\hat{A})$, also $y_0 \in D(B)$,
und $\binom{0}{y_0} = \hat{A}\binom{x_0}{0} = \binom{0}{B x_0}$,
also $B x_0 = y_0$.
Somit ist $B$ abgeschlossen.
\dremark{17.10.'08/1}%
}%

Die Menge $\Adjlor{X}$ der beschr"ankten,
von links adjungierbaren Multiplikatoren ist in $\Adjloru{X}$ enthalten:

\begin{bemerkung}\label{AdjlorSubseteqAdjloru}
Sei $X$ ein Operatorraum.
Dann gilt: $\Adjlor{X} \subseteq \Adjloru{X}$.
\dremark{Gilt $\Adjlor{X} = \{ A \in \Adjloru{X} \setfdg D(A) = X \}$?}%
\dremark{Die erzeugte Gruppe ist sogar normstetig.}%
\end{bemerkung}

In den Beweis dieser Proposition geht das folgende Lemma ein:

\begin{lemma}\label{NullTTs0InAdjlor}
Sei $X$ ein Operatorraum und $T \in \Adjlor{X}$.
Sei $a := \isoIMl^{-1}(T)$,
also $a \in \mIMls(X)$ und
$j(T x) = a \multis j(x)$ f"ur alle $x \in X$ (siehe Proposition~\ref{AdjlorIsomIMls}).
Dann ist $\hat{T} := \cmpmatrix{ 0 & T \\ T^* & 0 } \in \Adjlor{C_2(X)}$ selbstadjungiert.
Sei
\[ j_{2,1} : C_2(X) \to C_2(I(X)), x \mapsto \binom{j(x_1)}{j(x_2)}. \]
\index[S]{j21@$j_{2,1}$}%
Es gilt: $j_{2,1}(\hat{T} x) = \cmpmatrix{ 0 & a \\ a^* & 0 } \multis j_{2,1}(x)$
f"ur alle $x \in C_2(X)$.
\end{lemma}

\begin{proof}
Nach \eqref{eqTsxEqasx} gilt: $j(T^* y) = a^* \multis j(y)$ f"ur alle $y \in X$.
Also ergibt sich:
\[ j_{2,1}\left(\hat{T}x\right)
=  \cmpmatrix{j(T x_2) \\ j(T^* x_1)}
=  \cmpmatrix{a \multis j(x_2) \\ a^* \multis j(x_1)}
=  \cmpmatrix{0 & a \\ a^* & 0} \multis j_{2,1}(x) \]
f"ur alle $x \in C_2(X)$.
Es ist $\hat{a} := \cmpmatrix{0 & a \\ a^* & 0} \in M_2(I_{11}(X))$
mit $\hat{a}^* = \hat{a}$.
Da $j$ und somit auch $j_{2,1}$ vollst"andige Isometrien sind,
folgt $\hat{T} \in \Adjlor{C_2(X)}$ und $\hat{T}^* = \hat{T}$.
\dremark{$j$ ist vollst"andige Isometrie.
Es gilt
$\norm{j_{2,1}\cmpmatrix{x \\ y}}
=  \norm{j_2\cmpmatrix{x & 0 \\ y & 0}}
=  \norm{\cmpmatrix{x & 0 \\ y & 0}}
=  \norm{\cmpmatrix{x \\ y}}$.
Setze $\rho := j_{2,1}$.
Es gilt:
$\rho(\tilde{T}^* x) = \tilde{a}^* \rho(x) = \tilde{a} \rho(x) = \rho(\tilde{T} x)$
f"ur alle $x \in X$, also $\tilde{T}^* = \tilde{T}$.}%
\end{proof}

\begin{proof}[Beweis von Proposition~\ref{AdjlorSubseteqAdjloru}]
\dmarginpar{pr}
Sei $T \in \Adjlor{X}$.
Nach Lemma~\ref{NullTTs0InAdjlor} ist
$\hat{T} := \cmpmatrix{0 & T \\ T^* & 0} \in \Adjlor{C_2(X)}$ selbstadjungiert.
Es gilt f"ur jeden Operatorraum $Y$ nach Proposition~\ref{AdjlorSubseteqAdjlco}:
\[ \bigl\{ A \in \Adjlor{Y} \setfdg A^* = -A  \bigr\} \subseteq \Adjlco{Y}. \]
Da $\mri \hat{T}$ schiefadjungiert ist, folgt
$\mri \hat{T} \in \Adjlco{C_2(X)}$, also $T \in \Adjloru{X}$.
\dremark{Beweis ohne die Verwendung von \ref{NullTTs0InAdjlor} und \ref{AdjlorSubseteqAdjlco}:

Sei $T \in \Adjlor{X}$.
Man findet ein $a \in \mIMl(X)$ mit $j(Tx) = a \multis j(x)$ f"ur alle $x \in X$.
Es gilt
\[ j_{2,1}\Bigl(\hat{T}\cmpmatrix{x \\ y}\Bigr)
=  \cmpmatrix{j(T y) \\ j(T^* x)}
=  \cmpmatrix{a \multis j(y) \\ a^* \multis j(x)}
=  \cmpmatrix{0 & a \\ a^* & 0} \multis j_{2,1}\Bigl(\cmpmatrix{x \\ y}\Bigr) \]
f"ur alle $x,y \in X$.
Man hat $\hat{a} := \cmpmatrix{0 & a \\ a^* & 0} \in M_2(I_{11}(X))$
mit $\hat{a}^* = \hat{a}$.
Da $j_{2,1}$ eine vollst"andige Isometrie ist, folgt $\hat{T} \in \Adjlor{C_2(X)}$
und $\hat{T}^* = \hat{T}$.
\dremark{$j$ ist vollst"andige Isometrie.
Es gilt
$\norm{j_{2,1}\cmpmatrix{x \\ y}}
=  \norm{j_2\cmpmatrix{x & 0 \\ y & 0}}
=  \norm{\cmpmatrix{x & 0 \\ y & 0}}
=  \norm{\cmpmatrix{x \\ y}}$.
Setze $\rho := j_{2,1}$.
Es gilt:
$\rho(\hat{T}^* x) = \hat{a}^* \rho(x) = \hat{a} \rho(x) = \rho(\hat{T} x)$
f"ur alle $x \in X$, also $\hat{T}^* = \hat{T}$.}%

Weil $\Adjlor{C_2(X)}$ eine unitale \csalgebra ist,
wird durch $U_t := \mre^{t \mri \hat{T}}$ f"ur alle $\indtGr$ eine unit"are Gruppe
in $\Adjlor{C_2(X)}$ definiert.
\dremark{$U_t^* = (\mre^{t \mri \hat{T}})^*
= \sum \frac{1}{n!} (t \mri \hat{T})^*
= \sum \frac{1}{n!} t^n (-\mri)^n \hat{T^*}^n
= \sum \frac{1}{n!} (-t)^n \mri^n \hat{T}^n
= \mre^{-t \mri \hat{T}} = U_{-t}$,
$U_t^* U_t = U_{-t} U_t = U_0 = \Id_X$}%
Aus der Normstetigkeit von $(U_t)_\indtGr$ folgt die strikte Stetigkeit
von $(U_t)_\indtGr$.
\dremark{$U_t^* = U_{-t}$, also $\norm{U_t^* x - x} = \norm{U_{-t} x - x} \leq \norm{U_{-t} - \Id_X} \cdot \norm{x}$}%
}%
\dremark{3.7.'08/3}%
\end{proof}

In unitalen Operatorr"aumen (siehe Definition~\ref{defUnitalerOpraum}) gilt,
analog zu Proposition~\ref{AdjlorsaEqAdjlco}:

\begin{bemerkung}\label{XunitalAdjlorEqAdjloru}
Sei $X$ ein unitaler Operatorraum.
Dann gilt: $\Adjlor{X} = \Adjloru{X}$.
\end{bemerkung}

Zum Beweis notieren wir das folgende

\begin{lemma}\label{NullTS0AdjFolgtTAdj}
Sei $X$ ein unitaler Operatorraum.
Seien $T : D(T) \subseteq X \to X$ und $S : D(S) \subseteq X \to X$ Operatoren
mit $\hat{T} := \cmpmatrix{0 & T \\ S & 0} \in \Adjlor{C_2(X)}$ und $\hat{T}^* = \hat{T}$.
\begin{enumaufz}
\item Es gilt: $S,T \in \Adjlor{X}$ und $S^* = T$.
\item Es existiert ein $a \in \mIMls(X)$ mit der Eigenschaft:
  $j_{2,1}(\hat{T}x) = a \multis j_{2,1}(x)$ f"ur alle $x \in C_2(X)$.
\end{enumaufz}
\end{lemma}

\begin{proof}
Mit $e_X$ sei das Einselement von $X$ bezeichnet.
Weil die Definition eines unitalen Operatorraumes
unabh"angig von der Einbettung in eine spezielle \csalgebra ist
(Anmerkung~\ref{AnmUnitalerORunabh}),
findet man o.B.d.A. einen Hilbertraum $H$ mit $X \subseteq L(H)$ und $e_X = \Id_H$.

Betrachte $\isoIMl^{-1}\restring_{\Adjlor{C_2(X)}} : \Adjlor{C_2(X)} \to \mIMls(C_2(X))$,
und setze ferner $c := \isoIMl^{-1}(\hat{T}) \in \mIMls(C_2(X)) \subseteq I_{11}(C_2(X))$.
Es gilt: $I_{11}(C_2(X)) \cong I_{11}(M_2(X))$ (Proposition~\ref{I11C2XEqI11M2X}).
Weiter sind
\[ p_1 := \Id_H \mathop{\oplus} 0, p_2 := {(\Id_H \oplus \Id_H)} - p_1 \in I_{11}(M_2(X)) \]
Projektionen.
Bez"uglich $p_1$ und $p_2$ kann man $c$ schreiben als
$c = \cmpmatrix{ \alpha & \beta \\ \gamma & \delta }$.
Es gilt f"ur alle $x \in X$:
\begin{align*}
   \binom{j(Tx)}{0}
&= \dremarkm{j_{2,1}\left( \binom{Tx}{0} \right)
=}  j_{2,1}\left(\hat{T}\left(\binom{0}{x}\right)\right)
=  c \multis j_{2,1}\left(\binom{0}{x}\right)
\dremarkm{=  c \multis \binom{0}{j(x)}}
=  \binom{\beta \multis j(x)}{\delta \multis j(x)},
\dremarkm{  \\
   \binom{0}{j(Sx)}
&= j_{2,1}( \binom{0}{Sx} )
=  j( \hat{T}\binom{x}{0} )
=  c \multis j_{2,1}(\binom{x}{0})
=  c \multis \binom{j(x)}{0}
=  \binom{\alpha \multis j(x)}{\gamma \multis j(x)} }%
\end{align*}
also $j(Tx) = \beta \multis j(x)$ und $0 = \delta \multis j(x)$.
Analog ergibt sich:
$j(Sx) = \gamma \multis j(x)$ und $0 = \alpha \multis j(x)$ f"ur alle $x \in X$.
Es ergibt sich mit
\refb{\cite{BlecherLeMerdy04OpAlg}, Proposition 4.4.12,}{Proposition~\ref{Lcinj}}:
$\alpha = 0$ und $\delta = 0$.
\dremark{$\cmpmatrix{ \alpha & 0 \\ 0 & 0 } \multis j_{2,1}(x) = 0$ f"ur alle $x \in C_2(X)$.}%
Da $\mIMls(C_2(X))$ und $\Adjlor{C_2(X)}$ nach
\refb{\cite{Zarikian01Thesis}, Proposition 1.7.4,}{Proposition~\ref{AdjlorIsomIMls}}
unter $\isoIMl$ als unitale \csalgebren isomorph sind,
folgt:
\[ c^* = \isoIMl^{-1}(\hat{T}^*) = \isoIMl^{-1}(\hat{T}) = c. \]
Man hat:
$\beta^*
= (p_1 \multis c \multis p_2)^*
\dremarkm{= p_2 \multis c^* \multis p_1}
= p_2 \multis c \multis p_1
= \gamma$.
Somit gilt $\beta^* \multis j(x) = \gamma \multis j(x) \in j(X)$ f"ur alle $x \in X$,
also $T \in \Adjlor{X}$.
Weiter erh"alt man: $S \in \Adjlor{X}$ und
$  S^*
= \isoIMl(\gamma)^*
\dremarkm{= \isoIMl(\gamma^*)}
= \isoIMl(\beta)
= T$.
\dremark{$S^* = \isoIMl^{-1}(\beta)^* = \isoIMl^{-1}(\beta^*) = T$}%
\dremark{27.10.'08/1}%
\end{proof}

\begin{proof}[Beweis von Proposition~\ref{XunitalAdjlorEqAdjloru}]
\bewitemphqq{$\subseteq$} folgt mit Proposition~\ref{AdjlorSubseteqAdjloru}.

\bewitemphq{$\supseteq$} Sei $T \in \Adjloru{X}$.
Dann findet man ein $S : D(S) \subseteq X \to X$ derart,
\dass gilt: $\mri \cmpmatrix{0 & T \\ S & 0} \in \Adjlco{C_2(X)}$.
Setze $\hat{T} := \cmpmatrix{0 & T \\ S & 0}$.
F"ur jeden Operatorraum $Y$ hat man nach Proposition~\ref{AdjlorsaEqAdjlco}:
\[ \bigl\{ A \in \Adjlor{Y} \setfdg A^* = -A \bigr\} = \Adjlco{Y}. \]
Somit ergibt sich $(\mri \hat{T})^* = -\mri \hat{T} \in \Adjlor{C_2(X)}$,
also $\hat{T}^* = \hat{T}$.
Mit Lemma~\ref{NullTS0AdjFolgtTAdj} erh"alt man: $T \in \Adjlor{X}$.
\dremark{25.9.'08/3, 27.10.'08/1}%
\end{proof}

Der folgende wichtige Satz zeigt, \dass auf einem Hilbert-\cstern{}Modul
die unbeschr"ankten Multiplikatoren mit den regul"aren Operatoren
"ubereinstimmen.
Damit sind jene eine Verallgemeinerung der regul"aren Operatoren auf Operatorr"aume.

\begin{satz}\label{zshgRegOpUnbeschrMult}
Sei $E$ ein Hilbert-\cstern{}Modul.
\begin{enumaufz}
\item Es gilt: $\Adjloru{E} = \Multwor{E}$.
\item Sei $T \in \Adjloru{E}$.
Nach Definition findet man ein $S : D(S) \subseteq E \to E$ derart,
\dass $\mri \cmpmatrix{0 & T \\ S & 0} \in \Adjlco{C_2(E)}$ ist.
Dann gilt: $T^* = S \in \Multwor{E}$.
\dremark{Beweis "uber Multiplikation mit Matrizen m"oglich?
  Benutze \cmzB St"orungssatz von Damaville. (WW)}%
\end{enumaufz}
\end{satz}

\begin{proof}
\dremark{\bewitemphq{(i).$\subseteq$ und (ii)}}%
Wir zeigen zun"achst (ii) und \glqq$\subseteq$\grqq{} von (i).
Sei $T \in \Adjloru{E}$ und $\hat{T} := \cmpmatrix{0 & T \\ S & 0}$,
also $\mri \hat{T} \in \Adjlco{C_2(E)}$.
Es gilt: $C_2(E) \cong E \oplus E$ (Beispiel~\ref{HCsModAdjMultEqAdj}).
Mit Satz~\ref{MultworInAdjlco} erh"alt man:
\begin{equation}\label{eqAdjloruEEqMultwor}
   \Adjlco{C_2(E)}
\cong \Adjlco{E \oplus E}
=  \left\{ R \in \Multwor{E \oplus E} \setfdg R^* = -R \right\}.
\end{equation}
Somit ist $\mri \hat{T} \in \Multwor{E \oplus E}$ schiefadjungiert.\dremark{also $\hat{T}$ s. a.}
\dremark{Beweisidee: $\hat{T} \in \Multwor{E \oplus E}$, $T$ ist als Ecke auf $\Multwor{E}$.}%

Nach Proposition~\ref{TinAdjloruIstabg} sind $S$ und $T$ abgeschlossen und dicht definiert.
Nach Lemma~\ref{NullTS0sternEq} gilt:
$\cmpmatrix{0 & T \\ S & 0}^*  = \cmpmatrix{0 & S^* \\ T^* & 0}$.
Wegen $\hat{T}^* = \hat{T}$ folgt $S = T^*$,
\dremark{$\cmpmatrix{0 & T \\ S & 0} = \cmpmatrix{0 & T \\ S & 0}^*
  = \cmpmatrix{0 & S^* \\ T^* & 0}$}%
insbesondere ist $T^*$ dicht definiert.

Zu zeigen bleibt: $1+T^*T$ hat dichtes Bild.
Es gilt:
\begin{align*}
   (1 + \hat{T}^*\hat{T})(x,y)
\dremarkm{=  (x,y) + \hat{T}(Ty,T^*x)}
=  (x,y) + (T T^* x,T^* Ty)
=  \bigl((1+TT^*)x, (1+T^*T)y\bigr)
\end{align*}
f"ur alle $(x,y) \in D(1+\hat{T}^*\hat{T})$.
Weil $\norm{y} \leq \norm{(x,y)}$ nach \eqref{eqNormHCsModEPlusE} f"ur alle $x,y \in E$ gilt
und weil $1 + \hat{T}^*\hat{T}$ dichtes Bild in $E \oplus E$ hat,
folgt\dremark{mit \ref{dichtInHSModEF}}: $1+T^*T$ hat dichtes Bild in $E$.
Somit ist $T$ regul"ar.

Mit \cite{LanceHmod}, Corollary 9.6, folgt: $T^* \in \Multwor{E}$.
\smallskip

\dremark{\bewitemphq{(i).$\supseteq$}}%
Es bleibt \glqq$\supseteq$\grqq{} von (i) zu zeigen.
Sei dazu $T \in \Multwor{E}$.
Nach Proposition~\ref{NullTTs0Selbstadj} ist $h := \cmpmatrix{0 & T \\ T^* & 0}$ regul"ar und selbstadjungiert.
\dremark{Weiter gilt: $E \oplus E \cong C_2(E)$ (Beispiel~\ref{HCsModAdjMultEqAdj}).
Nach Satz~\ref{MultworInAdjlco} hat man:
\[ \left\{ R \in \Multwor{E \oplus E} \setfdg R^* = -R \right\} \subseteq \Adjlco{E \oplus E}. \]}%
Mit \eqref{eqAdjloruEEqMultwor} ergibt sich $\mri h \in \Adjlco{C_2(E)}$, also $T \in \Adjloru{E}$.
\dremark{Alternativer Beweis ohne Verwendung von \ref{AdjlorsaEqAdjlco}:
Mit \ref{erzeugerUnitaereGrWor} folgt,
\dass $\mri h$ eine unit"are
$C_0$-Gruppe $(U_t)_\indtGr$ auf $E \oplus E$ mit
$U_t \in \Adjhm{E \oplus E} \overset{\eqref{eqAdjEplusEeqAdjlor}}{\cong} \Adjlor{C_2(E)}$
f"ur alle $\indtGr$ erzeugt.}%
\dremark{3.7.'08/1}%
\end{proof}

Mit dem obigen Satz folgt,
\dass im allgemeinen nicht $\Adjlor{X} = \Adjloru{X}$ gilt:

\begin{beispiel}
Sei $\Omega$ ein lokalkompakter Hausdorffraum, der nicht kompakt ist.
Dann gilt:
\[ \Adjlor{C_0(\Omega)} \cong C_b(\Omega) \neq C(\Omega) \cong \Adjloru{C_0(\Omega)}. \]
\end{beispiel}

\begin{proof}
Da $X := C_0(\Omega)$ eine kommutative \csalgebra ist, ergibt sich mit
\refb{\cite{BlecherLeMerdy04OpAlg}, Corollary 8.4.2,}{Beispiel~\ref{bspMlCsAlgEqMlOR}}:
$\Adjlor{X} \cong \Multcs{X} \cong C_b(\Omega)$.
Nach \refb{\cite{Woronowicz91}, Example~2,}{Beispiel~\ref{bspRegOp}} gilt: $\Multwor{X} = C(\Omega)$.
Man erh"alt mit Satz~\ref{zshgRegOpUnbeschrMult}:
$\Adjloru{X} = \Multwor{X} \cong C(\Omega)$.
\end{proof}

\section{Zusammenhang zwischen $C_0$-Halbgruppen auf $X$ und auf $C_2(X)$}

Wir untersuchen in diesem Abschnitt, unter welchen Voraussetzungen die Operatoren $A$
und $\cmpmatrix{A & 0 \\ 0 & 0}$ auf einem Operatorraum $X$
(bzw. auf dem Spaltenoperatorraum $C_2(X)$)
eine $C_0$-Halbgruppe erzeugen und
welche Beziehung zwischen diesen beiden Operatoren besteht.
Als wichtiges Resultat werden wir eine intrinsische Charakterisierung der
Elemente aus $\Adjlco{X}$ und $\Adjloru{X}$ beweisen.
\skiptext

Zwischen Halbgruppen auf $X$ und $C_2(X)$ besteht der folgende Zusammenhang:

\begin{bemerkung}\label{halbgrInXundC2X}
Sei $X$ ein Operatorraum und $(T_t)_\indtHG$ eine Familie in $\mLinStet(X)$.
F"ur alle $t \in \indHG$ setze
$S_t := \matoplusc{T_t}{\Id_X} := \cmpmatrix{T_t & 0 \\ 0 & \Id_X}$.
Es gilt:
\index[S]{AoplusB@$A \matoplus B$}%
\begin{enumaufz}
\item $(T_t)_\indtHG$ ist genau dann eine Operatorhalbgruppe (bzw. $C_0$"=Halbgruppe) auf $X$,
wenn $(S_t)_\indtHG$ eine Operatorhalbgruppe (bzw. $C_0$"=Halbgruppe) auf $C_2(X)$ ist.

\item Sei $(T_t)_\indtHG$ eine $C_0$-Halbgruppe
  und $A : D(A) \subseteq X \to X$.
  Dann ist $A$ genau dann Erzeuger von $(T_t)_\indtHG$,
  wenn $\matnull{A} := \cmpmatrix{A & 0 \\ 0 & 0} : D(A) \cmtimes X \subseteq C_2(X) \to C_2(X)$
  Erzeuger von $(S_t)_\indtHG$ ist.
\index[S]{etaA@$\matnull{A}$}%
\dremark{Gilt analog auch f"ur Operatorgruppen.}%
\dremark{Aussagen "uber das Wachstumsverhalten?}%
\end{enumaufz}
\end{bemerkung}

\begin{proof}
\bewitemph{(i):}\dremark{\bewitemph{(I):}}
Offensichtlich ist $(T_t)_\indtHG$ genau dann eine Operatorhalbgruppe,
wenn $(S_t)_\indtHG$ eine Operatorhalbgruppe ist.

\dremark{\bewitemph{(II):}}%
Somit bleibt zu zeigen, \dass $(T_t)_\indtHG$ genau dann stark stetig ist,
wenn $(S_t)_\indtHG$ stark stetig ist.

\bewitemph{\glqq$\Rightarrow$\grqq:}
Sei $(T_t)_\indtHG$ eine $C_0$-Halbgruppe auf $X$.
Es gilt f"ur alle $\binom{x}{y} \in C_2(X)$:
\[ \normlr{ S_t \binom{x}{y} - \binom{x}{y} }
\dremarkm{=  \normlr{ \binom{T_t x - x}{0} }}
=  \norm{ T_t x - x } \longrightarrow 0 \qquad\text{f"ur } t \to 0. \]

\bewitemph{\glqq$\Leftarrow$\grqq} folgt analog zu \grqq$\Rightarrow$\grqq.
\dremark{folgt aus $\norm{T_t x - x} = \normlr{S_t \binom{x}{0} - \binom{x}{0}}$
f"ur alle $\indtHG$ und $x \in X$.}%

\bewitemph{(ii):} Wegen
\begin{align*}
   \normlr{ \frac{S_t\binom{x}{y} - \binom{x}{y}}{t} - \matnull{A}\left( \binom{x}{y} \right) }
\dremarkm{&= \normlr{ \frac{\binom{T_t x-x}{0}}{t} - \binom{Ax}{0} }  \\}
&= \normlr{ \binom{\frac{1}{t} ( T_t x-x ) - Ax}{0} }
=  \normlr{ \frac{T_t x-x}{t} - Ax }
\end{align*}
f"ur alle $x \in D(A)$, $y \in X$ und $t \in \mbbR_{>0}$ folgt die Behauptung.
\dremark{
Sei $x \in D(A)$, $y \in X$.
Es gilt
$ 0
= \lim \norm{\frac{T_t x - x}{t} - Ax}
= \lim \norm{ \frac{S_t \binom{x}{y} - \binom{x}{y}}{t} - \matnull{A}\binom{x}{y} }$,
also $\binom{x}{y} \in D(\matnull{A})$.
}%
\dremark{Bew.: Siehe Z. 11.4.'07}%
\end{proof}

\dremww{
\begin{lemma}\label{nullNotinResA0}
Sei $X$ ein Operatorraum, $A : D(A) \subseteq X \to X$ linear und abgeschlossen.
Dann gilt: $0 \notin \rho( \matnull{A} )$.
\end{lemma}

\begin{proof}
F"ur alle $\lambda \in \mbbC$ gilt
$  A_\lambda := \lambda \Id_{C_2(X)} - \matnull{A}
=  \cmpmatrix{\lambda \Id_X - A & 0 \\ 0 & \lambda \Id_X }$,
also ist $A_0$ nicht invertierbar.
Somit folgt:\dremark{Z. 3.8.'07/1}
\[ 0 \notin \rho( \matnull{A} )
=  \bigl\{ \lambda \in \mbbC \setfdg \lambda \Id_{C_2(X)} - \matnull{A} :
   D(\matnull{A}) \arrowbij C_2(X) \text{ ist bijektiv auf } C_2(X) \bigr\}.  \qedhere \]
\end{proof}
}%

Die Operatoren $A$ und $\eta_A = \begin{smallpmatrix} A & 0 \\ 0 & 0 \end{smallpmatrix}$
h"angen wie folgt zusammen:

\begin{lemma}\label{zshgAundA0}
Sei $X$ ein Operatorraum und $A : D(A) \subseteq X \to X$ linear.
\begin{enumaufz}
\item $D(A)$ ist genau dann dicht in $X$, wenn $D(A) \cmtimes X$ dicht in $C_2(X)$ ist,
  wobei $D(A) \cmtimes X$ mit der von $C_2(X)$ induzierten Norm versehen wird.
\item $A$ ist genau dann abgeschlossen, wenn $\matnull{A}$ abgeschlossen ist.
\item Sei $A$ abgeschlossen.
Dann gilt: $\rho(A) \setminus \{0\} = \rho( \matnull{A} )$.
Insbesondere hat man: $\mbbR_{>0} \subseteq \rho(A) \iff \mbbR_{>0} \subseteq \rho( \matnull{A} )$.
\item Sei $A$ abgeschlossen, $\lambda \in \rho(A) \setminus \{0\}$ und $n \in \mbbN$.
Dann gilt f"ur die Resolvente:
\begin{enumaufzB}
\item $\norm{R(\lambda,A)_n} \leq \norm{R(\lambda,\matnull{A})_n}$,
\item $\normlr{\lambda R(\lambda,\matnull{A})_n\binom{x}{y}} \leq \norm{\lambda R(\lambda,A)_n x} + \norm{y}$
f"ur alle $x,y \in M_n(X)$,
insbesondere $\norm{\lambda R(\lambda,\matnull{A})_n} \leq \norm{\lambda R(\lambda,A)_n} + 1$.
\end{enumaufzB}
\end{enumaufz}
\end{lemma}

\dfrage{Gibt es f"ur (iv).(II) eine bessere Absch"atzung? Gegenbeispiel?}%
\dremark{Man kann die Absch"atzung in (iv).(II) nicht wesentlich verbessern.
  Sei $X$ ein Banachraum, $A := \Id_X$.
  Es gilt:
$\normlr{ \lambda R(\lambda,\matnull{A}) \binom{0}{y} }
=  \normlr{ \lambda \binom{R(\lambda,A)0}{(\lambda \Id_X)^{-1}y} }
= \norm{y}$ und
\[ \norm{ \lambda R(\lambda,A) x}
=  \normlr{ \lambda \left( (\lambda - 1) \Id_X \right)^{-1} x }
=  \normlr{ \frac{\lambda}{\lambda-1} x }
=  \normlr{ \left(1 + \frac{1}{\lambda-1}\right) x }
\overset{\lambda \to 0}{\longrightarrow} 0. \]
}%

\begin{proof}
\bewitemph{(i)} und \bewitemph{(ii)} erh"alt man unmittelbar.

\bewitemph{(iii):} Es gilt:\dremark{\ref{nullNotinResA0}} $0 \notin \rho(\matnull{A})$.
F"ur alle $\lambda \in \mbbC \setminus \{0\}$ ist
$\lambda \Id_X - A$ genau dann bijektiv auf $X$,
wenn $\lambda \Id_{C_2(X)} - \matnull{A} = \cmpmatrix{\lambda \Id_X - A & 0 \\ 0 & \lambda \Id_X}$
bijektiv auf $C_2(X)$ ist.
Somit folgt (iii).

\bewitemph{(iv):}
\bewitemph{(I):}
Es gilt f"ur alle $x \in X$:
\[ \norm{R(\lambda,A)x}
\dremarkm{=  \normlr{ \binom{R(\lambda,A)x}{0} }}
=  \normlr{ R(\lambda,\matnull{A}) \binom{x}{0} }
\leq \normlr{ R(\lambda,\matnull{A}) } \cdot \norm{x}. \]

\bewitemph{(II):}
Es gilt f"ur alle $x,y \in X$:
\dremark{Z. 3.8.'07/1}%
\begin{align*}
   \normlr{ \lambda R\bigl(\lambda,\matnull{A}\bigr)\binom{x}{y} }
&=  \normlr{ \lambda \binom{R(\lambda,A)x}{(\lambda \Id_X)^{-1}y} }  \\
&\leq \norm{\lambda R(\lambda,A)x} + \norm{y}
\leq \bigl(\norm{\lambda R(\lambda,A)} + 1\bigr) \normlr{ \binom{x}{y} }.
\end{align*}
\end{proof}

F"ur $C_0$-Halbgruppen auf Operatorr"aumen definieren wir:

\begin{definitn}
Sei $X$ ein Operatorraum. 
Eine $C_0$-Halbgruppe $(T_t)_\indtHG$ (bzw. $C_0$-Gruppe $(T_t)_\indtGr$) auf $X$ hei"st
\defemph{vollständig kontraktiv},
  falls $T_t$ f"ur alle $\indtHG$ (bzw. $\indtGr$)
  vollst"andig kontraktiv ist.
\index[B]{vollständig kontraktive $C_0$-Halb\-grup\-pe}%
\end{definitn}

\dremark{
\begin{bemerkung}\label{A0ErzHGFolgtAErz}
Sei $X$ ein Operatorraum, $M \in \mbbR_{\geq 1}$ und $\omega \in \mbbR$.
Sei $A : D(A) \subseteq X \to X$ so,
\dass $\matnull{A} : D(\matnull{A}) \subseteq C_2(X) \to C_2(X)$ Erzeuger einer
$C_0$-Halbgruppe $(S_t)_\indtHG$ ist mit
$\normcb{S_t} \leq M \mre^{\omega t}$ f"ur alle $\indtHG$.
Dann ist $A$ Erzeuger einer $C_0$-Halbgruppe $(T_t)_\indtHG$ mit
$\normcb{T_t} \leq M \mre^{\omega t}$ f"ur alle $\indtHG$.
\end{bemerkung}

\dfrage{Sei $A$ Erzeuger einer (kontraktiven) $C_0$-Halbgruppe auf $X$.
  Unter welchen Voraussetzungen ist $\eta_A$ Erzeuger
  einer (kontraktiven) $C_0$-Halbgruppe auf $C_2(X)$?}%

\begin{proof}
Nach \ref{halbgrInXundC2X} ist $A$ Erzeuger einer $C_0$-Halbgruppe
$(T_t)_\indtHG$ auf $X$ und $\eta_A$ Erzeuger einer $C_0$-Halbgruppe
auf $C_2(X)$ der Gestalt $(S_t)_\indtHG$, wobei $S_t := \matoplusc{T_t}{\Id_X}$
f"ur alle $\indtHG$.
Es gilt f"ur alle $\indtHG$, $x \in X$
\[ \norm{T_t x}
\leq  \normlr{ \cmpmatrix{T_t & 0 \\ 0 & \Id_X}  \cmpmatrix{x \\ 0} }
\leq  \norm{S_t} \,\norm{x}
\leq  M \mre^{\omega t} \norm{x}, \]
also folgt: $\normcb{T_t} \leq M \mre^{\omega t}$.
\dremark{Z. 30.7.'07/1}%
\end{proof}
}%

Es sei daran erinnert, \dass auf einem Operatorraum $X$ mit $\Multlor{X}$
die Menge der Linksmultiplikatoren bezeichnet wird,
die mit der Multiplikatornorm $\normMl{\cdot}{X}$ versehen wird
(siehe Definition~\ref{defLinksmultOR}).
\medskip

Der folgende Satz stellt mit Hilfe von $\matnull{A}$ eine intrinsische Charakterisierung der
Elemente von $\Adjlco{X}$ bereit, die insbesondere nicht auf $\Adjlor{X}$ zur"uckgreift:

\begin{satz}\label{A0ErzvkC0GrFolgt}\label{charMultlcomatnull}\label{charAdjlcomatnull}
Sei $X$ ein Operatorraum und $A : D(A) \subseteq X \to X$.
\begin{enumaufz}
\item Die folgenden Aussagen sind "aquivalent:
\begin{enumaequiv}
\item $A$ erzeugt eine $C_0$-Halbgruppe $(T_t)_\indtHG$ auf $X$ mit
  $T_t$ aus der abgeschlossenen Kugel $\overline{\cmkug}_\Multlor{X}(0,1)$
f"ur alle $\indtHG$.
\item $\matnull{A}$ erzeugt eine \cmvk{e} $C_0$-Halbgruppe auf $C_2(X)$.
\dremark{Gilt etwas "Ahnliches zu (b)$\Rightarrow$(a)
  f"ur vollst"andig quasikontraktive \cohgn?}%
\end{enumaequiv}

\item Die folgenden Aussagen sind "aquivalent:
\begin{enumaequiv}
\item $A \in \Adjlco{X}$.
\item $\matnull{A}$ erzeugt eine \cmvk{e} $C_0$-Gruppe auf $C_2(X)$.
\end{enumaequiv}
\end{enumaufz}
\end{satz}

Zum Beweis verwenden wir das folgende Lemma, welches man leicht nachrechnet:

\begin{lemma}\label{vkC0GrIstVollstIsometr}
Sei $X$ ein Operatorraum und $(T_t)_\indtGr$ eine \cmvk{e} $C_0$-Gruppe auf $X$.
Dann ist $T_t$ eine vollst"andig isometrische Bijektion auf $X$
f"ur alle $\indtGr$.
\end{lemma}

\dremww{Beweis.
Sei $\indtGr$.
Dann gilt $\normcb{T_t} \leq 1$ und
$\normcb{T_t^{-1}} = \normcb{T_{-t}} \leq 1$.
Somit \dremark{\ref{normTTinvLeq1Isometrie}} ist $T_t$ eine vollst"andige Isometrie,
die wegen $T_t^{-1} = T_{-t}$ bijektiv auf $X$ ist.
}%

\begin{proof}[Beweis von Satz~\ref{A0ErzvkC0GrFolgt}]
\bewitemph{(i):}
\bewitemphq{(a)$\Rightarrow$(b)}
Es gelte (a).
Weiter ist $T_t \matoplus \Id_X$ nach Satz~\ref{charMultOR}
(Charakterisierung der Linksmultiplikatoren) \cmvk{} f"ur alle $\indtHG$.
Nach Proposition~\ref{halbgrInXundC2X} ist $(\matoplusc{T_t}{\Id_X})_\indtHG$
eine $C_0$"=Halbgruppe mit Erzeuger $\eta_A$.

\bewitemphq{(b)$\Rightarrow$(a)}
Sei $(S_t)_\indtHG$ die von $\matnull{A}$ erzeugte \cmvk{e} $C_0$"=Halbgruppe.
Mit Proposition~\ref{halbgrInXundC2X} erh"alt man, \dass $A$ eine $C_0$-Halbgruppe $(T_t)_\indtHG$
erzeugt und \dass gilt:
$S_t = \matoplusc{T_t}{\Id_X}$ f"ur alle $\indtHG$.
Nach Satz~\ref{charMultOR} ist $T_t \in \overline{\cmkug}_\Multlor{X}(0,1)$ f"ur alle $\indtHG$.

\bewitemph{(ii)} folgt analog unter Verwendung von Proposition~\ref{charAlXZar}
(Charakterisierung der unit"aren Elementen von $\Adjlor{C_2(X)}$)
und Lemma~\ref{vkC0GrIstVollstIsometr}.
\dremark{25.9.'08/2, 10.7.'08/1}%
\end{proof}

Ist $(T_t \matoplus \Id_X)_\indtHG$ auf $C_2(X)$ eine \cmvk{e} $C_0$-Halbgruppe,
so ist $(T_t)_\indtHG$ insbesondere vollst"andig kontraktiv.
Die Umkehrung gilt \iallg jedoch nicht:\dmarginpar{zutun}\dremark{Was gilt f"ur $H^{oh}$, $H^c$?}

\begin{beispiel}\label{bspNichtKontrInHr}
Setze $H := L^2(\mbbR)$. 
Sei $T_t : H^r \to H^r, f \mapsto \mre^{t \mri \Id_\mbbR} f$, f"ur alle $\indtHG$,
wobei mit $H^r$ der Zeilen-Hilbertoperatorraum bezeichnet wird (siehe Proposition~\ref{defHr}).
Dann gilt:
\begin{enumaufz}
\item $(T_t)_\indtHG$ ist eine $C_0$-Halbgruppe auf $H^r$
  mit der Eigenschaft: $\normcb{T_t} = 1$ f"ur alle $\indtHG$.
\item Sei $\lambda \in \mbbR_{>0}$.
  Es ist $(T_t \matoplus \Id_{H^r})_\indtHG$
  eine $C_0$-Halbgruppe auf $C_2(H^r)$ mit
  $\norm{\lambda T_1 \matoplus \Id_{H^r}} \geq \sqrt{1+\lambda^2}$.
  Insbesondere sieht man f"ur $\lambda = 1$, \dass diese Halbgruppe nicht kontraktiv ist.
\end{enumaufz}
\end{beispiel}

\begin{proof}
\bewitemph{(i):}
Offensichtlich ist $(T_t)_\indtHG$ eine Operatorhalbgruppe,
die nach \cite{EngelNagelSemigroups}, Proposition I.4.11, stark stetig ist.

Sei $\indtHG$.
F"ur alle $f \in X := H^r$ gilt:
\begin{equation}\label{eqBspNormHrTt}
   \norm{T_t f}^2
=  \skalprb{\mre^{t \mri \Id_\mbbR} f}
\dremarkm{=  \int_\mbbR \overline{\mre^{t\mri s}} \mre^{t\mri s} \overline{f(s)} f(s) \,ds
=  \int_\mbbR \mre^{-t\mri s} \mre^{t\mri s} \overline{f(s)} f(s) \,ds}
=  \norm{f}^2.
\end{equation}
Mit \cite{Paulsen02ComplBoundedMaps}, S. 200,
\dremark{Da $X$ nach \refb{\cite{Paulsen02ComplBoundedMaps}, S. 200,}{Beispiel~\ref{HrHomogen}} homogen ist,}%
folgt: $\normcb{T_t} = \norm{T_t} = 1$.

\bewitemph{(ii):}
Nach Proposition~\ref{halbgrInXundC2X} ist $(T_t \matoplus \Id_{X})_\indtHG$ eine $C_0$-Halbgruppe
auf $C_2(X)$.

F"ur alle $f,g \in X$ folgt mit
\refb{\cite{EffrosRuan00OperatorSpaces}, S. 55,}{Proposition~\ref{defHr}}:
\begin{equation}\label{eqBspNormHr}
   \normlr{ \cmpmatrix{f \\ g} }^2_{C_2(X)}
\dremarkm{=  \normlr{ \cmpmatrix{f & 0 \\ g & 0} }^2_{M_2(X)}}
\overset{\dremarkm{\ref{defHr}}}{=}
   \normlr{ \cmpmatrix{ \skalpr{f}{f} & \skalpr{f}{g} \\ \skalpr{g}{f} & \skalpr{g}{g} } }_{M_2}.
\end{equation}

Sei $f : \mbbR \to \mbbC, s \mapsto \mre^{-\mri s} \cdot 1_{[0,2\pi]}$ und
$g := 1_{[0,2\pi]}$, wobei $1_{[0,2\pi]}$ die charakteristische Funktion auf $[0,2\pi]$ bezeichne.
Es gilt: $\norm{f}^2 = 2\pi = \norm{g}^2$,\dremark{$\mre^{\mri s}$ ist $2\pi$-periodisch}
$\skalpr{f}{g}
\dremarkm{=  \int_\mbbR \overline{f(s)} g(s) \,ds
=  \int_0^{2\pi} \mre^{\mri s}}%
=  0$ und
$  \skalpr{g}{T_1 f}
=  \int_0^{2\pi} \mre^{\mri s} \mre^{-\mri s} \, ds
=  2\pi$.
Mit \eqref{eqBspNormHr} und \eqref{eqBspNormHrTt} ergibt sich:
\dremark{$(*)$: Eigenwerte: $0$, $2\pi(1+\lambda^2)$}%
\begin{align*}
   \normlr{ \cmpmatrix{f \\ g} }^2_{C_2(X)}
&\overset{\dremarkm{\eqref{eqBspNormHr}}}{=}
   \normlr{ \cmpmatrix{2\pi & 0 \\ 0 & 2\pi} }_{M_2}
=  2\pi \quad\text{und}  \\
   \normlr{ \cmpmatrix{\lambda T_1 f \\ g} }^2_{C_2(X)}
&\overset{\dremarkm{\eqref{eqBspNormHr}}}{=}
   \normlr{ \cmpmatrix{ \lambda^2 \skalpr{T_1 f}{T_1 f} & \lambda\skalpr{T_1 f}{g} \\
                        \lambda\skalpr{g}{T_1 f} & \skalpr{g}{g} } }_{M_2}  \\
&\overset{\dremarkm{\eqref{eqBspNormHrTt}}}{=}
   \normlr{ \cmpmatrix{ \lambda^2 2\pi & \lambda 2\pi \\ \lambda 2\pi & 2\pi } }_{M_2}
\overset{\dremarkm{(*)}}{=}  2\pi(1+\lambda^2).
\end{align*}
Man erh"alt: $\norm{\lambda T_1 \matoplus \Id_{X}} \geq \sqrt{1+\lambda^2}$.
\end{proof}

Mit Satz~\ref{A0ErzvkC0GrFolgt} erh"alt man eine intrinsische Charakterisierung
der Elemente von $\Adjloru{X}$:

\begin{satz}\label{charAdjloruIntr}
Sei $X$ ein Operatorraum und $A : D(A) \subseteq X \to X$.
Dann sind die folgenden Aussagen "aquivalent:
\begin{enumaequiv}
\item $A \in \Adjloru{X}$.
\item Es existiert eine Abbildung $B : D(B) \subseteq X \to X$ mit der Eigenschaft, \dass
$\mri \cmpmatrix{ 0 & A \\ B & 0 } \matoplus \cmpmatrix{ 0 & 0 \\ 0 & 0 }$
Erzeuger einer \cmvk{en} $C_0$-Gruppe auf $C_2(C_2(X)) \cong C_4(X)$ ist.
\dremark{Die Abbildung ist eigentlich nicht definiert.}%
\end{enumaequiv}
\end{satz}

\dremark{
Beweis.
Es gilt: $C_4(X) \cong C_2(C_2(X))$.

\bewitemphq{(a)$\Rightarrow$(b)}
Es existiert eine Abbildung $B : D(B) \subseteq X \to X$ derart,
\dass $\mri \cmpmatrix{ 0 & A \\ B & 0 }$ eine unit"are $C_0$-Gruppe in $\Adjlor{C_2(X)}$
erzeugt.
Mit \ref{A0ErzvkC0GrFolgt} folgt die Behauptung.

\bewitemphq{(b)$\Rightarrow$(a)}
Nach \ref{A0ErzvkC0GrFolgt} erzeugt $\cmpmatrix{0 & \mri A \\ \mri B & 0}$
eine unit"are $C_0$-Gruppe in $\Adjlor{C_2(X)}$,
also gilt: $A \in \Adjloru{X}$.
}%

\section{Die S"atze von Hille-Yosida und von Lu\-mer-Phil\-lips in Operatorr"aumen}

In Satz~\ref{A0ErzvkC0GrFolgt} und Satz~\ref{charAdjloruIntr} werden unbeschr"ankte
Multiplikatoren mit Hilfe von Erzeugern \cmvk{er} $C_0$-Halbgruppen bzw. $C_0$-Gruppen
charakterisiert.
Daher ist von Interesse, wann ein Operator Erzeuger einer
\cmvk{en} $C_0$-Halbgruppe (bzw. $C_0$-Gruppe) ist.

Charakterisierungen hierf"ur werden in diesem Abschnitt bereitgestellt.
So werden Analoga des Satzes von Hille-Yosida und
des Satzes von Lumer-Phillips in Operatorr"aumen gezeigt.
\dremark{Wof"ur? Erzeuger charakterisieren, insbesondere von \cmvk{en} $C_0$-Halbgruppen
  und $C_0$-Gruppen.
  Satz von Hille-Yosida geht in den Beweis von \ref{charMlC0X} ein.}%
\skiptext

Eine entsprechende Version von Satz~\ref{unglFuerResHochn} gilt in Operatorr"aumen:

\begin{satz}\label{eigErzvkHG}\label{eigErzHGOR}
Sei $X$ ein Operatorraum, $\omega \in \mbbR$ und $M \in \mbbR_{\geq 1}$.
Sei $A$ der Erzeuger einer $C_0$-Halb\-grup\-pe $(T_t)_\indtHG$ auf $X$ mit der Eigenschaft:
\[ \normcb{T_t} \leq M \mre^{\omega t} \qquad\text{f"ur alle } \indtHG. \]
Dann gilt f"ur alle $\lambda \in \mbbC$ mit $\Re \lambda > \omega$ und alle $k, n \in \mbbN$:
\begin{enumaufz}
\item $R(\lambda,A_n)^k x = \frac{1}{(k-1)!} \int_0^\infty s^{k-1} \mre^{-\lambda s} T_s^{(n)} x \,ds$
  f"ur alle $x \in M_n(X)$,
\item $\normlr{R(\lambda,A)^k}_\text{cb} \leq \frac{M}{(\Re(\lambda)  - \omega)^k}$.
\end{enumaufz}
\end{satz}

Zum Beweis formulieren wir die folgenden beiden Lemmata:

\begin{lemma}\label{AnErzeugtTtn}
Sei $X$ ein Operatorraum und $(T_t)_\indtHG$ eine $C_0$-Halbgruppe auf $X$ mit Erzeuger~$A$.
Dann ist $(T_t^{(n)})_\indtHG$ eine $C_0$-Halbgruppe auf $M_n(X)$ mit Erzeuger $A_n$
f"ur alle $n \in \mbbN$.\dremark{Es gilt: $D(A_n) = M_n(A)$.}%
\end{lemma}

\begin{proof}
Offensichtlich ist $(T_t^{(n)})_\indtHG$ eine Operatorhalbgruppe.
Die $C_0$-Ei\-gen\-schaft von $(T_t^{(n)})_\indtHG$ und die Eigen\-schaft, \dass $A_n$ Erzeuger ist,
erh"alt man mit der Ungleichung
$\norm{x}_n \leq \sum_{i,j=1}^n \norm{x_{ij}}$ f"ur alle $x \in M_n(X)$
(Lemma~\ref{normaLeqNormSumAij}).
\dremark{bzw. \ref{konvergenzInMnX}}%
\dremark{Ausf"uhrlicher Beweis: 27.4.'07, S. 1}
\end{proof}

Durch einfache Umformungen erh"alt man:

\begin{lemma}\label{RlAnEqRlAn}
Sei $X$ ein Operatorraum und $A : D(A) \subseteq X \to X$ linear.
Dann gilt f"ur alle $\lambda \in \rho(A)$ und $n \in \mbbN$:
$R(\lambda,A_n) = R(\lambda,A)_n$.
\end{lemma}

\dremark{Beweis:
Es gilt:
$  R(\lambda,A)_n
=  ( (\lambda-A)^{-1} )_n
\overset{\ref{TnInvEqTInvn}}{=}  ( (\lambda \Id_X - A)_n )^{-1}
=  ( \lambda \Id_{M_n(X)} - A_n )^{-1}
=  R(\lambda,A_n)$.
}%

\begin{proof}[Beweis von Satz~\ref{eigErzHGOR}]
\bewitemph{(i):}
Nach Lemma~\ref{AnErzeugtTtn} ist $A_n$ Erzeuger der $C_0$"=Halbgruppe $(T_t^{(n)})_\indtHG$.
Mit Satz~\ref{unglFuerResHochn} erh"alt man f"ur alle $x \in M_n(X)$:
\[  R(\lambda,A_n)^k x
= \frac{1}{(k-1)!} \int_0^\infty s^{k-1} \mre^{-\lambda s} T_s^{(n)} x \, ds. \]
\dremark{Bew. analog zu \cite{WernerFunkana3}, Satz VII.4.10, l"anger:
Da f"ur alle $i,j \in \haken{n}$ die Einbettung
$\iota_{ij} : X \to M_n(X), x \mapsto (\delta_{ik}\delta_{j\ell} x)_{k,\ell \in \haken{n}}$,
linear und stetig ist,
erh"alt man mit \ref{ResolventeVonHG}.(ii) f"ur alle $x \in M_n(X)$:
\begin{align*}
   R(\lambda,A)_n x
&= ( R(\lambda,A) x_{ij} )_{i,j}
\overset{\ref{ResolventeVonHG}.(ii)}{=}
   \left( \int_0^\infty \mre^{-\lambda s} T_s x_{ij} \,ds \right)_{i,j}  \\
&= \sum_{i,j=1}^n \iota_{ij} \left( \int_0^\infty \mre^{-\lambda s} T_s x_{ij} \,ds \right)
\overset{\dremarkm{\text{\cite{vGrudAna2Phys}, 12.25}}}{=}
   \sum_{i,j=1}^n \int_0^\infty \iota_{ij}\left( \mre^{-\lambda s} T_s x_{ij} \right) \,ds  \\
&= \int_0^\infty \sum_{i,j=1}^n \mre^{-\lambda s} \iota_{ij}\left( T_s x_{ij} \right) \,ds
=  \int_0^\infty \mre^{-\lambda s} T_s^{(n)} x \,ds.
\end{align*}}%

\bewitemph{(ii):}
\dremark{Beweis geht analog zu \cite{EngelNagelSemigroups}, Corollary II.1.11.}%
Mit Lemma~\ref{RlAnEqRlAn} und (i) ergibt sich f"ur alle $x \in M_n(X)$:
\dremark{$(*)$: Benutze partielle Integration.}%
\dremark{18.4.'07, S. 1}%
\begin{align*}
   \normlr{\left(R(\lambda,A)^k\right)_n x}
&\overset{\dremarkm{\ref{RlAnEqRlAn}}}{=} \normlr{R(\lambda,A_n)^k x}  \\
&\overset{\dremarkm{\text{(i)}}}{=}
   \frac{1}{(k-1)!} \normlr{ \int_0^\infty s^{k-1} \mre^{-\lambda s} T_s^{(n)} x \,ds }  \\
\dremarkm{&\leq   \frac{1}{(k-1)!}
   \int_0^\infty \abslr{s^{k-1} \mre^{-\lambda s}} \normlr{T_s^{(n)}} \norm{x} \,ds  \\}
&\leq  \frac{M}{(k-1)!} \int_0^\infty s^{k-1} \mre^{(\omega - \Re \lambda) s} \,ds \,\norm{x}  \\
\dremarkm{\\ &\overset{\dremarkm{(*)}}{=}
   \left[ \frac{M}{(\omega - \Re \lambda)^k} \mre^{(\omega - \Re \lambda)s} \right]_0^\infty
   \,\norm{x}}
&=  \frac{M}{(\Re(\lambda) - \omega)^k} \,\norm{x}.
\end{align*}
\end{proof}

Analog zum Satz von Hille-Yosida f"ur kontraktive $C_0$-Halbgruppen
(Satz \ref{satzHilleYoskHG}) gilt im Operatorraum:

\begin{satz}[Satz von Hille-Yosida f"ur \cmvk{e} $C_0$-Halb\-grup\-pen]\label{satzHilleYosvkHG}
Sei $X$ ein Operatorraum.
Ein Operator $A : D(A) \subseteq X \to X$ ist genau dann Erzeuger einer
\cmvk{en} $C_0$-Halbgruppe, wenn $A$ dicht definiert und abgeschlossen ist,
$\mbbR_{>0} \subseteq \rho(A)$ gilt und
\begin{equation*}
\normcb{ \lambda R(\lambda,A) } \leq 1
   \qquad \text{f"ur alle } \lambda \in \mbbR_{>0}.
\end{equation*}
\end{satz}

Im Beweis des Satzes verwenden wir das folgende Lemma.
Es beschreibt die Resolventenmengen der Amplifikationen und folgt durch direktes Nachrechnen.
\dremark{Evtl. etwas zum Beweis sagen. (WW)}%

\begin{lemma}\label{rhoAnEqrhoA}
Sei $X$ ein Operatorraum und $A : D(A) \subseteq X \to X$ linear.
Dann gilt: $\rho(A_n) = \rho(A)$ f"ur alle $n \in \mbbN$.
\end{lemma}

\dremww{
Beweis.
\bewitemphq{$\subseteq$}
Sei $\lambda \in \rho(A_n)$.
Da $\lambda - A_n$ bijektiv auf $M_n(X)$ ist,
ist $\lambda -A$ bijektiv auf $X$.
Wegen $(\lambda - A)^{-1} x = \pr_{11}((\lambda - A_n)^{-1} (x \oplus 0_{M_{n-1}(X)}))$
ist $(\lambda - A)^{-1} \in L(X)$.
Es folgt: $\lambda \in \rho(A)$.

\bewitemphq{$\supseteq$}
Sei $\lambda \in \rho(A)$.
Da $\lambda - A$ bijektiv auf $X$ ist,
ist $\lambda - A_n = (\lambda - A)_n : D(A_n) = M_n(D(A)) \to M_n(X)$ bijektiv auf $M_n(X)$.
Nach \ref{TstetigFolgtTnSt} ist $\lambda - A_n = (\lambda - A)_n$ stetig.
\dremark{17.11.'08/1, 2.12.'08/2}%
}%

\begin{proof}[Beweis von Satz~\ref{satzHilleYosvkHG}]
\bewitemph{\glqq$\Rightarrow$\grqq} folgt
mit dem Satz von Hille-Yosida (Satz~\ref{satzHilleYoskHG}) und
mit der Absch"atzung der Resolvente aus Satz~\ref{eigErzvkHG}.
\dremark{$\omega=0$, $M=1$, $k=1$}%

\bewitemph{\glqq$\Leftarrow$\grqq:}
Nach Satz~\ref{satzHilleYoskHG} erzeugt
$A$ eine kontraktive $C_0$-Halbgruppe $(T_t)_\indtHG$ auf $X$.
Sei $n \in \mbbN$.
Weiter ist $A_n$ abgeschlossen\dremark{\ref{AabgFolgtAnabg}} und
dicht definiert\dremark{\ref{DdichtFolgtMnDdicht}}.
Ferner gilt:
$R(\lambda,A)_n = R(\lambda,A_n)$ (Lemma~\ref{RlAnEqRlAn}).
Mit Lemma~\ref{rhoAnEqrhoA} erh"alt man:
$\mbbR_{>0} \subseteq \rho(A) \overset{\dremarkm{\ref{rhoAnEqrhoA}}}{=} \rho(A_n)$.
Daher folgt mit Satz~\ref{satzHilleYoskHG}: $\normbig{T_t^{(n)}} \leq 1$
f"ur alle $\indtHG$.
Also erzeugt $A$ eine \cmvk{e} \cohg.
\end{proof}

Entsprechend zu Satz~\ref{satzHilleYosvkHG} beweist man den Satz von Hille-Yosida
f"ur beliebige $C_0$-Halbgruppen in Operatorr"aumen:

\begin{satz}[Satz von Hille-Yosida im Operatorraum, allgemeiner Fall]\label{satzHilleYosHGOR}
Sei $X$ ein Operatorraum und $A : D(A) \subseteq X \to X$ linear.
Sei $\omega \in \mbbR$ und $M \in \mbbR_{\geq 1}$.
Dann sind die folgenden Aussagen "aquivalent:
\begin{enumaequiv}
\item $A$ erzeugt eine $C_0$-Halbgruppe $(T_t)_\indtHG$ mit
\[ \normcb{T_t} \leq M\mre^{\omega t} \qquad\text{f"ur alle } \indtHG. \]
\item $A$ ist dicht definiert, abgeschlossen mit $\mbbR_{>\omega} \subseteq \rho(A)$,
und f"ur jedes $\lambda \in \mbbR_{>\omega}$ gilt:
\[ \normbig{\bigl[(\lambda - \omega)R(\lambda,A)\bigr]^k}_\text{cb}
   \leq M \qquad\text{f"ur alle } k \in \mbbN. \]
\end{enumaequiv}
\end{satz}

\dremark{
Beweis.
\bewitemph{\glqq(a)$\Rightarrow$(b)\grqq} folgt
mit dem Satz von Hille-Yosida (Satz~\ref{satzHilleYosHG}) und \ref{eigErzvkHG}.

\bewitemph{\glqq(b)$\Rightarrow$(a)\grqq:}
Nach dem Satz von Hille-Yosida (Satz \ref{satzHilleYosHG}) erzeugt
$A$ eine $C_0$-Halbgruppe $(T_t)_\indtHG$ auf $X$.
Weiter ist $A_n$ dicht definiert,\dremark{\ref{DdichtFolgtMnDdicht}}
abgeschlossen,\dremark{\ref{AabgFolgtAnabg}}
und es gilt:
$R(\lambda,A)_n = R(\lambda,A_n)$.
Mit \ref{rhoAnEqrhoA} erh"alt man:
$\mbbR_{> \omega} \subseteq \rho(A) \overset{\ref{rhoAnEqrhoA}}{=} \rho(A_n)$.
Daher folgt mit \ref{satzHilleYosHG}: $\normlr{T_t^{(n)}} \leq M \mre^{\omega t}$
f"ur alle $\indtHG$.

\begin{bemerkung}
Sei $X$ ein Operatorraum, $A$ Erzeuger einer \cmvk{en} $C_0$-Halbgruppe in $X$.
Dann ist $A_n$ Erzeuger einer kontraktiven $C_0$-Halbgruppe in $M_n(X)$.
\dremark{M"u\cms{}te gelten, pr"ufen. Z. 17.11.'08}%
\end{bemerkung}
}%

Man erh"alt eine Version des Satzes von Hille-Yosida f"ur $C_0$-Gruppen
(Satz \ref{satzHilleYosidaGr}) auf Operatorr"aumen:

\begin{satz}[Satz von Hille-Yosida f"ur $C_0$-Gruppen im Operatorraum]\label{satzHilleYosidaGrOR}
Sei $X$ ein Operatorraum und $A : D(A) \subseteq X \to X$ linear.
Sei $\omega \in \mbbR$ und $M \in \mbbR_{\geq 1}$.
Dann sind die folgenden Aussagen "aquivalent:
\begin{enumaequiv}
\item $A$ erzeugt eine $C_0$-Gruppe $(T_t)_\indtGr$ mit der Eigenschaft:
\[ \normcb{T_t} \leq M \mre^{\omega \abs{t}} \qquad\text{f"ur alle } \indtGr. \]

\item $A$ und $-A$ erzeugen $C_0$-Halbgruppen
$(T^+_t)_\indtHG$ bzw. $(T^-_t)_\indtHG$, die folgendes erf"ullen:
\[ \normcb{T^+_t}, \normcb{T^-_t} \leq M \mre^{\omega t} \qquad\text{f"ur alle } \indtHG. \]

\item $A$ ist dicht definiert, abgeschlossen und
f"ur jedes $\lambda \in \mbbR$ mit $\abs{\lambda} > \omega$
gilt $\lambda \in \rho(A)$ und
\dremark{Wird im Beweis von \ref{stoerungAnBeschraenkt} benutzt.}%
\[ \normbig{ [(\abs{\lambda} - \omega) R(\lambda,A)]^k }_{\text{cb}} \leq M
\qquad\text{f"ur alle } k \in \mbbN. \]
\end{enumaequiv}
\end{satz}

\begin{proof}
\bewitemphqq{(a)$\Rightarrow$(b)} folgt mit dem Text vor Satz~\ref{satzHilleYosidaGr}
(Satz von Hille-Yosida f"ur $C_0$-Gruppen).
\smallskip

\bewitemphq{(b)$\Rightarrow$(c)}
Es gelte \bewitemph{(b)}.
Nach Satz~\ref{satzHilleYosidaGr} ist $A$ dicht definiert und abgeschlossen mit
$]{-\infty},-\omega[ \,\,\cup\,\, ]\omega,\infty[ \,\,\subseteq \rho(A)$.
Nach Lemma~\ref{AnErzeugtTtn} erzeugt $A_n$ (bzw. $-A_n$) die $C_0$"=Halbgruppe
$\left( (T_t^+)^{(n)} \right)_\indtHG$ (bzw. $\left( (T_t^-)^{(n)} \right)_\indtHG$)
mit der Eigenschaft:
\[ \normlr{(T_t^+)^{(n)}}, \normlr{(T_t^-)^{(n)}} \leq M \mre^{\omega t}
   \quad\text{f"ur alle } \indtHG. \]
Mit $R(\lambda,A)_n = R(\lambda,A_n)$ (Lemma~\ref{RlAnEqRlAn}),
$\rho(A) = \rho(A_n)$ (Lemma~\ref{rhoAnEqrhoA}) und Satz~\ref{satzHilleYosidaGr} folgt:
\begin{equation}\label{eqHilleYosGrOR}
   \normbig{ [(\abs{\lambda} - \omega) R(\lambda,A)_n]^k }
=  \normbig{ [(\abs{\lambda} - \omega) R(\lambda,A_n)]^k }
\leq M
\end{equation}
f"ur alle $\lambda \in \mbbR$ mit $\abs{\lambda} > \omega$ und alle $k \in \mbbN$.

\bewitemphq{(c)$\Rightarrow$(a)}
Es gelte \bewitemph{(c)}.
Nach Satz~\ref{satzHilleYosidaGr} erzeugt $A$ eine $C_0$-Gruppe $(T_t)_\indtGr$ auf $X$.
Weiter ist $A_n$ dicht definiert, abgeschlossen und
erzeugt $(T_t^{(n)})_\indtHG$ nach Lemma~\ref{AnErzeugtTtn}.
Nach Lemma~\ref{rhoAnEqrhoA} hat man
$\rho(A_n) = \rho(A)$, und es gilt \eqref{eqHilleYosGrOR}.
Mit Satz~\ref{satzHilleYosidaGr} folgt somit \bewitemph{(a)}.
\end{proof}

\dremww{Eher weg:
"Ahnlich wie in \ref{UmnormHG} f"ur Banachr"aume existiert
f"ur eine beschr"ankte $C_0$-Halb\-grup\-pe auf einem Operatorraum
eine Matrixnorm derart, \dass die $C_0$"=Halbgruppe bez"uglich dieser Norm \cmvk{} ist:

\begin{bemerkung}
Sei $X$ ein Operatorraum und $(T_t)_\indtHG$ eine beschr"ankte $C_0$-Halbgruppe auf $X$.
Somit existiert ein $M \in \mbbR_{\geq 1}$ mit $\norm{T_t} \leq M$ f"ur alle $\indtHG$.
Dann wird durch
\begin{equation*}
\norm{x}_{s,n} := \sup_\indtHG \normlr{T_t^{(n)} x}_n \qquad\text{f"ur alle } n \in \mbbN, x \in M_n(X)
\end{equation*}
eine Matrixnorm auf $X$ so definiert, \dass $X$, versehen mit dieser Norm, ein Operatorraum ist,
in dem die $C_0$-Halbgruppe $(T_t)_\indtHG$ \cmvk{} ist.
\dremark{Bew.: Z. 29.12.'07, S. 3/4 (Umnormieren)}%
\end{bemerkung}

\dfrage{Wird durch $\norm{\cdot}_s$ eine "aquivalente Operatorraumnorm auf $X$ definiert,
  \dheisst: $\exists C \in \mbbR_{\geq 1}\,\, \forall n \in \mbbN\,\, \forall x \in M_n(X):
    \norm{x}_n \leq \norm{x}_{s,n} \leq C \norm{x}_n$?
  Alternativ: Z.\,B. verlangen, \dass $\sup_{\indtHG} \normcb{T_t} < \infty$.}%
}%

Auch f"ur den Satz von Lumer-Phillips (Satz~\ref{SatzLumerPhillips}) beweisen wir eine
entsprechende Fassung in Operatorr"aumen.
Bevor wir diese formulieren, f"uhren wir zun"achst den folgenden Begriff ein:

\begin{definitn}
Sei $X$ ein Operatorraum.
Man nennt einen Operator $A : D(A) \subseteq X \to X$ \defemph{vollständig dissipativ},
falls $A_n$ dissipativ f"ur alle $n \in \mbbN$ ist, \dheisst
\index[B]{vollständig dissipativer Operator}%
\[ \normlr{(\lambda-A_n)x} \geq \lambda \norm{x}
   \quad\text{f"ur alle } \lambda \in \mbbR_{>0} \text{ und } x \in D(A_n). \]
\end{definitn}

\begin{satz}[Satz von Lumer-Phillips f"ur \cmvk{e} $C_0$-Halb\-grup\-pen]\label{SatzLumerPhillipsOR}
Sei $X$ ein Operatorraum und $A : D(A) \subseteq X \to X$ linear und dicht definiert.
Dann erzeugt $A$ genau dann eine \cmvk{e} $C_0$"=Halbgruppe, wenn $A$ vollst"andig dissipativ ist
und ein $\lambda \in \mbbR_{>0}$ so existiert, \dass $\lambda - A$ surjektiv auf $X$ ist.
\end{satz}

\begin{proof}
\bewitemph{\glqq$\Rightarrow$\grqq:}
$A$ erzeuge eine \cmvk{e} $C_0$-Halb\-grup\-pe.
Nach Satz~\ref{satzHilleYosHGOR} (Satz von Hille-Yosida) gilt $\mbbR_{>0} \subseteq \rho(A)$
und $\norm{\lambda  R(\lambda,A)}_\text{cb} \leq 1$
f"ur alle $\lambda \in \mbbR_{>0}$,
also folgt
\[ \norm{(\lambda-A)_n x} \geq \frac{1}{ \norm{(\lambda-A)_n^{-1}x} } \geq \lambda \norm{x} \]
f"ur alle $n \in \mbbN$ \text{ und } $x \in D(A_n)$.
\smallskip

\bewitemph{\glqq$\Leftarrow$\grqq} folgt mit dem Satz von Lumer-Phillips
(Satz~\ref{SatzLumerPhillips}).
\dremark{
Sei $A$ vollst"andig dissipativ und $\lambda \in \mbbR_{>0}$ so,
\dass $\lambda - A$ surjektiv ist.
Nach dem Satz von Lumer-Phillips (\ref{SatzLumerPhillips})
erzeugt $A_n$ eine kontraktive $C_0$-Halbgruppe.}%
\end{proof}

\dremww{
Mit dem obigen Satz\dremark{\ref{SatzLumerPhillipsOR}} erh"alt man:

\begin{satz}[Satz von Lumer-Phillips f"ur \cmvk{e} $C_0$-Grup\-pen]\label{SatzLumerPhillipsGrOR}
Sei $X$ ein Operatorraum und $A : D(A) \subseteq X \to X$ linear und dicht definiert.
Dann erzeugt $A$ genau dann eine \cmvk{e} $C_0$-Gruppe,
wenn $A$ und $-A$ vollst"andig dissipativ sind
und $\lambda, \mu \in \mbbR_{>0}$ so existieren, \dass $\lambda - A$ und $\mu + A$ surjektiv auf $X$ sind.
\end{satz}
}%

F"ur einen beliebigen Banachraum $X$
liefert der obige Satz unter Verwendung der minimalen Operatorraumstruktur $\min(X)$
(siehe Definitions"=Proposition~\ref{defminX}) gerade den Satz von Lumer-Phillips f"ur Banachr"aume
(Satz~\ref{SatzLumerPhillips}),
also ist Satz~\ref{SatzLumerPhillipsOR} eine Verallgemeinerung dieses Satzes.
Ebenso sind die hier formulierten S"atze von Hille-Yosida
f"ur Operatorr"aume Verallgemeinerungen der entsprechenden S"atze in Banachr"aumen.
\skiptext

Mit Hilfe des Satzes von Hille-Yosida f"ur \cmvk{e} $C_0$"=Halbgruppen
(Satz~\ref{satzHilleYosvkHG}) und des Satzes von Lumer-Phillips in Operatorr"aumen
(Satz~\ref{SatzLumerPhillipsOR}) erh"alt man die folgende Charakterisierung
von bestimmten $C_0$"=Linksmultiplikatoren:

\begin{bemerkung}
Sei $X$ ein Operatorraum und $A : D(A) \subseteq X \to X$ linear.
Dann sind die folgenden Aussagen "aquivalent:
\begin{enumaequiv}
\item $A$ erzeugt eine $C_0$-Halbgruppe $(T_t)_\indtHG$ auf $X$ mit
  $T_t$ aus der abgeschlossenen Kugel $\overline{\cmkug}_\Multlor{X}(0,1)$
f"ur alle $\indtHG$.
\item $A$ ist dicht definiert und abgeschlossen mit $\mbbR_{>0} \subseteq \rho(A)$ und
\[ \normcb{\lambda R(\lambda,\matnull{A})} \leq 1 \qquad\text{f"ur alle } \lambda \in \mbbR_{>0}. \]
\item $A$ ist dicht definiert, $\matnull{A}$ vollst"andig dissipativ,
  und es existiert ein $\lambda \in \mbbR_{>0}$ derart,
  \dass $\lambda - A$ surjektiv auf $X$ ist.
\end{enumaequiv}
\end{bemerkung}

\begin{proof}
\bewitemphq{(a)$\Rightarrow$(b)}
Es gelte (a).
Nach Satz~\ref{charMultlcomatnull} (Charakterisierung von $\Multlco{X}$ mittels $\matnull{A}$)
erzeugt $\matnull{A}$ eine vollst"andig kontraktive $C_0$"=Halbgruppe auf $C_2(X)$.
Nach dem Satz von Hille-Yosida f"ur \cmvk{e} $C_0$"=Halbgruppen
(Satz~\ref{satzHilleYosvkHG}) ist $\matnull{A}$ dicht definiert und abgeschlossen
mit $\mbbR_{>0} \subseteq \rho(\matnull{A})$ und
\[ \normcb{\lambda R(\lambda,\matnull{A})} \leq 1 \qquad\text{f"ur alle } \lambda \in \mbbR_{>0}. \]
Dann ist $A$ nach Lemma~\ref{zshgAundA0} dicht definiert und abgeschlossen,
und es gilt: $\mbbR_{>0} \subseteq \rho(\matnull{A}) = \rho(A) \setminus\{0\}$.

\bewitemphq{(b)$\Rightarrow$(c)}
Es gelte (b).
Nach Lemma~\ref{zshgAundA0} ist $\matnull{A}$ dicht definiert und abgeschlossen,
und es gilt: $\mbbR_{>0} \subseteq \rho(\matnull{A})$.
Es folgt mit Satz~\ref{satzHilleYosvkHG},
\dass $\matnull{A}$ Erzeuger einer \cmvk{en} $C_0$-Halbgruppe ist.
Mit dem Satz von Lumer-Phillips im Operatorraum (Satz~\ref{SatzLumerPhillipsOR}) erh"alt man,
\dass $\matnull{A}$ vollst"andig dissipativ ist.
Wegen $\mbbR_{>0} \subseteq \rho(A)$ folgt die Existenz des gesuchten $\lambda$.
\dremark{Existenz von $\lambda$: $\lambda \in \rho(A)$}%

\bewitemphq{(c)$\Rightarrow$(a)}
Es gelte (c).
Dann ist $\lambda - \matnull{A}$ surjektiv auf $C_2(X)$.
Nach Lemma~\ref{zshgAundA0} ist $\matnull{A}$ dicht definiert.
Es ergibt sich mit Satz~\ref{SatzLumerPhillipsOR},
\dass $\matnull{A}$ eine \cmvk{e} $C_0$-Halbgruppe auf $C_2(X)$ erzeugt.
Mit Satz~\ref{charMultlcomatnull} erh"alt man (a).
\end{proof}

Unter Verwendung des Satzes von Hille-Yosida f"ur $C_0$-Gruppen im Operatorraum
(Satz~\ref{satzHilleYosidaGrOR})\dremark{und \ref{SatzLumerPhillipsGrOR}}
erh"alt man eine entsprechende Aussage f"ur $\Adjlco{X}$:

\begin{bemerkung}\label{charAdjlormatnull2}
Sei $X$ ein Operatorraum und $A : D(A) \subseteq X \to X$ linear.
Dann sind die folgenden Aussagen "aquivalent:
\begin{enumaequiv}
\item $A \in \Adjlco{X}$.
\item $A$ ist dicht definiert und abgeschlossen mit $\mbbR\setminus\{0\} \subseteq \rho(A)$ und
\[ \normcblr{\abs{\lambda} R(\lambda,\matnull{A})} \leq 1 \qquad\text{f"ur alle } \lambda \in \mbbR\setminus\{0\}. \]
\item $A$ ist dicht definiert, $\matnull{A}$ und $-\matnull{A}$ sind vollst"andig dissipativ,
  und es existieren $\lambda, \mu \in \mbbR_{>0}$ derart,
  \dass $\lambda - A$ und $\mu + A$ surjektiv auf $X$ sind.
\end{enumaequiv}
\end{bemerkung}

Mit der obigen Proposition folgt zusammen mit Satz~\ref{charAdjlcomatnull}
(Charakterisierung von $\Adjlco{X}$ mittels $\matnull{A}$)
die folgende Charakterisierung der Elemente aus $\Adjloru{X}$:

\begin{bemerkung}
Sei $X$ ein Operatorraum und $A : D(A) \subseteq X \to X$ linear.
Dann sind die folgenden Aussagen "aquivalent:
\begin{enumaequiv}
\item $A \in \Adjloru{X}$.
\item Es existiert ein $B : D(B) \subseteq X \to X$ derart, \dass
  $\tilde{A} := \mri \cmsmallpmatrix{0 & A \\ B & 0}$ dicht definiert und abgeschlossen ist
  mit $\mbbR\setminus\{0\} \subseteq \rho(\tilde{A})$ und
\[ \normcblr{\abs{\lambda} R(\lambda,\matnull{\tilde{A}})} \leq 1 \qquad\text{f"ur alle } \lambda \in \mbbR\setminus\{0\}. \]
\item Es existiert ein $B : D(B) \subseteq X \to X$ mit den folgenden Eigenschaften:
  $\tilde{A} := \mri \cmsmallpmatrix{0 & A \\ B & 0}$ ist dicht definiert,
  $\matnull{\tilde{A}}$ und $-\matnull{\tilde{A}}$ sind vollst"andig dissipativ,
  und $\lambda, \mu \in \mbbR_{>0}$ existieren derart,
  \dass $\lambda - \tilde{A}$ und $\mu + \tilde{A}$ surjektiv auf $C_2(X)$ sind.
\end{enumaequiv}
\end{bemerkung}

\dremark{Beweis.
Beachte: $\mri \cmsmallpmatrix{0 & A \\ A^* & 0} \in \Adjlco{C_2(X)}$}%

\dremark{Evtl. weitere Aussagen von Banachr"aumen auf Operatorr"aume "ubertragen.}%

\section[Unbeschr"ankte Multiplikatoren von $X$ auf $\ternh{X}$ "uberf"uhren]{Unbeschr"ankte Multiplikatoren von $X$ auf\newlinef{}$\ternh{X}$ "uberf"uhren}

Wir zeigen in diesem Abschnitt,
wie man zu einem $C_0$-Linksmultiplikator auf einem Operatorraum $X$
einen $C_0$-Linksmultiplikator auf der tern"aren H"ulle $\ternh{X}$ erh"alt
und halten Eigenschaften des hierdurch entstehenden Operators fest.
Entsprechende Aussagen beweisen wir auch f"ur unbeschr"ankte schiefadjungierte Multiplikatoren
und f"ur unbeschr"ankte Multiplikatoren.
Au"serdem wird bewiesen, unter welchen Voraussetzungen man einen unbeschr"ankten Multiplikator
so auf einen Unteroperatorraum einschr"ankten kann,
\dass man wieder einen unbeschr"ankten Multiplikator erh"alt.
Ferner wird gezeigt, \dass man unbeschr"ankte Multiplikatoren auf $X$
als Einschr"ankungen von regul"aren Operatoren auf $\ternh{X}$ auf\/fassen kann.
\skiptext

Es sei an die folgende Definition erinnert:

\begin{definitn}
Sei $X$ ein Banachraum und $A : D(A) \subseteq X \to X$ linear.
Ein Unterraum $U$ von $D(A)$ hei"st
\defemphi{wesentlicher Bereich} f"ur $A$,
falls $U$ dicht in $D(A)$ bez"uglich der Graphennorm
\dliter{\cite{EngelNagelSemigroups}, Def. II.1.6}%
\[ \norm{x}_A = \norm{x} + \norm{Ax} \qquad\text{f"ur alle } x \in X \]
liegt.
\end{definitn}

\dmarginpar{pr}\dremark{Evtl. an Bezeichnungen erinnern ($j$,$\isoIMl$,\dots).}%

Der folgende Satz schildert, unter welchen Voraussetzungen man ein
$A \in \Multlco{X}$ auf die tern"are H"ulle $\ternh{X}$
(zur Definition siehe Definitions"=Satz \ref{defTernHuelle}) "uberf"uhren kann:

\begin{satz}\label{halbgrHochliftenAufTX}
Sei $X$ ein Operatorraum.
Sei $A \in \Multlco{X}$ derart,
\dass f"ur die von $A$ erzeugte $C_0$-Halbgruppe $(T_t)_\indtHG$ gilt:
\begin{equation}\label{eqNormatEndl}
\exists \delta \in \mbbR_{>0} : \sup_{t \in [0,\delta]} \normMl{T_t}{X} < \infty.
\end{equation}
Dann findet man ein eindeutig bestimmtes $B \in \Multlco{\ternh{X}}$ mit der Eigenschaft:
$j(D(A)) \subseteq D(B)$ und $j \circ A = B \circ j\restring_{D(A)}$.
\dremark{Zur Bedingung \eqref{eqNormatEndl}: (i) Sie ist f"ur eine unit"are $C_0$-Halbgruppe
  $(T_t)_t$ in einem Hilbert-\cstern{}Modul $X$ erf"ullt,
  denn es gilt $1 \overset{\ref{HCsModulIsometrieIstVollstIsom}}{=} \normcb{T_t} = \norm{a_t}$
  (2. \glqq$=$\grqq: nach \ref{NormcbLeqNormMl}).
  (ii) Falls $\tau^c_{T_t}$ \cmvk{} ist, ist $a_t \leq 1$ nach \ref{charMultOR}.
  (iii) Etwas zur Bedingung schreiben.}%
\dremark{Man kann $\ternh{X}$ nicht durch eine tern"are Erweiterung ersetzen,
  da das Produkt mit $a_t$ definiert sein mu\cms{}.}%
\dremark{Ist $A \in \Multlor{X}$, so folgt: $B \in \Multlor{\ternh{X}}$.}%
\end{satz}

Bevor wir zum Beweis der obigen Satzes kommen, halten wir fest:

\begin{bemerkung}\label{halbgrHochliftenAufTX2}
Seien $X$, $A$, $(T_t)_\indtHG$ und $B$ wie in Satz~\ref{halbgrHochliftenAufTX}.
\begin{enumaufz}
\item F"ur alle $\indtHG$ setze $a_t := \isoIMl^{-1}(T_t)$,
also $j(T_t x) = a_t \multis j(x)$ f"ur alle $x \in X$ (siehe Satz~\ref{IsomIMlAufMl}).
Definiere $S_t : \ternh{X} \to \ternh{X}, z \mapsto a_t \multis z$, f"ur alle $\indtHG$.
Dann ist $(S_t)_\indtHG$ die eindeutig bestimmte $C_0$-Halbgruppe auf $\ternh{X}$ mit
$S_t \in \Multlor{\ternh{X}}$ und $S_t \circ j = j \circ T_t$ f"ur alle $\indtHG$.
Au"serdem wird $(S_t)_\indtHG$ von $B$ erzeugt.

\item Es gilt:
\begin{align*}
   &\phantom{=}\,\,\, B\bigl( j(x_1) \multis j(x_2)^* \multis j(x_3) \multis j(x_4)^* \cdots j(x_{2k+1}) \bigr)  \\
&= j(A(x_1)) \multis j(x_2)^* \multis j(x_3) \multis j(x_4)^* \cdots j(x_{2k+1})
\end{align*}
f"ur alle $k \in \mbbN_0$, $x_1 \in D(A)$ und $x_2, \dots, x_{2k+1} \in X$.

\dremark{\item Es gilt: $a_t \multis j(D(A)) \subseteq j(D(A))$ und
  $S_t(W) \subseteq W$ f"ur alle $\indtHG$.}%

\item F"ur jeden Operatorraum $Y$ und $M \subseteq Y$ definiere
\begin{align*}
  W_{M,Y}
  &:= \cmlin\bigl\{ j(x_1) \multis j(x_2)^* \multis j(x_3) \multis j(x_4)^* \cdots j(x_{2k+1}) \setfdg  \\
  &\phantom{:= \cmlin\big\{}  k \in \mbbN_0, x_1 \in M, x_2, \dots, x_{2k+1} \in Y \bigr\}
  \subseteq \ternh{Y}.
\end{align*}
\index[S]{WMY@$W_{M,Y}$}%
Dann ist $W := W_{D(A),X}$
ein wesentlicher Bereich f"ur~$B$, das hei"st, $\overline{B\restring_W} = B$.
\dremark{$B$ ist die Abschlie"sung von $B\restring_W$
  (\cite{EngelNagelSemigroups}, Exerc. II.1.15.(2)).}%
\dremark{Idee: Auf $W$ verh"alt sich $B$ wie $A$.}%
\dremark{Frage: Gilt $W = D(B)$?}%
\end{enumaufz}
\end{bemerkung}

Die Bedingung \eqref{eqNormatEndl} wird benutzt, um mittels
der Charakterisierung der starken Stetigkeit von Halbgruppen
(Proposition~\ref{charHGIstC0}) zu zeigen,
\dass die Halbgruppe $(S_t)_\indtHG$ stark stetig ist.

\begin{proof}[Beweis von Satz~\ref{halbgrHochliftenAufTX} und Proposition~\ref{halbgrHochliftenAufTX2}]
\bewitemph{(i):} Mit Lemma~\ref{atIstOpHalbgr} folgt,
\dass $(S_t)_\indtHG$ eine Operatorhalbgruppe ist.
\dremark{Analog zu \ref{atIstOpHalbgr}.(iv) ($R_t$).}%

Zun"achst wird die starke Stetigkeit von $(S_t)_\indtHG$ gezeigt.

F"ur alle $x,y,z \in X$ gilt:
\begin{align*}
&  \norm{a_t \multis j(x) \multis j(y)^* \multis j(z) - j(x) \multis j(y)^* \multis j(z)}  \\
\leq\,\, &\norm{a_t \multis j(x) - j(x)} \cdot \norm{j(y)^* \multis j(z)}  \\
=\,\, &\norm{T_t x - x} \cdot \norm{j(y)^* \multis j(z)} \to 0 \qquad\text{f"ur } t \to 0.
\end{align*}
Analog folgt f"ur alle $n \in \mbbN$ und $x \in X^{2n+1}$ mit
$y := j(x_1) \multis j(x_2)^* \multis j(x_3) \multis j(x_4)^* \cdots j(x_{2n+1})$:
 $\norm{a_t \multis y - y} \to 0$ f"ur $t \to 0$.
Setze
\[ V := \cmlin\left\{ j(x_1) \multis j(x_2)^* \multis j(x_3) \multis j(x_4)^* \cdots j(x_{2n+1}) \setfdg
                n \in \mbbN_0, x \in X^{2n+1} \right\}. \]
Dann gilt:
$\lim_{t \downarrow 0} S_t z = \lim_{t \downarrow 0} a_t \multis z = z$ f"ur alle $z \in V$.
Nach \eqref{eqTXeqSpan} ist
\begin{equation}\label{eqTXeqSpanLiften}
\ternh{X} \cong
   \overline{\cmlin}\left\{x_1 \multis x_2^* \multis x_3 \multis x_4^* \cdots x_{2n+1} \setfdg
     n \in \mbbN_0, x \in j(X)^{2n+1} \right\}  \subseteq I(X),
\end{equation}
somit liegt $V$ dicht in $\ternh{X}$.
Nach Proposition~\ref{I11TXeqI11X} hat man: $I_{11}(\ternh{X}) \cong I_{11}(X)$.
Man erh"alt
\[ \norm{S_t z}_{\ternh{X}}
=  \norm{a_t \multis z}_{\ternh{X}}
\leq  \norm{a_t}_{I_{11}(\ternh{X})} \, \norm{z}_{\ternh{X}}
=  \norm{a_t}_{I_{11}(X)} \, \norm{z}_{\ternh{X}} \]
f"ur alle $\indtHG$ und $z \in \ternh{X}$,
also $\norm{S_t} \leq \norm{a_t} = \normMl{T_t}{X}$.
Zusammen mit \eqref{eqNormatEndl} und
der Charakterisierung der starken Stetigkeit von Halbgruppen (Proposition~\ref{charHGIstC0}) folgt,
\dass $(S_t)_\indtHG$ stark stetig ist.
\dremark{stark stetig: Z. 16.8.'07/1 und 15.8.'07/1}%
\dremark{Sei $\sigma : \ternh{X} \to I(\mcalS(X)), z \mapsto z$.
  Es gilt $\sigma(S_t z) = a_t \multis \sigma(z)$ f"ur alle $z \in \ternh{X}$,
  also $S_t \in \Multlor{\ternh{X}}$.}%
\smallskip

Sei $B$ der Erzeuger von $(S_t)_\indtHG$.
Da offensichtlich $S_t \in \Multlor{\ternh{X}}$ f"ur alle $\indtHG$ gilt,
\dremark{$S_t(z) = a_t \multis z$, $I(S(X))$ ist \csalgebra, die $a_t$ und $z$ enth"alt.}%
ist $B \in \Multlco{\ternh{X}}$.
Man hat f"ur alle $x \in D(A)$:
\begin{equation*}
   j(Ax)
=  j \left( \lim_{t \downarrow 0} \frac{T_t x - x}{t} \right)
\dremarkm{=  \lim_{t \downarrow 0} \frac{j(T_t x) - j(x)}{t}}
=  \lim_{t \downarrow 0} \frac{a_t \multis j(x) - j(x)}{t}
=  B(j(x)).
\end{equation*}
Insbesondere folgt: $j(D(A)) \subseteq D(B)$.
\smallskip

Um zu zeigen, \dass $(S_t)_\indtHG$
durch die aufgef"uhrten Eigenschaften eindeutig bestimmt ist,
geben wir uns eine weitere $C_0$-Halbgruppe $(R_t)_\indtHG$ auf $\ternh{X}$ vor mit
$R_t \in \Multlor{\ternh{X}}$ und $R_t \circ j = j \circ T_t$ f"ur alle $\indtHG$.
Sei $\indtHG$.
Man findet nach Satz~\ref{charMultOR} (Charakterisierung der Linksmultiplikatoren)
ein $d_t \in \mIMl(\ternh{X}) \subseteq I_{11}(\ternh{X})$ mit $R_t(z) = d_t \multis z$
f"ur alle $z \in \ternh{X}$.
Nach Proposition~\ref{I11TXeqI11X} gilt: $I_{11}(\ternh{X}) \cong I_{11}(X)$.
Somit folgt: $d_t \in I_{11}(X)$.
Man hat f"ur alle $x \in X$:
\[ (a_t - d_t)j(x)
=  S_t(j(x)) - R_t(j(x))
=  j(T_t x) - j(T_t x)
=  0. \]
Mit \refb{\cite{BlecherLeMerdy04OpAlg}, Proposition 4.4.12,}{Proposition~\ref{Lcinj}}
erh"alt man\dremark{wegen $d_t \in I_{11}(\ternh{X}) \cong I_{11}(X)$}:
$a_t = d_t$, also $S_t = R_t$.
\smallskip

\bewitemph{(ii) und (iii):}
\dremark{\bewitemph{($\alpha$):}}%
Da $D(A)$ dicht in $X$ liegt und \eqref{eqTXeqSpanLiften} gilt, ist $W$ dicht in $\ternh{X}$.
\dmarginpar{pr}\dremark{Nach Dennis Beweis notieren.}%
\dremark{Sei $z \in \ternh{X}$ und $\varepsilon \in \mbbR_{>0}$.
Nach \eqref{eqTXeqSpan} findet man $k \in \mbbN$, $\lambda \in \mbbC^k$, $m \in \mbbN_0^k$,
$x_1^{(i)}, \dots, x_{2m_i +1}^{(i)} \in X$ f"ur alle $i \in \haken{m_i}$ so,
\dass f"ur $y := \sum_{i=1}^k \lambda_i j(x_1^{(i)}) \cdots j(x_{2m_i +1}^{(i)})$ gilt:
$\norm{z - y} < \frac{\varepsilon}{2}$.
Da $D(A)$ dicht in $X$ liegt, findet man f"ur alle $i \in \haken{k}$ ein $d_i \in D(A)$ mit
\[ \norm{x_1^{(i)} - d_i}
<  \frac{\varepsilon}{2k}
   \frac{1}{\norm{ \lambda_i j(x_2^{(i)})^* \cdots j(x_{2m_i +1}^{(i)}) } +1}. \]
Setze $w := \sum_{i=1}^k \lambda_i j(d_i) \multis j(x_2^{(i)})^* \cdots j(x_{2m_i +1}^{(i)}) \in W$.
Es gilt:
\begin{align*}
  \norm{z - w}
&\leq  \norm{z-y} + \norm{y-w}  \\
&<  \frac{\varepsilon}{2} +
   \norm{ \sum_{i=1}^k \lambda_i \bigl( j(x_1^{(i)}) - j(d_i) \bigr) \multis
            j(x_2^{(i)})^* \cdots j(x_{2m_i +1}^{(i)}) }  \\
&\leq   \frac{\varepsilon}{2} +
   \sum_{i=1}^k \norm{ \lambda_i \bigl( j(x_1^{(i)}) - j(d_i) \bigr) \multis
            j(x_2^{(i)})^* \cdots j(x_{2m_i +1}^{(i)}) }  \\
&\leq   \frac{\varepsilon}{2} +
   \sum_{i=1}^k  \norm{x_1^{(i)} - d_i}\,
            \norm{\lambda_i j(x_2^{(i)})^* \cdots j(x_{2m_i +1}^{(i)})}
\leq  \frac{\varepsilon}{2} + k \frac{\varepsilon}{2k}
=  \varepsilon.
\end{align*}
}%

\dremark{Zeige: $S_s(W) \subseteq W$.}%
Sei $s \in \mbbR_{\geq 0}$ und $x \in D(A)$.
Es gilt:
\begin{align*}
   a_s \multis j(Ax)
&=  a_s \multis j\left( \lim_{t \downarrow 0} \frac{ T_t x - x }{t} \right)
\dremarkm{=  j\left( T_s \left(\lim_{t \downarrow 0} \frac{ T_t x - x }{t} \right) \right) } \\
&=  j\left( \lim_{t \downarrow 0} \frac{ T_s T_t x - T_s x }{t} \right)
=  j\left( \lim_{t \downarrow 0} \frac{ T_t(T_s x) - T_s x }{t} \right).
\end{align*}
Man erh"alt $T_s x \in D(A)$, also $a_s \multis j(x) = j(T_s x) \in j(D(A))$.
\dremark{Man hat somit:
\[ a_s \multis j(D(A)) \subseteq j(D(A)). \dremarkm{\qquad (*)} \]}%
Hiermit folgt: $S_s (W) \subseteq W$.
\dmarginpar{pr.: Bew. not.}\dremark{Beweis notieren.}%
\dremark{Sei $w \in W$.
Dann hat $w$ die Gestalt
$w = \sum_{i=1}^k \lambda_i j(x_1^{(i)}) \multis j(x_2^{(i)})^* \cdots j(x_{2m_i +1}^{(i)})$,
wobei $x_1^{(i)} \in D(A)$, $x_2^{(i)}, \dots, x_{2m_i +1}^{(i)} \in X$.
Es folgt:
\[ S_s w
=  a_s \multis
   \sum_{i=1}^k \lambda_i j(x_1^{(i)}) \multis j(x_2^{(i)})^* \cdots j(x_{2m_i +1}^{(i)})
=  \sum_{i=1}^k \lambda_i \underbrace{a_s \multis j(x_1^{(i)})}_{\in j(D(A)) \text{ nach } (*)}
                \multis j(x_2^{(i)})^* \cdots j(x_{2m_i +1}^{(i)}) \in W. \]
}%

\dremark{\bewitemph{($\beta$):}}
F"ur alle $k \in \mbbN_0$, $x_1 \in D(A)$ und $x_2, \dots, x_{2k+1} \in X$ gilt:
\begin{equation}\label{eqjAx1EqB}
\begin{split}
&\phantom{=}\,\,\,  j(A x_1) \multis \underbrace{j(x_2)^* \multis j(x_3) \multis j(x_4)^* \multis j(x_5) \cdots j(x_{2k+1})}_{=:z}  \\
&= \lim_{t \downarrow 0} \frac{a_t \multis j(x_1) - j(x_1)}{t} \multis z
=  B( j(x_1) \multis z ).
\end{split}
\end{equation}
\dremark{
\begin{align*}
&\phantom{=}\,\,\,  j(A x_1) \multis \underbrace{j(x_2)^* \multis j(x_3) \multis j(x_4)^* \multis j(x_5) \cdots j(x_{2k+1})}_{=:z}  \\
&= \lim_{t \downarrow 0} \frac{a_t \multis j(x_1) - j(x_1)}{t} \multis z  \dremarkm{\quad (\dagger)\\
&= \lim_{t \downarrow 0} \frac{a_t \multis j(x_1) \multis z - j(x_1) \multis z}{t}}
=  B( j(x_1) \multis z ).
\end{align*}}%
Somit ergibt sich: $W \subseteq D(B)$.
\dremark{Sei $w \in W$ wie oben.
Dann gilt:
\begin{align*}
   B(w)
&= B( \sum_{i=1}^k \lambda_i j(x_1^{(i)}) \multis j(x_2^{(i)})^* \cdots j(x_{2m_i +1}^{(i)}) )
=  \sum \lambda_i B( j(x_1^{(i)}) \multis j(x_2^{(i)})^* \cdots j(x_{2m_i +1}^{(i)}) )  \\
&\overset{(\dagger)}{=}
   \sum \lambda_i j(A x_1^{(i)}) \multis j(x_2^{(i)})^* \cdots j(x_{2m_i +1}^{(i)}).
\end{align*}
}%
Nach \cite{EngelNagelSemigroups}, Proposition~II.1.7, ist
$W$ ein wesentlicher Bereich f"ur $B$.

Um die Eindeutigkeit von $B$ zu zeigen, geben wir uns ein
$C \in \Multlco{\ternh{X}}$ vor mit der Eigenschaft:
$j(D(A)) \subseteq D(C)$ und $j \circ A = C \circ j\restring_{D(A)}$.
Dann gilt nach \eqref{eqjAx1EqB}:
\begin{align*}
   &\phantom{=}\,\,\, B\bigl( j(x_1) \multis j(x_2)^* \multis j(x_3) \multis j(x_4)^* \cdots j(x_{2k+1}) \bigr)  \\
&= j(A(x_1)) \multis j(x_2)^* \multis j(x_3) \multis j(x_4)^* \cdots j(x_{2k+1})  \\
&= C\bigl( j(x_1) \multis j(x_2)^* \multis j(x_3) \multis j(x_4)^* \cdots j(x_{2k+1}) \bigr)
\end{align*}
f"ur alle $k \in \mbbN_0$, $x_1 \in D(A)$ und $x_2, \dots, x_{2k+1} \in X$.
Also folgt: $B\restring_W = C\restring_W$.
Da $W$ ein wesentlicher Bereich f"ur $B$ ist,
erh"alt man mit \cite{WernerFunkana6}, VII.5.35,
\dass $B$ gleich der Abschlie"sung von $B\restring_W$ ist,
also gleich dem Abschlu\cms{} des Graphen von $B\restring_W$,
notiert als $\overline{B\restring_W}$.
Somit ergibt sich:
\dremark{29.9.'08/1, 16.10.'08/1}%
\[ B = \overline{B\restring_W} = \overline{C\restring_W} = C.  \qedhere \]
\end{proof}

Im folgenden formulieren wir Versionen der vorherigen zwei Aussagen
\dremark{\ref{halbgrHochliftenAufTX}, \ref{halbgrHochliftenAufTX2}}%
f"ur $\Adjlco{X}$.
Hierbei wird die Bedingung \eqref{eqNormatEndl} nicht mehr ben"otigt,
da sie wegen $\normMl{T}{X} = \norm{T}$
f"ur alle $T \in \Adjlor{X}$ (Proposition~\ref{inAlXNormTEqNormcb}) automatisch erf"ullt ist.

\begin{satz}\label{grHochliftenAufTX}
Sei $X$ ein Operatorraum und $A \in \Adjlco{X}$.
Dann findet man ein eindeutig bestimmtes $B \in \Adjlco{\ternh{X}}$ mit
$j(D(A)) \subseteq D(B)$ und $j \circ A = B \circ j\restring_{D(A)}$.
\dremark{Beachte: $\ternh{X}$ ist \hcsmodul.
  Man kann somit den Satz von Stone anwenden.}%
\end{satz}

Wir notieren zun"achst die folgende Proposition und beweisen
anschlie"send den obigen Satz.

\begin{bemerkung}\label{grHochliftenAufTX2}
Seien $X$, $A$ und $B$ wie in Satz~\ref{grHochliftenAufTX}.
\begin{enumaufz}
\item Sei $(T_t)_\indtGr$ die von $A$ erzeugte $C_0$-Gruppe.
Setze $a_t := \isoIMl^{-1}(T_t) \in \mIMl^*(X)$ f"ur jedes $\indtGr$.
Es gilt:
\[ \exists \delta \in \mbbR_{>0} :
   \sup_{t \in [-\delta,\delta]} \norm{a_t}
=  \sup_{t \in [-\delta,\delta]} \normMl{T_t}{X} < \infty. \]
Definiere $S_t : \ternh{X} \to \ternh{X}, z \mapsto a_t \multis z$, f"ur alle $\indtGr$.
Dann ist $(S_t)_\indtGr$ die eindeutig bestimmte $C_0$-Gruppe auf $\ternh{X}$ mit
$S_t \in \Adjlor{\ternh{X}}$ unit"ar und $S_t \circ j = j \circ T_t$ f"ur alle $\indtGr$.
Au"serdem wird $(S_t)_\indtGr$ von $B$ erzeugt.

\item Die Aussagen (ii) und (iii) aus Proposition~\ref{halbgrHochliftenAufTX2}
gelten entsprechend.
\end{enumaufz}
\end{bemerkung}

\begin{proof}[Beweis von Satz~\ref{grHochliftenAufTX} und   Proposition~\ref{grHochliftenAufTX2}]
\bewitemph{(i):}
Mit Lemma~\ref{atIstOpGr} folgt, \dass $(S_t)_\indtGr$ eine Operatorgruppe ist.
\dremark{Analog zu (i).}%

Da $(T_t)_\indtHG$ und $(T_{-t})_\indtHG$ $C_0$-Halbgruppen sind,
findet man nach der Charakterisierung der starken Stetigkeit von Halbgruppen
(Proposition~\ref{charHGIstC0})
ein $\delta \in \mbbR_{>0}$ mit: $\sup_{t \in [-\delta,\delta]} \norm{T_t} < \infty$.
Es gilt nach Proposition~\ref{inAlXNormTEqNormcb}: $\normMl{R}{X} = \norm{R}$ f"ur alle $R \in \Adjlor{X}$.
Also hat man: $\norm{a_t} = \normMl{T_t}{X} = \norm{T_t}$ f"ur alle $\indtGr$.
Es ergibt sich:
\[ \sup_{t \in [-\delta,\delta]} \norm{a_t}
=  \sup_{t \in [-\delta,\delta]} \norm{T_t}
<  \infty. \]
Somit folgt (i) analog zum Beweis von Satz~\ref{halbgrHochliftenAufTX}
und Proposition~\ref{halbgrHochliftenAufTX2}.
Da $T_t$ f"ur alle $\indtGr$ unit"ar ist
und $\mIMl(X)$ und $\Adjlor{X}$ nach Proposition~\ref{AdjlorIsomIMls} als unitale \csalgebren isomorph sind,
ist $a_t$ unit"ar und somit auch $S_t$.
Da nach Proposition~\ref{halbgrHochliftenAufTX2} $(S_t)_\indtHG$ (bzw. $(S_{-t})_\indtHG$)
eine $C_0$-Halbgruppe mit Erzeuger $B$ (bzw. $-B$) ist,
ist nach Satz~\ref{satzHilleYosidaGr} (Satz von Hille-Yosida f"ur Gruppen)
$B$ Erzeuger der $C_0$-Gruppe $(S_t)_\indtGr$.

Die restlichen Aussagen ergeben sich mit Proposition~\ref{halbgrHochliftenAufTX2}.
\end{proof}

Mit Satz~\ref{grHochliftenAufTX} erh"alt man:

\begin{satz}\label{AdjloruHochliftenAufTX}
Sei $X$ ein Operatorraum und $A \in \Adjloru{X}$.
Dann findet man ein eindeutig bestimmtes $B \in \Adjloru{\ternh{X}}$ mit
$j(D(A)) \subseteq D(B)$ und $j \circ A = B \circ j\restring_{D(A)}$.
\dremark{Gelten weitere Aussagen analog zu Aussagen "uber $\Adjlco{X}$?}%
\end{satz}

\begin{bemerkung}\label{AdjloruHochliftenAufTX2}
Seien $X$, $A$ und $B$ wie in Satz~\ref{AdjloruHochliftenAufTX}.
\begin{enumaufz}
\item Es gilt:
\begin{align*}
   &\phantom{=}\,\,\, B\bigl( j(x_1) \multis j(x_2)^* \multis j(x_3) \multis j(x_4)^* \cdots j(x_{2k+1}) \bigr)  \\
&= j(A(x_1)) \multis j(x_2)^* \multis j(x_3) \multis j(x_4)^* \cdots j(x_{2k+1})
\end{align*}
f"ur alle $k \in \mbbN_0$, $x_1 \in D(A)$ und $x_2, \dots, x_{2k+1} \in X$.

\item Es ist $W := W_{D(A),X}$
ein wesentlicher Bereich f"ur~$B$, \dass hei"st, $\overline{B\restring_W} = B$.
\end{enumaufz}
\end{bemerkung}

Zum Beweis formulieren wir das folgende

\begin{lemma}\label{lemmaAdjloruSIstEind}
Sei $X$ ein Operatorraum und $T \in \Adjloru{X}$.
Sei $S : D(S) \subseteq X \to X$ derart,
\dass gilt: $\tilde{T} := \mri \cmpmatrix{0 & T \\ S & 0} \in \Adjlco{C_2(X)}$.
Dann findet man ein eindeutig bestimmtes $C \in \Adjloru{\ternh{X}} = \Multwor{\ternh{X}}$ mit
den Eigenschaften:
\begin{gather*}
\tilde{C} := \mri \cmpmatrix{ 0 & C \\ C^* & 0 } \in \Adjlco{\ternh{C_2(X)}},  \\
j(D(\tilde{T})) \subseteq D(\tilde{C}) \quad\text{und}\quad
j \circ \tilde{T} = \tilde{C} \circ j\restring_{D(\tilde{T})}.
\end{gather*}
\end{lemma}

\begin{proof}
\dmarginpar{pr}%
Nach Satz~\ref{grHochliftenAufTX} findet man
ein eindeutig bestimmtes Element $\tilde{C} \in \Adjlco{\ternh{C_2(X)}}$ mit
$j(D(\tilde{T})) \subseteq D(\tilde{C})$ und
$j \circ \tilde{T} = \tilde{C} \circ j\restring_{D(\tilde{T})}$.
Nach Proposition~\ref{grHochliftenAufTX2} gilt
\begin{align*}
&\phantom{=\,}\,\,\,  \tilde{C}\bigl( j(x_1) \multis j(x_2)^* \multis j(x_3) \multis j(x_4)^*
   \cdots j(x_{2k+1}) \bigr)  \\
&=\,  j(\tilde{T}(x_1)) \multis j(x_2)^* \multis j(x_3) \multis j(x_4)^* \cdots j(x_{2k+1})
\end{align*}
f"ur alle $k \in \mbbN_0$, $x_1 \in D(\tilde{T})$ und $x_2, \dots, x_{2k+1} \in C_2(X)$,
und au"serdem ist $W := W_{D(\tilde{T}),C_2(X)}$
ein wesentlicher Bereich f"ur $\tilde{C}$.
Es gilt nach \cite{BlecherLeMerdy04OpAlg}, 8.3.12.(4): $\ternh{C_2(X)} \cong C_2(\ternh{X})$.
Somit folgt, \dass man Operatoren
\[ C_1 : D(C_1) \subseteq \ternh{X} \to \ternh{X} \quad\text{und}\quad
   C_2 : D(C_2) \subseteq \ternh{X} \to \ternh{X} \]
findet mit $\tilde{C} = \mri \cmpmatrix{ 0 & C_1 \\ C_2 & 0}$.
\dremark{Z. 22.10.'08/1}%
Es ist $\ternh{X}$ ein \hcsmodul und
$-\mri\tilde{C} \in \Adjlco{C_2(\ternh{X})}$,
also $C_1 \in \Adjloru{\ternh{X}} = \Multwor{\ternh{X}}$ nach Satz~\ref{zshgRegOpUnbeschrMult}.
Mit diesem Satz folgt au"serdem: $C_2 = C_1^*$.

Aus der Eindeutigkeit von $\tilde{C}$ folgt die Eindeutigkeit von $C_1$.
\end{proof}

\begin{proof}[Beweis von Satz~\ref{AdjloruHochliftenAufTX} und   Proposition~\ref{AdjloruHochliftenAufTX2}]
\dmarginpar{pr}%
Man findet eine Abbildung $S : D(S) \subseteq X \to X$ derart, \dass gilt:
$\tilde{A} := \mri \cmpmatrix{0 & A \\ S & 0} \in \Adjlco{C_2(X)}$.
Nach Lemma~\ref{lemmaAdjloruSIstEind} gibt es
ein eindeutig bestimmtes $B \in \Adjloru{\ternh{X}} = \Multwor{\ternh{X}}$ mit
den Eigenschaften:
\begin{gather*}
\tilde{B} := \mri \cmpmatrix{ 0 & B \\ B^* & 0 } \in \Adjlco{\ternh{C_2(X)}},  \\
j(D(\tilde{A})) \subseteq D(\tilde{B}) \quad\text{und}\quad
j \circ \tilde{A} = \tilde{B} \circ j\restring_{D(\tilde{A})}.
\end{gather*}
Es folgt: $j(D(A)) \subseteq D(B)$
\dremark{denn $D(\tilde{A}) = D(S) \times D(A)$, $D(\tilde{B}) = D(B^*) \times D(B)$}%
und $j(Ax) = B(j(x))$ f"ur alle $x \in D(A)$.
\dremark{denn $j_{2,1}\binom{\mri Ay}{\mri Sx}
  =  j \circ \tilde{A} \binom{x}{y}
  =  \tilde{B}\left( j\binom{x}{y} \right)
  =  \tilde{B}\left( \binom{j(x)}{j(y)} \right)
  =  \binom{\mri B(j(y))}{\mri B^*(j(x))}.$}%

Weiter gilt nach Proposition~\ref{grHochliftenAufTX2}:
\begin{equation*}
\begin{split}
&\phantom{=\,}\,\,\,  \tilde{B}\bigl( j(x_1) \multis j(x_2)^* \multis j(x_3) \multis j(x_4)^*
   \cdots j(x_{2k+1}) \bigr)  \\
&=\,  j(\tilde{A}(x_1)) \multis j(x_2)^* \multis j(x_3) \multis j(x_4)^* \cdots j(x_{2k+1})
\end{split}
\end{equation*}
f"ur alle $k \in \mbbN_0$, $x_1 \in D(\tilde{A})$ und $x_2, \dots, x_{2k+1} \in C_2(X)$.
Hiermit erh"alt man:
\begin{equation}\label{eqAdjloruHochliftBjx1}
\begin{split}
&\phantom{=\,}   \,\,\,B\bigl( j(x_1) \multis j(x_2)^* \multis j(x_3) \multis j(x_4)^* \cdots j(x_{2k+1}) \bigr)  \\
&=\, j(A(x_1)) \multis j(x_2)^* \multis j(x_3) \multis j(x_4)^* \cdots j(x_{2k+1})
\end{split}
\end{equation}
f"ur alle $k \in \mbbN_0$, $x_1 \in D(A)$ und $x_2, \dots, x_{2k+1} \in X$.

Da $W_{D(\tilde{A}),C_2(X)}$
nach Proposition~\ref{grHochliftenAufTX2} ein wesentlicher Bereich f"ur $\tilde{B}$ ist,
ist $W_{D(A),X}$ ein wesentlicher Bereich f"ur $B$.
\dremark{Z. 3.2.'09/1}%

Um die Eindeutigkeit von $B$ zu zeigen, geben wir uns
ein $C \in \Adjloru{\ternh{X}}$ vor mit
$j(D(A)) \subseteq D(C)$ und $j \circ A = C \circ j\restring_{D(A)}$.
Mit \eqref{eqAdjloruHochliftBjx1} folgt:
\[ B\restring_W = C\restring_W. \]
Es ist $W$ ein wesentlicher Bereich f"ur $B$ und f"ur $C$.
Weiter ist nach \cite{WernerFunkana6}, VII.5.35,
$B$ gleich der Abschlie"sung von $B\restring_W$, also gleich dem Abschlu\cms{}
des Graphen von $B\restring_W$, notiert als $\overline{B\restring_W}$.
Man erh"alt:
\[ B = \overline{B\restring_W} = \overline{C\restring_W} = C.  \qedhere \]
\end{proof}

Sei $T \in \Adjloru{X}$.
Dann findet man nach Definition~\ref{defAdjloru} ein $S : D(S) \subseteq X \to X$ derart,
\dass gilt: $\mri \cmpmatrix{0 & T \\ S & 0} \in \Adjlco{C_2(X)}$.
Die Abbildung $S$ ist eindeutig:

\begin{bemerkung}\label{AdjloruSIstEind}
Sei $X$ ein Operatorraum und $T \in \Adjloru{X}$.
Seien $S_1 : D(S_1) \subseteq X \to X$ und $S_2 : D(S_2) \subseteq X \to X$ mit der Eigenschaft,
\dass $\mri \cmpmatrix{0 & T \\ S_1 & 0}, \mri \cmpmatrix{0 & T \\ S_2 & 0} \in \Adjlco{C_2(X)}$ sind.
Dann gilt: $S_1 = S_2$.
\end{bemerkung}

\begin{proof}
\dmarginpar{pr}%
Setze $\hat{S}_1 := \cmpmatrix{0 & T \\ S_1 & 0}$
und $\hat{S}_2 := \cmpmatrix{0 & T \\ S_2 & 0}$.
Es gilt nach \cite{BlecherLeMerdy04OpAlg}, 8.3.12.(4): $\ternh{C_2(X)} \cong C_2(\ternh{X})$.
Nach Lemma~\ref{lemmaAdjloruSIstEind} findet man  $C_1, D_1 \in \Multwor{\ternh{X}}$ mit
den Eigenschaften:
\begin{gather*}
\tilde{C} := \mri \cmpmatrix{ 0 & C_1 \\ C_1^* & 0 },
 \tilde{D} := \mri \cmpmatrix{ 0 & D_1 \\ D_1^* & 0 } \in \Adjlco{\ternh{C_2(X)}}, \\
j(D(\hat{S}_1)) \subseteq D(\tilde{C})  \quad\text{und}\quad
j(D(\hat{S}_2)) \subseteq D(\tilde{D}),  \\
j \circ \mri \hat{S}_1 = \tilde{C} \circ j\restring_{D(\hat{S}_1)}  \quad\text{und}\quad
j \circ \mri \hat{S}_2 = \tilde{D} \circ j\restring_{D(\hat{S}_2)}.
\end{gather*}
Setze $\hat{C} := -\mri\tilde{C}$ und $\hat{D} := -\mri\tilde{D}$.

Es gilt nach Proposition~\ref{grHochliftenAufTX2}:
\begin{equation}\label{eqEindShatCEq}
\begin{split}
&\phantom{=\,}\,\,\,  \mri\hat{C}\bigl( j(x_1) \multis j(x_2)^* \multis j(x_3) \multis j(x_4)^*
   \cdots j(x_{2k+1}) \bigr)  \\
&=\,  j(\mri\hat{S}_1(x_1)) \multis j(x_2)^* \multis j(x_3) \multis j(x_4)^* \cdots j(x_{2k+1})
\end{split}
\end{equation}
f"ur alle $k \in \mbbN_0$, $x_1 \in D(\hat{S}_1)$ und $x_2, \dots, x_{2k+1} \in C_2(X)$.
Mit \eqref{eqEindShatCEq} erh"alt man:
\begin{align*}
&\phantom{=\,}   \,\,\,C_1\bigl( j(x_1) \multis j(x_2)^* \multis j(x_3) \multis j(x_4)^* \cdots j(x_{2k+1}) \bigr)  \\
&=\, j(T(x_1)) \multis j(x_2)^* \multis j(x_3) \multis j(x_4)^* \cdots j(x_{2k+1})  \\
&=\, D_1\bigl( j(x_1) \multis j(x_2)^* \multis j(x_3) \multis j(x_4)^* \cdots j(x_{2k+1}) \bigr)
\end{align*}
f"ur alle $k \in \mbbN_0$, $x_1 \in D(T)$ und $x_2, \dots, x_{2k+1} \in X$.

Es folgt:
\[ C_1\restring_{W_{D(T),X}} = D_1\restring_{W_{D(T),X}}. \]
Da $W_{D(\hat{S}_1),C_2(X)}$ (bzw. $W_{D(\hat{S}_2),C_2(X)}$)
nach Proposition~\ref{grHochliftenAufTX2} ein wesentlicher Bereich f"ur $\hat{C}$ (bzw. $\hat{D}$) ist,
ist $W := W_{D(T),X}$ ein wesentlicher Bereich f"ur $C_1$ und f"ur $D_1$.
\dremark{Z. 3.2.'09/1}%
Weiter ist nach \cite{WernerFunkana6},\dremark{Aufgabe} VII.5.35, $C_1$ gleich der Abschlie"sung
$\overline{C_1\restring_W}$ von $C_1\restring_W$.
Hiermit erh"alt man:
\dremark{Alternativ kann man dies auch
  mit \ref{AdjloruHochliftenAufTX} und \ref{AdjloruHochliftenAufTX2} beweisen}%
\[ C_1 = \overline{C_1\restring_W} = \overline{D_1\restring_W} = D_1. \]
Damit folgt $C_2 = C_1^* = D_1^* = D_2$, also $\hat{C} = \hat{D}$.

Sei $(U_t)_\indtGr$ (bzw. $(V_t)_\indtGr$) die von $\mri \hat{S}_1$ (bzw. $\mri \hat{S}_2$)
erzeugte $C_0$-Gruppe auf $C_2(X)$ und
$(S_t)_\indtGr$ die von $\mri \hat{C} = \mri \hat{D}$
erzeugte $C_0$-Gruppe auf $\ternh{C_2(X)}$.
Setze $a_t := \isoIMl^{-1}(U_t)$ und $b_t := \isoIMl^{-1}(V_t)$ f"ur alle $\indtGr$.
Es ergibt sich mit Proposition~\ref{grHochliftenAufTX2}:
\[ a_t \multis j(x)
=  j(U_t x)
\overset{\dremarkm{\ref{grHochliftenAufTX2}}}{=}  S_t(j(x))
\overset{\dremarkm{\ref{grHochliftenAufTX2}}}{=}  j(V_t x)
=  b_t \multis j(x) \]
f"ur alle $x \in X$.
Mit \refb{\cite{BlecherLeMerdy04OpAlg}, Proposition 4.4.12,}{Proposition~\ref{Lcinj}}
folgt: $a_t = b_t$ f"ur alle $\indtGr$.
Man erh"alt\dremark{$U_t = V_t$, also} $\hat{S}_1 = \hat{S}_2$, also $S_1 = S_2$.
%
\dremark{Z. 4.11.'08/1}%
\end{proof}

Nach der obigen Proposition\dremark{\ref{AdjloruSIstEind}} ist die Abbildung $S_1$ eindeutig.
Ist $X$ ein \hcsmodul, so gilt: $S_1 = T^*$ (Satz~\ref{zshgRegOpUnbeschrMult}).
Daher kann man $S_1$ als die Adjungierte von $T$ ansehen
und die entsprechende Bezeichnung einf"uhren:

\begin{definitn}
Sei $X$ ein Operatorraum und $T \in \Adjloru{X}$.
Dann findet man einen eindeutig bestimmten Operator $T^* : D(T^*) \subseteq X \to X$,
genannt \defemphi{Adjungierte} von $T$, derart,
\dass gilt: $\mri \cmpmatrix{ 0 & T \\ T^* & 0 } \in \Adjlco{C_2(X)}$.
\index[S]{Ts@$T^*$}%
\end{definitn}

Die obige Bezeichnung $T^*$ f"uhrt nicht zu Bezeichnungskollisionen:

\begin{bemerkung}
Sei $X$ ein Operatorraum und $T \in \Adjlor{X}$.
Es gilt nach Proposition~\ref{AdjlorSubseteqAdjloru}: $\Adjlor{X} \subseteq \Adjloru{X}$.
Dann stimmt $T^*$ aus der obigen Definition
mit $T^*$, der Linksadjungierten in $\Adjlor{X}$ (siehe S. \pageref{linksadjEind}), "uberein.
\end{bemerkung}

\begin{proof}
Mit $T^\star$ sei die Adjungierte aus der obigen Definition bezeichnet und
mit $T^*$ die Linksadjungierte in $\Adjlor{X}$.
Mit Lemma~\ref{NullTTs0InAdjlor} erh"alt man,
\dass $\tilde{T} := \mri\cmsmallpmatrix{0 & T \\ T^* & 0} \in \Adjlor{C_2(X)}$ schiefadjungiert ist.
Nach Proposition~\ref{AdjlorSubseteqAdjlco} gilt:
$\tilde{T} \in \Adjlor{C_2(X)} \subseteq \Adjlco{C_2(X)}$.
Mit der Eindeutigkeitsaussage aus Proposition~\ref{AdjloruSIstEind} ergibt sich:
$T^* = T^\star$.
\end{proof}

Wir haben gezeigt, wie man unbeschr"ankte Multiplikatoren von $X$ auf die
tern"are H"ulle $\ternh{X}$, also auf einen gr"o"seren Operatorraum, "uberf"uhren kann.
Im folgenden untersuchen wir,
unter welchen Voraussetzungen ein unbeschr"ankter Multiplikator
so auf einen Unteroperatorraum eingeschr"ankt werden kann,
\dass man wieder einen unbeschr"ankten Multiplikator erh"alt.

\begin{bemerkung}\label{ArestringXInMultlco}
Seien $X$, $Y$ Operatorr"aume mit $X \subseteq Y$.
Ferner sei $A \in \Multlco{Y}$ (bzw. $A \in \Adjlco{Y}$) derart,
\dass die von $A$ erzeugte $C_0$"=Halbgruppe (bzw. $C_0$-Gruppe) $X$ invariant l"a\cms{}t.
Dann folgt: $A\restring_X  \in \Multlco{X}$ (bzw. $A\restring_X \in \Adjlco{X}$).
\dremark{Anwendung: $Y$ ist \hcsmodul}%
\end{bemerkung}

\begin{proof}
\bewitemph{(i):}
Sei $(T_t)_\indtHG$ die von $A \in \Multlco{Y}$ erzeugte $C_0$-Halbgruppe.
Nach \cite{EngelNagelSemigroups}, Corollary II.2.3, erzeugt $A\restring_X$
die $C_0$-Halbgruppe $(T_t\restring_X)_\indtHG$.
Betrachte $j = j_Y : Y \to I(Y) \subseteq I(\mcalS(Y))$.
Dann ist $\sigma := j_Y\restring_X$ eine vollst"andige Isometrie in die \csalgebra $I(\mcalS(Y))$.
Sei $\indtHG$.
Setze $a_t := \isoIMl^{-1}(T_t) \in I_{11}(Y) \subseteq I(\mcalS(Y))$,
also $j_Y(T_t x) = a_t \multis j_Y(x)$ f"ur alle $x \in X$.
Es gilt:
\begin{equation}\label{eqsigmaTtXEqatx}
\sigma(T_t\restring_X (x)) = j_Y(T_t x) = a_t \multis j_Y(x) = a_t \multis \sigma(x)
\end{equation}
f"ur alle $x \in X$.
Somit folgt: $T_t\restring_X \in \Multlor{X}$.

\bewitemph{(ii):}
Sei $(T_t)_\indtGr$ die von $A \in \Adjlco{Y}$ erzeugte $C_0$-Gruppe
von unit"aren Elementen aus $\Adjlor{Y}$.
Analog zu (i) ergibt sich, \dass $A\restring_X$ die $C_0$-Gruppe $(T_t\restring_X)_\indtGr$
erzeugt.
Sei $\indtGr$.
Setze $a_t := \isoIMl^{-1}(T_t)$.
Da $T_t$ unit"ar ist, gilt: $T_t^* \dremarkm{=T_t^{-1}} = T_{-t}$.
Somit folgt:
\[ a_t^* \multis \sigma(x)
=  j_Y(T_t^* x)
=  j_Y(T_{-t} x)
\in \dremarkm{j(X) = } \sigma(X) \]
f"ur alle $x \in X$.
Mit \eqref{eqsigmaTtXEqatx} und \cite{Zarikian01Thesis}, Proposition 1.7.7,
erh"alt man: $T_t\restring_X \in \Adjlor{X}$.
Da $a_t$ unit"ar ist
und $e_{I_{11}(Y)} \multis j_Y(x) = j_Y(x)$ f"ur alle $x \in X$ gilt (siehe \eqref{eqeI11xEqx}),
ist auch $T_t\restring_X$ unit"ar.
\dremark{$j((T_t\restring_X)^* T_t\restring_X (x))
  = a_t^* \multis a_t \multis j(x)
  = e_{I_{11}(Y)} \multis j_Y(x)
  \overset{\eqref{eqeI11xEqx}}{=}  j(x)$.
  Da $j$ injektiv ist, folgt: $(T_t\restring_X)^* T_t\restring_X = \Id_X$.}%
Insgesamt folgt: $A\restring_X \in \Adjlco{X}$.
\dremark{Z. 3.4.'09}%
\end{proof}

Mit der obigen Proposition\dremark{\ref{ArestringXInMultlco}} erh"alt man
ein entsprechendes Resultat f"ur die unbeschr"ankten Multiplikatoren:

\begin{bemerkung}\label{ArestringXInAdjloru}
Seien $X$, $Y$ Operatorr"aume mit $X \subseteq Y$.
Weiter sei $A \in \Adjloru{Y}$ derart,
\dass die von $\mri \cmsmallpmatrix{0 & A \\ A^* & 0}$
erzeugte $C_0$-Gruppe $C_2(X)$ invariant l"a\cms{}t.
Dann folgt: $A\restring_X \in \Adjloru{X}$.
\end{bemerkung}

\dremww{Beweis.
Es gilt: $\tilde{A} := \mri \cmsmallpmatrix{0 & A \\ A^* & 0} \in \Adjlco{C_2(Y)}$.
Mit Proposition~\ref{ArestringXInMultlco} folgt:
$\tilde{A}\restring_{C_2(X)} \in \Adjlco{C_2(X)}$.
Man erh"alt: $A\restring_X \in \Adjloru{X}$.
}%

In den beiden vorherigen Resultaten wurde gefordert,
\dass die erzeugte \cohg den Unterraum $X$ invariant l"a\cms{}t.
Auf diese Forderung kann man verzichten,
wenn man zus"atzliche Voraussetzungen an den Erzeuger stellt:

\begin{bemerkung}\label{ArestringXInMultlcoLP}
Seien  $X$, $Y$ Operatorr"aume mit $X \subseteq Y$.
Ferner sei $A \in \Multlco{Y}$ derart,
\dass die von $A$ erzeugte $C_0$-Halbgruppe $(T_t)_\indtHG$ erf"ullt:
$\normMl{T_t}{Y} \leq 1$ f"ur alle $\indtHG$.
Es existiere ein $\lambda \in \mbbR_{>0}$ so, \dass $\lambda - A\restring_X$ surjektiv auf $X$ ist.
Dann gilt $A\restring_X \dremarkm{= \erzur{A}{X}} \in \Multlco{X}$, und
$A\restring_X$ erzeugt die $C_0$"=Halbgruppe $(T_t\restring_X)_\indtHG$.
\dremark{Mit \cite{EngelNagelSemigroups}, Corollary II.2.3, folgt:
  $\erzur{A}{X} = A\restring_X$.}%
\end{bemerkung}

\begin{proof}
Nach dem Charakterisierungssatz von $\Multlco{Y}$ mittels $\matnull{A}$
(Satz \ref{charMultlcomatnull}) erzeugt $\matnull{A}$ eine \cmvk{e}
$C_0$-Halbgruppe auf $C_2(Y)$.
Definiere
\[ \erzur{A}{X} := \{ (x,y) \setfdg x \in D(A) \cap X, y = Ax \in X \}. \]
\index[S]{AX@$\erzur{A}{X}$}%
\dremark{Beachte: $A\restring_X$ surjektiv auf $X$ $\iff$ $\erzur{A}{X}$ surjektiv auf $X$.}%
Es sind $\matnull{A}$ und damit auch $\matnull{\erzur{A}{X}}$
nach dem Satz von Lumer-Phillips in Operatorr"aumen (Satz~\ref{SatzLumerPhillipsOR})
vollst"andig dissipativ.
Da $\lambda - \matnull{\erzur{A}{X}}$ surjektiv auf $C_2(X)$ ist,
folgt mit Satz~\ref{SatzLumerPhillipsOR},
\dass $\matnull{\erzur{A}{X}}$ eine \cmvk{e} $C_0$"=Halbgruppe erzeugt.
Es ergibt sich mit Satz~\ref{charMultlcomatnull}: $\erzur{A}{X} \in \Multlco{X}$.

Sei $(T_t)_\indtHG$ die von $A$ erzeugte $C_0$"=Halbgruppe.
Weil die von $\erzur{A}{X}$ erzeugte $C_0$"=Halbgruppe\dremark{$(S_t)_\indtHG$}
mit $(T_t\restring_X)_\indtHG$ "ubereinstimmt,
erh"alt man mit \cite{EngelNagelSemigroups}, Corollary II.2.3:
$\erzur{A}{X} = A\restring_X$.
\dremark{Es gilt $\erzur{A}{X} \subseteq A$,
  also $\lambda - \erzur{A}{X} \subseteq \lambda - A$.
  Mit \cite{UllmannDiplom}, Satz 2.8.(a), folgt:
  $(\lambda - \erzur{A}{X})^{-1} \subseteq (\lambda - A)^{-1}$.
  F"ur die Yosida"=Approximation folgt:
  $ (\erzur{A}{X})_\lambda
  = \lambda^2 (\lambda - \erzur{A}{X})^{-1} - \lambda
  \subseteq \lambda^2 (\lambda - A)^{-1} - \lambda
  =  A_\lambda$.
  Es folgt
  $ T_t x
  = \lim_{\lambda \to \infty} \mre^{t A_\lambda} x
  = \lim_{\lambda \to \infty} \mre^{t (\erzur{A}{X})_\lambda} x
  = S_t x$ f"ur alle $x \in X$,
  also $S_t = T_t\restring_X$.
  Somit l"a\cms{}t $T_t$ $X$ invariant.}%
\dremark{Z. 11.4.'09/2}%
\end{proof}

Eine entsprechende Version der obigen Proposition\dremark{\ref{ArestringXInMultlcoLP}}
f"ur $\Adjlco{X}$ beweist man analog und erh"alt:

\begin{bemerkung}\label{ArestringXInAdjlcoLP}
Seien  $X$, $Y$ Operatorr"aume mit $X \subseteq Y$.
Weiter sei $A \in \Adjlco{Y}$.
Es gelte: Es existieren $\lambda,\mu \in \mbbR_{>0}$ so,
\dass $\lambda - A\restring_X$ und $\mu + A\restring_X$ surjektiv auf $X$ sind.
Dann gilt: $A\restring_X \dremarkm{= \erzur{A}{X}} \in \Adjlco{X}$.
\end{bemerkung}

Mit der obigen Aussage\dremark{\ref{ArestringXInAdjlcoLP}} ergibt sich:

\begin{bemerkung}
Seien  $X$, $Y$ Operatorr"aume mit $X \subseteq Y$.
Ferner sei $A \in \Adjloru{Y}$.
Setze $\tilde{A} := \mri \cmsmallpmatrix{0 & A \\ A^* & 0}$.
Es gelte: Es existieren $\lambda,\mu \in \mbbR_{>0}$ so,
\dass $\lambda - \tilde{A}\restring_{C_2(X)}$ und
$\mu + \tilde{A}\restring_{C_2(X)}$ surjektiv auf $C_2(X)$ sind.
Dann gilt: $A\restring_X \in \Adjloru{X}$.
\end{bemerkung}

\dremark{Beweis.
Mit \ref{ArestringXInAdjlcoLP} folgt:
$ \mri \cmpmatrix{0 & A\restring_X \\ A^*\restring_X & 0}
= \tilde{A}\restring_{C_2(X)}
= \erzur{\tilde{A}}{C_2(X)}
\in \Adjlco{C_2(X)}$,
also $A\restring_X \in \Adjloru{X}$.
}%

Man kann jeden unbeschr"ankten schiefadjungierten Multiplikator
als Einschr"ankung eines geeigneten regul"aren Operators auf\/fassen:

\begin{satz}\label{AdjlcoEqRegOpAufTX}
Sei $X$ ein Operatorraum und $A : D(A) \subseteq X \to X$ linear.
Dann sind die folgenden Aussagen "aquivalent:
\begin{enumaequiv}
\item $A \in \Adjlco{X}$.

\item Es existiert ein $B \in \Multwor{\ternh{X}}$ schiefadjungiert mit den Eigenschaften:
\begin{enumaufz}
\item $j(D(A)) \subseteq D(B)$,
\item $j \circ A = B \circ j\restring_{D(A)}$ und
\item es existieren $\lambda, \mu \in \mbbR_{>0}$ so,
  \dass $\lambda - B\restring_{j(X)}$ und $\mu + B\restring_{j(X)}$ surjektiv auf $j(X)$ sind.
\end{enumaufz}
\end{enumaequiv}
\end{satz}

Im obigen Satz kann man Punkt (iii) unter (b) auch durch eine der folgenden beiden
Bedingungen ersetzen:
\begin{enumerate}
\item[$(iii)'$] Es existieren $\lambda, \mu \in \mbbR_{>0}$ so,
  \dass $\lambda - A$ und $\mu + A$ surjektiv auf $X$ sind.
\item[$(iii)''$] Die von $B$ erzeugte $C_0$-Gruppe l"a\cms{}t $j(X)$ invariant.
\dremark{Folgt mit Proposition~\ref{ArestringXInMultlco}.}%
\end{enumerate}

\begin{proof}
\bewitemphq{(a)$\Rightarrow$(b)}
Es gelte (a).
Man findet nach Satz~\ref{grHochliftenAufTX} einen Operator $B \in \Adjlco{\ternh{X}}$ mit
$j(D(A)) \subseteq D(B)$ und $j \circ A = B \circ j\restring_{D(A)}$.
Mit Satz~\ref{MultworInAdjlco} erh"alt man: $B \in \Multwor{\ternh{X}}$ mit $B^* = -B$.

Nach Proposition~\ref{inAlXNormTEqNormcb} gilt:
$\normMl{T}{X} = \norm{T}$ f"ur alle $T \in \Adjlor{X}$.
Da $A$ Erzeuger einer $C_0$-Gruppe von unit"aren Elementen von $\Adjlor{X}$ ist,
folgt somit, \dass $A$ Erzeuger einer kontraktiven $C_0$-Gruppe ist.
Da $A$ und $-A$ jeweils Erzeuger einer kontraktiven $C_0$-Halbgruppe sind,\dremark{\ref{satzHilleYosidaGr}}
findet man nach dem Satz von Lumer-Phillips (Satz~\ref{SatzLumerPhillips})
$\lambda,\mu \in \mbbR_{>0}$ derart,
\dass $\lambda - A$ und $\mu + A$ surjektiv auf $X$ sind.
Somit ergibt sich Punkt (iii) von (b).
\smallskip

\bewitemphq{(b)$\Rightarrow$(a)}
Es gelte (b).
Mit Satz~\ref{MultworInAdjlco} erh"alt man: $B \in \Adjlco{\ternh{X}}$.
Da $\lambda - B\restring_{j(X)}$ und $\mu + B\restring_{j(X)}$ surjektiv auf $j(X)$ sind,
ergibt sich mit Proposition~\ref{ArestringXInAdjlcoLP}:
$B\restring_{j(X)} \in \Adjlco{j(X)}$.
Es folgt: $A \in \Adjlco{X}$.
\end{proof}

Analog zum obigen Satz beweist man ein Analogon f"ur $\Adjloru{X}$:

\pagebreak
\begin{satz}
Sei $X$ ein Operatorraum und $A : D(A) \subseteq X \to X$ linear.
Dann sind die folgenden Aussagen "aquivalent:
\begin{enumaequiv}
\item $A \in \Adjloru{X}$.

\item Es existiert ein $B \in \Multwor{\ternh{X}}$ mit den Eigenschaften:
\begin{enumaufz}
\item $j(D(A)) \subseteq D(B)$,
\item $j \circ A = B \circ j\restring_{D(A)}$ und
\item es existieren $\lambda, \mu \in \mbbR_{>0}$ derart,
  \dass $\lambda - \mri\cmsmallpmatrix{0 & B \\ B^* & 0}\restring_{j_{2,1}(C_2(X))}$ und
  $\mu + \mri\cmsmallpmatrix{0 & B \\ B^* & 0}\restring_{j_{2,1}(C_2(X))}$
  surjektiv auf $j_{2,1}(C_2(X))$ sind.
\end{enumaufz}
\end{enumaequiv}
\end{satz}

Im obigen Satz kann man Punkt (iii) unter (b) auch durch eine der folgenden beiden
Bedingungen ersetzen:
\begin{enumerate}
\item[$(iii)'$] Es existieren $\lambda, \mu \in \mbbR_{>0}$
  und ein $C : D(C) \subseteq X \to X$ so,
  \dass $\lambda - \mri\cmsmallpmatrix{0 & A \\ C & 0}$ und
  $\mu + \mri\cmsmallpmatrix{0 & A \\ C & 0}$
  surjektiv auf $C_2(X)$ sind.
\item[$(iii)''$] Die von $\mri\cmsmallpmatrix{0 & B \\ B^* & 0}$
  erzeugte $C_0$-Gruppe l"a\cms{}t $j_{2,1}(C_2(X))$ invariant.
\end{enumerate}

Mit Hilfe von Proposition~\ref{ArestringXInAdjlcoLP} zeigt man,
\dass die Adjungierte eines Elementes aus $\Adjloru{X}$ die folgenden Eigenschaften besitzt:

\begin{bemerkung}
Sei $X$ ein Operatorraum und $A \in \Adjloru{X}$.
Dann gilt: $A^* \in \Adjloru{X}$ und ${(A^*)}^* = A$.
\end{bemerkung}

\begin{proof}
Es gilt: $\tilde{A} := \mri \cmpmatrix{0 & A \\ A^* & 0} \in \Adjlco{C_2(X)}$.
Nach Lemma~\ref{lemmaAdjloruSIstEind} findet man ein
$C \in \Adjloru{\ternh{X}} = \Multwor{\ternh{X}}$ mit:
\begin{gather}
\tilde{C} := \mri \cmpmatrix{ 0 & C \\ C^* & 0 } \in \Adjlco{\ternh{C_2(X)}},  \nonumber \\
j(D(\tilde{A})) \subseteq D(\tilde{C}) \quad\text{und}\quad
j \circ \tilde{A} = \tilde{C} \circ j\restring_{D(\tilde{A})}.  \label{eqAssEqAjAEq}
\end{gather}
Weiter hat man nach \cite{LanceHmod}, Corollary 9.6: $C^* \in \Multwor{\ternh{X}}$.
Da nach \refb{\cite{LanceHmod}, Corollary 9.4,}{Proposition \ref{TssEqT}} gilt $C^{**} = C$,
erh"alt man mit Proposition~\ref{NullTTs0Selbstadj}, \dass
$\tilde{B} := \mri \cmpmatrix{0 & C^* \\ C & 0} \in \Multwor{\ternh{X} \oplus \ternh{X}}$
schiefadjungiert ist.
Mit Proposition~\ref{normHCsModEoplusEEqC2E} ergibt sich:
$\ternh{X} \oplus \ternh{X} \cong C_2(\ternh{X}$.
Es folgt mit Satz~\ref{MultworInAdjlco}: $\tilde{B} \in \Adjlco{C_2(\ternh{X})}$.

Nach Proposition~\ref{charAdjlormatnull2} findet man $\lambda, \mu \in \mbbR_{>0}$
so, \dass $\lambda - \tilde{A}$ und $\mu + \tilde{A}$ surjektiv auf $C_2(X)$ sind.
Setze $Y := j_{2,1}(C_2(X))$.
Somit sind $\lambda - \tilde{C}\restring_Y$ und
$\mu + \tilde{C}\restring_Y$ nach \eqref{eqAssEqAjAEq}
surjektiv auf $Y$,
also auch $\lambda - \tilde{B}\restring_Y$ und
$\mu + \tilde{B}\restring_Y$.
\dremark{denn
  $\tilde{C}\restring_{j(D\tilde{A})} = j \circ \tilde{A} \circ j\restring_{j(D(\tilde{A}))}^{-1}$,
  also $\mri C\restring_{j(D(A))} = j \circ \mri A \circ j\restring_{j(D(A))}^{-1}$,
  setze $\check{A} := \mri \cmpmatrix{0 & A^* \\ A & 0}$,
  also $\tilde{B}\restring_{j(D(\check{A}))} = j \circ \check{A} \circ j\restring_{j(D(\check{A}))}^{-1}$.}%
Mit Proposition~\ref{ArestringXInAdjlcoLP} erh"alt man
$\tilde{B}\restring_Y \in \Adjlco{Y}$,
also $\mri \cmpmatrix{0 & A^* \\ A & 0} \in \Adjlco{C_2(X)}$.
Es folgt: $A^* \in \Adjloru{X}$ und ${(A^*)}^* = A$.\dremark{mit \ref{AdjloruSIstEind}}
\dremark{Z. 24.4.'09/1}%
\end{proof}

\dremww{
\begin{bemerkung}
Sei $X$ ein Operatorraum, $D$ eine dichte Teilmenge von $X$, $H$ ein Hilbertraum und
$\sigma : X \to \mLinStet(H)$ linear und stetig mit $\sigma(X)H$ dicht in $H$.
Dann liegt $\sigma(D)(H)$ dicht in $H$.
\dremark{Typische Anwendung:
  $A : D(A) \subseteq X \to X$ linear und dicht definiert, $D = D(A)$.}%
\end{bemerkung}

\begin{bemerkung}
Sei $X$ ein Operatorraum, $D$ eine dichte Teilmenge von $X$.
Dann gilt: $\troerzi{j(M)}{\ternh{X}} = \ternh{X}$,
wobei mit $\troerzi{j(M)}{\ternh{X}}$ der von $j(M)$ in $\ternh{X}$
erzeugte TRO bezeichnet sei.\dremark{23.10.'07}
\end{bemerkung}

\begin{anmerkung}
Anderer Ansatz zur Definition eines unbeschr"ankten Multiplikators auf $L^p(\mcalM)$
f"ur eine beliebige semiendliche Von-Neumann-Algebra $(\mcalM, \tau)$,
die auf einem Hilbertraum $H$ operiert
(siehe \cite{JungeLeMerdyXu06HinfFuncCalc}, S. 69 oben und (8.3)):

F"ur jedes $a \in \mcalM$ definiere den beschr"ankten Operator
\begin{equation}\label{eqDefLa}
\mcalL_a : L^p(\mcalM) \to L^p(\mcalM), x \mapsto ax.
\end{equation}
Sei $a : D(a) \subseteq H \to H$ ein dicht definierter, abgeschlossener Operator auf $H$.
Es wird nicht \eqref{eqDefLa} direkt benutzt, um $\mcalL_a$ zu definieren,
da die Multiplikation mit dem unbeschr"ankten Operator $a$ zu technischen Problemen f"uhrt.
Stattdessen wird mit der Resolvente multipliziert.

Sei $z \in \rho(a)$.
Definiere
\[ \mcalL_a := z - \mcalL^{-1}_{R(z,a)} : \mcalL_{R(z,a)}(L^p(\mcalM)) \to L^p(\mcalM). \]
\end{anmerkung}
}%

\section{Unbeschr"ankte Multiplikatoren auf einen Hilbertraum "uberf"uhren}

In diesem Abschnitt beweisen wir, wie man f"ur einen beliebigen Operatorraum
$X$ Elemente aus $\Multlco{X}$ (bzw. $\Adjlco{X}$)
und die von ihnen erzeugte $C_0$-Halbgruppen (bzw. $C_0$-Gruppen)
auf einen Hilbertraum $H$ "uberf"uhren kann.
Dies erm"oglicht es, eine Charakterisierung
der Elemente von $\Multlco{X}$, $\Adjlco{X}$ und $\Adjloru{X}$ zu zeigen.
\dmarginpar{"U}\dremark{Evtl. an einige Begriffe erinnern ($j$, $\isoIMl$, \dots)}%
\skiptext

Das folgende Beispiel dient der Motivation und legt dar,
wie man eine bestimmte Art von $C_0$-Halbgruppen, sogenannte Multiplikationshalbgruppen,
von dem Operatorraum $C_0(\Omega)$ auf den Hilbertraum $\ell^2(\Omega)$ "uberf"uhren kann:

\begin{beispiel}
Sei $\Omega$ ein lokalkompakter Hausdorffraum und $f \in C(\Omega)$ mit
$\sup_{\omega \in \Omega} \Re f(\omega) < \infty$.
Setze $X := C_0(\Omega)$, $H := \ell^2(\Omega)$,
$T_t : X \to X, g \mapsto \mre^{tf} \cdot g$ und
$S_t : H \to H, h \mapsto \mre^{tf} \cdot h$ f"ur alle $\indtHG$.
\dremark{$f$ ist me\cms{}bar, also ist $(S_t)_\indtHG$ nach
  \cite{EngelNagelSemigroups}, Prop. I.4.11, eine $C_0$-Halbgruppe.}%
Dann gilt:
\begin{enumaufz}
\item $(T_t)_\indtHG$ (bzw. $(S_t)_\indtHG$) ist eine $C_0$-Halbgruppe auf $X$ (bzw. $H$),
die sogenannte \defemphi{Multiplikationshalbgruppe}, mit Erzeuger
\begin{align*}
&M_f : D(M_f) \subseteq X \to X, g \mapsto fg  \\
(\text{bzw. } &\tilde{M}_f : D(\tilde{M}_f) \subseteq H \to H, h \mapsto fh).
\end{align*}

\item $\sigma : X \to \mLinStet(H), g \mapsto (h \mapsto gh)$, ist eine nicht-ausgeartete
\sterns{}Darstellung von $X$ mit der Eigenschaft, \dass
$[\sigma(D(M_f))H]$ dicht in $H$ liegt.
\item $\sigma(M_f(g)) = \tilde{M}_f \circ (\sigma(g))$ f"ur alle $g \in D(M_f)$.
\dmarginpar{zutun}\dremark{Was hat das Beispiel mit dem "Uberf"uhren zu tun?}%
\dremark{Beispiel zur Motivation}%
\end{enumaufz}
\end{beispiel}

\begin{proof}
\bewitemph{(i)} folgt mit \cite{EngelNagelSemigroups}, Lemma II.2.9.

\bewitemph{(ii) und (iii)} rechnet man leicht nach.
\dmarginpar{pr}\dremark{Beweis notieren}%
\dremark{\bewitemph{(ii):} Sei $g \in X$.
Zeige:
\[ \Bild(\sigma(g)) \subseteq H. \dremarkm{\qquad(\dagger)} \]
Sei $h \in H$, $\varepsilon \in \mbbR_{>0}$.
Setze $M := \norm{g}_\infty$.
Es gilt
$  \sum_{\omega \in \Omega} \abs{g(\omega) h(\omega)}^2
\leq  \sum_{\omega \in \Omega} M^2 \abs{h(\omega)}^2
< \infty$,
also $gh \in H$.

Zeige: $\sigma$ ist \sterns{}Darstellung.
Es gilt:
$  \skalpr{\sigma(g)^* h}{\xi}
=  \skalpr{h}{\sigma(g) \xi}
=  \sum h(\omega) \overline{g}(\omega) \overline{\xi}(\omega)
=  \skalpr{\overline{g} h}{\xi}
=  \skalpr{\sigma(\overline{g}) h}{\xi}$.
Da $\skalpr{\cdot}{\cdot\cdot}$ nicht-ausgeartet ist,
folgt: $\sigma(g)^* = \sigma(\overline{g})$.

Zeige: $\erzl{\sigma(D(M_f))H}$ dicht in $H$.
Nach \cite{EngelNagelSemigroups}, Beweis von Proposition I.4.2,
ist $C_c(\Omega)$ in $D(M_f)$ enthalten.
Sei $g \in H$, $\varepsilon \in \mbbR_{>0}$.
Nach \cite{PedersenAnalysisNowRev}, Proposition 6.4.11, ist das Bild von $C_c(\Omega)$
dicht in $\ell^2(\Omega)$.
Somit findet man ein $h \in C_c(\Omega)$ mit: $\norm{h - g}_2 < \frac{\varepsilon}{\norm{g}_2}$.
Es gilt: $\norm{hg - g^2}_2 = \norm{g}_2 \norm{h-g}_2 < \varepsilon$.
Somit ist $C_c(\Omega) H$ dicht in $H \cdot H$.
Da $H \cdot H$ dicht in $H$ liegt, folgt die Behauptung.

\bewitemph{(iii):} Sei $g \in D(M_f)$ und $h \in H$.
Es gilt $fg \in X$, also nach $(\dagger)$: $fgh \in H$.
Es folgt:
\dremark{Z. 7.2.'09, S. 1/2; Z. 1.10.'07, S. 2/3}%
\[ \sigma(M_f(g))(h)
=  \sigma(fg)(h)
=  fgh
=  \tilde{M}_f(gh)
=  \tilde{M}_f(\sigma(g)(h))
=  (\tilde{M}_f \circ (\sigma(g))(h).  \qedhere \]
}%
\end{proof}

Im folgenden Lemma halten wir eine speziell an $\ternh{X}$
angepa\cms{}te Einbettung fest, die in den darauf\/folgenden S"atzen verwendet wird:

\pagebreak
\begin{lemma}
Sei $X$ ein Operatorraum.
Dann findet man einen Hilbertraum $K$,
eine unitale, injektive \sterns{}Darstellung $\pi : I(\mcalS(X)) \to \mLinStet(K)$
und abgeschlossene Unterr"aume $K_1$, $K_2$ von $K$ derart,
\dass gilt:
\begin{enumaufz}
\item $[\pi(\ternh{X})K_1]$ liegt dicht in $K_2$,
\item $[\pi(\ternh{X})K_1^\perp] = \{0\}$ und
\item $\erzl{\pi(j(X)) K} \subseteq K_2$.
\end{enumaufz}
Eine solches Quadrupel $(\pi, K, K_1, K_2)$ wird \defemph{$\ternh{X}$-Einbettung} von $X$ genannt.
\dremark{Evtl. etwas dazu sagen, wie man diese Hilbertr"aume findet.}%
\index[B]{TXEinbettung@$\ternh{X}$-Einbettung}%
\end{lemma}

\begin{proof}
Man findet einen Hilbertraum $K$ und
eine unitale, injektive \sterns{}Dar\-stel\-lung $\pi : I(\mcalS(X)) \to \mLinStet(K)$.
Nach Proposition~\ref{TROnichtAusgeartet} (nicht-ausgeartete Einbettung eines TRO)
findet man abgeschlossene Unterr"aume $K_1$, $K_2$ von $K$
so, \dass $\erzl{\pi(\ternh{X}) K_1}$ dicht in $K_2$ liegt
und $\erzl{\pi(\ternh{X}) K_1^\perp} = \{0\}$ gilt.
\dremark{$\pi$ ist eine vollst"andige Isometrie und ein tern"arer Morphismus.
  Somit ist $\pi(\ternh{X})$ ein TRO.}%
Wegen $\erzl{\pi(\ternh{X}) K_1^\perp} = \{0\}$ gilt:
$\erzl{\pi(\ternh{X}) K} = \erzl{\pi(\ternh{X}) K_1} \subseteq K_2$.
\end{proof}

\begin{satz}\label{halbgrHochliftenAufHR}
Sei $X$ ein Operatorraum und
$(\pi,K,K_1,K_2)$ eine $\ternh{X}$"=Einbettung von $X$.
Sei $A \in \Multlco{X}$ derart,
\dass f"ur die von $A$ erzeugte $C_0$-Halbgruppe $(T_t)_\indtHG$ gilt:
\begin{equation}\label{eqNormatEndl2}
\exists \delta \in \mbbR_{>0} : \sup_{t \in [0,\delta]} \normMl{T_t}{X} < \infty.
\end{equation}
Sei $B \in \Multlco{\ternh{X}}$ mit der Eigenschaft:
$j(D(A)) \subseteq D(B)$ und $j \circ A = B \circ j\restring_{D(A)}$.
\dremark{Kurz: Seien $X$, $A$, $(T_t)_\indtHG$ und $B$ wie in Satz~\ref{halbgrHochliftenAufTX}.}%
Setze $b_t := \pi(\isoIMl^{-1}(T_t))\restring_{K_2}$ f"ur alle $\indtHG$.

Dann ist $(b_t)_\indtHG$ eine $C_0$"=Halbgruppe auf $K_2$.
Sei $C$ der Erzeuger von $(b_t)_\indtHG$.
Es gilt:
\begin{enumaufz}
\item $\erzl{\pi(D(B))K_1} \subseteq D(C)$,
\dremark{Nach Dennis ist $\erzl{\cdot}$ hier "uberfl"ussig.}%
\item $\pi(B(z))(\xi) = C(\pi(z)(\xi))\dremarkm{ = L_C(\pi(z))(\xi)}$ f"ur alle $z \in D(B)$ und $\xi \in K_1$,
\item $\lim_{t \downarrow 0} b_t \circ y = y$ f"ur alle $y \in \pi(j(X))$,
\item $b_t \circ \pi(j(x)) \in \pi(j(X))$ f"ur alle $t \in \mbbR_{>0}$ und $x \in X$.
\dremark{Frage: Sind $(S_t)_t$ (bzw. $B$)oder $(b_t)_t$ (bzw. $C$) eindeutig?}%
\dremark{Frage: Kommt man von (ii) nach (i)? Von (iii) nach (ii)?}%
\dremark{Mit Erzeugern formulieren}%
\dremark{Idee von Walther: Weniger voraussetzen, um R"uckrichtung zu erhalten,
  \cmzB{} nur $\forall x \in X \,\forall \xi \in H \,:\,
  \lim_{t \to 0} T_t(x(\xi)) = x(\xi)$ in (i)}%
\dremark{Man stellt $I(\mcalS(X))$ und nicht $\ternh{X}$ dar,
  da die Multiplikation (\cmzB von $a_t j(x)$) in $I(\mcalS(X))$ stattfindet
  und es nicht klar ist, ob man einen Hilbertraum $\tilde{H}$
  mit $\ternh{X} \subseteq L(\tilde{H})$ so findet,
  \dass die Einbettung mit der Multiplikation vertr"aglich ist.}%
\dremark{Es gelte: $A \in \Multlor{X}$ mit $j(Ax) = \tilde{a} \multis j(x)$ f"ur alle $x \in X$.
  Nach dem Beweis von \ref{MultlorSubseteqMultlco} gilt dann
  $a_t j(x) = j(T_t x) = \mre^{t \tilde{a}} \multis j(x)$ f"ur alle $x \in X$,
  also folgt: $a_t = \mre^{t \tilde{a}}$.
  Man erh"alt f"ur alle $z \in Z$:
  $S_t z = \mre^{t \tilde{a}} \multis z = \mre^{t L_{\tilde{a}}}(z)$.
  Es folgt: $B = L_{\tilde{a}} \in L(\ternh{X})$.
  Weiter ergibt sich:
$  b_t
=  \pi(a_t)\restring_{K_2}
=  \pi(\mre^{t \tilde{a}})\restring_{K_2}
=  \pi(\sum_n \frac{1}{n!} t^n \tilde{a}^n)\restring_{K_2}
=  \sum_n \frac{1}{n!} t^n \pi(\tilde{a})^n \restring_{K_2}
=  \mre^{t \pi(\tilde{a})}\restring_{K_2}$.
  Es folgt: $C = \pi(\tilde{a})\restring_{K_2} \in L(K_2)$.}%
\end{enumaufz}
\end{satz}

\begin{proof}
Setze $a_t := \isoIMl^{-1}(T_t)$ und $c_t := \pi(a_t)$ f"ur alle $\indtHG$.

Da $\erzl{\pi(\ternh{X}) K_1}$ dicht in $K_2$ ist und
$c_t \cdot \erzl{\pi(\ternh{X}) K_1} \subseteq K_2$ f"ur alle $\indtHG$ gilt,
\dremark{$c_t \pi(z)\xi = \pi(a_t \multis z) \xi \in \pi(\ternh{X}) K_1$}%
ist $c_t(K_2) \subseteq K_2$, also auch $b_t(K_2) \dremarkm{=c_t(K_2)} \subseteq K_2$.
\dremark{Sei $\xi \in K_2$.
  Dann findet man $z \in \ternh{X}^\mbbN$ und $\zeta \in K_1^\mbbN$
  mit $\lim \pi(z_n)(\zeta_n) = \xi$
  (eigentlich: $\xi = \lim \sum_i \pi(z_n^{(i)})(\zeta_n^{(i)})$).
  Es gilt $c_t(\xi) = c_t (\lim \pi(z_n)(\zeta_n))
    = \lim \pi(a_t \multis z_n)(\zeta_n) \in K_2$.}%
Mit Lemma~\ref{atIstOpHalbgr} folgt,
\dass $(b_t)_\indtHG$ eine Operatorhalbgruppe auf $K_2$ ist.
\dremark{$c_{s+t} = \pi(a_{s+t}) = \pi(a_s \multis a_t) = c_s b_t$, somit
  $b_{s+t} = c_{s+t}\restring_{K_2} = c_s c_t\restring_{K_2}
   = c_s\restring_{K_2} c_t\restring_{K_2} = b_s b_t$}%
Mit $\erzl{\pi(\ternh{X}) K} \subseteq K_2$ erh"alt man:
\begin{equation}\label{eqctpizEqbt}
c_t \pi(z) = b_t \pi(z)  \qquad\text{f"ur alle } \indtHG \text{ und } z \in \ternh{X}.
\end{equation}
Definiere $S_t : \ternh{X} \to \ternh{X}, z \mapsto a_t \multis z$
f"ur alle $\indtHG$.
Da $(S_t)_\indtHG$ nach Proposition~\ref{halbgrHochliftenAufTX2}
eine $C_0$-Halbgruppe auf $\ternh{X}$ ist,
ergibt sich mit \eqref{eqctpizEqbt} f"ur alle $z \in \ternh{X}$ und $\xi \in K_1$:
\begin{align}
   \norm{ b_t(\pi(z)(\xi)) - \pi(z)(\xi) }
\dremarkm{&= \norm{ (b_t \pi(z) - \pi(z))(\xi) }  \nonumber\\}
&\overset{\dremarkm{\eqref{eqctpizEqbt}}}{=}  \norm{ (c_t \pi(z) - \pi(z))(\xi) }  \nonumber\\
&= \norm{ \pi(a_t \multis z - z)(\xi) }
\leq  \norm{ a_t \multis z - z } \, \norm{\xi}  \label{eqHGHochhebenbtpiz}\\
&= \norm{S_t z - z} \, \norm{\xi}
\to  0  \qquad\text{f"ur } t \downarrow 0.  \nonumber
\end{align}
Des weiteren hat man f"ur alle $\indtHG$:
$\norm{b_t} \leq \norm{c_t} \dremarkm{= \norm{\pi(a_t)}} = \norm{a_t}$.
Mit \eqref{eqNormatEndl2} und Proposition~\ref{charHGIstC0}
(Charakterisierung der starken Stetigkeit von Halbgruppen) erh"alt man,
\dass $(b_t)_\indtHG$ eine $C_0$-Halbgruppe auf $K_2$ ist.

\bewitemph{(i) und (ii):}
Es gilt mit \eqref{eqHGHochhebenbtpiz} f"ur alle $z \in D(B)$ und $\xi \in K_1$:
\begin{align*}
   \pi(B(z))(\xi)
&= \pi\left( \lim_{t \downarrow 0} \frac{a_t \multis z - z}{t} \right)(\xi)  \\
\dremarkm{&=  \lim_{t \downarrow 0} \frac{\pi(a_t \multis z)(\xi) - \pi(z)(\xi)}{t}  \\
&= \lim_{t \downarrow 0} \frac{c_t(\pi(z)(\xi)) - \pi(z)(\xi)}{t} }%
&\overset{\dremarkm{\eqref{eqHGHochhebenbtpiz}}}{=}
   \lim_{t \downarrow 0} \frac{b_t(\pi(z)(\xi)) - \pi(z)(\xi)}{t}
=  C(\pi(z)(\xi)).
\end{align*}
Insbesondere folgt: $\erzl{\pi(D(B))K_1} \subseteq D(C)$.

\bewitemph{(iii):}
Da $(T_t)_\indtHG$ stark stetig ist, gilt f"ur alle $x \in X$:
\[ \norm{b_t \pi(j(x)) - \pi(j(x))}
\dremarkm{=  \norm{ \pi(a_t \multis j(x)) - \pi(j(x)) }}
=  \norm{ a_t \multis j(x) - j(x) }
=  \norm{ j(T_t x) - j(x) }
\to 0 \]
f"ur $t \to 0$.

\bewitemph{(iv):}
Mit \eqref{eqctpizEqbt} ergibt sich f"ur alle $\indtHG$ und $x \in X$:
\dremark{19.8.'08/2, 16.8.'07/2, 3.9.'07/1}%
\[ b_t \pi(j(x))
\overset{\dremarkm{\eqref{eqctpizEqbt}}}{=}  c_t \pi(j(x))
\dremarkm{=  \pi(a_t) \pi(j(x))}
=  \pi(a_t \multis j(x)) \in \pi(j(X)).  \qedhere \]
\end{proof}

Ein Analogon von Satz~\ref{halbgrHochliftenAufHR} f"ur $\Adjlco{X}$ lautet:

\begin{satz}\label{grHochliftenAufHR}
Sei $X$ ein Operatorraum und
$(\pi,K,K_1,K_2)$ eine $\ternh{X}$"=Einbettung von $X$.
Sei $A \in \Adjlco{X}$ und $(T_t)_\indtGr$ die von $A$ erzeugte $C_0$-Gruppe.
Dann findet man ein eindeutig bestimmtes $B \in \Adjlco{\ternh{X}}$ mit der Eigenschaft:
$j(D(A)) \subseteq D(B)$ und $j \circ A = B \circ j\restring_{D(A)}$.
\dremark{Kurz: Seien $X$, $A$, $(T_t)_\indtGr$ und $B$ wie in Proposition~\ref{grHochliftenAufTX2}.}%
Setze $b_t := \pi(\isoIMl^{-1}(T_t))\restring_{K_2}$ f"ur alle $\indtGr$.

Dann ist $(b_t)_\indtGr$ eine $C_0$-Gruppe von unit"aren Elementen aus $\mLinStet(K_2)$.
Sei $C$ der Erzeuger von $(b_t)_\indtGr$.
Es gilt:
\begin{enumaufz}
\item $\erzl{\pi(D(B))K_1} \subseteq D(C)$,
\item $\pi(B(z))(\xi) = C(\pi(z)(\xi)) \dremarkm{= L_C(z)(\xi)}$
  f"ur alle $z \in D(B)$ und $\xi \in K_1$,
\item $\lim_{t \to 0} b_t \circ y = y$ f"ur alle $y \in \pi(j(X))$,
\item $b_t \circ \pi(j(x)), b_t^* \circ \pi(j(x)) \in \pi(j(X))$
  f"ur alle $x \in X$ und $t \in \mbbR$.
\end{enumaufz}
\end{satz}

\begin{proof}
Setze $a_t := \isoIMl^{-1}(T_t)$ und $c_t := \pi(a_t)$ f"ur alle $\indtHG$.
Analog zum Beweis von Satz~\ref{halbgrHochliftenAufHR} folgt mit
Satz~\ref{satzHilleYosidaGr} (Satz von Hille-Yosida f"ur $C_0$"=Gruppen),
\dass $(b_t)_\indtGr$ eine $C_0$-Gruppe auf $K_2$ ist.
Nach Lemma~\ref{atIstOpGr} ist $(a_t)_\indtGr$ eine Gruppe von unit"aren Elementen
aus $\mIMls(X)$.
F"ur alle $\indtGr$, $z \in \ternh{X}$ und $\xi \in K_1$ gilt
\dremark{$(*)$: Nach \ref{adjVertauschtRestr} gilt:
  $b_t^* = (c_t\restring_{K_2})^* = c_t^*\restring_{K_2}$,
  insbesondere folgt $b_t^* \pi(j(x)) \in \pi(j(X))$.}%
\begin{align*}
   b_t^* b_t(\pi(z)(\xi))
\overset{\dremarkm{(*)}}{=}  c_t^* c_t(\pi(z)(\xi))
=  \pi(a_t^* \multis a_t \multis z)(\xi)
=  \pi(z)(\xi)
\dremarkm{\overset{\dremarkm{(*)}}{=}  c_t c_t^* (\pi(z)(\xi))}
=  b_t b_t^*(\pi(z)(\xi)),
\end{align*}
somit folgt:\dremark{da $\erzl{\pi(\ternh{X} K_1}$ dicht in $K_2$}
$b_t$ ist unit"ar.
Die restlichen Aussagen erh"alt man mit Satz~\ref{halbgrHochliftenAufHR}.
\end{proof}

Mit Hilfe von Satz~\ref{halbgrHochliftenAufHR} kann man die Elemente von $\Multlco{X}$
charakterisieren:

\begin{satz}[Charakterisierungssatz f"ur $\Multlco{X}$]\label{charMlC0X}
Sei $X$ ein Operatorraum und
$(\pi,K,K_1,K_2)$ eine $\ternh{X}$"=Einbettung von $X$.
Setze $\sigma := \pi \circ j : X \to \mLinStet(K)$.
Sei $A : D(A) \subseteq X \to X$ linear.
Dann sind die folgenden Aussagen "aquivalent:
\begin{enumaequiv}
\item $A \in \Multlco{X}$,
und f"ur die von $A$ erzeugte $C_0$-Halbgruppe $(T_t)_\indtHG$ gilt:
\begin{equation}\label{eqNormMlTtEndl}
\exists \delta \in \mbbR_{>0} : \sup_{t \in [0,\delta]} \normMl{T_t}{X} < \infty.
\end{equation}

\item Man findet eine \cohg $(b_t)_\indtHG$ auf $K_2$ derart, \dass gilt:
\begin{enumaufz}
\item $\sigma(Ax) = \lim_{t \downarrow 0} \frac{b_t \sigma(x) - \sigma(x)}{t}$
  f"ur alle $x \in D(A)$,
\item $\lim_{t \downarrow 0} b_t \circ y = y$ f"ur alle $y \in \sigma(X)$,
\item $b_t \circ y \in \sigma(X)$ f"ur alle $y \in \sigma(X)$ und $t \in \mbbR_{>0}$,
\item $A$ ist dicht definiert mit $\rho(A) \neq \emptyset$.
\end{enumaufz}
\end{enumaequiv}
\end{satz}

\begin{proof}

\bewitemphq{(a)$\Rightarrow$(b)} Es gelte \bewitemph{(a)}.
F"ur alle $\indtHG$ setze $a_t := \isoIMl^{-1}(T_t)$,
also $j(T_t x) = a_t \multis j(x)$ f"ur alle $x \in X$,
und $b_t := \pi(a_t)\restring_{K_2}$.
Nach Satz~\ref{halbgrHochliftenAufHR} ist $(b_t)_\indtHG$ eine $C_0$-Halbgruppe auf $K_2$.

\bewitemph{(ii) und (iii)} folgen mit Satz~\ref{halbgrHochliftenAufHR}.
\dremark{\bewitemph{(o):} Es gilt $j(D(A)) \subseteq D(B)$,
  also folgt $\pi(j(D(A))) \subseteq \pi(D(B))$.
  Mit $\erzl{\pi(D(B)) K_1} \subseteq D(C)$ erh"alt man:
  $\erzl{\pi(j(D(A))) K_1} \subseteq \erzl{\pi(D(B)) K_1} \subseteq C$.
  \bewitemph{(iv):} Es gilt $a_t \multis j(x) \in j(X)$ f"ur alle $x \in X$,
  also $b_t \pi(j(x)) = c_t \pi(j(x)) = \pi(a_t \multis j(x)) \in \pi(j(X)) = Y$.}%

\bewitemph{(iv)} folgt mit dem Satz von Hille-Yosida, allgemeiner Fall (Satz~\ref{satzHilleYosHG}).

\bewitemph{(i):}
Man hat f"ur alle $x \in D(A)$:
\dremark{$\pi(a_t) \pi(z) = b_t \pi(z)$ f"ur alle $z \in \ternh{X}$.}%
\begin{align*}
   \sigma(Ax)
&= \sigma \left( \lim_{t \downarrow 0} \frac{T_t x - x}{t} \right)
\dremarkm{=  \lim_{t \downarrow 0} \frac{\sigma(T_t x) - \sigma(x)}{t}  }
=  \lim_{t \downarrow 0} \frac{\pi(a_t) \sigma(x) - \sigma(x)}{t}
=  \lim_{t \downarrow 0} \frac{b_t \sigma(x) - \sigma(x)}{t}.
\end{align*}

\bewitemphq{(b)$\Rightarrow$(a)} Es gelte \bewitemph{(b)}.
Definiere $c_t : K \to K, \xi \mapsto b_t(p_{K_2}(\xi))$,
wobei mit $p_{K_2}$ die Projektion von $K$ auf $K_2$ bezeichnet werde.
Man erh"alt:
\begin{equation}\label{eqctsigmaxEqbt}
c_t \sigma(x) = b_t \sigma(x) \qquad\text{f"ur alle } x \in X.
\end{equation}
Setze $Y := \sigma(X)$.
Sei
\[ S_t := L_{c_t} : Y \to Y, y \mapsto c_t y  \quad\text{und}\quad
   T_t := \sigma^{-1} \circ S_t \circ \sigma : X \to X \]
f"ur alle $\indtHG$.
Nach (iii) gilt $S_t \in \Multlor{Y}$ und $T_t \in \Multlor{X}$ f"ur alle $\indtHG$,
\dremark{$\sigma(X) \subseteq \pi(\ternh{X})$,
  $\sigma(T_t x) = S_t(\sigma(x)) = c_t \sigma(x) = b_t \sigma(x) \in Y$}%
und es ergibt sich:
\dremark{$X$ ist vollst"andig isometrisch isomorph zu $Y$.}
\begin{equation}\label{eqCharNormTt}
\normMl{T_t}{X} = \normMl{S_t}{Y} \leq \norm{c_t} = \norm{b_t}.
\end{equation}
Da $(b_t)_\indtHG$ eine $C_0$-Halbgruppe ist, findet man nach der
Charakterisierung der starken Stetigkeit von Halbgruppen (Proposition~\ref{charHGIstC0})
ein $\delta \in \mbbR_{>0}$ mit der Eigenschaft: $\sup_{t \in [0,\delta]} \norm{b_t} < \infty$.
Mit \eqref{eqCharNormTt} erh"alt man:
\[ \sup_{t \in [0,\delta]} \normMl{T_t}{X}
\overset{\dremarkm{\eqref{eqCharNormTt}}}{\leq}  \sup_{t \in [0,\delta]} \norm{b_t}
< \infty. \]
Weiter gilt nach \eqref{eqctsigmaxEqbt} und (ii):
\begin{align*}
   \norm{T_t x - x}
&= \dremarkm{\norm{\sigma(T_t x) - \sigma(x)}
= } \norm{c_t \sigma(x) - \sigma(x)}  \\
&\overset{\dremarkm{\eqref{eqctsigmaxEqbt}}}{=}  \norm{b_t \sigma(x) - \sigma(x)}
\to 0 \qquad\text{f"ur } t \to 0.
\end{align*}
Also ist $(T_t)_\indtHG$ eine $C_0$-Halbgruppe auf $X$.
Sei $\tilde{A}$ der Erzeuger dieser Halbgruppe.
Mit \eqref{eqctsigmaxEqbt} folgt f"ur alle $x \in D(\tilde{A})$:
\[  \sigma(\tilde{A}x)
=  \sigma\left( \lim_{t \downarrow 0} \frac{T_t x - x}{t} \right)
=  \lim_{t \downarrow 0} \frac{c_t \sigma(x) - \sigma(x)}{t}
\overset{\dremarkm{\eqref{eqctsigmaxEqbt}}}{=}
   \lim_{t \downarrow 0} \frac{b_t \sigma(x) - \sigma(x)}{t}. \]
Da $\tilde{A}$ der Erzeuger von $(T_t)_\indtHG$ ist\dremark{also $D(\tilde{A})$ maximal}
und (i) gilt,
erh"alt man $A \subseteq \tilde{A}$.
Wegen (iv) ergibt sich mit Lemma~\ref{BSubseteqAFolgtGleich}: $A = \tilde{A}$.
\dremark{Z. 3.9.'08, S. 1/2}%
\end{proof}

Eine entsprechende Version von Satz~\ref{charMlC0X} f"ur $\Adjlco{X}$ ist der folgende

\begin{satz}[Charakterisierungssatz f"ur $\Adjlco{X}$]\label{charAdjlcoHR}
Sei $X$ ein Operatorraum und
$(\pi,K,K_1,K_2)$ eine $\ternh{X}$"=Einbettung von $X$.
Setze $\sigma := \pi \circ j : X \to \mLinStet(K)$.
Sei $A : D(A) \subseteq X \to X$ linear.
Dann sind die folgenden Aussagen "aquivalent:
\begin{enumaequiv}
\item $A \in \Adjlco{X}$,
\dheisst $A$ ist Erzeuger einer $C_0$-Gruppe $(T_t)_\indtGr$ auf $X$
mit $T_t \in \Adjlor{X}$ unit"ar f"ur alle $\indtGr$.

\item Man findet eine $C_0$-Gruppe $(b_t)_\indtHG$
von unit"aren Elementen aus $\mLinStet(K_2)$\dremark{Satz von Stone} derart,
\dass gilt:
\begin{enumaufz}
\item $\sigma(Ax) = \lim_{t \to 0} \frac{b_t \sigma(x) - \sigma(x)}{t}$
  f"ur alle $x \in D(A)$,
\item $\lim_{t \to 0} b_t \circ y = y$ f"ur alle $y \in \sigma(X)$,
\item $b_t \circ y, b_t^* \circ y \in \sigma(X)$ f"ur alle $y \in \sigma(X)$ und $t \in \mbbR$,
\item $A$ ist dicht definiert mit $\rho(A) \neq \emptyset$.
\dremark{Man h"atte gern einen analogen Charakterisierungssatz f"ur $\Adjloru{X}$ (WW).
  Problem hierbei: Bei (b)$\Rightarrow$(a) ist nicht klar,
  ob der Erzeuger die Gestalt $\cmpmatrix{0 & A \\ B & 0}$ hat.}%
\end{enumaufz}
\end{enumaequiv}
\end{satz}

\begin{proof}\bewitemphq{(a)$\Rightarrow$(b)} Es gelte \bewitemph{(a)}.
Mit Proposition~\ref{grHochliftenAufTX2} folgt:
\[ \exists \delta \in \mbbR_{>0} :
   \sup_{t \in [0,\delta]} \normMl{T_t}{X}
<  \infty. \]

Analog zum Beweis von Satz~\ref{charMlC0X} folgen (i) bis (iv),
wobei $(b_t)_\indtGr$ nach Satz~\ref{grHochliftenAufHR}
eine $C_0$-Gruppe von unit"aren Elementen ist.

\bewitemphq{(b)$\Rightarrow$(a)} Es gelte \bewitemph{(b)}.
Man erh"alt analog zum Beweis von Satz~\ref{charMlC0X},
\dass $A$ Erzeuger einer $C_0$-Gruppe $(T_t)_\indtGr$ auf $X$ mit
$T_t \in \Multlor{X}$ f"ur alle $\indtGr$ ist.
Mit (iii) folgt: $T_t \in \Adjlor{X}$ f"ur alle $\indtGr$.\dremark{denn $b_t^* \sigma(y) \in \sigma(X)$}
F"ur alle $\indtGr$ und $x \in X$ gilt
\[ \sigma(T_t^* T_t x) = b_t^* b_t \sigma(x) = \sigma(x) = \sigma(T_t T_t^* x), \]
also ist $T_t$ unit"ar.\dremark{denn $j$ ist injektiv}
\end{proof}

Mit dem obigen Satz erh"alt man:

\begin{satz}[Charakterisierungssatz f"ur $\Adjloru{X}$]\label{charAdjloruHR}
Sei $X$ ein Operatorraum und
$(\pi,K,K_1,K_2)$ eine $\ternh{C_2(X)}$-Einbettung von $C_2(X)$.
Sei $A : D(A) \subseteq X \to X$ linear.
Dann sind die folgenden Aussagen "aquivalent:
\begin{enumaequiv}
\item $A \in \Adjloru{X}$.

\item Man findet eine Abbildung $B : D(B) \subseteq X \to X$
und eine $C_0$-Gruppe $(b_t)_\indtHG$
von unit"aren Elementen aus $K_2$ derart, \dass gilt:
\begin{enumaufz}
\item $\sigma(\tilde{A}x) = \lim_{t \to 0} \frac{b_t \sigma(x) - \sigma(x)}{t}$
  f"ur alle $x \in D(\tilde{A})$,
  wobei $\tilde{A} := \mri \cmsmallpmatrix{ 0 & A \\ B & 0}$ ist,
\item $\lim_{t \to 0} b_t \circ y = y$ f"ur alle $y \in \sigma(C_2(X))$,
\item $b_t \circ y, b_t^* \circ y \in \sigma(C_2(X))$
  f"ur alle $y \in \sigma(C_2(X))$ und $t \in \mbbR$,
\item $\tilde{A}$ ist dicht definiert mit $\rho(\tilde{A}) \neq \emptyset$.
\dremark{Man h"atte gern einen analogen Charakterisierungssatz f"ur $\Adjloru{X}$ (WW).
  Problem hierbei: Bei (b)$\Rightarrow$(a) ist nicht klar,
  ob der Erzeuger die Gestalt $\cmpmatrix{0 & A \\ B & 0}$ hat.}%
\end{enumaufz}
\end{enumaequiv}
\end{satz}

\dremark{
Beweis.
\bewitemphq{(a)$\Rightarrow$(b)}
Es gelte (a).
Dann findet man ein $B : D(B) \subseteq X \to X$ mit
$\tilde{A} := \mri \cmpmatrix{ 0 & A \\ B & 0 } \in \Adjlco{C_2(X)}$.
Mit Satz~\ref{charAdjlcoHR} folgt (b).

\bewitemphq{(b)$\Rightarrow$(a)}
Es gelte (b).
Mit Satz~\ref{charAdjlcoHR} erh"alt man: $\tilde{A} \in \Adjlco{C_2(X)}$.
Es folgt: $A \in \Adjloru{X}$.
}%

\section{Die strikte $X$-Topologie auf $\mLinStet(H)$}
\label{secStrikteXTop}

Im vorherigen Abschnitt wurde gezeigt, wie man eine $C_0$-Halbgruppe
von einem Operatorraum auf einen Hilbertraum "uberf"uhren kann.
In diesem Abschnitt wird eine $C_0$-Halbgruppe $(S_t)_\indtHG$
von einem Hilbertraum $H$ auf einen Operatorraum $X \subseteq L(H)$
"ubertragen.
Daf"ur gen"ugt es nicht, von $(S_t)_\indtHG$ die starke Stetigkeit
auf $H$ zu verlangen, man ben"otigt
\[ \lim_{t \downarrow 0} S_t \circ x = x \qquad\text{f"ur alle } x \in X, \]
wie man in Satz~\ref{halbgrHochliftenAufHR}.(iii) sieht.
Um dies abstrakter zu fassen, wird die sogenannte strikte $X$-Topologie eingef"uhrt.
Des weiteren untersuchen wir, welche Eigenschaften diese Topologie besitzt,
insbesondere Zusammenh"ange zur Normtopologie und zur
starken Topologie.

\begin{definitn}
Sei $X \subseteq \mLinStet(H)$ ein Operatorraum.
F"ur alle $x \in X$ definiere die Halbnorm
\[ p_x : \mLinStet(H) \to \mbbR, T \mapsto \norm{T \circ x}. \]
Die von $(p_x)_{x \in X}$ erzeugte lokalkonvexe Topologie auf $\mLinStet(H)$
wird \defemphi{strikte $X$-Topologie} genannt.
\dremark{Wof"ur? Siehe Gespr"ach 17.1.'08/1}
\end{definitn}

In den folgenden beiden Propositionen werden Eigenschaften der strikten $X$-Topologie beleuchtet.

\begin{bemerkung}
Sei $X \subseteq \mLinStet(H)$ ein Operatorraum.
\begin{enumaufz}
\item Die strikte $X$-Topologie auf $\mLinStet(H)$ ist gr"ober als die Normtopologie.
\item Ist $X$ unital, so gilt in (i) Gleichheit.
\item Die strikte $X$-Topologie ist genau dann hausdorffsch, wenn gilt:
\[ \forall T \in L(H)\setminus\{0\} \,\,\exists x \in X : \norm{T \circ x} \neq 0. \]
\item Ist $\erzl{XH}$ dicht in $H$, so ist die strikte $X$-Topologie hausdorffsch.
\end{enumaufz}
\end{bemerkung}

\begin{proof}
\bewitemph{(i):}
Sei $(\xi_\lambda)_\lambda$ ein Netz in $L(H)$, welches gegen ein $\xi_0 \in L(H)$ konvergiert.
Dann gilt f"ur alle $x \in X$:
\[ \norm{\xi_\lambda \circ x - \xi_0 \circ x}
\leq  \norm{\xi_\lambda - \xi_0} \, \norm{x}
\overset{\lambda}{\longrightarrow}  0. \]

\bewitemph{(ii)} folgt unmittelbar.\dremark{$p_{\Id_X} = \norm{\cdot}$}

\bewitemph{(iii)} erh"alt man mit \cite{WernerFunkana6}, Lemma VIII.1.4.

\bewitemph{(iv):}
Sei $\erzl{X H}$ dicht in $H$ und $T \in L(H) \setminus \{0\}$.
Dann findet man ein $\xi_0 \in H$ mit $T \xi_0 \neq 0$.
Setze $\varepsilon := \normlr{T \xi_0}$.
Da $T$ stetig ist, gibt es ein $\delta \in \mbbR_{> 0}$
derart, \dass gilt:
\[ \forall \eta \in H : \normlr{\eta - \xi_0} < \delta
  \Rightarrow  \normlr{T \eta - T \xi_0} < \frac{\varepsilon}{2}. \]
Sei $\eta \in \cmkug_H(\xi_0, \delta) =: M$.
Es folgt
\[ \frac{\varepsilon}{2}
>  \normlr{T \eta - T \xi_0}
\geq  \dremarkm{\abslr{ \norm{T \eta} - \normlr{T \xi_0} }
=}  \abslr{ \norm{T \eta} - \varepsilon }, \]
\dremark{also $\norm{T \eta} - \varepsilon > - \frac{\varepsilon}{2}$,
  also $\norm{T \eta} > \frac{\varepsilon}{2}$.}%
also $\norm{T \eta} \neq 0$.

Da $\erzl{X H}$ dicht in $H$ liegt, ist $M \cap \erzl{X H} \neq \emptyset$.
Somit findet man $n \in \mbbN$, $x \in X^n$ und $\xi \in H^n$ mit
$\sum_{i=1}^n x_i(\xi_i) \in M \cap \erzl{X H}$.
Es gilt: $0 \neq T\left(\sum_{i=1}^n x_i(\xi_i)\right) = \sum_{i=1}^n T\left(x_i(\xi_i)\right)$.
Daher gibt es ein $j \in \haken{n}$ mit $T(x_j(\xi_j)) \neq 0$.
Mit (iii) folgt die Behauptung.
\dremark{\bewitemph{(iv):} Z. 28.11.'08/1}%
\dremark{Bew.: Z. 23.1.'08}%
\end{proof}

\begin{bemerkung}
Sei $X \subseteq \mLinStet(H)$ ein Operatorraum.
\begin{enumaufz}
\item Ist $\erzl{X H} = H$, so ist die strikte $X$-Topologie auf $\mLinStet(H)$
  feiner als die starke Topologie.
\dremark{Man ben"otigt die Voraussetzung $\erzl{X H} = H$,
  da Netze \iallg{} nicht beschr"ankt sind.}%
\item Sei $\erzl{X H}$ dicht in $H$.
Sei $(T_\lambda)_{\lambda \in \Lambda}$ ein beschr"anktes Netz in $\mLinStet(H)$ und
$T_0 \in \mLinStet(H)$.
Falls $T_\lambda \to T_0$ in der strikten $X$-Topologie, dann folgt: $T_\lambda \to T_0$ stark.
\dremark{Beispiel f"ur nicht-beschr"anktes Netz, f"ur das (ii) nicht gilt?}%
\end{enumaufz}
\end{bemerkung}

\begin{proof}
\bewitemph{(i):}
Es gelte $\erzl{X H} = H$.
Sei $(T_\lambda)_\lambda$ ein Netz in $L(H)$,
welches in der strikten $X$-Topologie gegen ein $T_0 \in L(H)$ konvergiert.
Sei $\xi \in H$.
Dann findet man $n \in \mbbN$, $x \in X^n$ und $\eta \in H^n$ mit
$\xi = \sum_{i=1}^n x_i(\eta_i)$.
Somit folgt:
\begin{align*}
   \norm{ T_\lambda(\xi) - T_0(\xi) }
\dremarkm{&=  \normlr{ T_\lambda\left(\sum_{i=1}^n x_i(\eta_i)\right) -
                       T_0\left(\sum_{i=1}^n x_i(\eta_i)\right) }  \\}
&\leq  \sum_{i=1}^n \norm{ T_\lambda \circ x_i - T_0 \circ x_i } \, \norm{\eta_i}
\overset{\lambda}{\longrightarrow}  0.
\end{align*}

\bewitemph{(ii):}
Es gelte: $T_\lambda \to T_0$ in der strikten $X$-Topologie.
Sei $\xi \in H$ und $\varepsilon \in \mbbR_{>0}$.
Setze $M := \sup_{\lambda \in \Lambda} \norm{T_\lambda} + 1$.
Dann gibt es $n \in \mbbN$, $x \in X^n$ und $\eta \in H^n$ mit
$\norm{ \xi - \sum_{i=1}^n x_i(\eta_i) } < \frac{\varepsilon}{2 ( M + \norm{T_0} )}$.
Weiter findet man ein $\lambda_0 \in \Lambda$ mit der Eigenschaft:
\[ \forall \lambda > \lambda_0\,\,  \forall i \in \haken{n} :
   \norm{T_\lambda \circ x_i - T_0 \circ x_i} < \frac{\varepsilon}{2 n\norm{\eta_i}+1}. \]
F"ur alle $\lambda > \lambda_0$ gilt:\dremark{Bew.: Z. 23.1.'08}%
\begin{align*}
&   \norm{T_\lambda(\xi) - T_0(\xi)}  \\
\dremarkm{\leq&  \normbig{T_\lambda(\xi) - T_\lambda\Bigl(\sum_{i=1}^n x_i(\eta_i)\Bigr)} +
      \normbig{T_\lambda\Bigl(\sum_i x_i(\eta_i)\Bigr) - T_0\Bigl(\sum_i x_i(\eta_i)\Bigr)} +  \\
&     \normbig{T_0\Bigl(\sum_i x_i(\eta_i)\Bigr) - T_0(\xi)}  \\}
\leq&  \left( \norm{T_\lambda} + \norm{T_0} \right) \normlr{\xi - \sum_{i=1}^n x_i(\eta_i)} +
      \sum_{i=1}^n \normlr{T_\lambda(x_i(\eta_i)) - T_0(x_i(\eta_i))}  \\
\leq&  \left( M + \norm{T_0} \right) \normlr{\xi - \sum_{i=1}^n x_i(\eta_i)} +
      \sum_{i=1}^n \norm{T_\lambda \circ x_i - T_0 \circ x_i} \, \norm{\eta_i}
<  \varepsilon.  \qedhere
\end{align*}
\end{proof}

Mit Lemma~\ref{KHStarkKonvFolgtKonv} erh"alt man:

\begin{beispiel}
Sei $H$ ein Hilbertraum.
Sei $T_0 \in \mLinStet(H)$ und
$(T_\lambda)_\lambda$ ein gleichm"a"sig beschr"anktes Netz in $\mLinStet(H)$
mit $T_\lambda \to T_0$ stark.
Dann gilt: $T_\lambda \to T_0$ in der strikten $\kptOp(H)$-Topologie.
\end{beispiel}

Im folgenden wird f"ur einen beliebigen Operatorraum $X \subseteq L(H)$
eine $C_0$"=Halbgruppe von $H$ auf $X$ "uberf"uhrt:

\begin{bemerkung}\label{C0HGvonHRzuOR}
Sei $X \subseteq \mLinStet(H)$ ein Operatorraum.
Sei $(S_t)_\indtHG$ eine $C_0$-Halbgruppe auf $H$ mit Erzeuger $C$ und den Eigenschaften:
\begin{enumaufz}
\item $(S_t)_\indtHG$ konvergiert in der strikten $X$-Topologie gegen $\Id_H$ f"ur $t \to 0$,
\item $S_t \circ x \in X$ f"ur alle $t \in \mbbR_{>0}$ und $x \in X$.
\end{enumaufz}
Definiere $T_t := L_{S_t} : X \to X, x \mapsto S_t \circ x$, f"ur alle $\indtHG$.
Dann ist $(T_t)_\indtHG$ eine $C_0$-Halbgruppe auf $X$ mit
$T_t \in \Multlor{X}$ f"ur alle $\indtHG$ und der Eigenschaft:
\[ \exists \delta \in \mbbR_{>0} : \sup_{t \in [0,\delta]} \normMl{T_t}{X} < \infty. \]
Sei $A$ der zugeh"orige Erzeuger.
Dann gilt $A \in \Multlco{X}$, $\erzl{D(A)H} \subseteq D(C)$ und
\[ A(x)(\xi) = C(x(\xi)) \dremarkm{= L_C(x)(\xi)}
  \qquad\text{f"ur alle } x \in D(A) \text{ und } \xi \in H. \]

\dremark{W.W.: Mit Erzeuger formulieren}%
\dremark{Motivation: In \ref{halbgrHochliftenAufHR} gilt
  f"ur die $C_0$-Halbgruppe $(a_t\restring_K)_t$ auf $K$ $a_t \circ j(x) \in j(X)$ und
  $\lim_{t \downarrow 0} a_t \circ j(x) = j(x)$ f"ur alle $x \in X$.
  Wenn man dies f"ur $(S_t)_t$ fordert, erh"alt man eine $C_0$-Halbgruppe.}
\dremark{Die Bedingung (i) ist "aquivalent zu:
  $\lim_{t \downarrow 0} S_t \circ x = x$ f"ur alle $x \in X$.}%
\end{bemerkung}

\begin{proof}
Offensichtlich ist $(T_t)_\indtHG$ eine $C_0$-Halbgruppe auf $X$ mit
$T_t \in \Multlor{X}$ f"ur alle $\indtHG$.\dmarginpar{Bew. not.?}\dremark{Evtl. beweisen.}
\dremark{(i) $T_0(x) = L_{S_0}(x) = \Id_H \circ x = x$, also $T_0 = \Id_X$.
  (ii) $T_s T_t (x) = S_s \circ S_t \circ x = S_{s+t} \circ x = T_{s+t}(x)$.
  (iii) $\lim_t T_t x = \lim_t S_t \circ x = x$.
  (iv) $\sigma : X \to L(H), x \mapsto x$ ist vollst"andig isometrisch.
  Es gilt: $\sigma(T_t x) = S_t \circ \sigma(x) \in X$.}%
Da $(S_t)_\indtHG$ eine $C_0$-Halbgruppe  ist, findet man nach der
Charakterisierung der starken Stetigkeit von Halbgruppen (Proposition~\ref{charHGIstC0})
ein $\delta \in \mbbR_{>0}$ mit: $\sup_{t \in [0,\delta]} \norm{S_t} < \infty$.
Damit folgt:
\[ \sup_{t \in [0,\delta]} \normMl{T_t}{X}
\leq  \sup_{t \in [0,\delta]} \norm{S_t}
<   \infty. \]

F"ur alle $x \in D(A)$ und $\xi \in H$ gilt
\begin{align*}
   A(x)(\xi)
&= \left(\lim_{t \downarrow 0} \frac{T_t x - x}{t} \right)(\xi)
=  \dremarkm{\lim_{t \downarrow 0} \frac{S_t \circ x - x}{t}(\xi)  \\
&=} \lim_{t \downarrow 0} \frac{S_t(x(\xi)) - x(\xi)}{t}
=  C(x(\xi))
\dremarkm{=  L_C(x)(\xi)},
\end{align*}
insbesondere erh"alt man: $\erzl{D(A)H} \subseteq D(C)$.
\dremark{Z. 14.1.'08/1, 5.2.'08/1, evtl. noch weiterer Zettel}%
\end{proof}

Eine Proposition~\ref{C0HGvonHRzuOR} entsprechende Version f"ur $\Adjlor{X}$ lautet:
\pagebreak

\begin{bemerkung}
Sei $X \subseteq \mLinStet(H)$ ein Operatorraum.
Sei $(S_t)_\indtGr$ eine $C_0$-Gruppe auf $H$ mit Erzeuger $C$ und den Eigenschaften:
\begin{enumaufz}
\item $(S_t)_\indtGr$ konvergiert in der strikten $X$-Topologie gegen $\Id_H$ f"ur $t \to 0$,
\item $S_t \circ x, S_t^* \circ x \in X$ f"ur alle $t \in \mbbR$ und $x \in X$.
\end{enumaufz}
Definiere $T_t := L_{S_t} : X \to X, x \mapsto S_t \circ x$, f"ur alle $\indtGr$.
Dann ist $(T_t)_\indtGr$ eine $C_0$-Gruppe auf $X$ mit
$T_t \in \Adjlor{X}$ f"ur alle $\indtGr$.
Sei $A$ der zugeh"orige Erzeuger.
Dann gilt $\erzl{D(A)H} \subseteq D(C)$ und
\dmarginpar{Frage}\dremark{Frage: Gilt dies f"ur unit"are $C_0$-Gruppen?}%
\[ A(x)(\xi) = C(x(\xi)) \dremarkm{= L_C(x)(\xi)}
  \qquad\text{f"ur alle } x \in D(A) \text{ und } \xi \in H. \]
\end{bemerkung}

Die obige Proposition beweist man analog zu Proposition~\ref{C0HGvonHRzuOR}.

\dremark{Evtl. hinzuf"ugen: Konvergenz bzgl. strikter $X$-Topologie in $\kptOp(H)$.}%

\dremww{
\begin{anmerkung}[Literatur zur strikten Topologie]
\begin{enumaufz}
\item \cite{BlackadarOpAlg}, §I.3.1, §I.3.2, nach I.8.6.3, II.7.2.9, §III.7.3
\item \cite{Wegge-OlsenKTheory}, §2.3, §2.4 (hier findet man am meisten)
\item \cite{LanceHmod}, S. 11, S. 18, S. 76--77
\item \cite{Busby68DoubleCentralizersExt}, 3.4ff
\end{enumaufz}
\end{anmerkung}
}%

\section{Störungstheorie}

In diesem Abschnitt formulieren wir einige Resultate aus der
St"orungstheorie auf Banachr"aumen f"ur Operatorr"aume.
So wird die Summe eines unbeschr"ankten und eines beschr"ankten
Multiplikators und die Summe zweier unbeschr"ankter Multiplikatoren untersucht.
\skiptext

Bei Damaville (\cite{Damaville04RegulariteDesOp}, Proposition 2.1.(1),
siehe auch \cite{Damaville07RegulariteDOpReg}) findet man
das folgende St"orungsresultat f"ur Operatoren auf \hcsmoduln:

\begin{bemerkung}\label{damavilleSumStoer}
Seien $E$, $F$ \hcsmodul{}n "uber $\mfrakA$.
F"ur alle $A \in \Multwor{E,F}$ und $B \in \Adjhm{E,F}$ gilt: $A + B \in \Multwor{E,F}$.
\dliter{\cite{Damaville04RegulariteDesOp}, Proposition 2.1.(1)}%
\end{bemerkung}

Um eine Verallgemeinerung der obigen Proposition (f"ur den Fall $E = F$)
auf Operatorr"aumen zu beweisen,
sei zun"achst an die folgende Aussage erinnert:

\begin{satz}[\cite{EngelNagelSemigroups}, Theorem III.1.10]\label{EngNagTh3.1.10}
Sei $X$ ein Banachraum und $A$ Erzeuger einer $C_0$-Halbgruppe $(T_t)_\indtHG$ auf $X$.
Sei $B \in \mLinStet(X)$.
Dann erzeugt $C := A + B$ eine $C_0$-Halbgruppe $(S_t)_\indtHG$,
die wie folgt definiert ist:
F"ur alle $\indtHG$ setze
\[ S_t := \sum_{n=0}^\infty S_{t,n}, \]
wobei f"ur alle $n \in \mbbN_0$ rekursiv
$S_{t,0} = T_t$ und
\[ S_{t,n+1} = \int_0^t T_{t-s} B S_{s,n} \, ds \]
definiert sei.
\dremark{Nach \cite{EngelNagelSemigroups}, Theorem III.1.3 gilt:
  Sei $M \in \mbbR_{\geq 1}$ und $\omega \in \mbbR$ mit
  $\norm{T_t} \leq M \mre^{\omega t}$ f"ur alle $\indtHG$.
  Dann folgt: $\norm{S_t} \leq M \mre^{(\omega + M\norm{B})t}$ f"ur alle $\indtHG$.}%
\end{satz}

Hiermit erh"alt man:

\begin{satz}
Sei $X$ ein Operatorraum.
\begin{enumaufz}
\item Ist $A \in \Multlco{X}$ und $B \in \Multlor{X}$, so ist $A + B \in \Multlco{X}$.

\item Ist $A \in \Adjlco{X}$ und $B \in \Adjlor{X}$ mit $B^* = -B$,
so ist $A+B \in \Adjlco{X}$.

\item Ist $A \in \Adjloru{X}$ und $B \in \Adjlor{X}$,
so ist $A+B \in \Adjloru{X}$.
\dremark{Alt:
\item Ist $A \in \Adjlco{X}$ und $B \in \Adjlor{X}$, so
erzeugt $A + B$ eine $C_0$-Gruppe $(S_t)_\indtGr$
mit $S_t \in \Adjlor{X}$ f"ur alle $\indtGr$.\dremark{$S_t$ unit"ar? Wohl nicht.}

\item Ist $A \in \Adjloru{X}$ und $B \in \Adjlor{X}$,
so existiert ein $S : D(S) \subseteq X \to X$ so,
\dass $\mri \cmpmatrix{0 & A + B \\ S & 0}$ Erzeuger einer
$C_0$-Gruppe $(U_t)_\indtGr$ ist mit $U_t \in \Adjlor{C_2(X)}$
f"ur alle $\indtGr$.
}%
\dremark{Wenn man beweisen k"onnte, \dass $U_t$ f"ur alle $\indtGr$ unit"ar ist,
  w"are (iii) eine Verallgemeinerung von \cite{Damaville04RegulariteDesOp}, Prop. 2.1.(1).
  Dies evtl. erw"ahnen.}%
\end{enumaufz}
\end{satz}

\begin{proof}
\bewitemph{(i):}
Sei $A \in \Multlco{X}$ und $B \in \Multlor{X}$.
Sei $(T_t)_\indtHG$ die von $A$ erzeugte $C_0$-Halbgruppe.
Setze $C := A+B$.
Definiere f"ur alle $\indtHG$ und $n \in \mbbN_0$ rekursiv
$S_{t,0} = T_t$ und
\[ S_{t,n+1} = \int_0^t T_{t-s} B S_{s,n} \, ds. \]
F"ur alle $\indtHG$ setze $S_t := \sum_{n=0}^\infty S_{t,n}$.
Nach Satz~\ref{EngNagTh3.1.10} erzeugt $C$ die $C_0$"=Halbgruppe $(S_t)_\indtHG$.
Per Induktion zeigt man f"ur alle $n \in \mbbN_0$ und $\indtHG$: $S_{t,n} \in \Multlor{X}$.
\dremark{Es gilt $S_{t,0} = T_t \in \Multlor{X}$.
  Setze $a_t := \isoIMl^{-1}(T_t)$, $b := \isoIMl^{-1}(B)$ und
  $c_t^{(n)} := \isoIMl^{-1}(S_{t,n})$ f"ur alle $\indtHG$, $n \in \mbbN_0$.
  Es gilt:\dremark{$(*)$: $\isoIMl^{-1}$ ist linear und stetig}
\begin{align*}
   \isoIMl^{-1}(S_{t,n+1})
&= \isoIMl^{-1}\left( \int_0^t T_{t-s} B S_{s,n} \, ds \right)
\overset{(*)}{=}
   \int_0^t \isoIMl^{-1} \left( T_{t-s} B S_{s,n} \right) \, ds  \\
&= \int_0^t \isoIMl^{-1}(T_{t-s}) \isoIMl^{-1}(B) \isoIMl^{-1}(S_{s,n}) \, ds
=  \int_0^t a_{t-s} b c^{(n)}_s \, ds
\in \mIMl(X).
\end{align*}
Somit folgt: $S_{t,n+1} \in \Multlor{X}$.
}%
Somit erh"alt man: $S_t \in \Multlor{X}$ f"ur alle $\indtHG$.
\smallskip

\bewitemph{(ii):}
Sei $A \in \Adjlco{X}$ und $B \in \Adjlor{X}$ mit $B^* = -B$.
Man findet nach Satz~\ref{grHochliftenAufTX} ein $\check{A} \in \Adjlco{\ternh{X}}$
mit $j(D(A)) \subseteq D(\check{A})$ und $j \circ A = \check{A} \circ j\restring_{D(A)}$.
Es ist $b := \isoIMl^{-1}_X(B) \in \mIMls(X) \subseteq I_{11}(X)$ schiefadjungiert.
Nach Proposition~\ref{I11TXeqI11X} gilt: $I_{11}(X) \cong I_{11}(\ternh{X})$.
Also kann man $b$ als ein Element von $\mIMls(\ternh{X}) \subseteq I_{11}(\ternh{X})$ auf\/fassen.
\dremark{denn $I(\mcalS(\ternh{X})) \cong I(\mcalS(X))$,
  also gilt:
  $b \multis_{I(\mcalS(\ternh{X}))} j(x)
  = b \multis_{I(\mcalS(X))} j(x) \in j(X)$.}%
Dann ist $\check{B} := \isoIMl_{\ternh{X}}(b) \in \Adjlor{\ternh{X}}$ schiefadjungiert.

Nach Satz~\ref{MultworInAdjlco} ist $\check{A} \in \Multwor{\ternh{X}}$.
Nach Beispiel~\ref{HCsModAdjMultEqAdj} gilt: $\Adjlor{\ternh{X}} \cong \Adjhm{\ternh{X}}$.
Somit hat man: $\check{B} \in \Adjhm{\ternh{X}}$.
Mit dem St"orungsresultat von Damaville (Proposition~\ref{damavilleSumStoer})
erh"alt man: $\check{A} + \check{B} \in \Multwor{\ternh{X}}$.
Da $\check{A} + \check{B}$ schiefadjungiert ist,
folgt wiederum mit Satz~\ref{MultworInAdjlco}:
$\check{A} + \check{B} \in \Adjlco{\ternh{X}}$.

Weil die von $A$ erzeugte $C_0$-Halbgruppe $X$ invariant l"a\cms{}t,
l"a\cms{}t die von $\check{A}\restring_{j(X)}$
erzeugte $C_0$-Halbgruppe $j(X)$ invariant.
Des weiteren gilt: $B(j(X)) \subseteq j(X)$.
Da die von $\check{A} + \check{B}$ erzeugte $C_0$"=Halbgruppe $(\check{S}_t^+)_\indtHG$
die in Satz~\ref{EngNagTh3.1.10} angegebene Gestalt hat,
l"a\cms{}t $(\check{S}_t^+)_\indtHG$ $j(X)$ invariant.
Analog erh"alt man, \dass die von $-(\check{A} + \check{B})$ erzeugte
$C_0$"=Halbgruppe $j(X)$ invariant l"a\cms{}t.
Daher l"a\cms{}t die von $\check{A} + \check{B}$ erzeugte $C_0$"=Gruppe $j(X)$ invariant.
Mit Proposition~\ref{ArestringXInMultlco} erh"alt man:
$(\check{A} + \check{B})\restring_{j(X)} \in \Adjlco{j(X)}$.
Also folgt: $A + B \in \Adjlco{X}$.
\dremark{Es gilt: $A = j_X^{-1} \circ \check{A} \circ j_X\restring_{D(A)}$.
  Weiter gilt $j_X(Bx)
  = b \multis_{I(\mcalS(X))} j_X(x)
  = b \multis_{I(\mcalS(\ternh{X}))} j_X(x)
  = \check{B}(j_X(x))$,
  also $B = j_X^{-1} \circ \check{A} \circ j_X$.}%

\dremark{Alt:
Analog zu (i) erh"alt man, \dass $A+B$ eine $C_0$-Halbgruppe $(S_t^{+})_\indtHG$
erzeugt mit $S^+_t \in \Adjlor{X}$ f"ur alle $\indtHG$.
Weiter erzeugt $-A$ nach dem Satz von Hille-Yosida f"ur Gruppen (Satz~\ref{satzHilleYosidaGr})
eine $C_0$-Halbgruppe, also $-A \in \Multlco{X}$.
Somit folgt mit (i), \dass $-A + (-B) = -(A+B)$ eine $C_0$-Halbgruppe
$(S^-_t)_\indtHG$ erzeugt mit $S^-_t \in \Adjlor{X}$ f"ur alle $\indtHG$.
Mit Satz~\ref{satzHilleYosidaGr} erh"alt man,
\dass $A+B$ eine $C_0$-Gruppe $(S_t)_\indtGr$ erzeugt
mit $S_t \in \Adjlor{X}$ f"ur alle $\indtGr$.
\dremark{$S_t = S^+_t$, falls $t \geq 0$, $S_t = S^-_t$, falls $t < 0$.}%
}%
\smallskip

\bewitemph{(iii):}
Sei $A \in \Adjloru{X}$ und $B \in \Adjlor{X}$.
Man findet ein $C : D(C) \subseteq X \to X$ mit:
$\mri \cdot \underbrace{\cmpmatrix{0 & A \\ C & 0}}_{=: \hat{A}} \in \Adjlco{C_2(X)}$.
Nach Lemma~\ref{NullTTs0InAdjlor} gilt:
$\hat{B} := \cmpmatrix{0 & B \\ B^* & 0} \in \Adjlor{C_2(X)}$.
\dremark{also $\mri \hat{B} \in \Adjlor{C_2(X)}$}%
Da $\mri \hat{B}$ nach Lemma~\ref{NullTTs0InAdjlor} schiefadjungiert ist, folgt mit (ii):
\[ \mri \cmpmatrix{0 & A+B \\ C+B^* & 0} = \mri \hat{A} + \mri \hat{B} \in \Adjlco{C_2(X)}. \]

\dremark{Alt: \bewitemph{(iii):}
Man findet $C : D(C) \subseteq X \to X$ mit:
$\mri \cdot \underbrace{\cmpmatrix{0 & A \\ C & 0}}_{=: \hat{A}} \in \Adjlco{C_2(X)}$.
Nach Lemma~\ref{NullTTs0InAdjlor} gilt:
$\hat{B} := \cmpmatrix{0 & B \\ B^* & 0} \in \Adjlor{C_2(X)}$.
\dremark{also $\mri \hat{B} \in \Adjlor{C_2(X)}$}%
Mit (ii) folgt, \dass
$\mri \cmpmatrix{0 & A+B \\ C+B^* & 0} = \mri \hat{A} + \mri \hat{B}$
Erzeuger einer $C_0$-Gruppe $(S_t)_\indtGr$ auf $C_2(X)$ ist mit $S_t \in \Adjlor{C_2(X)}$
f"ur alle $\indtGr$.}%
\dremark{Z. 5.11.'08, S. 3/4}%
\end{proof}

\dremark{\textbf{(ii) und (iii) sind so falsch:}
\begin{bemerkung}
Sei $X$ ein Operatorraum.
\begin{enumaufz}
\item Ist $A \in \Multlco{X}$ und $B \in \Multlor{X}$, so ist $A + B \in \Multlco{X}$.

\item Ist $A \in \Adjlco{X}$ und $B \in \Adjlor{X}$, so gilt: $A + B \in \Adjlco{X}$.
\dremark{Es ist seltsam, \dass $B$ nicht schiefadjungiert sein mu\cms{}.(WW)}%

\item Ist $A \in \Adjloru{X}$ und $B \in \Adjlor{X}$, so ist $A + B \in \Adjloru{X}$.
\end{enumaufz}
\end{bemerkung}

\begin{proof}
\dmarginpar{Bew. pr}%
\bewitemph{(i):}
Sei $A \in \Multlco{X}$ und $B \in \Multlor{X}$.
Sei $(T_t)_\indtHG$ die von $A$ erzeugte $C_0$-Halbgruppe.
Setze $C := A+B$.
Definiere f"ur alle $\indtHG$ und $n \in \mbbN_0$ rekursiv
$S_{t,0} = T_t$ und
\[ S_{t,n+1} = \int_0^t T_{t-s} B S_{s,n} \, ds. \]
F"ur alle $\indtHG$ setze $S_t := \sum_{n=0}^\infty S_{t,n}$.
Nach Satz~\ref{EngNagTh3.1.10} erzeugt $C$ die $C_0$"=Halbgruppe $(S_t)_\indtHG$.
Per Induktion zeigt man f"ur alle $n \in \mbbN_0$ und $\indtHG$: $S_{t,n} \in \Multlor{X}$.
\dremark{Es gilt $S_{t,0} = T_t \in \Multlor{X}$.
  Setze $a_t := \isoIMl^{-1}(T_t)$, $b := \isoIMl^{-1}(B)$ und
  $c_t^{(n)} := \isoIMl^{-1}(S_{t,n}$ f"ur alle $\indtHG$, $n \in \mbbN_0$.
  Es gilt:\dremark{$(*)$: $\isoIMl^{-1}$ ist linear und stetig}
\begin{align*}
   \isoIMl^{-1}(S_{t,n+1})
&= \isoIMl^{-1}\left( \int_0^t T_{t-s} B S_{s,n} \, ds \right)
\overset{(*)}{=}
   \int_0^t \isoIMl^{-1} \left( T_{t-s} B S_{s,n} \right) \, ds  \\
&= \int_0^t \isoIMl^{-1}(T_{t-s}) \isoIMl^{-1}(B) \isoIMl^{-1}(S_{s,n}) \, ds
=  \int_0^t a_{t-s} b c^{(n)}_s \, ds
\in \mIMl(X).
\end{align*}
Somit folgt: $S_{t,n+1} \in \Multlor{X}$.
}%
Somit erh"alt man: $S_t \in \Multlor{X}$ f"ur alle $\indtHG$.

\bewitemph{(ii):}
Sei $A \in \Adjlco{X}$ und $B \in \Adjlor{X}$.
Wegen $A \in \Adjlco{X} \subseteq \Multlco{X}$ und $B \in \Adjlor{X} \subseteq \Multlor{X}$
ergibt sich mit (i): $A + B \in \Multlco{X}$.
Daher erzeugt $\matnull{A+B} = \begin{smallpmatrix} A+B & 0 \\ 0 & 0 \end{smallpmatrix}$
nach dem Charakterisierungssatz f"ur $\Multlco{X}$
(Satz~\ref{charMultlcomatnull}) eine \cmvk{e} $C_0$-Halbgruppe.
\dmarginfz\dremark{\textbf{\normalsize{}Achtung:} So nicht richtig!
  Die von $A+B$ erzeugte Halbgruppe $(S_t^+)_\indtHG$
  mu"s erf"ullen: $\normMl{S_t^+}{X} \leq 1$ f"ur alle $\indtHG$.}%

Nach Satz~\ref{charAdjlcomatnull} erzeugt
$\matnull{A}$ eine \cmvk{e} $C_0$-Gruppe.
Nach dem Satz von Hille-Yosida f"ur $C_0$-Gruppen im Operatorraum
(Satz~\ref{satzHilleYosidaGrOR}) erzeugt $-\matnull{A} = \matnull{-A}$
eine \cmvk{e} $C_0$-Halbgruppe.
Somit gilt nach Satz~\ref{charMultlcomatnull}: $-A \in \Multlco{X}$.
Mit (i) folgt: $-(A+B) = -A + (-B) \in \Multlco{X}$.
Also erzeugt $\matnull{-(A+B)} = -\matnull{A+B}$ nach Satz~\ref{charMultlcomatnull}
eine \cmvk{e} $C_0$-Halbgruppe.
Mit Satz~\ref{satzHilleYosidaGrOR} folgt: $\matnull{A+B}$ erzeugt eine \cmvk{e} $C_0$-Gruppe.
Es ergibt sich mit Satz~\ref{charAdjlcomatnull}: $A+B \in \Adjlco{X}$.
\dremark{3.4.'09/1}%
\dremark{Alter Beweis:
Analog zu (i) erh"alt man, \dass $A+B$ eine $C_0$-Halbgruppe $(S_t^{+})_\indtHG$
erzeugt mit $S^+_t \in \Adjlor{X}$ f"ur alle $\indtHG$.
Weiter erzeugt $-A$ nach dem Satz von Hille-Yosida f"ur Gruppen (Satz~\ref{satzHilleYosidaGr})
eine $C_0$-Halbgruppe, also $-A \in \Multlco{X}$.
Somit folgt mit (i), \dass $-A + (-B) = -(A+B)$ eine $C_0$-Halbgruppe
$(S^-_t)_\indtHG$ erzeugt mit $S^-_t \in \Adjlor{X}$ f"ur alle $\indtHG$.
Mit Satz~\ref{satzHilleYosidaGr} erh"alt man,
\dass $A+B$ eine $C_0$-Gruppe $(S_t)_\indtGr$ erzeugt
mit $S_t \in \Adjlor{X}$ f"ur alle $\indtGr$.
\dremark{$S_t = S^+_t$, falls $t \geq 0$, $S_t = S^-_t$, falls $t < 0$.}%
}%

\bewitemph{(iii):}
Sei $A \in \Adjloru{X}$ und $B \in \Adjlor{X}$.
Man findet ein $C : D(C) \subseteq X \to X$ mit:
$\mri \cdot \underbrace{\cmpmatrix{0 & A \\ C & 0}}_{=: \hat{A}} \in \Adjlco{C_2(X)}$.
Nach Lemma~\ref{NullTTs0InAdjlor} gilt:
$\hat{B} := \cmpmatrix{0 & B \\ B^* & 0} \in \Adjlor{C_2(X)}$.
\dremark{also $\mri \hat{B} \in \Adjlor{C_2(X)}$}%
Mit (ii) folgt:
$\mri \cmpmatrix{0 & A+B \\ C+B^* & 0} = \mri \hat{A} + \mri \hat{B} \in \Adjlco{C_2(X)}$.
\dremark{Z. 5.11.'08, S. 3/4}%
\end{proof}
}%

Wie man leicht sieht, gilt analog zum Fall von $C_0$-Halbgruppen in Banachr"aumen
(\cite{EngelNagelSemigroups}, Abschnitt II.2.7):

\begin{bemerkung}
Sei $X$ ein Operatorraum.
Seien $A,B \in \Multlco{X}$ derart, \dass die von $A$ und $B$
erzeugten $C_0$-Halbgruppen $(S_t)_\indtHG$ (bzw. $(T_t)_\indtHG$) kommutieren.
Dann wird durch $R_t := S_t T_t$ f"ur alle $\indtHG$ eine $C_0$-Halb\-grup\-pe auf $X$ definiert
mit Erzeuger $C \in \Multlco{X}$.
Ferner ist $D(A) \cap D(B)$ ein wesentlicher Bereich f"ur $C$
und\dremark{5.11.'08/1}%
\[ Cx = Ax + Bx \qquad\text{f"ur alle } x \in D(A) \cap D(B). \]
\end{bemerkung}

\dremark{Beweis.
Da $\Multlor{X}$ nach Satz~\ref{MlXunitaleBA} eine unitale Banachalgebra ist,
folgt: $R_t = S_t T_t \in \Multlor{X}$ f"ur alle $\indtHG$.
}%

Eine entsprechende Version dieser Proposition gilt f"ur $\Adjlco{X}$
anstelle von $\Multlco{X}$.
\dremark{
\begin{bemerkung}
Sei $X$ ein Operatorraum.
Seien $A,B \in \Adjlco{X}$ derart, \dass die von $A$ und $B$
erzeugten $C_0$-Gruppen $(S_t)_\indtGr$ (bzw. $(T_t)_\indtGr$) von unit"aren Elementen kommutieren.
Dann wird durch $R_t := S_t T_t$ f"ur alle $\indtGr$ eine unit"are $C_0$-Grup\-pe
auf $X$ definiert mit Erzeuger $C$.
Es gilt $C \in \Adjlco{X}$.
Ferner ist $D(A) \cap D(B)$ ein wesentlicher Bereich f"ur $C$
und\dremark{5.11.'08/1}%
\[ Cx = Ax + Bx \qquad\text{f"ur alle } x \in D(A) \cap D(B). \]
\end{bemerkung}

Beweis.
Sei $\indtHG$.
Da $\Adjlor{X}$ nach Satz~\ref{AlXunitaleCsAlg} eine unitale \csalgebra ist,
folgt: $R_t = S_t T_t \in \Adjlor{X}$.
Weiter gilt:
$R_t^* R_t
= (S_t T_t)^* S_t T_t
= T_t^* S_t^* S_t T_t
= T_t^* \Id_X T_t
= T_t^* T_t
= \Id_X
= R_t R_t^*$.

}%
\skiptext

F"ur die Formulierung des n"achsten Satzes
sei an die folgende Definition erinnert:

\begin{definitn}
Sei $X$ ein Banachraum und
$A : D(A) \subseteq X \to X$ linear.
Ein Operator $B : D(B) \subseteq X \to X$ hei"st
\defemph{$A$-beschr"ankt}, falls $D(A) \subseteq D(B)$ gilt und
falls $a,b \in \mbbR_{>0}$ existieren mit
\begin{equation}\label{eqAbeschraenkt}
\norm{Bx} \leq a \norm{Ax} + b \norm{x}
\end{equation}
f"ur alle $x \in D(A)$.
\index[B]{Abeschraenkt@$A$-beschränkt}%
Die \defemph{$A$-Schranke} von $B$ ist
\index[B]{ASchranke@$A$-Schranke}%
\index[S]{a0AB@$a_0(A,B)$}%
\dliter{\cite{EngelNagelSemigroups}, Def. III.2.1}%
\[ a_0(A,B)
:=  \inf \bigl\{ a \in \mbbR_{\geq 0} \setfdg \text{es existiert } b \in \mbbR_{>0}
\text{ derart, \dass{} } \eqref{eqAbeschraenkt} \text{ gilt} \bigr\}. \]
\end{definitn}

Im folgenden Satz ist der St"orungsoperator $B : D(B) \subseteq X \to X$ dissipativ, \dheisst
\[ \normlr{(\lambda-B)x} \geq \lambda \norm{x}
   \quad\text{f"ur alle } \lambda \in \mbbR_{>0} \text{ und } x \in D(B). \]

\begin{satz}[\cite{EngelNagelSemigroups}, Theorem III.2.7]\label{stoerungMitDissOp}
Sei $X$ ein Banachraum und $A$ Erzeuger
einer kontraktiven $C_0$-Halbgruppe auf $X$.
Sei $B : D(B) \subseteq X \to X$
dissipativ und $A$-beschr"ankt mit $A$-Schranke $a_0(A,B) < 1$.
Dann ist $A+B$ mit Definitionsbereich $D(A)$
Erzeuger einer kontraktiven $C_0$-Halb\-grup\-pe auf~$X$.
\dremark{Ist nach \cite{Davies80OneParameterSemigroups}, S. 73 unten,
  auch f"ur $C_0$-Gruppen n"utzlich.}%
\end{satz}

Wie man in \cite{Goldstein85SemigroupsAndApplicat}, S. 39 unten,
sieht, gilt der obige Satz f"ur $a_0(A,B) = 1$ im allgemeinen nicht.
\skiptext

Die folgende Proposition ist ein Analogon des obigen Satzes\dremark{\ref{stoerungMitDissOp}}
f"ur Operatorr"aume und folgt mit dem obigen Satz\dremark{\ref{stoerungMitDissOp}} und
Satz~\ref{satzHilleYosidaGrOR} (Satz von Hille-Yosida f"ur $C_0$-Gruppen im Operatorraum):

\begin{bemerkung}\label{stoerungAnBeschraenkt}
Sei $X$ ein Operatorraum, $A$ Erzeuger
einer \cmvk{en} $C_0$-Halbgruppe (bzw. $C_0$-Gruppe) auf $X$ und $B : D(B) \subseteq X \to X$
(bzw. $B$ und $-B$) vollst"andig dissipativ.
Es sei $B_n$ $A_n$-beschr"ankt mit $A_n$-Schranke
\[ a_0(A_n,B_n) < 1 \qquad\text{f"ur alle } n \in \mbbN. \]
Dann ist $A+B$ mit Definitionsbereich $D(A)$ Erzeuger einer \cmvk{en}
$C_0$-Halb\-grup\-pe (bzw. $C_0$-Gruppe) auf $X$.
\dmarginpar{DB}\dremark{Vorschlag von Dennis: Gegenbsp. auff"uhren.}%
\dremark{Anwenden auf Multiplikatoren.}%
\dremark{Z. 17.11.'08/2}%
\end{bemerkung}

\begin{proof}
\bewitemph{(i):}
Sei $n \in \mbbN$.
Es ist $A_n$ Erzeuger einer kontraktiven $C_0$"=Halbgruppe und $B_n$ dissipativ
mit $a_0(A_n,B_n) < 1$.
Nach Satz~\ref{stoerungMitDissOp} ist somit $A_n + B_n = (A+B)_n$ Erzeuger
einer kontraktiven $C_0$-Halbgruppe.

\bewitemph{(ii):}
Nach dem Satz von Hille-Yosida f"ur $C_0$-Gruppen im Operatorraum
(Satz~\ref{satzHilleYosidaGrOR})
sind $A$ und $-A$ Erzeuger einer \cmvk{en} $C_0$"=Halbgruppe auf $X$.
Weiter ist $-B_n$ $(-A_n)$-beschr"ankt mit $(-A_n)$"=Schranke
$a_0(-A_n, -B_n) = a_0(A_n, B_n)$ f"ur alle $n \in \mbbN$.
Mit (i) folgt, \dass $A+B$ und $(-A)+(-B) = -(A+B)$ Erzeuger
einer \cmvk{en} $C_0$"=Halbgruppe auf $X$ sind.
Somit erzeugt $A+B$ nach Satz~\ref{satzHilleYosidaGrOR} eine \cmvk{e} $C_0$-Gruppe auf $X$.
\end{proof}

Mit Hilfe dieser Aussage und des Charakterisierungssatzes
(Satz~\ref{charMultlcomatnull}) ergibt sich:

\begin{bemerkung}
Sei $X$ ein Operatorraum und $A \in \Multlco{X}$ so,
\dass $A$ eine $C_0$-Halbgruppe $(T_t)_\indtHG$ auf $X$
erzeugt mit $T_t \in \overline{\cmkug}_{\Multlor{X}}(0,1)$ f"ur alle $\indtHG$.
Sei $B : D(B) \subseteq X \to X$ derart,
\dass $\matnull{B} = \begin{smallpmatrix} B & 0 \\ 0 & 0 \end{smallpmatrix}$ vollst"andig dissipativ ist und
\dass $B_n$ $A_n$-beschr"ankt ist mit $A_n$-Schranke
\[ a_0(A_n,B_n) < 1 \qquad\text{f"ur alle } n \in \mbbN. \]
Dann gilt $A + B \in \Multlco{X}$,
genauer erzeugt $A+B$ eine $C_0$-Halbgruppe $(S_t)_\indtHG$
mit $S_t \in \overline{\cmkug}_{\Multlor{X}}(0,1)$ f"ur alle $\indtHG$.
\end{bemerkung}

Beispielsweise ist $\matnull{B}$ vollst"andig dissipativ,
falls $B \in \Multlco{X}$ ist.

\begin{proof}
Nach Satz~\ref{charMultlcomatnull} (Charakterisierung von $\Multlco{X}$ mittels $\matnull{A}$)
erzeugt $\matnull{A}$ eine \cmvk{e} $C_0$-Halbgruppe auf $C_2(X)$.
Man findet ein $b \in \mbbR_{>0}$ derart, \dass  gilt:
\[ \norm{Bx} \leq a_0(A,B) \norm{Ax} + b\norm{x} \]
f"ur alle $x \in D(A)$.
Es folgt:
\begin{align*}
   \norm{\matnull{B} x}
=  \norm{B x_1}
\leq  a_0(A,B) \norm{A x_1} + b \norm{x_1}
\leq  a_0(A,B) \norm{\matnull{A} x} + b \norm{x}
\end{align*}
f"ur alle $x \in D(\matnull{A})$.
\dremark{
$  \norm{(\matnull{B})_n x}
=  \norm{\matnull{B_n} \tilde{x}}
=  \norm{B_n \tilde{x}_1}
\leq  a \norm{A_n \tilde{x}_1} + b\norm{x}
=  a \norm{\matnull{A_n} \tilde{x}} + b\norm{x}
=  a \norm{(\matnull{A})_n x} + b\norm{x}$}%
Somit ergibt sich, \dass $(\matnull{B})_n$ $(\matnull{A})_n$-beschr"ankt ist
mit $a_0\bigl( (\matnull{A})_n, (\matnull{B})_n \bigr) < 1$ f"ur alle $n \in \mbbN$.
Nach Proposition~\ref{stoerungAnBeschraenkt} ist $\matnull{A} + \matnull{B} = \matnull{A+B}$
Erzeuger einer \cmvk{en} $C_0$-Halbgruppe auf $C_2(X)$.
Also erh"alt man mit Satz~\ref{charMultlcomatnull} die Behauptung.
\dremark{10.2.'09/2}%
\end{proof}

Eine entsprechende Version der obigen Proposition f"ur $\Adjlco{X}$ lautet:

\begin{bemerkung}
Sei $X$ ein Operatorraum und $A \in \Adjlco{X}$.
Sei $B : D(B) \subseteq X \to X$ derart,
\dass $\matnull{B}$ und $-\matnull{B}$ vollst"andig dissipativ sind und
\dass $B_n$ $A_n$-beschr"ankt ist mit $A_n$-Schranke
\[ a_0(A_n,B_n) < 1 \qquad\text{f"ur alle } n \in \mbbN. \]
Dann gilt: $A + B \in \Adjlco{X}$.
\dremark{Es ist unklar, wie man bei einer Version f"ur $\Adjloru{X}$
  eine \glqq sch"one\grqq{} Bedingung f"ur die Beschr"anktheit formuliert.}%
\end{bemerkung}

\dremark{
\begin{lemma}
Ist $B$ $A$-beschr"ankt, so ist $B_n$ $A_n$-beschr"ankt.
\end{lemma}

Beweis.
Wegen
\begin{align*}
   \normlr{B_n x}_n
\leq  \sum_{i,j=1}^n \normlr{B x_{ij}}
\leq  \sum_{i,j=1}^n \left( a \normlr{A x_{ij}} + b \normlr{x_{ij}} \right)
\leq  n^2 \left( a \norm{Ax} + b \norm{x} \right)
\end{align*}
f"ur alle $x \in D(A_n)$ ist $B_n$ $A_n$-beschr"ankt.

\begin{anmerkung}\label{stoerungsthUebertragen}
Die folgenden Aussagen kann man evtl. auf Operatorr"aume "ubertragen:
\begin{enumaufz}
\item \cite{Pazy83SemigroupsApplicatPDE}, Theorem 3.4.3 (bereits versucht).
\item \cite{EngelNagelSemigroups}, Theorem III.4.8, Theorem III.4.9
  (Approximationss"atze von Kato-Trotter).
\item \cite{KatoPerturbationTh2ed}, Theorem IX.2.11.
\item \cite{Damaville04RegulariteDesOp}, Prop. 2.1.(1), (2) und (3)
  (multiplikative St"orung):
  (1) (additive St"orung).
  Beweisideen:
\begin{enumaufzB}
\item \cite{EngelNagelSemigroups}, III.3.18 ff (Zusammenhang zwischen
    additiven und multiplikativen St"orungen) benutzen.
    Problem: $B \in L(X_1,X) = L(D(A),X)$, also $B \notin \Multlor{X}$.
\item Direkt beweisen. Problem: Wie sieht die erzeugte $C_0$-Halbgruppe aus?
\end{enumaufzB}
\end{enumaufz}
\end{anmerkung}

Evtl. auf Operatorr"aume "ubertragen.
Problem f"ur $\Adjloru{X}$:
Es ist unklar, warum die $C_0$-Gruppe unit"ar ist.

\begin{bemerkung}[Vgl. \cite{Damaville04RegulariteDesOp}, Proposition 2.1.(1)]
Seien $E$, $F$ \hcsmodul{}n "uber $\mfrakA$.
Ist $A \in \Multwor{E,F}$, $B \in \Adjhm{E,F}$, dann ist $A + B \in \Multwor{E,F}$.
\end{bemerkung}
}%

\appendix

\makeatletter \@addtoreset {equation}{section}
\makeatother
\renewcommand {\theequation}{\Alph{chapter}.\arabic{section}.\arabic{equation}}

\chapter{Beschr"ankte Multiplikatoren auf Operatorr"aumen}

Wir wiederholen in diesem Kapitel bekannte Tatsachen "uber Linksmultiplikatoren
und von links adjungierbare Multiplikatoren auf Operatorr"aumen.

\dmarginpar{pr}\dremark{Anhang evtl. deutlicher im Inhaltsverzeichnis kennzeichnen}%

\section{Multiplikatoren}

In diesem Abschnitt erinnern wir an einige Fakten zu Linksmultiplikatoren
auf Operatorr"aumen.
Diese verallgemeinern beispielsweise Linksmultiplikatoren auf \csalgebren
(siehe Beispiel~\ref{bspMlCsAlgEqMlOR}).
Wir notieren eine alternative Beschreibung und
eine intrinsische Charakterisierung der Linksmultiplikatoren.

\begin{definitn}\label{defLinksmultOR}
Sei $X$ ein Operatorraum und $T : X \to X$ linear.
\begin{enumaufz}
\item $T$ hei"st \defemphi{Linksmultiplikator} (bzw. \defemphi{Rechtsmultiplikator}) auf $X$,
falls ein Hilbertraum $H$,
eine vollst"andige Isometrie $\sigma : X \to \mLinStet(H)$
und ein $a \in \mLinStet(H)$ so existieren, \dass $\sigma(Tx) = a\sigma(x)$
(bzw. $\sigma(Tx) = \sigma(x) a$)
f"ur alle $x \in X$ gilt.
\item Mit $\Multlor{X}$ (bzw. $\Multror{X}$) sei die Menge aller Linksmultiplikatoren
(bzw. Rechtsmultiplikatoren) auf $X$ bezeichnet.
\index[S]{Mscrl@$\Multlor{X}$ (Linksmultiplikatoren)}%
\index[S]{Mscrr@$\Multror{X}$ (Rechtsmultiplikatoren)}%
\dremark{\glqq OR\grqq{} f"ur Operatorraum}%
\dliter{\cite{Zarikian01Thesis}, 1.6.2}%
\dremww{
\item Sei $H$ ein Hilbertraum und $\sigma : X \to \mLinStet(H)$ eine vollst"andige Isometrie.
Setze
\[ M_\ell^\sigma(X) := \{ T : X \to X \setfdg T \text{ linear},
  \exists a \in \mLinStet(H)\,\forall x \in X : \sigma(Tx) = a \sigma(x) \}. \]
}%
\item F"ur alle $T \in \Multlor{X}$ definiere die Multiplikatornorm
\dliter{\cite{Zarikian01Thesis}, Def. 1.6.1}%
\index[S]{normMl@$\normMl{\cdot}{X}$ (Multiplikatornorm)}%
\begin{align*}
   \normMl{T}{X}
   &:=  \inf\bigl\{ \norm{a} \setfdg \text{es existieren ein Hilbertraum } H,\, a \in \mLinStet(H) \text{ und } \\
   &\phantom{:=  \inf\bigl\{ \norm{a} \setfdg} \,\,\,\text{eine vollst"andige Isometrie } \sigma : X \to \mLinStet(H) \text{ mit: }\\
   &\phantom{:=  \inf\bigl\{ \norm{a} \setfdg} \,\,\,\forall x \in X: \sigma(Tx) = a\sigma(x) \bigr\}.
\end{align*}
\end{enumaufz}
\end{definitn}

In der Definition von $\Multlor{X}$ kann man $\mLinStet(H)$ durch eine beliebige \csalgebra ersetzen.
Dadurch wird die Menge $\Multlor{X}$ nicht vergr"o"sert, auch die Norm $\normMl{\cdot}{X}$
"andert sich nicht. Genauer gilt:

\begin{bemerkung}[\cite{BlecherLeMerdy04OpAlg}, 4.5.1]
Sei $X$ ein Operatorraum. Dann gilt:
\begin{align*}
\Multlor{X} = \bigl\{ &T : X \to X \setfdg T \text{ linear, es existiert eine \csalgebra{} } \mfrakA, \\
  &a \in A \text{ und eine vollst"andige Isometrie } \sigma : X \to \mfrakA \text{ mit: }  \\
  &\forall x \in X : \sigma(Tx) = a \sigma(x) \bigr\}.
\end{align*}
\end{bemerkung}

\dremark{Beweis.
Sei $T : X \to X$ linear.
  Sei $\mfrakA$ eine \csalgebra, $\sigma : X \to \mfrakA$ vollst"andig isometrisch und
  $a \in \mfrakA$ mit $\sigma(Tx) = a\sigma(x)$ f"ur alle $x \in X$.
  Man findet eine isometrische \sterns{}Darstellung $\pi$ von $\mfrakA$ auf
  einem Hilbertraum $K$.
  Insbesondere ist $\pi$ vollst"andig isometrisch.
  Somit ist $\eta := \pi \circ \sigma : X \to \mLinStet(K)$ vollst"andig isometrisch,
  und es gilt $\eta(Tx) = \pi(a \sigma(x)) = \pi(a) \pi(\sigma(x)) = \pi(a) \eta(x)$
  f"ur alle $x \in X$ und $\norm{a} = \norm{\pi(a)}$.
  Mit der letzten Gleichung folgt:
\begin{align*}
   \norm{T}_{M_\ell(X)} =
   \inf\{ \norm{a} \setfdg& \text{es ex. \csalgebra{} } \mfrakA, a \in \mfrakA \text{ und
   eine vollst"andige Isometrie} \\
   &\sigma : X \to \mfrakA \text{ mit: }
   \forall x \in X: \sigma(Tx) = a\sigma(x) \}.  \qedhere
\end{align*}
}%

\dremark{Ausblick: Sei $\mfrakA$ \csalgebra.
  Dann ist $M_\ell(\mfrakA)$ gleich der "ublichen Algebra der Linksmultiplikatoren in $\mfrakA$}%
\dremark{Nach der Definition eines Linksmultiplikators ist nicht klar, ob die Summe
zweier Linksmultiplikatoren wieder ein Linksmultiplikator ist.
Die Multiplikatornorm erf"ullt offensichtlich alle Normeigenschaften bis
auf die Dreiecksungleichung.}%
\dremark{Um dieses zu zeigen, wird eine andere Konstruktion der injektiven H"ulle aufgef"uhrt.}

Wie man leicht sieht, gilt die folgende

\begin{bemerkung}\label{NormcbLeqNormMl}
Sei $X$ ein Operatorraum und $T \in \Multlor{X}$.
Dann gilt: $\norm{T}_\text{cb} \leq \normMl{T}{X}$.
\end{bemerkung}

\dremww{Beweis.
Sei $j : X \to I(X)$ die kanonische Einbettung.
Nach \ref{IsomIMlAufMl} findet man ein $a \in \mIMl(X)$ mit $T = \tilde{L}_a$
und $\norm{a}_{I_11} = \normMl{T}{X}$.
F"ur alle $n \in \mbbN$, $x \in M_n(X)$ folgt:
\begin{align*}
   \norm{T_n x}
&= \norm{\sigma_n(T_n x)}
=  \norm{(\sigma(T x_{ij}))_{ij}}
=  \norm{(a \sigma(x_{ij}))_{ij}}
= \norm{\diag_n(a,\dots,a) \cdot \sigma_n(x)} \\
&\leq \norm{\diag_n(a,\dots,a)} \cdot \norm{\sigma_n(x)}
=  \norm{a} \cdot \norm{x}
=  \normMl{T}{X} \cdot \norm{x}.  \qedhere
\end{align*}
}%

Falls $X$ ein unitaler Operatorraum oder ein Hilbert-\cstern{}Modul ist,
\dremark{oder von der Form $C_n(\mfrakA)$, wobei $\mfrakA$ eine Operatoralgebra mit c.a.i.,}%
kann man $\Multlor{X}$ vollst"andig isometrisch in $\mCB(X)$ einbetten
(\cite{Blecher01ShilovBoundary}, S. 20).
Im allgemeinen ist $\Multlor{X}$ allerdings nicht isometrisch in $\mCB(X)$
eingebettet (\cite{Blecher01ShilovBoundary}, Example 4.4).
\dremark{$\Adjlor{X}$ hingegen schon (\cite{BlecherLeMerdy04OpAlg}, Prop. 4.5.8.(4)).}%
\skiptext

Im folgenden halten wir eine wichtige Konstruktion fest,
mit der man unter anderem die injektive H"ulle eines Operatorraumes erh"alt:

\begin{bemerkung}[Vgl. \cite{BlecherLeMerdy04OpAlg}, 4.4.2]\label{BleLeM4.4.2}
Sei $X \subseteq \mLinStet(H)$ ein Operatorraum.
Betrachte das Paulsen"=System $\mcalS(X) \subseteq M_2(\mLinStet(H))$.
Man findet eine \cmvk{e} Abbildung $\Phi : M_2(\mLinStet(H)) \to M_2(\mLinStet(H))$,
deren Bild eine injektive H"ulle $I(\mcalS(X))$ von $\mcalS(X)$ ist.
\dremark{Man findet nach \ref{exInjHuelle} eine minimale $\mcalS(X)$-Projektion $\Phi$
auf $M_2(\mLinStet(H))$,}
Ferner ist $I(\mcalS(X))$ mit der Multiplikation $x \multis y := x \cdot_\Phi y := \Phi(xy)$
\dremark{nach \ref{injHuelleInjOpsysIstCsalg}}%
eine unitale \csalgebra.
Seien $p_1 := \Id_H \oplus \,\,0$ und $p_2 := 0 \oplus \Id_H$ die kanonischen Projektionen in $I(\mcalS(X))$.
Wegen $\Phi(p_i) = p_i$ f"ur alle $i \in \haken{2}$ ist $\Phi$ eckenerhaltend.
\dremark{D. h. jede Ecke wird in sich abgebildet}%
Nach Definition von $\cdot_\Phi$ sind $p_1$ und $p_2$
orthogonale Projektionen in $I(\mcalS(X))$.
Somit wird durch
\[ I(\mcalS(X)) \to M_2(I(\mcalS(X))),
   a \mapsto \begin{pmatrix}
               p_1 \cdot_\Phi a \cdot_\Phi p_1 & p_1 \cdot_\Phi a \cdot_\Phi p_2 \\
               p_2 \cdot_\Phi a \cdot_\Phi p_1 & p_2 \cdot_\Phi a \cdot_\Phi p_2
             \end{pmatrix}, \]
ein injektiver \sterns{}Homomorphismus definiert.
Man kann also $I(\mcalS(X))$ als aus $2 \times 2$-Matrizen bestehend auffassen.
Wir verwenden $I_{k\ell}(X)$ oder kurz $I_{k\ell}$ als Notation f"ur die Ecke $(k,\ell)$ von $I(\mcalS(X))$
f"ur alle $k,\ell \in \haken{2}$,
also hat man: $I_{k\ell}(X) \subseteq I(\mcalS(X)) \subseteq M_2(\mLinStet(H))$.
\index[S]{Ikl@$I_{k\ell}$}%
\index[S]{IklX@$I_{k\ell}(X)$}%
\dremark{$I_{k\ell} = p_k \cdot_\Phi I(\mcalS(X)) \cdot_\Phi p_\ell$}%
Es gilt somit:
\[ I(\mcalS(X)) \cong \begin{pmatrix} I_{11} & I_{12} \\ I_{21} & I_{22} \end{pmatrix}
   \subseteq M_2(I(\mcalS(X))). \]
Weiter sind $I_{11}$ und $I_{22}$ unitale \csalgebren.
Sei $j$ die kanonische Abbildung von $X$ nach $I_{12}$.
Dann ist $(I_{12}, j)$ eine injektive H"ulle von~$X$.
\dremark{$X \mapsto \mcalS(X), x \mapsto \cmpmatrix{0 & x \\ 0 & 0}$,
$\mcalS(X) \to I(\mcalS(X)), x \mapsto x (=\Phi(x))$,
$I(\mcalS(X)) \cong \cmpmatrix{ I_{11} & I_{12} \\ I_{21} & I_{22} },
x \mapsto (p_i \cdot_\Phi x \cdot_\Phi p_j)_{i,j}$,
$j(x) = \cmpmatrix{0 & x \\ 0 & 0} =: \hat{x}$,
denn nach \cite{EffrosRuan00OperatorSpaces}, (6.1.5), S. 107, gilt
\begin{equation}\label{eqabcEqPhiabc}
a \cdot_\Phi (b \cdot_\Phi c) = \Phi(abc) = (a \cdot_\Phi b) \cdot_\Phi c
\end{equation}
f"ur alle $a,b,c \in I(\mcalS(X))$, also folgt:
$p_1 \cdot_\Phi \hat{x} \cdot_\Phi p_2
=  \Phi(p_1 \hat{x} p_2)
=  \Phi(\hat{x})
=  \hat{x}$.
}%

Sei mit $P$ die Ecke $(1,2)$ von $\Phi$ bezeichnet.
Dann ist $P$ idempotent.\dremark{$P \circ P = P$}
Es ist $I_{12}$ ein TRO\dremark{\ref{bspTROCsAlg}}
mit dem Tripelprodukt $[\cdot,\cdot\cdot,\cdots]$,
welches von dem Produkt der \csalgebra $I(\mcalS(X))$ herr"uhrt.
\dremark{Also $[x,y,z] = x \cdot_\Phi y^* \cdot_\Phi z$.}%
Man hat f"ur alle $x,y,z \in I_{12}$:
\dremark{$I(\mcalS(X))$ ist nicht definiert.}%
\dremark{Weitere Eigenschaften von $I_{11}$ und $I_{22}$, aber auch zur injektiven H"ulle,
  in \cite{BlecherPaulsen00MultOfOpSpaces}.
  Multiplikation in $I_{11}$: $a \multis b := a \cdot_\Phi b$.}%
\dliter{\cite{BlecherLeMerdy04OpAlg}, 4.4.2, \cite{Zarikian01Thesis}, Abschnitt 1.4.2}%
\[ [x,y,z] = P(xy^*z) \]
und
\begin{equation}\label{eqeI11xEqx}
e_{I_{11}(X)} \multis x = x.
\end{equation}
\end{bemerkung}

\dremark{Beweis:
Zeige: $I_{12}$ ist injektiv.
Seien $X, Y$ Operatorr"aume mit $X \subseteq Y$.
Setze $\mfrakA := I(\mcalS(X))$ und $Z := p_1 \mfrakA p_2$.
Sei $\varphi : X \to Z$ \cmvk{}.
Da $\mfrakA$ injektiv ist, findet man eine \cmvk{e} Fortsetzung
$\Psi : Y \to \mfrakA$ von $\varphi$.
\[
\xymatrix{
Y \ar@{-->}[r]^\Psi & \mfrakA \\
X \ar[u]^\subseteq \ar[r]_\varphi & Z \ar[u]^\subseteq }
\]
Dann ist $\Phi : Y \to Z, x \mapsto p_1 \Phi(x) p_2$
eine \cmvk{e} Fortsetzung von $\varphi$.

Es ist $xy^*z \in I_{12}$.
Mit \eqref{eqabcEqPhiabc} erh"alt man:
\[ [x,y,z]
=  x \cdot_\Phi y^* \cdot_\Phi z
\overset{\eqref{eqabcEqPhiabc}}{=}  \Phi(xy^*z)
=  P(xy^*z). \]
Beachte:
$I_{ij}(X) \subseteq I(\mcalS(X)) \subseteq M_2(\mLinStet(H)) \cong \mLinStet(H \oplus H)$.
Also gilt $e_{I_{11}} = p_1 \Id_{H \oplus H} p_1 = p_1 = \Id_H \oplus 0$.
Sei $x \in I_{12}$.
Dann findet man $y \in I(\mcalS(X))$ mit $x = p_1 \cdot_\Phi y \cdot_\Phi p_2$.
Es gilt:
$e_{I_{11}} \cdot x = p_1 \multis p_1 \multis y \multis p_2
\overset{\eqref{eqabcEqPhiabc}}{=}  p_1 \multis y \multis p_2
=  x$.
}%

\dremark{Beachte:
$I_{ij}(X) \subseteq I(\mcalS(X)) \subseteq M_2(\mLinStet(H)) \cong \mLinStet(H \oplus H)$.
Es gilt:
$I_{11} \hookrightarrow \begin{pmatrix} \mLinStet(H) & 0 \\ 0 & 0 \end{pmatrix}
  \cong \mLinStet(H)$.
Ist $X$ unital, so gilt
\begin{equation}\label{eqeIXEq0IdH00}
e_{I(X)} = j(e_X) = j(\Id_H) = \begin{pmatrix} 0 & \Id_H \\ 0 & 0 \end{pmatrix},
\end{equation}
denn $j : X \to I(X)$ ist unital und $e_X = \Id_H$.
\dremark{Siehe auch Z. 16.8.'07/3}}%

Die \csalgebren $I_{11}(X)$, $I_{22}(X)$ und $I(\mcalS(X))$ h"angen nicht von der Einbettung
$X \subseteq \mLinStet(H)$ ab. Genauer gilt:
Sei $u : X \to \mLinStet(K)$ vollst"andig isometrisch und $Y := u(X)$.
Dann findet man einen eckenerhaltenden, unitalen \sterns{}Isomorphismus
$\Phi : I(\mcalS(X)) \to I(\mcalS(Y))$,
der durch Einschr"anken einen unitalen \sterns{}Isomorphismus von $I_{11}(X)$ auf $I_{11}(Y)$
(bzw. von $I_{22}(X)$ auf $I_{22}(Y)$) induziert.
\dliter{\cite{BlecherLeMerdy04OpAlg}, 4.4.5.(2)}%
\skiptext

\dremww{
\begin{beispiel}
Sei $H$ ein Hilbertraum.
Dann gilt $I(\mcalS(H^c)) = \begin{pmatrix} \mLinStet(H) & H^c \\ (H^c)^* & \mbbC \end{pmatrix}$.
\dremark{Z. 21.8.'07/1}%
\end{beispiel}

\begin{beispiel}
Sei $H$ ein Hilbertraum und $a \in \mLinStet(H)$ selbstadjungiert.
Dann ist $X := \mbbC \Id_H + \mbbC a$ ein unitales Operatorsystem,
welches isomorph zur \csalgebra{} $\mbbC^2$ ist.
\dremark{Beispiel daf"ur, \dass i.a. $I(\mcalS(X)) \neq \mLinStet(H)$.}%
\end{beispiel}
}%

Man kann $I_{11}(X)$ (bzw. $I_{22}(X)$) mit Hilfe der Multiplikatoren
einer bestimmten \csalgebra beschreiben:

\begin{bemerkung}[Vgl. \cite{BlecherLeMerdy04OpAlg}, 4.4.11, und \cite{BlecherPaulsen00MultOfOpSpaces}, Corollary 1.8]\label{I11EqMCX}
Sei $X$ ein Operatorraum.
\begin{enumaufz}
\item Es sind $\mcalC(X) := I(X) I(X)^*$ und $\mcalD(X) := I(X)^* I(X)$
\csalgebren, wobei das Produkt in der \csalgebra $I(\mcalS(X))$ gebildet wird.
Weiter ist $\mcalC(X)$ (bzw. $\mcalD(X)$) eine \cstern{}Unteralgebra von $I_{11}(X)$
(bzw. $I_{22}(X)$).
\index[S]{CX@$\mcalC(X)$ ($=I(X) I(X)^*$)}%
\index[S]{DX@$\mcalD(X)$ ($=I(X)^* I(X)$)}%

\item Es gilt:
\begin{align*}
I_{11}(X) &= \Multcs{\mcalC(X)} = \Multlcs{\mcalC(X)} \dremarkm{= \Multrcs{\mcalC(X)}}
\quad\text{und}  \\
I_{22}(X) &= \Multcs{\mcalD(X)} = \Multlcs{\mcalD(X)} \dremarkm{= \Multrcs{\mcalD(X)}}.
\end{align*}
\end{enumaufz}
\end{bemerkung}


\begin{bemerkung}\label{I11C2XEqI11M2X}\label{I11TXeqI11X}\label{ISTXeqISX}
Sei $X$ ein Operatorraum.
Dann gilt:
\begin{enumaufz}
\item $\mcalC(C_2(X)) \cong \mcalC(M_2(X))$ und $\mcalD(C_2(X)) \cong \mcalD(X)$.
\item $I_{11}(C_2(X)) \cong I_{11}(M_2(X))$ und $I_{22}(C_2(X)) \cong I_{22}(X)$.
\item $\mcalC(\ternh{X}) \cong \mcalC(X)$ und $\mcalD(\ternh{X}) \cong \mcalD(X)$.
\item $I_{11}(\ternh{X}) \cong I_{11}(X)$ und $I_{22}(\ternh{X}) \cong I_{22}(X)$.
\item $I(\mcalS(\ternh{X})) \cong I(\mcalS(X))$.
\end{enumaufz}
\end{bemerkung}

\begin{proof}
\bewitemph{(i):}
Da f"ur jeden Operatorraum $Y$ gilt $M_{m,n}(I(Y)) \cong I(M_{m,n}(Y))$
(Proposition~\ref{injHuelleMnX}), erh"alt man:
\begin{align*}
   \mcalC(C_2(X))
&= I(C_2(X)) I(C_2(X))^*
\overset{\dremarkm{\ref{injHuelleMnX}}}{\cong}
   C_2(I(X)) C_2(I(X))^*  \\
&\cong  \cmpmatrix{ I(X) & 0 \\ I(X) & 0 } \cmpmatrix{ I(X)^* & I(X)^* \\ 0 & 0 }  \\
&= \cmpmatrix{ I(X)I(X)^* & I(X)I(X)^* \\ I(X)I(X)^* & I(X)I(X)^* }  \\
&= M_2(I(X)) M_2(I(X))^*
\overset{\dremarkm{\ref{injHuelleMnX}}}{\cong}
   I(M_2(X)) I(M_2(X))  \\
&= \mcalC(M_2(X)) 
\end{align*}
und
%
\begin{align*}
   \mcalD(C_2(X))
&= I(C_2(X))^* I(C_2(X))
\overset{\dremarkm{\ref{injHuelleMnX}}}{\cong}
   C_2(I(X))^* C_2(I(X))  \\
&\cong  \cmpmatrix{ I(X)^* & I(X)^* \\ 0 & 0 } \cmpmatrix{ I(X) & 0 \\ I(X) & 0 }  \\
&= \cmpmatrix{ I(X)^*I(X) & 0 \\ 0 & 0 }
\cong  I(X)^*I(X)
=  \mcalD(X).
\end{align*}

\bewitemph{(ii):}
Es ergibt sich mit Proposition~\ref{I11EqMCX} und (i):
\dremark{8.10.'08/1}%
\begin{align*}
   I_{11}(C_2(X))
&\overset{\dremarkm{\ref{I11EqMCX}}}{=}  \Multcs{\mcalC(C_2(X))}
\overset{\dremarkm{(i)}}{\cong}  \Multcs{\mcalC(M_2(X))}
\overset{\dremarkm{\ref{I11EqMCX}}}{=}  I_{11}(M_2(X))  \quad\text{und}\\
   I_{22}(C_2(X))
&\overset{\dremarkm{\ref{I11EqMCX}}}{=}  \Multcs{\mcalD(C_2(X))}
\overset{\dremarkm{(i)}}{\cong}  \Multcs{\mcalD(X)}
\overset{\dremarkm{\ref{I11EqMCX}}}{=}  I_{22}(X).
\end{align*}

\bewitemph{(iii):}
Nach \eqref{ITXeqIX} gilt: $I(\ternh{X}) \cong I(X)$.
Somit folgt:
\[ \mcalC(\ternh{X}) = I(\ternh{X})I(\ternh{X})^* \cong I(X)I(X)^* = \mcalC(X). \]
Analog erh"alt man: $\mcalD(\ternh{X}) \cong \mcalD(X)$.

\bewitemph{(iv)} ergibt sich mit Proposition~\ref{I11EqMCX} und (iii).

\bewitemph{(v):}
Nach \eqref{ITXeqIX} gilt: $I(\ternh{X}) \cong I(X) \cong I_{12}(X)$.
Zusammen mit Proposition~\ref{BleLeM4.4.2} und~(iv) erh"alt man:
\begin{align*}
   I(\mcalS(\ternh{X}))
&\cong  \begin{pmatrix} I_{11}(\ternh{X}) & I(\ternh{X}) \\
                       I(\ternh{X})^* & I_{22}(\ternh{X}) \end{pmatrix}  \\
&\cong  \begin{pmatrix} I_{11}(X) & I(X) \\ I(X)^* & I_{22}(X) \end{pmatrix}
\cong   I(\mcalS(X)).
\end{align*}
\end{proof}

\begin{bemerkung}[\cite{BlecherLeMerdy04OpAlg}, Proposition 4.4.12]\label{Lcinj}
Sei $X$ ein Operatorraum und $a \in I_{11}(X)$.
Falls $a \multis j(x) = 0$ f"ur alle $x \in X$ gilt, dann folgt: $a=0$.
\dliter{\cite{PaulsenComplBoundedMapsPS}, 16.3, \cite{BlecherLeMerdy04OpAlg}, 4.4.12}%
\end{bemerkung}

\dremww{
\begin{bemerkung}
Sei $X$ ein Operatorraum, dessen injektive H"ulle eine \csalgebra{} $\mfrakA$ ist.
Dann ist $M_2(\mfrakA)$ eine injektive H"ulle f"ur das Paulsen-System $\mcalS(X)$.
Ferner kann man $I_{11}(X) = I_{22}(X) = \mfrakA$ erreichen.
\dliter{\cite{BlecherLeMerdy04OpAlg}, Prop. 4.4.13}%
\end{bemerkung}
}%

Um eine alternative Beschreibung von $\Multlor{X}$ zu erhalten,
definieren wir:

\begin{definitn}
Sei $X$ ein Operatorraum.
\begin{enumaufz}
\item Setze $\mIMl(X) := \bigl\{ a \in I_{11}(X) \setfdg a \multis j(X) \subseteq j(X) \bigr\}$
  und $\mIMl^*(X) := \mIMl(X) \cap \mIMl(X)^*$.
\index[S]{IMlX@$\mIMl(X)$}%
\index[S]{IMlXs@$\mIMl^*(X)$}%
\dremark{Bezeichnung ist mi\cms{}verst"andlich, \glqq I\grqq{} erinnert an injektive H"ulle}%
\item F"ur alle $a \in \mIMl(X)$ definiere den Linksmultiplikator
\index[S]{La@$\cmlmult_a$}%
\[ \cmlmult_a : X \to X, x \mapsto j^{-1}(a \multis j(x)). \]
\end{enumaufz}
\end{definitn}

Es ist $\mIMl(X)$ eine unitale Operatoralgebra und $\mIMl^*(X)$ eine unitale \csalgebra
(\cite{Zarikian01Thesis}, S. 12).
\dremark{Def. Operatoralgebra: abgeschlossene Unteralgebra von $\mLinStet(H)$,
  wobei $H$ ein Hilbertraum sei.}%

\begin{satz}[Vgl. \cite{Zarikian01Thesis}, Theorem 1.6.2]\label{IsomIMlAufMl}\label{MlXunitaleBA}
Sei $X$ ein Operatorraum.
Dann ist
\[ \isoIMl : \mIMl(X) \arrowbij \Multlor{X}, a \mapsto \cmlmult_a, \]
ein unitaler, isometrischer Isomorphismus von $\mIMl(X)$ auf
die unitale Banachalgebra $(\Multlor{X},\normMl{\cdot}{X})$.
\index[S]{Theta@$\isoIMl$}%
\end{satz}

Insbesondere gilt: F"ur jedes $T \in \Multlor{X}$ existiert genau ein $a \in \mIMl(X)$ mit
$T = \cmlmult_a$, und man hat $\normMl{T}{X} = \norm{a}_{I_{11}(X)}$.

\dremark{Beweisskizze.
Setze $e := e_{I_{11}(X)}$.
Unital: $\cmlmult_e(x) = j^{-1}(e \multis j(x)) \overset{\eqref{eqeI11xEqx}}{=} j^{-1}(j(x)) = x$,
also $\isoIMl(e) = \cmlmult_e = \Id_X$.

Nach Definition gilt $\normMl{L_c}{X} \leq \norm{c}$.
Nach \ref{Lcinj} ist $L_c$ injektiv.
$\Phi$ surjektiv: Nicht so leicht, andere Aussage geht ein.
}%

Man kann somit $\Multlor{X}$ so mit einer Operatorraumstruktur versehen,
\dass $\Multlor{X}$ vollst"andig isometrisch isomorph zu $\mIMl(X)$ ist.
\dremark{Ziel: Intrinsische Charakterisierung der Linksmultiplikatoren}%
\skiptext

Die folgende wichtige Aussage stellt eine intrinsische Charakterisierung
der Linksmultiplikatoren bereit.
Eine Variante dieses Satzes f"ur sogenannte matrixgeordnete Operatorr"aume
wurde von W. Werner in \cite{WWerner99SmallKgr} (siehe auch \cite{WWerner03}, Theorem 3.10)
bewiesen.
Aus dieser Variante folgt die "Aquivalenz der Punkte (a) und (b) des folgenden Satzes:

\dremark{
\begin{satz}
Sei $X$ ein Operatorraum und $T : X \to X$ linear.
Dann sind die folgenden Aussagen "aquivalent:
\dliter{\cite{Zarikian01Thesis}, 1.6.4}%
\begin{enumaequiv}
\item $T \in \Multlor{X}$ mit $\normMl{T}{X} \leq 1$.
\item Die Abbildung
\[ \tau^c_T : C_2(X) \to C_2(X), \binom{x}{y} \mapsto \binom{Tx}{y}, \]
  ist vollst"andig kontraktiv.
\item Es existiert ein $a \in I_{11}(X)$ mit $\norm{a} \leq 1$ so, \dass gilt:
\[ j(Tx) = a \multis j(x) \qquad\text{f"ur alle } x \in X \]
(vgl. \cite{Zarikian01Thesis}, Proposition 1.6.4).
\end{enumaequiv}
\end{satz}

\dremark{Beweisskizze:
\bewitemph{\glqq(a)$\Rightarrow$(b)\grqq:} Siehe \cite{Zarikian01Thesis}.

\bewitemph{\glqq(b)$\Rightarrow$(c)\grqq:} Benutzt wird: Lemma von Paulsen, Lemma von Choi.

\bewitemph{\glqq(c)$\Rightarrow$(a)\grqq:} Klar.
}%

Indem man \ref{charMultOR} auf $\lambda T$ anstelle von $T$ anwendet, erh"alt man:
}%

\begin{satz}[Vgl. \cite{Zarikian01Thesis}, Proposition 1.6.4]\label{charMultORmitLambda}\label{charMultOR}
Sei $X$ ein Operatorraum, $T : X \to X$ linear und $\lambda \in \mbbR_{>0}$.
Dann sind die folgenden Aussagen "aquivalent:
\begin{enumaequiv}
\item $T \in \Multlor{X}$ mit $\normMl{T}{X} \leq \frac{1}{\lambda}$.
\item Die Abbildung
\[ \tau^c_{\lambda T} : C_2(X) \to C_2(X), \binom{x}{y} \mapsto \binom{\lambda Tx}{y}, \]
  ist vollst"andig kontraktiv.
\item Es existiert ein $a \in I_{11}(X)$ mit $\norm{a} \leq \frac{1}{\lambda}$ so, \dass gilt:
\[ j(Tx) = a \multis j(x) \qquad\text{f"ur alle } x \in X. \]
\dremark{ $\exists b \in I_{11}(X) : \norm{b} \leq 1 \wedge j(\lambda Tx) = b j(x)$.
  Setze $a := \frac{1}{\lambda} b$. }%
\end{enumaequiv}
\end{satz}

\dremark{
\begin{folgerung}
Sei $X$ ein Operatorraum und $T : X \to X$ linear.
Dann gilt $T \in \Multlor{X}$ genau dann, wenn ein $\lambda \in \mbbR_{>0}$ so existert,
\dass die Abbildung $\tau^c_{\lambda T}$ vollst"andig kontraktiv ist.
\dliter{\cite{WWerner03}, Cor. 3.13}%
\end{folgerung}

\begin{proof}
\bewitemph{\glqq$\Rightarrow$\grqq:}
Sei $T \in \Multlor{X}$.
Dann findet man eine \csalgebra $\mfrakA$, eine vollst"andige Isometrie $\sigma : X \to \mfrakA$
und ein $b \in \mfrakA$ so, \dass gilt: $\sigma(Tx) = b \sigma(x)$ f"ur alle $x \in X$.
Setze $\lambda :=
\begin{cases}
\frac{1}{\norm{b}}, & \text{falls } b \neq 0, \\
1, & \text{sonst}
\end{cases}.$
Dann gilt:
\begin{align*}
   \normlr{ \begin{pmatrix} \lambda Tx \\ y \end{pmatrix} }
&= \normlr{ \begin{pmatrix} \lambda Tx & 0 \\ y & 0 \end{pmatrix} }
=  \normlr{ \sigma_2\Bigl(\begin{pmatrix} \lambda Tx & 0 \\ y & 0 \end{pmatrix}\Bigr) }
=  \normlr{ \begin{pmatrix} \sigma(\lambda Tx) & 0 \\ \sigma(y) & 0 \end{pmatrix} }  \\
&= \normlr{ \begin{pmatrix} \lambda b \sigma(x) & 0 \\ \sigma(y) & 0 \end{pmatrix} }
=  \normlr{ \begin{pmatrix} \lambda b & 0 \\ 0 & 1 \end{pmatrix}
            \begin{pmatrix} \sigma(x) & 0 \\ \sigma(y) & 0 \end{pmatrix} }  \\
&\leq \normlr{ \begin{pmatrix} \lambda b & 0 \\ 0 & e \end{pmatrix} } \cdot
     \normlr{ \begin{pmatrix} \sigma(x) & 0 \\ \sigma(y) & 0 \end{pmatrix} }
=  \max\{ \abs{\lambda} \cdot \norm{b}, \norm{e} \} \normlr{ \begin{pmatrix} x \\ y \end{pmatrix} }
=  \normlr{ \begin{pmatrix} x \\ y \end{pmatrix} }.
\end{align*}
\bewitemph{Fall $n \geq 2$:} Analog.

\bewitemph{\glqq$\Leftarrow$\grqq:} Siehe \ref{charMultORmitLambda}\dremark{oder Artikel}.
\end{proof}}%

\section{Von links adjungierbare Multiplikatoren}

In diesem Abschnitt wird an die von links adjungierbaren Multiplikatoren auf
einem beliebigen Operatorraum $X$, als Menge $\Adjlor{X}$, erinnert,
die auf einem \hcsmodul mit den adjungierbaren Operatoren "ubereinstimmen.
Des weiteren wird eine intrinsische Charakterisierung unit"arer Elemente aus $\Adjlor{X}$
angegeben.

\begin{definitn}\label{defAdjAbbOR}
Sei $X$ ein Operatorraum und $T : X \to X$.
\begin{enumaufz}
\item $T$ wird \defemph{von links adjungierbarer Multiplikator} genannt,
falls ein Hilbertraum $H$, eine lineare, vollst"andige Isometrie $\sigma : X \to \mLinStet(H)$
und eine Abbildung $S : X \to X$ so existieren, \dass gilt:
\index[B]{von links adjungierbarer Multiplikator}%
\[ \sigma(Tx)^*\sigma(y) = \sigma(x)^*\sigma(Sy) \qquad\text{f"ur alle } x,y \in X. \]
\item Die Menge der von links adjungierbaren Multiplikatoren auf $X$ wird mit $\Adjlor{X}$ bezeichnet.
\index[S]{AdjlaX@$\Adjlor{X}$}%
\dliter{\cite{Zarikian01Thesis}, 1.7.1}%
\end{enumaufz}
\end{definitn}

Es folgt, \dass jedes Element aus $\Adjlor{X}$ linear und nach dem Satz vom abgeschlossenen
Graphen beschr"ankt ist.\dremark{Unklar}
\dremark{Eine Charakterisierung der selbstadjungierten Elementen
(bzw. positiven Elementen, bzw. der Projektionen)
aus $\Adjlor{X}$ steht in \cite{Blecher01ShilovBoundary}, Th. 4.10.
Eine Charakterisierung der selbstadjungierten (bzw. unit"aren) Elementen aus $\Adjlor{X}$
findet man in \cite{Zarikian01Thesis}, Abschnitt 1.7.2}%

\begin{bemerkung}[\cite{Zarikian01Thesis}, Corollary 1.7.3]\label{linksadjEind}
Sei $X$ ein Operatorraum und $T : X \to X$.
F"ur alle $i \in \haken{2}$ sei $H_i$ ein Hilbertraum, $\sigma_i : X \to L(H_i)$
eine vollst"andige Isometrie und $S_i : X \to X$ mit der Eigenschaft:
\[ \sigma_i(Tx)^* \sigma_i(y) = \sigma_i(x)^* \sigma_i(S_i y) \]
f"ur alle $x,y \in X$.
Dann gilt: $S_1 = S_2$.
\end{bemerkung}

Wegen dieser Proposition ist es sinnvoll, von \emph{der} Linksadjungierten $T^*$
eines Elementes $T \in \Adjlor{X}$ zu sprechen.
\index[S]{Ts@$T^*$}%
Es gilt $T^* \in \Adjlor{X}$ und $T^{**} = T$.
\skiptext

Schr"ankt man den isometrischen Isomorphismus
\[ \isoIMl : \mIMl(X) \arrowbij \Multlor{X}, a \mapsto \cmlmult_a, \]
aus Satz~\ref{IsomIMlAufMl} auf $\mIMls(X)$ ein, so erh"alt man:

\begin{bemerkung}[Vgl. \cite{Zarikian01Thesis}, Proposition 1.7.4, 1.7.5 und 1.7.6]\label{AlXunitaleCsAlg}\label{AdjlorIsomIMls}
Sei $X$ ein Operatorraum.
Dann gilt: $\Adjlor{X} \subseteq \Multlor{X}$.
Ferner ist
\[ \isoIMl\restring_{\mIMls(X)} : \mIMls(X) \arrowbij \Adjlor{X}, a \mapsto \tilde{L}_a, \]
ein unitaler \sterns{}Iso\-mor\-phis\-mus auf die unitale \csalgebra{} $\Adjlor{X}$.
\dremark{Isomorphie: 1.7.6}%
\end{bemerkung}

Insbesondere folgt f"ur ein beliebiges $T \in \Adjlor{X}$:
Sei $a \in \mIMls(X)$ mit $j(Tx) = a \multis j(x)$ f"ur alle $x \in X$.
Dann gilt:
\begin{equation}\label{eqTsxEqasx}
j(T^* x) = a^* \multis j(x) \quad\text{f"ur alle } x \in X.
\end{equation}

\dremww{
\begin{satz}
Sei $X$ ein Operatorraum, $T : X \to X$ linear.
Die folgenden Aussagen sind "aquivalent:
\begin{enumaequiv}
\item $T \in \Adjlor{X}$.
\item Es existiert eine Abbildung $S : X \to X$ derart, \dass gilt:
\dremark{Das innere Produkt ist das Shilov innere Produkt f"ur $X$,
  siehe \cite{BlecherLeMerdy04OpAlg}, 8.3.14}%
\[ \skalpr{Tx}{y} = \skalpr{x}{Ry} \qquad\text{f"ur alle } x,y \in X. \]
\end{enumaequiv}
\end{satz}
}%

Im Gegensatz zu $\Multlor{X}$ kann man $\Adjlor{X}$
stets isometrisch in $\mCB(X)$ einbetten:

\begin{bemerkung}[Vgl. \cite{BlecherLeMerdy04OpAlg}, Proposition 4.5.8.(4)]\label{inAlXNormTEqNormcb}
Sei $X$ ein Operatorraum.
Dann gilt: $\norm{T} = \normcb{T} = \normMl{T}{X}$ f"ur alle $T \in \Adjlor{X}$.
\dliter{\cite{Zarikian01Thesis}, S. 29}%
\end{bemerkung}

Die unit"aren von links adjungierbaren Multiplikatoren kann man wie folgt charakterisieren:

\begin{bemerkung}[\cite{Zarikian01Thesis}, Proposition 1.7.8]\label{charAlXZar}
Sei $X$ ein Operatorraum und $U : X \to X$ linear.
Dann sind die folgenden Aussagen "aquivalent:
\begin{enumaequiv}
\item $U$ ist ein unit"ares Element von $\Adjlor{X}$.
\item $\tau_U^c : C_2(X) \arrowbij C_2(X), \binom{x}{y} \mapsto \binom{Ux}{y},$
  ist ein vollst"andig isometrischer Isomorphismus auf $C_2(X)$.
\item Es existiert ein Hilbertraum, eine vollst"andige Isometrie $\sigma : X \to L(H)$
  und ein $u \in L(H)$ unit"ar derart, \dass
  $\sigma(Ux) = u \sigma(x)$ f"ur alle $x \in X$ gilt.
\end{enumaequiv}
\end{bemerkung}

Die Linksmultiplikatoren auf einem Operatorraum sind eine
Verallgemeinerung der Linksmultiplikatoren einer \csalgebra:

\begin{beispiel}[Vgl. \cite{Blecher04MultCsModAlgStructure}, Example 5.6]\label{bspMlCsAlgEqMlOR}\label{HCsModAdjMultEqAdj}
\begin{enumaufz}
\item Es sei $\mfrakA$ eine \csalgebra.
Mit $X$ sei $\mfrakA$ aufgefa\cms{}t als Operatorraum bezeichnet.
Dann ist $\Multlor{X}$ isomorph zu der Menge $\Multlcs{\mfrakA}$ der
Linksmultiplikatoren auf der \csalgebra $\mfrakA$ und
$\Adjlor{X}$ isomorph zu der Menge $\Multcs{\mfrakA}$ der
Multiplikatoren auf $\mfrakA$.
\dremark{Kurz: $\Multlor{X} \cong \Multlcs{\mfrakA}$
\dremark{$ := \{ x \in L(H) \setfdg x\pi(\mfrakA) \subseteq \pi(\mfrakA) \}$,
  wobei $H$ ein Hilbertraum und
  $\pi : \mfrakA \to L(H)$ eine nicht-ausgeartete \sterns{}Darstellung von $\mfrakA$}%
und $\Adjlor{X} \cong \Multcs{\mfrakA}$.}%
\dliter{\cite{BlecherLeMerdy04OpAlg}, Cor. 8.4.2 oder Z. 8.2.'08/1 (zu Fu"s)}%

\item Sei $H$ ein Hilbertraum.
Dann gilt: $\Multror{H^r} = \Multlor{H^c} = \Adjlor{H^c} = \mLinStet(H)$ und
$\Multror{H^c} = \Multlor{H^r} = \Adjlor{H^r} \cong \mbbC$.
\dliter{\cite{Blecher04MultCsModAlgStructure}, Example 5.6, S. 15,
  \cite{BlecherLeMerdy04OpAlg}, 4.6.8 (ohne Beweis), Z. 15.6.'07, S. 1/2, Z. 21.8.'07/1}%
\dremark{Es gilt: $\Multlor{\operatorname{Max}(\ell^1_n)} = \mbbC$.}%
\dliter{\cite{BlecherLeMerdy04OpAlg}, 4.6.8}%

\item Sei $E$ ein \hcsmodul.
Dann ist $\Adjlor{E}$ gleich der Menge $\Adjhm{E}$ der adjungierbaren Abbildungen auf $E$.
\dliter{\cite{BlecherLeMerdy04OpAlg}, Cor. 8.4.2, auch vor 4.5.1, S. 167,
  s. auch \cite{Blecher04MultCsModAlgStructure}, Example 5.6, S. 15}%
\dremark{
\item Sei $X$ ein unitaler Operatorraum,
$j : X \to C^*_e(X)$ die (unitale) Einbettung von $X$ in $C^*_e(X)$.
Dann gilt:\dliter{\cite{Blecher01MultAndDualOpAlg}, S. 505 unten,
  \cite{Blecher01ShilovBoundary}, Prop. 4.3}
\[ \Multlor{X} \cong \{ a \in C^*_e(X) \setfdg a \multis j(X) \subseteq j(X) \}. \]

\item Sei $\mfrakA$ eine approximativ unitale Operatoralgebra.
Dann ist $\Multlor{\mfrakA}$ die gew"ohnliche Algebra der Linksmultiplikatoren $LM(\mfrakA)$.
\dliter{\cite{Blecher01ShilovBoundary}, 4.17}%

\item Sei $E$ Banachraum.
Dann stimmt $\Multlor{\text{Min}(E)}$ mit der klassischen Funktionen"=Multiplikatoralgebra
$\mathscr{M}(E)$ "uberein und $\Adjlor{\text{Min}(E)}$ mit der
klassischen Zentralisatoralgebra darin.

\item Sei $E$ ein Hilbert-\cstern{}Rechtsmodul.
Dann ist $\Multlor{E}$ (bzw. $\Adjlor{E}$) gleich der Algebra der beschr"ankten
Rechtsmodulabbildungen (bzw. der adjungierbaren Abbildungen) auf $E$.
\dliter{\cite{BlecherLeMerdy04OpAlg}, Cor. 8.4.2, auch vor 4.5.1, S. 167,
  s. auch \cite{Blecher04MultCsModAlgStructure}, Example 5.6, S. 15}%
}%
\end{enumaufz}
\end{beispiel}

\dremark{Beweis:
\bewitemph{(iii):} $\mathbb{B}_B(Z)$ (in \cite{BlecherLeMerdy04OpAlg})
ist nach \cite{BlecherLeMerdy04OpAlg}, 8.1.7, gleich $\Adjhm{Z}$.}%

\dremww{
\section{Beispiele f"ur Multiplikatoren}

\begin{satz}
Sei $H$ ein homogener hilbertscher Operatorraum mit $\dim(H) \geq 2$.
Dann sind die folgenden Aussagen "aquivalent:
\begin{enumaequiv}
\item $H$ ist ein Hilbert-Spaltenraum.
\item $H$ hat einen nichttrivialen linksadjungierbaren Multiplikator,
  \dheisst $\mcalA_\ell(H) \neq \mbbC$.
\dliter{\cite{BlecherZarikian06CalculusOfMIdeals}, Th. 4.4}%
\end{enumaequiv}
\end{satz}

\begin{satz}
Sei $X$ ein Operatorraum.
Dann sind die folgenden Aussagen "aquivalent:
\begin{enumaequiv}
\item $X$ ist vollst"andig isometrisch zu einem Hilbert-Spaltenraum.
\item $\mLinStet(X) = \mcalA_\ell(X)$,
  \dheisst jeder beschr"ankte Operator auf $X$ ist ein linksadjungierbarer Multiplikator.
\item $\mCB(X) = \mcalA_\ell(X)$,
  \dheisst jeder vollst"andig beschr"ankte Operator auf $X$
  ist ein linksadjungierbarer Multiplikator.
\item $\mLinStet(X) = \Multlor{X}$ isometrisch.
\item $\mCB(X) = \Multlor{X}$ isometrisch.
\dliter{\cite{BlecherZarikian06CalculusOfMIdeals}, Th. 4.5}%
\end{enumaequiv}
\end{satz}

\begin{beispiel}
Sei $E$ die Menge der Matrizen in $M_3$, bei denen nur die Eintr"age $(1,2)$ und $(1,3)$
ungleich 0 sind, also Matrizen der Form
$\begin{smallpmatrix}
0 & * & * \\ 0 & 0 & 0 \\ 0 & 0 & 0
\end{smallpmatrix}$.
Dann ist $A := \mbbC E_3 + E$ eine Unteralgebra von $M_3$.
Setze $Q := \begin{smallpmatrix}
3 & 1 & 1 \\ 1 & 3 & 1 \\ 1 & 1 & 3
\end{smallpmatrix}$, $P := Q^{\frac{1}{2}}$ und $X := AP$.
Dann gilt $C^*(X) = M_3$, $\ternh{X} = M_3$ und $\Multlor{X}=A$.
\dliter{\cite{Blecher01ShilovBoundary}, 4.4, S. 20}%
\end{beispiel}

\begin{beispiel}
Sei $X$ ein zweidimensionaler Operatorraum.
Dann ist $\mcalA_\ell(X)$ entweder gleich $\mbbC$, $\ell^\infty_2 = \mbbC \oplus \mbbC$
oder $M_2$. Es gilt:
\begin{enumaufz}
\item $\dim \mcalA_\ell(X) = 4$ genau dann, wenn $X \cong C_2$ vollst"andig isometrisch.
\item $\dim \mcalA_\ell(X) = 2 = \dim \mcalA_r(X)$ genau dann,
  wenn $X \cong \ell^\infty_2$ vollst"andig isometrisch.
\dliter{\cite{BlecherZarikian06CalculusOfMIdeals}, Prop. 4.3}%
\end{enumaufz}
\end{beispiel}
}%

\dremark{Weitere Beispiele:
\begin{enumaufz}
\item Verweis auf Literatur mit weiteren Beispielen zu $\Multlor{E}$:
\cite{BlecherLeMerdy04OpAlg}, S. $191_{12}$.
\item Siehe \cite{BlecherZarikian06CalculusOfMIdeals}, Kapitel 4 (Examples).
\item Siehe \cite{Blecher01ShilovBoundary}, Example 6.8, S. 40.
\end{enumaufz}}%

\chapter{$C_0$-Halbgruppen}

\dremark{Evtl. straffen}%

In diesem Kapitel stellen wir bekannte Tatsachen zu $C_0$-Halbgruppen vor.
Elementare Definitionen und Aussagen zu $C_0$-Halbgruppen werden im
ersten Abschnitt behandelt.
Im darauf\/folgenden Abschnitt werden die Erzeuger von kontraktiven $C_0$-Halbgruppen
in den S"atzen von Hille-Yosida und von Lumer-Phillips charakterisiert.
Im letzten Abschnitt wird kurz auf die Erzeuger von $C_0$-Gruppen eingegangen.


\section{Grundlagen}

\begin{definitn}\label{defC0Halbgr}
Sei $X$ ein Banachraum.
\begin{enumaufz}
\item Eine Familie $(T_t)_\indtHG$ (bzw. $(T_t)_\indtGr$) von Operatoren aus $\mLinStet(X)$ hei"st
  \defemph{(Operator-)Halbgruppe} (bzw. \defemph{(Operator-)Gruppe}) auf $X$, falls gilt:
\begin{enumaufzB}
\item $T_0 = \Id_X$,
\item $T_{s+t} = T_s T_t$ f"ur alle $s,\indtHG$ (bzw. $s,\indtGr$).
\end{enumaufzB}

\item Unter einer \defemph{$C_0$-Halbgruppe} (bzw. \defemph{$C_0$-Gruppe})
(oder auch \defemph{stark stetigen Halbgruppe} (bzw. \defemph{stark stetigen Gruppe}))
  auf $X$ versteht man eine Operatorhalbgruppe $(T_t)_\indtHG$
  (bzw. Operatorgruppe $(T_t)_\indtGr$), f"ur die gilt:
\index[B]{C0Halbgruppe@$C_0$-Halbgruppe}%
\index[B]{stark stetige Halbgruppe}%
\index[B]{C0Gruppe@$C_0$-Gruppe}%
\index[B]{stark stetige Gruppe}%
\begin{equation}\label{eqDefStarkStHG}
\lim_{t \downarrow 0} T_t x = x \qquad(\text{bzw. }
\lim_{t \to 0} T_t x = x) \qquad\text{f"ur alle } x \in X.
\end{equation}

\item Gilt statt \eqref{eqDefStarkStHG} die st"arkere Forderung
\begin{equation}\label{eqDefNormstHG}
\lim_{t \downarrow 0} \norm{T_t - \Id_X} = 0,
\end{equation}
so spricht man von einer \defemph{normstetigen Halbgruppe}.
\index[B]{normstetige Halbgruppe}%
\dremark{Motivation: \cite{WernerFunkana3}, S. 229/230}%
\end{enumaufz}
\end{definitn}

\eqref{eqDefStarkStHG} ist nichts anderes als die Stetigkeit von $t \mapsto T_t$ bei $t=0$
in der starken Operatortopologie auf $\mLinStet(X)$.
Analog bedeutet \eqref{eqDefNormstHG} gerade die Stetigkeit von $t \mapsto T_t$ bei $t=0$
in der Operatornormtopologie.

\begin{beispiel}[\cite{WernerFunkana6}, S. 358--360]  
\begin{enumaufz}
\item Sei $X$ ein Banachraum und $A \in \mLinStet(X)$.
Setze
\[ T_t := \mre^{tA} = \exp(tA) = \sum_{n=0}^\infty \frac{1}{n!} (tA)^n \]
f"ur alle $\indtHG$.\dremark{Mit $\mre^{tA}$ bezeichnen einige Autoren die zu
  $A$ geh"orige $C_0$-Halbgruppe, daher wird dies hier ausf"uhrlich notiert.}%
Dann ist $(T_t)_\indtHG$ eine normstetige Halbgruppe.

\item Sei $p \in [1,\infty[$.
Sei $X = C_0(\mbbR)$ oder $X = L^p(\mbbR)$.
Betrachte die durch
\[ (T_t f)(x) = f(x+t) \]
f"ur alle $\indtHG$, $f \in X$ und $x \in \mbbR$
gegebene \defemphi{Translationshalbgruppe}.
Dann ist $(T_t)_\indtHG$ eine $C_0$-Halb\-grup\-pe.
\dliter{\cite{WernerFunkana3}, S. 330/331}%

\item Sei $p \in [1,\infty[$ und $d \in \mbbN$.
Definiere
\[ \gamma_t : \mbbR^d \to \mbbR,
  x \mapsto \frac{1}{(4\pi t)^{d/2}} \exp\left( -\frac{\abs{x}^2}{4t} \right), \]
f"ur alle $t \in \mbbR_{>0}$.
Die \defemphi{Wärmeleitungshalbgruppe} auf $X := L^p(\mbbR^d)$
ist durch $T_0 = \Id_X$ bzw. $T_t f = \gamma_t * f$
f"ur alle $t \in \mbbR_{>0}$ und $f \in X$ definiert.
\dremark{Man sieht sehr sch"on, wie $\gamma_t$ f"ur $t \to \infty$ immer flacher wird.}%
\dliter{\cite{WernerFunkana3}, S.3 11}%
\end{enumaufz}
\end{beispiel}

\dremark{F"ur eine Halbgruppe $(\mre^{tA})_\indtHG$ wie im obigen Beispiel
s"ahe man $A$ als \glqq Erzeuger\grqq{} der Halbgruppe an.
Um $A$ aus der Exponentialfunktion $e^{tA}$ zur"uckzuerhalten, mu\cms{} man diese differenzieren.}%

Die folgende Proposition ist n"utzlich, um zu "uberpr"ufen, ob eine Halbgruppe stark stetig ist:

\begin{bemerkung}[\cite{EngelNagelSemigroups}, Proposition I.5.3]\label{charHGIstC0}
Sei $X$ ein Banachraum und $(T_t)_\indtHG$ eine Halbgruppe auf $X$.
Dann sind die folgenden Aussagen "aquivalent:
\begin{enumaequiv}
\item $(T_t)_\indtHG$ ist ein $C_0$-Halbgruppe.
\item Die Abbildung $u_x : \mbbR_{\geq 0} \to X, t \mapsto T_t x$,
  ist stetig f"ur alle $x \in X$.
\item Es existiert ein $\delta \in \mbbR_{>0}$ und eine dichte Teilmenge $D$ in $X$ mit
\begin{enumaufz}
\item $\sup_{t \in [0,\delta]} \norm{T_t} < \infty$,
\item $\lim_{t \downarrow 0} T_t x = x$ f"ur alle $x \in D$.
\end{enumaufz}
\end{enumaequiv}
\end{bemerkung}

\begin{bemerkung}[\cite{EngelNagelSemigroups}, Proposition I.5.5]\label{TtLeqMeomegat}
Sei $X$ ein Banachraum und $(T_t)_\indtHG$ eine $C_0$-Halb\-grup\-pe auf $X$.
Dann existieren $M \in \mbbR_{\geq 1}$ und $\omega \in \mbbR$ mit der Eigenschaft:
\[ \norm{T_t} \leq M \mre^{\omega t} \qquad \text{f"ur alle } \indtHG. \]
\end{bemerkung}

Es gibt $C_0$-Halbgruppen, f"ur die $M>1$ gilt, siehe beispielsweise
\cite{EngelNagelSemigroups}, Example I.5.7.(iii) (Translationshalbgruppe mit Sprung), oder
  \cite{Davies80OneParameterSemigroups}, Example 1.19.
\dremark{Translationshalbgruppe auf $C_0(\mbbR)$, versehen mit einer anderen Norm}%
\dmarginpar{Bsp. not.}\dremark{Vorschlag DB: Beispiel eintippen}%
\skiptext

Um das Wachstumsverhalten einer $C_0$-Halbgruppe zu beschreiben,
f"uhrt man die folgenden Begriffe ein:

\begin{definitn}
Sei $M \in \mbbR_{\geq 1}$ und $\omega \in \mbbR$.
Eine $C_0$-Halbgruppe $(T_t)_\indtHG$ hei"st
\begin{enumaufz}
\item \defemph{beschränkt (mit Schranke $M$)},
  falls man $\omega=0$ in Proposition~\ref{TtLeqMeomegat} w"ahlen kann,
\item \defemph{quasikontraktiv (zum Parameter $\omega$)},
  falls man $M=1$ in Proposition~\ref{TtLeqMeomegat} w"ahlen kann, und
\item \defemph{kontraktiv},
  falls man $M=1$ und $\omega=0$ in Proposition~\ref{TtLeqMeomegat} w"ahlen kann.
\dliter{\cite{EngelNagelSemigroups}, Def. I.5.6}%
\index[B]{beschränkte $C_0$-Halbgruppe}%
\index[B]{quasikontraktive $C_0$-Halb\-grup\-pe}%
\index[B]{kontraktive $C_0$-Halbgruppe}%
\end{enumaufz}
\end{definitn}

\dremark{
\begin{bemerkung}
Sei $\Omega$ ein lokalkompakter Hausdorffraum.
Jede Multiplikationshalbgruppe auf $C_0(\Omega)$ (bzw. $L^p(\Omega,\mu)$) ist quasikontraktiv.
\dremark{Literaturverweis fehlt, existiert wohl nicht.}%
\dremark{Etwas "Ahnliches gilt f"ur Banachverb"ande.}%
\end{bemerkung}

Beweis.
$C_0(\Omega)$ (F"ur $L^p(\Omega,\mu)$ analog):
Man findet ein $q \in C(\Omega)$ mit $T_t f = \mre^{tq} f$ f"ur alle $f \in C_0(\Omega)$
und $s := \sup_{\omega \in \Omega} \Re q(\omega) < \infty$.
Dann gilt: $\norm{T_t} = \mre^{ts}$.
\dremark{16.11.'07}%
}%

\dremww{
\begin{bemerkung}[Vgl. \cite{EngelNagelSemigroups}, Lemma II.3.10]\label{UmnormHG}
Sei $X$ ein Banachraum und
$(T_t)_\indtHG$ eine beschr"ankte $C_0$-Halbgruppe auf $X$.
Dann existiert ein $M \in \mbbR_{\geq 1}$ mit $\norm{T_t} \leq M$ f"ur alle $\indtHG$.
Es wird durch
\begin{equation*}
\norm{x}_s := \sup_\indtHG \norm{T_t x}
\end{equation*}
eine Norm auf $X$ definiert, die wegen $\norm{x} \leq \norm{x}_s \leq M \norm{x}$
f"ur alle $x \in X$ zur urspr"unglichen Norm auf $X$ "aquivalent ist.
Ferner ist $(T_t)_\indtHG$ auf $(X, \norm{\cdot}_s)$ eine kontraktive $C_0$-Halbgruppe.
\dremark{$\norm{T_t x}_s = \sup_{s \geq 0} \norm{T_{s+t} x} \leq
  \sup_{t \geq 0} \norm{T_t x} = \norm{x}_s$}%
\end{bemerkung}
}%

Man m"ochte einer $C_0$-Halbgruppe einen eindeutig bestimmten,
im allgemeinen unbeschr"ankten Operator, ihren Erzeuger, zuordnen.
Eine Hauptaufgabe der Theorie der Halbgruppen ist es, den Zusammenhang
zwischen $C_0$-Halbgruppe und Erzeuger zu studieren.
\dmarginpar{zutun}\dremark{Besser erkl"aren, was dies soll (DB).}%

\begin{definitn}
Sei $X$ ein Banachraum und $(T_t)_\indtHG$ eine $C_0$-Halbgruppe.
Der \defemph{Erzeuger} von $(T_t)_\indtHG$
ist der Operator $A$, der durch
\[ Ax = \lim_{t \downarrow 0} \frac{T_t x - x}{t} \]
auf dem Definitionsbereich
\index[B]{Erzeuger einer $C_0$-Halbgruppe}%
\[ D(A) = \left\{ x \in X \setfdg \lim_{t \downarrow 0} \frac{T_t x - x}{t} \text{ existiert} \right\} \]
definiert wird.
\end{definitn}

Definiere $u_x : \mbbR_{\geq 0} \to X, t \mapsto T_t x$.
Dann ist $Ax$ gleich dem Wert der (rechtsseitigen) Ableitung
von $u_x$ in $0$ f"ur alle $x \in D(A)$.

\begin{beispiel}[\cite{WernerFunkana6}, S. 363, \cite{EngelNagelSemigroups}, Proposition 1 in II.2.10] 
\rule{0.1ex}{0mm}\\[-1\baselineskip]
\begin{enumaufz}
\item Sei $X$ ein Banachraum und $A \in \mLinStet(X)$.
Der Erzeuger der Halbgruppe $(\mre^{tA})_\indtHG$ ist $A$ selbst, denn es gilt:
\begin{align*}
   \normlr{\frac{\mre^{tA}-\Id_X}{t} - A}
&= \dremarkm{ \left( \sum_{n=1}^\infty \frac{1}{n!} t^{n-1} A^n \right) - A
=} \normlr{\sum_{n=2}^\infty \frac{t^{n-1}A^n}{n!}}
\leq  \sum_{n=2}^\infty \frac{t^{n-1}\norm{A}^n}{n!}  \\
\dremarkm{&= \sum_{k=0}^\infty \frac{t^{k+1}\norm{A}^{k+2}}{(k+2)!} }
&\leq  t\norm{A}^2\mre^{t\norm{A}} \to 0 \qquad\text{f"ur } t \downarrow 0.
\end{align*}

\item Sei $p \in [1,\infty[$.
Wir betrachten die Translationshalbgruppe $(T_t)_\indtHG$
auf $X = C_0(\mbbR)$ oder $X = L^p(\mbbR)$.
F"ur alle $f \in X$ und $x \in \mbbR$ gilt
\[ \lim_{t \downarrow 0} \frac{T_t f(x) - f(x)}{t}
=  \lim_{t \downarrow 0} \frac{f(x+t) - f(x)}{t}
\dremarkm{=  \frac{d^+f}{dx}(x)}, \]
falls die Grenzwerte existieren.
Somit ist der Erzeuger $A$ von $(T_t)_\indtHG$ formal der Ableitungsoperator $f \mapsto f'$.
\dremark{formal, denn der Definitionsbereich ist unklar}%
Im Falle $X = C_0(\mbbR)$ erh"alt man als Definitionsbereich
\[ D(A) = \bigl\{ f \in C_0(\mbbR) \setfdg f' \text{ existiert und } f' \in C_0(\mbbR) \bigr\}, \]
im Falle $X = L^p(\mbbR)$\dliter{Beweis in \cite{EngelNagelSemigroups}, Prop. II.2.10}
\[ D(A) = \bigl\{ f \in L^p(\mbbR) \setfdg f \text{ ist absolutstetig und } f' \in L^p(\mbbR) \bigr\}. \]
\end{enumaufz}
\end{beispiel}

\dremark{
\begin{beispiel}
Die $n$-dimensionale Diffusionshalbgruppe $(T_t)_\indtHG$ auf $L^p(\mbbR^n)$, $p \in [1,\infty[$:
Die Fouriertransformation $\mcalF$ transformiert $(T_t\restring_\mathscr{S}(\mbbR^n))_\indtHG$
in eine Multiplikationshalbgruppe auf $\mathscr{S}(\mbbR^n)$ mit Erzeuger
\[ B(\mcalF f)(x) = -\abs{x}^2 (\mcalF f)(x) \qquad\text{f"ur alle } f \in \mathscr{S}(\mbbR^n). \]
\dliter{\cite{ArendtEtAl86OneParamSemigrOfPosOp}, A-I, 2.8}%
\end{beispiel}
}

Die folgenden beiden S"atze beschreiben Eigenschaften des Erzeugers einer
$C_0$-Halbgruppe.

\begin{satz}[\cite{WernerFunkana6}, Satz VII.4.6]
Der Erzeuger einer $C_0$-Halbgruppe ist dicht definiert und abgeschlossen.
\end{satz}

\begin{satz}[\cite{WernerFunkana6}, Korollar VII.4.8]\label{HGDerselbeErzFolgtEq}
Zwei $C_0$-Halbgruppen mit demselben Erzeuger stimmen "uberein.
\end{satz}

\section{Die S"atze von Hille-Yosida und von Lumer-Phillips}

Um den nachfolgenden Satz formulieren zu k"onnen, ben"otigen wir
die folgende

\begin{definitn}
Sei $X$ ein Banachraum und
$A : D(A) \subseteq X \to X$ linear.
\begin{enumaufz}
\item Man nennt
\begin{align*}
\rho(A) :=
   \bigl\{ \lambda \in \mbbC \setfdg &\lambda - A : D(A) \arrowbij X
      \text{ ist bijektiv auf } X  \\
   & \text{mit } (\lambda - A)^{-1} \in L(X) \bigr\}
\end{align*}
\index[S]{rhoA@$\rho(A)$ (Resolventenmenge)}%
die \defemphi{Resolventenmenge} von $A$ und
das Komplement $\sigma(A) := \mbbC \setminus \rho(A)$ das \defemph{Spektrum} von $A$.
\dremark{In \cite{WernerFunkana3}, Def. VII.2.14, wird $A$ dicht definiert vorausgesetzt.
  In \cite{EngelNagelSemigroups}, Def. IV.1.1, wird $A$ abgeschlossen vorausgesetzt,
  damit man eine sinnvolle Spektraltheorie erh"alt.
  Ich ben"otige aber $\rho(A)$ f"ur Operatoren,
  die nicht abgeschlossen sind (siehe \ref{BSubseteqAFolgtGleich}).
  Vgl. auch \cite{UllmannDiplom}, 2.63, S. 35}%
\index[B]{Spektrum eines Operators}%
\index[S]{sigmaA@$\sigma(A)$ (Spektrum)}%

\item F"ur alle $\lambda \in \rho(A)$ hei"st
\[ R(\lambda,A) := (\lambda - A)^{-1} = (\lambda \Id_X - A)^{-1} \]
\defemphi{Resolvente} von $A$ im Punkte $\lambda$.
\dremark{$\lambda \mapsto R(\lambda,A)$ bezeichnet man als Resolventenabbildung.}%
\index[S]{RlambdaA@$R(\lambda,A)$ ($:= (\lambda - A)^{-1}$)}%
\end{enumaufz}
\end{definitn}

Ist $A$ abgeschlossen und $\lambda - A$ bijektiv auf $X$,
so ist $(\lambda - A)^{-1}$ nach dem Satz
vom abgeschlossenen Graphen stetig.\dliter{\cite{WernerFunkana3}, Satz IV.4.4}
\skiptext

Der folgende Satz stellt eine wichtige Formel bereit,
die eine $C_0$"=Halbgruppe mit der Resolvente ihres Erzeugers in Verbindung bringt:

\begin{satz}[Vgl. \cite{EngelNagelSemigroups}, Theorem II.1.10 und Corollary II.1.11]\label{HGErzEnthRGromega}\label{ResolventeVonHG}\label{unglFuerResHochn}
Sei $X$ ein Banachraum und
$A$ der Erzeuger der $C_0$-Halbgruppe $(T_t)_\indtHG$ auf $X$.
Nach \refb{\cite{EngelNagelSemigroups}, Proposition I.5.5,}{Proposition~\ref{TtLeqMeomegat}}
findet man $\omega \in \mbbR$ und $M \in \mbbR_{\geq 1}$ mit
\[ \norm{T_t} \leq M \mre^{\omega t} \qquad\text{f"ur alle } \indtHG. \]
Dann gilt:
\begin{enumaufz}
\item $\bigl\{ \lambda \in \mbbC \setfdg \Re \lambda > \omega \bigr\} \subseteq \rho(A)$.
\item $R(\lambda,A)^k x = \frac{1}{(k-1)!} \int_0^\infty s^{k-1} \mre^{-\lambda s} T_s x \, ds$
  f"ur alle $x \in X$, $k \in \mbbN$ und $\lambda \in \mbbC$ mit $\Re \lambda > \omega$.
\item $\norm{R(\lambda,A)^k} \leq \frac{M}{(Re(\lambda) - \omega)^k}$
  f"ur alle $k \in \mbbN$ und $\lambda \in \mbbC$ mit $\Re \lambda > \omega$.
\end{enumaufz}
\end{satz}

\dremark{
\begin{satz}
Sei $X$ ein Banachraum und
$A$ der Erzeuger der $C_0$-Halbgruppe $(T_t)_\indtHG$ auf $X$.
Nach \refb{\cite{EngelNagelSemigroups}, Proposition I.5.5,}{\ref{TtLeqMeomegat}}
findet man $\omega \in \mbbR$ und $M \in \mbbR_{\geq 1}$ mit
\[ \norm{T_t} \leq M \mre^{\omega t} \qquad\text{f"ur alle } \indtHG. \]
Dann gilt:
\begin{enumaufz}
\item $\bigl\{ \lambda \in \mbbC \setfdg \Re \lambda > \omega \bigr\} \subseteq \rho(A)$.
\item $R(\lambda,A) x = \int_0^\infty \mre^{-\lambda s} T_s x \,ds$
  f"ur alle $x \in X$, $\lambda \in \mbbC$ mit $\Re \lambda > \omega$.
\item $\norm{R(\lambda,A)} \leq \frac{M}{\Re(\lambda) - \omega}$ f"ur alle $\lambda \in \mbbC$
  mit $\Re \lambda > \omega$
(vgl. \cite{EngelNagelSemigroups}, Theorem II.1.10).
\end{enumaufz}
\end{satz}
}%

Der folgende fundamentale Satz charakterisiert die Erzeuger kontraktiver $C_0$-Halbgruppen
(\cite{WernerFunkana6}, Theorem VII.4.11):

\begin{satz}[Satz von Hille-Yosida f"ur kontraktive $C_0$-Halbgruppen]\label{satzHilleYoskHG}
Sei $X$ ein Banachraum.
Ein Operator $A : D(A) \subseteq X \to X$ ist genau dann Erzeuger einer
kontraktiven $C_0$-Halbgruppe, wenn $A$ dicht definiert und abgeschlossen ist,
$\mbbR_{>0} \subseteq \rho(A)$ gilt und
\begin{equation*}
\norm{ \lambda R(\lambda,A) } \leq 1
   \qquad \text {f"ur alle } \lambda \in \mbbR_{>0}.
\end{equation*}
\end{satz}

Falls eine $C_0$-Halbgruppe $(T_t)_\indtHG$ quasikontraktiv ist, es also ein $\omega \in \mbbR$ mit
\[ \norm{T_t} \leq \mre^{\omega t} \qquad\text{f"ur alle } \indtHG \]
gibt, kann man die obige Charakterisierung auf die umskalierte kontraktive Halbgruppe,
die gegeben ist durch
\[ S_t := e^{-\omega t} T_t \qquad\text{f"ur alle } \indtHG, \]
anwenden.
Da der Erzeuger von $(S_t)_\indtHG$ die Gestalt $B = A - \omega$ hat,
erh"alt Satz~\ref{satzHilleYoskHG} die folgende Form:

\begin{folgerung}\label{charQuasikontrHG}
Sei $X$ ein Banachraum und $\omega \in \mbbR$.
F"ur einen Operator $A : D(A) \subseteq X \to X$ sind die folgenden
Aussagen "aquivalent:
\begin{enumaequiv}
\item $A$ erzeugt eine quasikontraktive Halbgruppe $(T_t)_\indtHG$ mit
\[ \norm{T_t} \leq \mre^{\omega t} \qquad\text{f"ur alle } \indtHG. \]
\item $A$ ist dicht definiert, abgeschlossen mit $\mbbR_{>\omega} \subseteq \rho(A)$,
und f"ur jedes $\lambda \in \mbbR_{>\omega}$ gilt:
\dliter{\cite{EngelNagelSemigroups}, Cor. II.3.6}%
\[ \norm{(\lambda - \omega)R(\lambda,A)} \leq 1. \]
\end{enumaequiv}
\end{folgerung}

\dremww{
\begin{proof}
\bewitemph{\glqq$\Rightarrow$\grqq:} Es gelte \bewitemph{(a)}.
Setze
\[ S_t := \mre^{-\omega t} T_t \qquad\text{f"ur alle } \indtHG. \]
Der Erzeuger der kontraktiven $C_0$-Halbgruppe $(S_t)_\indtHG$ ist $B := A-\omega$.
Nach Satz~\ref{satzHilleYoskHG} gilt $\mbbR_{>0} \subseteq \rho(B)$ und
$\norm{\lambda R(\lambda,B)} \leq 1$ f"ur alle $\lambda \in \mbbR_{>0}$.
Wegen $\sigma(B) = \sigma(A) - \omega$\dremark{\cite{EngelNagelSemigroups}, II.2.2}
folgt
\dremark{$\rho(B) = \mbbC \setminus \sigma(A) = \mbbC \setminus (\sigma(A)-\omega)
  = (\mbbC \setminus \sigma(A)) - \omega = \rho(A) - \omega$}%
\begin{equation}\label{eqQuasikontrRhoAB}
\mbbR_{>\omega} \subseteq \rho(B) + \omega = \rho(A).
\end{equation}
Weiter gilt f"ur alle $\mu \in \mbbR_{>\omega}$, wobei $\lambda := \mu-\omega>0$,
wegen $R(\lambda+\omega,A) = R(\lambda,B)$:\dremark{\cite{EngelNagelSemigroups}, II.2.2}
\begin{equation}\label{eqResQuasikontrAB}
   \norm{(\mu-\omega) R(\mu,A)}
=  \norm{\lambda R(\lambda+\omega,A)}
=  \norm{\lambda R(\lambda,B)}
\leq 1.
\end{equation}

\bewitemph{\glqq$\Leftarrow$\grqq:} Es gelte \bewitemph{(b)}.
Wegen $D(B) = D(A)$ ist $B$ dicht definiert und abgeschlossen.
Mit \eqref{eqQuasikontrRhoAB} und \eqref{eqResQuasikontrAB} folgt
$\mbbR_{>0} \subseteq \rho(B)$ und $\norm{\mu R(\mu,B)} \leq 1$ f"ur alle $\mu \in \mbbR_{>0}$.
Somit ist $B$ nach Satz~\ref{satzHilleYoskHG} Erzeuger
der kontraktiven $C_0$-Halbgruppe $(S_t)_\indtHG$.
Es folgt \bewitemph{(a)}.
\end{proof}
}%

Man kann den Erzeuger einer beliebigen $C_0$"=Halbgruppe charakterisieren,
ben"otigt dann allerdings Normabsch"atzungen f"ur \emph{alle} Potenzen der Resolvente
(vgl. \cite{EngelNagelSemigroups}, Theorem II.3.8):

\begin{satz}[Satz von Hille-Yosida, allgemeiner Fall]\label{satzHilleYosHG}
Sei $X$ ein Banachraum und $A : D(A) \subseteq X \to X$ linear.
Sei $\omega \in \mbbR$ und $M \in \mbbR_{\geq 1}$.
Dann sind die folgenden Aussagen "aquivalent:
\begin{enumaequiv}
\item $A$ erzeugt eine $C_0$-Halbgruppe $(T_t)_\indtHG$ mit
\[ \norm{T_t} \leq M\mre^{\omega t} \qquad\text{f"ur alle } \indtHG. \]
\item $A$ ist dicht definiert, abgeschlossen mit $\mbbR_{>\omega} \subseteq \rho(A)$,
und f"ur jedes $\lambda \in \mbbR_{>\omega}$ gilt:
\[ \normbig{\bigl[(\lambda - \omega)R(\lambda,A)\bigr]^k} \leq M \qquad\text{f"ur alle } k \in \mbbN. \]
\end{enumaequiv}
\end{satz}

\begin{lemma}\label{BSubseteqAFolgtGleich}
Sei $X$ ein Banachraum.
Seien $A : D(A) \subseteq X \to X$ und $B : D(B) \subseteq X \to X$ linear
mit $B \subseteq A$ und $\rho(B) \neq \emptyset$.
Dann gilt: $B = A$.
\dremark{Benutzt in: \ref{charMlC0X}}%
\end{lemma}

\begin{proof}
Sei $\lambda \in \rho(B)$.
Es gilt: $\lambda - B \subseteq \lambda - A$.
Somit folgt:\dremark{\cite{UllmannDiplom}, Satz 2.8.(a)}
$(\lambda - B)^{-1} \subseteq (\lambda - A)^{-1}$.
Da $(\lambda - B)^{-1} \in \mLinStet(X)$ ist, ergibt sich:
$(\lambda - B)^{-1} = (\lambda - A)^{-1}$.\dremark{Klar, sonst: \cite{UllmannDiplom}, 2.8.(b)}
Man erh"alt:\dremark{$(\lambda - B) = (\lambda - A)$, also } $B = A$.
\dremark{17.10.'07, S. 1}%
\dremark{Bew. gepr.}%
\end{proof}

Als einfache Folgerung erh"alt man:

\begin{folgerung}
Sei $X$ ein Banachraum, $A$ (bzw. $B$) Erzeuger einer $C_0$-Halbgruppe $(S_t)_\indtHG$
(bzw. $(T_t)_\indtHG$) auf $X$ mit $B \subseteq A$.
Dann gilt: $B = A$.
\end{folgerung}

\begin{proof}
Nach dem Satz von Hille-Yosida (Satz~\ref{satzHilleYosHG}) gilt: $\rho(B) \neq \emptyset$.
Mit Lemma~\ref{BSubseteqAFolgtGleich} erh"alt man: $B = A$.
\end{proof}

In der Praxis sind die Versionen des Satzes von Hille-Yosida
(Satz \ref{satzHilleYoskHG}, Folgerung~\ref{charQuasikontrHG} und Satz~\ref{satzHilleYosHG})
h"aufig schwierig anzuwenden,
da man eine explizite Kenntnis der Resolventenmenge ben"otigt.
Deshalb werden wir mit Hilfe dissipativer Operatoren
eine weitere Charakterisierung der Erzeuger kontraktiver $C_0$-Halb\-grup\-pen auf\/f"uhren.

\begin{definitn}
Sei $X$ ein Banachraum.
Ein Operator $A : D(A) \subseteq X \to X$
hei"st \defemph{dissipativ}, falls gilt:
\index[B]{dissipativer Operator}%
\[ \normlr{(\lambda-A)x} \geq \lambda \norm{x}
   \quad\text{f"ur alle } \lambda \in \mbbR_{>0} \text{ und } x \in D(A). \]
\end{definitn}

Mit Hilfe der folgenden Definition erh"alt man eine alternative
Charakterisierung dissipativer Operatoren:

\begin{definitn}
Sei $X$ ein Banachraum.
Definiere f"ur alle $x \in X$ die Menge
\index[S]{Jx@$\mcalJ(x)$ (Dualitätsmenge)}%
\[ \mcalJ(x) := \bigl\{ \varphi \in \dualr{X} \setfdg \varphi(x) = \norm{x}^2 = \norm{\varphi}^2 \bigr\}, \]
genannt \defemphi{Dualitätsmenge} im Punkte $x$,
wobei mit $\dualr{X}$\dremark{$=\mLinStet(X,\mbbK)$} der Dualraum von $X$ bezeichnet wird.
\end{definitn}

Nach dem Satz von Hahn-Banach ist die Dualit"atsmenge nicht-leer.
\skiptext

Dissipative Operatoren kann man wie folgt charakterisieren:

\begin{bemerkung}[\cite{EngelNagelSemigroups}, Proposition II.3.23]
Sei $X$ ein Banachraum.
Ein Operator $A : D(A) \subseteq X \to X$ ist genau dann dissipativ,
wenn f"ur jedes $x \in D(A)$ ein $\varphi \in \mcalJ(x)$ existiert mit
\begin{equation}\label{eqIffDissipativ}
\Re \varphi(Ax) \leq 0.
\end{equation}
Falls $A$ Erzeuger einer kontraktiven $C_0$-Halbgruppe ist,
dann gilt \eqref{eqIffDissipativ} f"ur alle $x \in D(A)$ und beliebige $\varphi \in \mcalJ(x)$.
\dliter{\cite{WernerFunkana3}, VII.4.15}%
\dremark{F"ur einige Funktionenr"aume ($C_0(\Omega)$, $L^p(\Omega,\mu)$) und Hilbertr"aume
  kann man die Menge $\mcalJ(x)$ leicht bestimmen und so pr"ufen, ob ein Operator dissipativ ist
  (\cite{EngelNagelSemigroups}, II.3.26, \cite{WernerFunkana3}, nach VII.4.14).}%
\end{bemerkung}

Im einem beliebigen Hilbertraum $H$ ist, wenn man $H$ kanonisch mit seinem Dualraum $\dualr{H}$
identifiziert,
\[ \mcalJ(x) = \{x\} \]
f"ur alle $x \in H$. Somit ist ein Operator $A$ in $H$ genau dann dissipativ,
wenn gilt:
\[ \Re \skalpr{Ax}{x} \leq 0 \]
f"ur alle $x \in D(A)$.
Dies gilt insbesondere f"ur selbstadjungierte Operatoren mit negativem Spektrum.
\dliter{\cite{WernerFunkana3}, nach VII.4.14}%
\skiptext

Es folgt die angek"undigte Charakterisierung der Erzeuger kontraktiver
$C_0$"=Halbgruppen (\cite{WernerFunkana6}, Theorem VII.4.16):

\begin{satz}[Satz von Lumer-Phillips]\label{SatzLumerPhillips}
Sei $X$ ein Banachraum und $A : D(A) \subseteq X \to X$ linear und dicht definiert.
Dann erzeugt $A$ genau dann eine kontraktive $C_0$-Halbgruppe, wenn $A$ dissipativ ist
und ein $\lambda \in \mbbR_{>0}$ so existiert, \dass $\lambda - A$ surjektiv auf $X$ ist.
\end{satz}

\dremww{
Es ist nicht sinnvoll, auf jedem Raum $C_0$-Halbgruppen zu betrachten:

\begin{satz}
Sei $X$ einer der folgenden Banachr"aume:
\begin{enumaufz}
\item $L^\infty(\Omega, \Sigma, \mu)$, wobei $(\Omega, \Sigma, \mu)$ ein positiver Ma"sraum sei,
\item $B(S, \Sigma)$, wobei $\Sigma$ eine $\sigma$-Algebra von Teilmengen von $S$ sei,
\item $C(\Omega)$, wobei $\Omega$ eine kompakter $\sigma$-stonescher Raum sei,
\item $C(\Omega)$, wobei $\Omega$ ein kompakter $F$-Raum sei,
\item ein Banachraum, der isomorph zu einem komplementierten Unterraum eines
  der R"aume in (i)--(iv) ist.
\end{enumaufz}
Dann ist jede $C_0$-Halbgruppe von beschr"ankten linearen Operatoren auf $X$
bereits gleichm"a"sig stetig.
\dremark{Es stellt sich die Frage, welche Teilr"aume komplementiert sind.}%
\dliter{\cite{Lotz85UniformConvOfOp}, S. 207;
  s. auch \cite{ArendtEtAl01LaplaceTransf}, Cor. 4.3.19, S. 275
  (hier nach A. Ullmann etwas allgemeiner);
  \cite{ArendtEtAl86OneParamSemigrOfPosOp}, A-II.3}%
\end{satz}

\begin{anmerkung}
In \cite{ArendtEtAl86OneParamSemigrOfPosOp}, B-II werden Erzeuger
von positiven Halbgruppen (\dheisst wohl $C_0$-Halbgruppen von positiven Operatoren)
in $C_0(\Omega)$ charakterisiert, wobei $\Omega$ ein lokalkompakter Hausdorffraum sei.
\end{anmerkung}

\begin{bemerkung}
F"ur einen Operator $(A, D(A))$ auf dem Banachraum $C_0(\Omega)$
oder $L^p(\Omega,\mu)$, $p \in [1,\infty[$, sind "aquivalent:
\dliter{\cite{EngelNagelSemigroups}, II.2.9}%
\begin{enumaequiv}
\item $(A,D(A))$ ist Erzeuger einer stark stetigen Multiplikationshalbgruppe.
\item $(A,D(A))$ ist ein Multiplikationsoperator derart, \dass
\[ \exists w \in \mbbR : \{ \lambda \in \mbbC \setfdg \Re \lambda > w \} \subseteq \rho(A). \]
\end{enumaequiv}
\end{bemerkung}

Es ist bemerkenswert, \dass bereits Bedingung \bewitemph{(b)} die Existenz einer $C_0$-Halbgruppe
gew"ahrleistet.\dremark{Siehe Z. 20.9.'07/2}

Man kann Multiplikationshalbgruppen nicht nur in $C_0(\Omega)$ und
$L^p(\Omega,\mu)$, sondern auch auf Banachverb"anden charakterisieren:

\begin{satz}
Sei $A$ der Erzeuger einer $C_0$-Halbgruppe $(T_t)_\indtHG$ auf einem
$\sigma$-ordnungsvollst"andigen reellen oder komplexen Banachverband $X$.
Dann sind die folgenden Aussagen "aquivalent:
\begin{enumaequiv}
\item $(T_t)_\indtHG$ ist eine Multiplikationshalbgruppe,
  \dheisst{} $T_t$ ist ein Multiplikationsoperator (\dheisst{} ein Element aus dem Zentrum)
  f"ur alle $t \in \mbbR_{>0}$.
\item Es existiert ein $\lambda \in \rho(A)$ derart, \dass $R(\lambda, A)$
  ein Multiplikationsoperator ist.
\item $R(\lambda, A)$ ist ein Multiplikationsoperator f"ur alle $t \in \rho(A)$.
\item $A$ ist lokal und $D(A)$ ein Ideal in $X$.
\item Ist $f \in D(A)$, dann gilt $Pf \in D(A)$ f"ur jede Band-Projektion $P$ auf $X$
und $APf = PAf$.
\dliter{\cite{ArendtEtAl86OneParamSemigrOfPosOp}, C-II, Th. 5.13}%
\end{enumaequiv}
\end{satz}

\begin{definitn}
Sei $\mcalM$ eine Von-Neumann-Algebra.
Ein $T \in \mLinStet(\mcalM)$ hei"st \defemphi{Schwarz-Abbildung},
falls gilt:
\dliter{\cite{ArendtEtAl86OneParamSemigrOfPosOp}, S. 370}%
\[ T(x) T(x)^* \leq T(xx^*) \qquad \text{f"ur alle } x \in \mcalM. \]
\end{definitn}

\begin{bemerkung}
Sei $\mcalM$ eine Von-Neumann-Algebra.
Jede Schwarz-Abbildung auf $M$ ist eine Kontraktion.
\dliter{\cite{ArendtEtAl86OneParamSemigrOfPosOp}, S. 370}%
\end{bemerkung}

\begin{proof}
Sei $x \in \mcalM$.
Es gilt $0 \leq T(x) T(x)^* \leq T(xx^*)$.
Mit \cite{Wegge-OlsenKTheory}, S. 308, folgt:
\[ \norm{Tx}^2 = \norm{T(x) T(x)^*} \leq \norm{T(xx^*)} \leq \norm{T} \cdot \norm{x}^2. \]
Also gilt $\norm{Tx} \leq \sqrt{\norm{T}} \cdot \norm{x}$ und $\norm{Tx} \leq \norm{T} \cdot \norm{x}$.
Da $\norm{T}$ die kleinste Lipschitz-Konstante ist, folgt: $\norm{T} \leq 1$.
\end{proof}

Zu kontraktiven $C_0$-Halbgruppen auf einem Hilbertraum:

\begin{satz}
Sei $(T_t)_\indtHG$ eine kontraktive $C_0$-Halbgruppe auf einem Hilbertraum $H$.
Dann existiert ein Hilbertraum $K$ und eine unit"are Gruppe $(U_t)_\indtGr$ auf $K$
derart, \dass
\[ J T_t = P U_t J \]
f"ur alle $\indtHG$,
wobei $J$ eine isometrische Einbettung ist und $P : K \to J(H)$
die orthogonale Projektion.
\dliter{\cite{KunstmannWeis04MaxLpRegForParEq}, Theorem 10.11, S. 214}%
\end{satz}

Literatur zum Ableitungsoperator auf verschiedenen R"aumen:
\begin{enumaufz}
\item \cite{WernerFunkana3}, S. 335
\item \cite{EngelNagelSemigroups}, S. 33--, S. 66-- (Erzeuger)
\item \cite{ReedSimonMethMathPhysIRevEd}, S. 251
\item \cite{RudinFuncAna:73}, Example 13.4 und Exercise 13.19 (S. 365)
\end{enumaufz}
}%

\section{Erzeuger von $C_0$-Gruppen}

Man kann ein dem Satz von Hille-Yosida (Satz~\ref{satzHilleYosHG}) entsprechendes Resultat
auch f"ur $C_0$-Gruppen formulieren.
Hierf"ur wird zun"achst definiert, was man unter einem Erzeuger einer $C_0$-Gruppe versteht:

\begin{definitn}
Der \defemph{Erzeuger} einer
$C_0$-Gruppe $(T_t)_\indtGr$ auf einem Banachraum $X$ ist der Operator $A : D(A) \subseteq X \to X$,
der durch
\[ Ax = \lim_{t \to 0} \frac{T_t x - x}{t} \]
f"ur alle $x$ aus seinem Definitionsbereich
\index[B]{Erzeuger einer $C_0$-Gruppe}%
\dliter{\cite{EngelNagelSemigroups}, Def. in II.3.11}
\[ D(A) = \left\{ x \in X \setfdg \lim_{t \to 0} \frac{T_t x - x}{t} \text{ existiert} \right\} \]
definiert wird.
\end{definitn}

Sei $(T_t)_\indtGr$ eine $C_0$-Gruppe mit Erzeuger $A$.
Dann kann man $T^{+}_t := T_t$ und $T^{-}_t := T_{-t}$ f"ur alle $\indtHG$ definieren.
Es folgt aus der obigen Definition, \dass $(T^{+}_t)_\indtHG$ und $(T^{-}_t)_\indtHG$
$C_0$-Halbgruppen mit Erzeuger $A$ bzw. $-A$ sind.
Falls ein Operator $A$ also Erzeuger einer $C_0$-Gruppe ist,
sind $A$ und $-A$ Erzeuger einer $C_0$-Halbgruppe.
Der folgende Satz zeigt, \dass auch die Umkehrung gilt
(vgl. \cite{EngelNagelSemigroups}, Abschnitt II.3.11):

\begin{satz}[Satz von Hille-Yosida f"ur $C_0$-Gruppen]\label{satzHilleYosidaGr}
Sei $X$ ein Banachraum und $A : D(A) \subseteq X \to X$ linear.
Sei $\omega \in \mbbR$ und $M \in \mbbR_{\geq 1}$.
Dann sind die folgenden Aussagen "aquivalent:
\begin{enumaequiv}
\item $A$ erzeugt eine $C_0$-Gruppe $(T_t)_\indtGr$ mit der Eigenschaft:
\[ \norm{T_t} \leq M \mre^{\omega \abs{t}} \qquad\text{f"ur alle } \indtGr. \]

\item $A$ und $-A$ erzeugen $C_0$-Halbgruppen
$(T^+_t)_\indtHG$ bzw. $(T^-_t)_\indtHG$, die folgendes erf"ullen:
\[ \norm{T^+_t}, \norm{T^-_t} \leq M \mre^{\omega t} \qquad\text{f"ur alle } \indtHG. \]

\item $A$ ist dicht definiert, abgeschlossen und
f"ur jedes $\lambda \in \mbbR$ mit $\abs{\lambda} > \omega$
hat man $\lambda \in \rho(A)$ und
\[ \normbig{ [(\abs{\lambda} - \omega) R(\lambda,A)]^k } \leq M
\qquad\text{f"ur alle } k \in \mbbN. \]
\end{enumaequiv}
\end{satz}

\dremww{
\section{Spektraltheorie f"ur $C_0$-Halbgruppen}

\begin{bemerkung}
Sei $\omega_0$ die Wachstumsschranke der $C_0$-Halbgruppe $(T_t)_\indtHG$ auf einem Banachraum $X$.
Dann gilt
\[ r(T_t) = \mre^{\omega_0 t} \qquad\text{f"ur alle } \indtHG. \]
\dliter{\cite{EngelNagelSemigroups}, IV.2.2,
  \cite{ArendtEtAl86OneParamSemigrOfPosOp}, A-III, Prop. 1.1}%
\end{bemerkung}

\begin{satz}
Sei $(T_t)_\indtHG$ eine schlie"slich normstetige Halbgruppe mit Erzeuger $A$
auf einem Banachraum $X$.
Dann gilt der Spektralabbildungssatz:\dliter{\cite{EngelNagelSemigroups}, IV.3.10}
\[ \sigma(T_t) \setminus \{0\} = \mre^{t \sigma(A)} \qquad\text{f"ur alle } \indtHG. \]
\end{satz}

\begin{folgerung}
F"ur eine schlie"slich normstetige Halbgruppe $(T_t)_\indtHG$ mit Erzeuger $A$
auf einem Banachraum $X$ gilt:\dliter{\cite{EngelNagelSemigroups}, IV.3.11}
\[ s(A) = \omega_0. \]
\end{folgerung}

Da viele verschiedene Regularit"atsbedingungen von Halbgruppen die Eigenschaft
\glqq schlie"slich normstetig\grqq implizieren, erh"alt man:

\begin{folgerung}
Der Spektralabbildungssatz
\[ \sigma(T_t) \setminus \{0\} = \mre^{t \sigma(A)} \qquad\text{f"ur alle } \indtHG \]
und die folgende Bedingung
\[ s(A) = \omega_0 \]
gelden f"ur die folgenden Klassen von $C_0$-Halbgruppen:
\begin{enumaufz}
\item schlie"slich kompakte Halbgruppen,
\item schlie"slich differenzierbare Halbgruppen,
\item analytische Halbgruppen,
\item normstetige Halbgruppen.
\end{enumaufz}
\dliter{\cite{EngelNagelSemigroups}, IV.3.12}
\end{folgerung}

\begin{anmerkung}
\begin{itemize}
\item Literatur zu $C_0$-Halbgruppen (von A. Ullmann empfohlen):
\cite{EngelNagelSemigroups}, \cite{Pazy83SemigroupsApplicatPDE},
\cite{Davies80OneParameterSemigroups}, \cite{Goldstein85SemigroupsAndApplicat}.

\item Weitere Literatur zu $C_0$-Halbgruppen:
  \cite{BratteliRobinson87OpAlgebras1}: Hier werden viele Aussagen "uber Operatorhalbgruppen
  gemacht, die $\sigma(X,F)$-stetig sind.
  Ferner werden Halbgruppen und analytische Elemente behandelt (siehe \cmzB{} Th. 3.1.18.).

\item Literatur zu Kontraktionshalbgruppen auf Hilbertr"aumen:
  \cite{Davies80OneParameterSemigroups}, Chapter 6.
\end{itemize}
\end{anmerkung}
}%

\dremark{
\chapter{Problemliste}

\textbf{zutun:}

\ref{bspNichtKontrInHr}: Was gilt f"ur $OH$? F"ur $H^c$? Erw"ahnen. (WW)

\begin{enumerate}[(1)]
\item Gilt $\Adjlor{X} \subseteq \Multlco{X}$?
  Gilt dies zumindest f"ur selbstadjungierte Elemente? (WW)
\item Gibt es eine Variante des Charakterisierungssatzes \ref{charAdjlcoHR}
  f"ur $\Adjloru{X}$ (Charakterisierung mittels Hilbertraum),
  die besser angepa\cms{}t ist?
  (Idee: Der Erzeuger auf dem Hilbertraum hat die Gestalt: $\cmpmatrix{0 & C \\ C^* & 0}$)
  (WW)
\item In \cite{Damaville04RegulariteDesOp}, Proposition 1.2.(iv), steht eine
  Bedingung daf"ur, wann ein Operator regul"ar ist.
  Gibt es eine entsprechende Bedingung f"ur unbeschr"ankte Multiplikatoren? (WW)
\item Sei $X$ ein \cmsa Operatorraum.
  F"ur alle $T \in \Multlor{X}$ mit $T^* \in \Multlor{X}$ gilt: $T \in \Adjlor{X}$.
  Gilt etwas "Ahnliches f"ur unbeschr"ankte Multiplikatoren? (WW)
\item Folgt Beispiel~\ref{bspMultlcoKH} ($\Multlco{\kptOp(H)} = \dots$)
  auch mit einem Charakterisierungssatz? (WW)
\item Beweis von \ref{zshgRegOpUnbeschrMult} ($\Adjloru{E} = \Multwor{E}$)
  "uber Multiplikation mit Matrizen m"oglich?
  Benutze \cmzB St"orungssatz von Damaville. (WW)
\item M"ogliche alternative Charakterisierung eines unbeschr"ankten Multiplikators $T$:
  $T$ ist Einschr"ankung eines regul"aren Operators $S$ auf einem
  geeigneten \hcsmodul $E$ (evtl. $E = \ternh{X}$) auf $X$
  (siehe \ref{halbgrHochliftenAufTX}).
  So etwas ist bewiesen, siehe \ref{AdjlcoEqRegOpAufTX}. (WW)
\item Gilt $\Adjlor{X} = \{ A \in \Adjloru{X} \setfdg D(A) = X \}$
  (siehe \ref{AdjlorSubseteqAdjloru})?
\item Wie h"angen die tern"are \sterns{}H"ulle und die injektive \sterns{}H"ulle zusammen?
\end{enumerate}

St"orungstherorie:

\begin{enumerate}[(1)]
\item Kann man weitere Resultate der St"orungstheorie auf Operatorr"aume "ubertragen
  und auf unbeschr"ankte Multiplikatoren anwenden?
  Siehe \ref{stoerungsthUebertragen},
  insbesondere \cite{Pazy83SemigroupsApplicatPDE}, Theorem 3.4.3, und
  \cite{Damaville04RegulariteDesOp}, Prop. 2.1 ((WW): (1)).
\end{enumerate}

\begin{enumerate}[(1)]
\item Sei $H$ ein Hilbertraum.
$H$ sei so mit einer Operatorraumnorm ausgestattet, \dass die Fouriertransformation $\mathscr{F}$ eine Isometrie ist.
Sei $\varphi$ die nat"urliche Einbettung des Operatorraumes $H$ in $\mLinStet(K)$,
wobei $K$ ein Hilbertraum sei.
$\mathscr{F}$ bildet einen Differentialoperator auf einen Multiplikator ab.
Wird dann dieser Multiplikator durch $\varphi$ auf einen unbeschr"ankten Multiplikator im
\cstern{}Algebren-Sinne abgebildet? (18.5.'05)

\item (Yosida-Approximation) Ansatz, um unbeschr"ankte Operatoren auf beschr"ankte zur"uckzuf"uhren:
Sei $X$ ein Banachraum und $A$ der Erzeuger einer Kontraktionshalbgruppe
(evtl. gen"ugen schw"achere Voraussetzungen).
F"ur alle $\lambda \in \mbbR_{>0}$ sei
\[ A_\lambda := \lambda A(\lambda - A)^{-1}
  = \lambda^2 (\lambda - A)^{-1} \in \mLinStet(X), \]
genannt die \defemphi{Yosida-Approximation}.
Es gilt:\dliter{\cite{WernerFunkana3}, VII.4.11}
\[ \lim_{\lambda \rightarrow \infty} A_\lambda x = Ax \qquad \text{f"ur alle } x \in D(A). \]

\item Gilt \ref{injHuelleMnX} auch f"ur unitale Operatorsysteme?
\end{enumerate}
}%



\ifthenelse{\boolean{finalVersion}}
{\cleardoublepage}%
{\clearpage}%

\addcontentsline{toc}{chapter}{Literaturverzeichnis}
\bibliographystyle{geralpha}
\bibliography{Diverses}

\dremark{N"utzliche Literatur zu \jbstripel{}n:
H. Upmeier.
\textsl{Symmetric Banach Manifolds and Jordan C*-Algebras}, North-Holland, 1985.}%

\renewcommand{\indexname}{Symbolverzeichnis}
\ifthenelse{\boolean{schreibidx}}%
{\ifthenelse{\boolean{finalVersion}}
{\cleardoublepage}%
{\clearpage}%
 \addcontentsline{toc}{chapter}{Symbolverzeichnis}
 \printindex[S]%
\ifthenelse{\boolean{finalVersion}}
{\cleardoublepage}%
{\clearpage}%
 \addcontentsline{toc}{chapter}{Stichwortverzeichnis}
 \printindex[B]%
}
{}

\dremark{
\vspace*{3ex}
\textbf{Schreibweise von Begriffen:}
Algebrenhomomorphismus,
Algebrenisomorphismus,
approximierendes Einselement,
Banachalgebra,
beidseitiges Ideal (anstelle von zweiseitiges Ideal),
Darstellung auf \emph{einem} Vektorraum (anstelle von \glqq{}auf einen\grqq,
    denn: operiert auf dem VR),
Disk-Algebra,
endlich-dimensional,
Funktional auf (nicht in),
selbstadjungiert (nicht hermitesch),
Gelfandtransformation,
Hilbertraum,
Ideal von $A$ (nicht in $A$),
kompakt(quasikompakt und hausdorffsch),
Linksideal,
lokalkompakt(lokalquasikompakt und hausdorffsch),
nicht-ausgeartet,
nicht-leer,
nicht-trivial,
Prähilbertraum,
Quotiengenalgebra (nicht Faktor-),
Unteralgebra (nicht Teil-),

\textbf{Schreibweise von R"aumen und Elementen:}
$L^p(\Omega)$,
Algebrenelemente $a,b,c$ (nicht $x,y,z$),
Banachraum $X,Y,Z$ (Elemente $x,y,z$),
\cstern{}Algebra $\mfrakA,\mfrakB,\mfrakC$ (Elemente $a,b,c$),
Elemente von $C_0(\Omega)$: f,g,h,
Hilbert-\cstern{}Modul $E,F,G$ (Elemente $x,y,z$),
Hilbertraum $H,K,L$ (Elemente $\xi,\eta,\zeta$),
Linksideal $L$,
komplexe Matrizen $\alpha, \beta, \gamma$,
lokalkompakter Hausdorffraum $\Omega$,
Operatorraum $X,Y,Z$ (Elemente $x,y,z$)
TRO $Z$, $Z_1$ (Elemente $x,y,z$)
Vektorraum $V,W,U$ (Elemente $v,w,u$),

\textbf{Schreibweise von Abbildungen:}
$C_0$-Halbgruppen: $(T_t)_\indtHG, (S_t)_\indtHG, (R_t)_\indtHG$,
Erzeuger von $C_0$-Halbgruppen: $A,B,C$,
Isomorphismen: $\varphi, \psi$,
Operatoren: $T,S,R$,

\textbf{Diverses:}
$M \subsetneq N$ (subsetneq) (nicht $\subsetneqq$),
dicht definiert, abgeschlossen (in dieser Reihenfolge)
das Monoid

\textbf{Konventionen:}
Hilbertraum: das Skalarprodukt ist linear in der 2. Variablen.
Hilbert-\cstern{}Modul: stets ein Rechts-Modul, dessen Skalarprodukt linear in der 2. Variablen ist.
Es ist $f : \mbbR \to \mbbR, x \mapsto x^2$, stetig.
}%


\end{document}

%% file: MathCommands.tex
\renewcommand{\proofname}{Beweis}
\newcommand{\abs}[1]{\ensuremath{\lvert#1\rvert}}
\newcommand{\bigobot}{\ensuremath{\mathbin{\text{\textcircled{$\bot$}}}}}
\newcommand{\chara}[1]{\ensuremath{\operatorname{char}(#1)}}  
\newcommand{\codim}{\operatorname{codim}}%
\newcommand{\complI}{\ensuremath{i}}          
\newcommand{\kard}[1]{\ensuremath{\lvert#1\rvert}}
\newcommand{\norm}[1]{\ensuremath{\lVert#1\rVert}}
\newcommand{\obGausskl}[1]{\ensuremath{\lceil#1\rceil}}   
\newcommand{\QuatAlg}{\ensuremath{\text{I\hspace*{-0.20em}H}}}
\newcommand{\restring}{\ensuremath{\rvert}}   
\newcommand{\restringa}[1]{\ensuremath{\rvert_{#1}}}   
\newcommand{\setfdg}{\;;\;}  
\newcommand{\untGausskl}[1]{\ensuremath{\lfloor#1\rfloor}}   
\newcommand{\schlietmcMathbb}[1]{\mathbbm{#1}}
\newcommand{\schlietmceMathbb}[1]{\ensuremath{\schlietmcMathbb{#1}}}
\newcommand{\mbbC}{\schlietmceMathbb{C}}
\newcommand{\mbbD}{\schlietmceMathbb{D}}
\newcommand{\mbbF}{\schlietmceMathbb{F}}
\newcommand{\mbbH}{\schlietmceMathbb{H}}
\newcommand{\mbbK}{\schlietmceMathbb{K}}
\newcommand{\mbbN}{\schlietmceMathbb{N}}
\newcommand{\mbbO}{\schlietmceMathbb{O}}
\newcommand{\mbbP}{\schlietmceMathbb{P}}
\newcommand{\mbbQ}{\schlietmceMathbb{Q}}
\newcommand{\mbbR}{\schlietmceMathbb{R}}
\newcommand{\mbbT}{\schlietmceMathbb{T}}
\newcommand{\mbbV}{\schlietmceMathbb{V}}
\newcommand{\mbbZ}{\schlietmceMathbb{Z}}
\newcommand{\CDach}{\ensuremath{\hat{\schlietmcMathbb{C}}}}   
\newcommand{\mcalo}{\ensuremath{\mathcal{o}}}
\newcommand{\mcalA}{\ensuremath{\mathcal{A}}}
\newcommand{\mcalB}{\ensuremath{\mathcal{B}}}
\newcommand{\mcalC}{\ensuremath{\mathcal{C}}}
\newcommand{\mcalD}{\ensuremath{\mathcal{D}}}
\newcommand{\mcalE}{\ensuremath{\mathcal{E}}}
\newcommand{\mcalF}{\ensuremath{\mathcal{F}}}
\newcommand{\mcalG}{\ensuremath{\mathcal{G}}}
\newcommand{\mcalH}{\ensuremath{\mathcal{H}}}   
\newcommand{\mcalI}{\ensuremath{\mathcal{I}}}
\newcommand{\mcalJ}{\ensuremath{\mathcal{J}}}
\newcommand{\mcalK}{\ensuremath{\mathcal{K}}}
\newcommand{\mcalL}{\ensuremath{\mathcal{L}}}
\newcommand{\mcalM}{\ensuremath{\mathcal{M}}}
\newcommand{\mcalN}{\ensuremath{\mathcal{N}}}
\newcommand{\mcalO}{\ensuremath{\mathcal{O}}}
\newcommand{\mcalP}{\ensuremath{\mathcal{P}}}
\newcommand{\mcalQ}{\ensuremath{\mathcal{Q}}}
\newcommand{\mcalR}{\ensuremath{\mathcal{R}}}
\newcommand{\mcalS}{\ensuremath{\mathcal{S}}}
\newcommand{\mcalT}{\ensuremath{\mathcal{T}}}
\newcommand{\mcalU}{\ensuremath{\mathcal{U}}}
\newcommand{\mcalV}{\ensuremath{\mathcal{V}}}
\newcommand{\mcalW}{\ensuremath{\mathcal{W}}}
\newcommand{\mcalX}{\ensuremath{\mathcal{X}}}
\newcommand{\mcalY}{\ensuremath{\mathcal{Y}}}
\newcommand{\mcalZ}{\ensuremath{\mathcal{Z}}}
%
\def\Dhaken#1{\setbox0=\hbox{$\displaystyle #1$\hskip 0.1em}\hbox{\vtop
{\copy0\vskip0.2pt\hrule}\vrule height 0.6\ht0 width 0.5pt}}
\def\Thaken#1{\hbox{\vtop{\hbox{$\textstyle #1$\hskip 0.05em}\vskip 0.3ex
\hrule}\vrule height 0.5ex width 0.5pt}}
\def\Shaken#1{\hbox{\vtop{\hbox{$\scriptstyle #1$\hskip 0.025em}\vskip
0.2ex\hrule}\vrule height 0.5ex}}
\def\SShaken#1{\hbox{\vtop{\hbox{$\scriptscriptstyle #1$\hskip
0.01em}\vskip 0.2ex\hrule}\vrule height 0.5ex}}
\def\haken#1{\mathchoice{\Dhaken{#1}}{\Thaken{#1}}{\Shaken{#1}}{\SShaken{#1}}}